\documentclass{memo-l}

\usepackage{amscd}
\usepackage[arrow,matrix,curve,cmtip,ps]{xy}
\usepackage{amssymb}

\theoremstyle{plain}   
\newtheorem{thm}{Theorem}[chapter]
\newtheorem{cor}[thm]{Corollary} 
\newtheorem{lem}[thm]{Lemma}    
\newtheorem{prop}[thm]{Proposition}

\theoremstyle{definition}
\newtheorem{defn}[thm]{Definition}

\theoremstyle{remark}
\newtheorem{rem}[thm]{Remark}    
\newtheorem{ex}[thm]{Example}   

\numberwithin{equation}{chapter}

\newcommand{\upref}[1]{\textup{\ref{#1}}}
\newcommand{\upeqref}[1]{\textup{\eqref{#1}}}
\newcommand{\chapref}[1]{Chapter~\upref{#1}}
\newcommand{\appxref}[1]{Appendix~\upref{#1}}
\newcommand{\secref}[1]{Section~\upref{#1}}
\newcommand{\thmref}[1]{Theorem~\upref{#1}}
\newcommand{\corref}[1]{Corollary~\upref{#1}}
\newcommand{\lemref}[1]{Lemma~\upref{#1}}
\newcommand{\propref}[1]{Proposition~\upref{#1}}
\newcommand{\defnref}[1]{Definition~\upref{#1}}
\newcommand{\remref}[1]{Remark~\upref{#1}}
\newcommand{\exref}[1]{Example~\upref{#1}}
\newcommand{\diagref}[1]{Diagram~\upeqref{#1}}
\newcommand{\eqeqref}[1]{Equation~\upeqref{#1}}

\DeclareMathOperator{\ad}{Ad}
         \newcommand{\Ad}{\ad}
\DeclareMathOperator{\aut}{Aut}
         \newcommand{\Aut}{\aut}

\DeclareMathOperator{\id}{id}

\DeclareMathOperator{\ind}{Ind}
         \newcommand{\Ind}{\ind}
\DeclareMathOperator{\infl}{Inf}
         \newcommand{\Infl}{\infl}
\DeclareMathOperator{\mor}{mor}
         \newcommand{\Mor}{\mor}

\DeclareMathOperator{\rep}{Rep}
         \newcommand{\Rep}{\rep}
\DeclareMathOperator{\spn}{span}
\DeclareMathOperator{\supp}{supp}

\renewcommand{\c}[1]{{\mathcal{#1}}}
\newcommand{\A}{\c{A}}
\newcommand{\B}{\c{B}}
\newcommand{\C}{\c{C}}
\newcommand{\D}{\c{D}}

\renewcommand{\H}{\c{H}}

\newcommand{\K}{\c{K}}
\renewcommand{\L}{\c{L}}

\renewcommand{\b}[1]{{\mathbb{#1}}}
\newcommand{\bbC}{\b{C}}
\newcommand{\R}{\b{R}}

\newcommand{\CC}{\bbC} 
\newcommand{\RR}{\R}
\newcommand{\FF}{\b{F}}

\newcommand{\what}[1]{\widehat{#1}}

\newcommand{\smtx}[1]{\left(\begin{smallmatrix} #1\end{smallmatrix}\right)}
\newcommand{\mtx}[1]{\begin{pmatrix} #1 \end{pmatrix}}

\newcommand{\rev}[1]{\widetilde{#1}}
\newcommand{\righttext}[1]{\qquad\text{#1 }}
\newcommand{\midtext}[1]{\quad\text{#1}\quad}
\newcommand{\mathbox}[1]{\mbox{\ensuremath{#1}}}

\newcommand{\deq}{=} 

\newcommand{\lk}{\left\langle}
\newcommand{\rk}{\right\rangle}
\newcommand{\pproj}{\smtx{1&0\\0&0}}
\newcommand{\qproj}{\smtx{0&0\\0&1}}
\renewcommand{\:}{\colon}
\newcommand{\<}{\langle}
\renewcommand{\>}{\rangle}
\newcommand{\iso}{\xrightarrow{\cong}}
\newcommand{\sptilde}{\tilde{\ }}
\newcommand{\sphat}{\hat{\ }}
\newcommand{\eps}{\epsilon}
\newcommand{\Chi}{\raisebox{2pt}{\ensuremath{\chi}}} 
\renewcommand{\phi}{\varphi}
\newcommand{\ip}{\text{imp}}
\renewcommand{\d}{\cdot} 
\newcommand{\dec}{^{\operatorname{dec}}}
\newcommand{\deltahat}{{\hat{\delta}}}
\newcommand{\deltahathat}{{\hat{\hat{\delta}}}}

\newcommand{\cotimes}{\mathrel{\sharp}}

\newcommand{\dashind}{\!\operatorname{-Ind}}

\newcommand{\pause}{\renewcommand{\qed}{}\end{proof}}
\newcommand{\resume}[1]{\begin{proof}[Back to the proof of #1]}

\newcommand{\ie}{\emph{i.e.}}
\newcommand{\eg}{\emph{e.g.}}
\newcommand{\cf}{\emph{cf.}}
\newcommand{\cstar}{\ensuremath{C^*}-}

\makeindex

\begin{document}

\frontmatter

\title[A Categorical Approach to Imprimitivity]%
{A Categorical Approach to Imprimitivity Theorems
for $C^*$-Dynamical Systems}

\author[Echterhoff]{Siegfried Echterhoff}
\address{Westf\"alische Wilhelms-Universit\"at\\
Mathematisches Institut\\
Einsteinstr. 62\\
D-48149 M\"unster\\
Germany}
\email{echters@math.uni-muenster.de}

\author[Kaliszewski]{S. Kaliszewski}
\address{Department of Mathematics\\
Arizona State University\\
Tempe, Arizona 85287}
\email{kaliszewski@asu.edu}

\author[Quigg]{John Quigg}
\address{Department of Mathematics\\
Arizona State University\\
Tempe, Arizona 85287}
\email{quigg@math.la.asu.edu}

\author[Raeburn]{Iain Raeburn}
\address{Department of Mathematics\\
University of Newcastle\\
Newcastle, New South Wales 2308
\\Australia}
\email{iain@maths.newcastle.edu.au}

\thanks{This research was partially supported by
National Science Foundation Grant \#DMS9401253,
the Australian Research Council, 
the Deutsche Forschungsgesellschaft (SFB 478), 
and the European Union.}



\date{May 30, 2002}

\begin{abstract}

Imprimitivity theorems provide a fundamental tool for studying the
representation theory and structure of crossed-product $C^*$-algebras.  
In this work, we show that the Imprimitivity Theorem for
induced algebras, Green's Imprimitivity Theorem for actions of groups,
and Mansfield's Imprimitivity Theorem for coactions of groups can all be
viewed as natural equivalences between various crossed-product functors
among certain equivariant categories.

The categories involved have $C^*$-algebras with actions or coactions
(or both) of a fixed locally compact group $G$ as their objects, and
equivariant equivalence classes of right-Hilbert bimodules as their
morphisms.  Composition is given by the balanced tensor product of
bimodules. 

The functors involved arise from taking crossed products; 
restricting, inflating, and decomposing actions and coactions; 
inducing actions; and various combinations of these.  

Several applications of this categorical approach are also presented,
including some intriguing relationships between the Green and Mansfield
bimodules, and between restriction and induction of representations.

\end{abstract}

\maketitle

\setcounter{page}{6}
\tableofcontents


\mainmatter

%
%

\chapter*{Introduction}
\label{intro-chap}

Given a dynamical system $(A,G,\alpha)$ in which a locally
compact group $G$ acts by automorphisms of a $C^*$-algebra $A$,
Mackey and Takesaki's \emph{induction} process allows us to construct
representations of $(A,G,\alpha)$ from representations of the system
$(A,H,\alpha|_H)$ associated to any closed subgroup $H$ of $G$.  Much is
known about induction: there are imprimitivity theorems which allow us
to recognize induced representations, and the process is functorial
with respect to intertwining operators.

In the modern framework of Rieffel, one introduces the crossed product
$A\times_\alpha G$, which is a $C^*$-algebra encapsulating the
representation theory of $(A,G,\alpha)$, and induces instead from
$A\times_\alpha H$ to $A\times_\alpha G$; induction of representations
from one $C^*$-algebra $D$ to another $C$ is achieved by tensoring the
underlying Hilbert space with a Hilbert bimodule $_CX_D$, which has a
$D$-valued inner product and in which the left action of $C$ is by
adjointable operators.  An \emph{imprimitivity theorem} tells us how
to expand the left action of $C$ to one of a larger algebra $E$ in
such a way that $_EX_D$ is an imprimitivity bimodule --- that is,
reversible.  The theorem then says that a representation of $C$ is
equivalent to one induced from $D$ if and only if there is a
compatible representation of $E$.

\emph{Duality} tells us how to recover a dynamical system
$(A,G,\alpha)$ from its crossed product $A\times_\alpha G$.  When $G$
is abelian, the crossed product carries a canonical dual action
$\widehat \alpha$ of the dual group $\widehat G$, and the
Takesaki-Takai Duality Theorem says that the double dual system
$\big((A\times_\alpha G)\times_{\widehat\alpha}\widehat
G,G,\widehat{\widehat\alpha}\big)$ is Morita equivalent to the original one. 
For nonabelian $G$, one has to use instead the dual coaction of $G$,
and recover the system from the crossed product by this dual coaction. 
For duality to be a useful tool, one has to understand these coactions
and their crossed products, and a good deal of progress has been made
in the past 15 years.  (An overview of this area has been provided in
an Appendix; see also \cite{RaeCPC} for a recent survey.)  Crucial
for us is Mansfield's theory of induction for crossed products by
coactions: he provides a Hilbert bimodule which allows us to induce
representations from crossed products by quotient groups, and an
imprimitivity theorem which characterizes these induced
representations.

Induction and duality interact in deep and mysterious ways.  One
general principle appears to be that duality swaps induction of
representations with restriction of representations.  This is
enormously appealing: restriction of representations (for example,
passing from a representation $U$ of $G$ to the representation $U|_H$
of a subgroup $H$) is ostensibly a trivial process.  Theorems making
this induction-restriction duality precise have been proved, first for
abelian groups in \cite{EchDI}, and later for arbitrary groups in
\cite{kqr:resind,ekr}.  We have gradually learned that it is best to
prove such theorems by manipulating the Hilbert bimodules which
implement the various induction and restriction processes;
however, the bimodules involved are hard to work with ---
especially Mansfield's --- and the results can safely be described as
``technically challenging''.  To make things worse, applications
frequently require that various isomorphisms and equivalences are
equivariant, and one is continually having to construct compatible
coactions on bimodules and check that they carry through complicated
arguments.  So it is definitely of interest to find a more systematic
approach.

Our goal here is to provide such a systematic approach
and to use it to complete our program of induction-restriction duality.  
We shall show that
many of the key technical problems in this area amount to asking for
functoriality of some construction or naturality of some equivalence
between functors.  Asking for equivalences to be equivariant amounts
to asking for an equivalence in a different category, one which
includes coactions or actions in its objects and morphisms.  We have
found that functoriality of the various crossed-product constructions
encompasses many results of the kind ``Morita equivalent systems have
Morita equivalent crossed products'', and naturality of the
equivalences many results of the kind ``induction is compatible with
Morita equivalence''.

To help see how our approach works, we consider one of our main
theorems.  It concerns the generalization of Green's Imprimitivity
Theorem to crossed products of induced algebras, which is, loosely
speaking, the analogue of the imprimitivity theorem for actions
$\alpha$ of a subgroup $H$ which do not extend to actions of $G$.  The
induced algebra $\Ind_H^G(A,\alpha)$ is a subalgebra of $C_b(G,A)$
which carries a left action $\tau$ of $G$ by translation, and the
generalization says that the crossed product
$\Ind_H^G(A,\alpha)\times_\tau G$ is Morita equivalent to
$A\times_\alpha H$.  We shall prove that this equivalence is natural,
and that it is equivariant for the dual coaction $\widehat \tau$ of
$G$ on $\Ind_H^G(A,\alpha)\times_\tau G$ and the inflation
$\Infl\widehat\alpha$ to $G$ of the dual coaction on $A\times_\alpha
H$.  To make this precise, we have to set up categories $\mathcal{C}$
of $C^*$-algebras, $\mathcal{A}(G)$ of dynamical systems
$(A,G,\alpha)$, and $\mathcal{C}(G)$ of cosystems $(A,G,\delta)$ in
which $\delta$ is a coaction of $G$ on $A$.  We then prove that
$(A,G,\alpha)\mapsto (\Ind_H^G(A,\alpha)\times_\tau G,\widehat \tau)$
and $(A,G,\alpha)\mapsto (A\times_\alpha H,\Infl\widehat \alpha)$ are
the object maps for functors from $\mathcal{A}(G)$ to
$\mathcal{C}(G)$, so that it makes sense to say that they are
naturally equivalent.

When we assert that, for example, $(A,G,\alpha)\mapsto (A\times_\alpha
G,\widehat \alpha)$ is a functor, we are completely ignoring the
morphisms, and we cannot appreciate what naturality means until we deal
with them too: a natural equivalence $T$ between two functors
$F,G\:\mathcal{A}\to\mathcal{B}$ assigns to each object $A$ of
$\mathcal{A}$ an equivalence $T(A)\:F(A)\to G(A)$ (that is, an
invertible morphism $T(A)$ in the category $\mathcal{B}$) such that,
for each morphism $\phi\:A\to B$ in $\mathcal{A}$, the diagram
\begin{equation*}
\xymatrix{
{F(A)}
\ar[r]^{T(A)}
\ar[d]_{F(\phi)}
&{G(A)}
\ar[d]^{G(\phi)}
\\
{F(B)}
\ar[r]_{T(B)}
&{G(B)}
}
\end{equation*}
commutes in $\mathcal{B}$.  In our categories, the morphisms will be
based on Hilbert bimodules; in $\mathcal{A}(G)$, for example, a
morphism from $(A,G,\alpha)$ to $(B,G,\beta)$ will be given by a
Hilbert bimodule $_AX_B$ with a compatible action $\gamma$ of $G$. 
The composition of morphisms will be based on the balanced tensor
product of bimodules, so that a diagram
\begin{equation*}
\xymatrix{
{A}
\ar[r]^{X}
\ar[d]_{Y}
&{B}
\ar[d]^{W}
\\
{C}
\ar[r]_{Z}
&{D}
}
\end{equation*}
of Hilbert bimodules commutes if $Y\otimes_C Z\cong X\otimes_B W$ as
Hilbert $A-D$ bimodules; in $\mathcal{A}(G)$ or $\mathcal{C}(G)$ this
isomorphism has to be appropriately equivariant.  The equivalences in
these categories are the morphisms which are given by imprimitivity
bimodules, so to prove that two of our functors $F,G$ are naturally
equivalent amounts to finding imprimitivity bimodules
$_{F(A)}X(A)_{G(A)}$ such that
\[
X(A)\otimes_{G(A)}G(Y)\cong F(Y)\otimes_{F(B)} X(B)
\]
as Hilbert $F(A)-G(B)$ bimodules for each Hilbert bimodule
${}_AY_B$. The modules are the usual ones, but
many of the details needed to establish these isomorphisms and their
properties are new.

This paper, like any other in which coactions appear, involves some
gritty technical arguments.  We will therefore begin by outlining the
main new issues which we face in this program, and how we have dealt
with them. Those who are interested in seeing how the categorical ideas
impact when there are no coactions around are encouraged to read our
previous paper \cite{ekqr:green} first.  Indeed, this might help even
those who are already coaction-compliant!

\section*{Outline}

We begin in \chapref{hilbert-chap} with a detailed discussion of the
Hilbert bimodules on which our morphisms are based.  The axioms are
intrinsically asymmetric; to see why, note that a homomorphism
$\phi\:A\to B$ gives $B$ the structure of an $A$-module, but not the
other way round.  Our modules $_AX_B$ will be right Hilbert
$B$-modules with a left action of $A$ given by a nondegenerate
homomorphism $\kappa$ of $A$ into the $C^*$-algebra $\L(X_B)$ of
adjointable operators on $X$.  As in \cite{kqr:resind}, we shall call
these \emph{right-Hilbert bimodules} to emphasize that the
Hilbert-module structure is on the right; we have stuck with this name
because the alternatives ($C^*$-correspondences or Hilbert bimodules)
do not carry the same sense of direction.  The theory of right-Hilbert
bimodules is similar to that of imprimitivity bimodules, but there
seem to be enough subtle differences to warrant a detailed discussion.

The first section contains the basic facts about multiplier bimodules and
homomorphisms between bimodules.  These are used repeatedly: a coaction
on a bimodule $X$, for example, is by definition a homomorphism of $X$
into the multiplier bimodule $M(X\otimes C^*(G))$.  Our treatment
is similar to that of imprimitivity bimodules in \cite{er:mult}.
\secref{hilbert-chap}.\ref{sec-tensor} is about the balanced tensor
products which are used to define the composition of morphisms; we
need to know in particular how this process extends to multipliers.
We also discuss external tensor products, which are crucial for the
definition of coactions on Hilbert bimodules.  The last section of
\chapref{hilbert-chap} is about linking algebras.  These are used
primarily as a technical tool in the proofs of naturality (an idea
lifted from \cite{er:stab}, and expounded in an easier setting in
\cite{ekqr:green}).

In \chapref{categories-chap} we describe the categories in which we
work.  The basic category $\mathcal{C}$ of $C^*$-algebras 
appears in \cite{ekqr:green};
we review the main facts in
Section~2.1.  The objects are $C^*$-algebras and the morphisms from
$A$ to $B$ are the isomorphism classes of right-Hilbert $A-B$
bimodules: we have to pass to isomorphism classes to ensure that the
composition law $[_CY_B]\circ[_AX_B]\deq[_A(X\otimes_C Y)_B]$ has the
required properties.  The other categories $\mathcal{A}(G)$,
$\mathcal{C}(G)$ and $\mathcal{AC}(G)$ are associated to a fixed
locally compact group $G$, and are obtained by adding, respectively,
actions of $G$, coactions of $G$, and both actions and coactions to the
objects and morphisms of $\mathcal{C}$.  
Adding actions is relatively
routine, but (as will be no surprise to those familiar with them)
adding coactions is a little harder.  
(Coactions on Hilbert bimodules first appeared in \cite{BS-CH}.)
We show that in each of these
categories, the equivalences (that is, the invertible morphisms) are
the morphisms in which the underlying bimodules are imprimitivity
bimodules, and then that every morphism is a composition of a morphism
coming from a nondegenerate homomorphism $\phi\:A\to M(C)$ and a
morphism based on an imprimitivity bimodule $_CX_B$.

In \chapref{functors-chap} we show that the various crossed products
appearing in our theorems define functors between appropriate
categories.  There are two main problems.  The first is to define
suitable crossed products.  We are interested here in coactions and
nonabelian duality, which is basically a theory about reduced crossed
products, so we have decided to give in gracefully and use reduced
crossed products throughout.  (This is definitely a choice: we have
already proved the naturality of Green's Imprimitivity Theorem for
full crossed products in \cite{ekqr:green}, and providing we were
willing to omit all statements about the coactions, we could
presumably do the same here.)  But because the objects in our
categories are $C^*$-algebras rather than isomorphism classes of
$C^*$-algebras, it is important that we don't just \emph{choose} a
regular representation willy-nilly.  So we shall discuss a specific
realization of the reduced product.  The second main problem is to
define crossed products of the Hilbert bimodules which define the
morphisms.
We do this differently for actions and coactions; for actions we make 
heavy use of the convenience of $C_c$-functions, and for coactions we 
realize the crossed product inside a certain multiplier bimodule.
For imprimitivity bimodules, it is handy to recognize that if $L(X)$ is
the linking algebra of $X$, then the 
bimodule crossed product $X\times G$ embeds as the top right corner 
of $L(X)\times G$,
and we have the important
relation $L(X)\times G=L(X\times G)$ almost by definition. We should
mention that defining these crossed products and establishing their
properties has been done before; see 
\cite{BS-CH}, \cite{bui:act}, \cite{BuiME}, \cite{er:mult},
and \cite{kas}.

We gather all the necessary functors in \chapref{functors-chap}; even
though some are easy, it is convenient to deal with them all at once. 
The key difficulty is the same in each case: it is not obvious that
crossed products preserve composition.  This amounts to proving things
like
\[
(X\otimes_B Y)\times G\cong(X\times G)\otimes_{B\times G}(Y\times G),
\]
and again our techniques are different for actions and coactions.

Our main theorems are in \chapref{natural-chap}.  We have already
discussed the first, which is about crossed products of induced
algebras, and which we prove in Section~4.1.  The proofs of this and
our other main theorems follow the same general pattern.  We factor
each morphism $_AX_B$ as a composition of a nondegenerate homomorphism
$\phi\:A\to M(C)$ and an imprimitivity bimodule ${}_CY_B$.  To prove
that
\[
\xymatrix{
{F(A)}
\ar[r]^{T(A)}
\ar[d]_{F(\phi)}
&{G(A)}
\ar[d]^{G(\phi)}
\\
{F(C)}
\ar[r]_{T(C)}
&{G(C)}
}
\]
commutes, we extend the homomorphisms $(F(\phi),G(\phi))$ to a
homomorphism of imprimitivity bimodules $T(A)\to M(T(C))$, and use a
general lemma which says this suffices.  To prove that
\[
\xymatrix{
{F(C)}
\ar[r]^{T(C)}
\ar[d]_{F(Y)}
&{G(C)}
\ar[d]^{G(Y)}
\\
{F(B)}
\ar[r]_{T(B)}
&{G(B)}
}
\]
commutes, we realize $T(C)$ and $T(D)$ as the diagonal corners in an
imprimitivity bimodule $Z$ over the linking algebras $L(F(Y))$ and
$L(G(Y))$, and use a general lemma from \cite{er:stab} which identifies
both $F(Y)\otimes_{F(B)}T(B)$ and $T(C)\otimes_{G(C)}G(Y)$ with the
top off-diagonal corner in $Z$.  The hard part in both halves is to
build the compatible coaction.

It may be known that this theorem about crossed products of induced
algebras is a generalization of Green's Imprimitivity Theorem, but it
does not appear to be well-documented.  We therefore give a careful
derivation,
which could be of some independent interest
(see the discussion preceding \thmref{thm-green-imp} in
\appxref{imprim-chap}).
We
then use this to deduce our second main theorem, which is a natural
and equivariant version of the Imprimitivity Theorem itself.  There
are many possible variations on this theme, depending on choices of
full and reduced crossed products and on whether or not the subgroup
is normal.  Here we have already decided to use reduced crossed
products, and we have further chosen to discuss what happens for
normal subgroups.  We have made this choice because in this case there
are several more actions and coactions in play, and the theorem has
something to say about all of them.  To see what is happening here,
recall that if $N$ is normal, we can view the imprimitivity algebra
$(A\otimes C_0(G/N))\times_{\alpha\otimes \tau} G$ in Green's theorem
as the crossed product $(A \times_\alpha
G)\times_{\widehat\alpha|}G/N$ by the restriction of the dual
coaction; thus this imprimitivity algebra carries a dual action
$(\widehat\alpha|)\,\widehat{}\;$ of $G/N$ as well as a dual coaction
$(\alpha\otimes\tau)\,\widehat{}\;$ of $G$.  Our theorem says that
Green's imprimitivity bimodule matches
$(\widehat\alpha|)\,\widehat{}\;$ with the so-called decomposition
action of $G$ on $A\times_\alpha N$ and
$(\alpha\otimes\tau)\,\widehat{}\;$ with the inflation to $G$ of the
dual coaction of $N$ on $A\times_\alpha N$.  This observation seems to
be new.  Indeed, we believe that the equivariance and the naturality
are both potentially important new pieces of information about Green's
theorem.

Our third main theorem is a version of Mansfield's Imprimitivity
Theorem.  This has all the same features as the version of Green's
theorem which we have just discussed: Mansfield's Morita equivalence
of $(A\times_\delta G)\times_{\widehat\delta}N$ with
$A\times_{\delta|}G/N$ is natural and equivariant for canonical
actions and coactions on the crossed products.  For this theorem, the
difficult part of the proof is establishing the naturality with
respect to ordinary homomorphisms $\phi\:A\to M(C)$; we have to work
hard to build compatible homomorphisms on Mansfield's bimodule.

In \chapref{apps-chap} we give some applications to our motivating
problem of understanding the relationships between induction and
duality.  In Section 5.1, we uncover some new and very intriguing
relationships between Green and Mansfield induction.  Important
special cases of these results say that the
Green bimodules $X_{\{e\}}^G(A)$ and Mansfield bimodules
$Y_{G/G}^G(A)$ are in duality:
\[
X_{\{e\}}^G(A)\times G\cong Y_{G/G}^G(A\times G)
\]
and
\[
X_{\{e\}}^G(A\times G)\cong Y_{G/G}^G(A)\times G.
\]
Results of this type require several applications of our main
theorems, and it is vital that we know everything is appropriately
equivariant.  
Our main new application to induction-restriction duality 
is \thmref{act-res-ind}, which completes the program of
\cite{EchDI, kqr:resind, ekr} by handling the restriction of
representations from $A\times_\alpha G$ to $A\times_{\alpha|}N$. 
We close with a new application of linking-algebra techniques to the
Symmetric Imprimitivity Theorem of \cite{RaeIC}.

Since this project is intrinsically involved with nonabelian duality,
we have necessarily made heavy use of coactions and their crossed
products.  There are several different sets of definitions available:
the subject is stabilizing, but some key questions of a fundamental
nature remain unresolved, and hence this is taking longer than one
might have wished.  So we have included as an appendix a survey of the
area, which outlines what we believe to be the most satisfactory
approach and describes how this approach relates to the others in the
literature.

A second appendix collects the precise versions of the imprimitivity
theorems we need; various formulations appear in the literature, so we
felt it would be handy to record exactly what we want.

Finally, the third appendix contains some technical results on
function spaces with values in locally convex spaces which are used
throughout the text to construct multipliers of bimodules.  In
applications, the locally convex spaces will be multiplier algebras or
bimodules with the strict topology: we need to know, for example, that
strictly continuous functions of compact support from $G$ to $M(X)$
define multipliers of $X\times G$, and that they do so in an orderly
fashion.

\section*{Epilogue}

Although this paper has turned out much longer than we intended, we
have made all sorts of simplifying assumptions to keep the length
down, and these are probably logically unnecessary.  First of all, we
have deliberately excised twisted crossed products, though some
residual traces remain in the presence of the decomposition actions
and coactions.  Any serious application of these ideas to the Mackey
machine --- which was, after all, our original motivation \cite{EchDI}
--- will require that we can handle twisted crossed products.  Second,
we have used reduced crossed products throughout.  For our present
applications involving nonabelian duality and crossed products by
coactions, this makes sense: the current duality theorems all factor
through the reduced crossed product.  But for applications to ordinary
crossed products this is not necessarily desirable, and there are
surely versions of \thmref{ind-act natural} and \thmref{G-act natural}
for full crossed products.  We have already described a version of
Green's theorem in \cite{ekqr:green}, but we neglected questions of
equivariance there.  Third, we have considered only some of the
important Morita equivalences.  The others, such as the Symmetric
Imprimitivity Theorem and the Stabilization Trick for twisted crossed
products, should be natural too.
(Working in the context of general locally compact quantum groups was
not even an issue, since there are currently no imprimitivity theorems
available in that generality!)

On the other hand, we have taken the liberty of treating actions
separately from coactions --- rather than viewing actions of $G$ 
as coactions of $C_0(G)$ --- although this would have led to a much
shorter exposition.  Our main reasons for this are that we think that
actions are much easier to understand than coactions, and that we feel
there may be more general interest in the action case than in the
coaction case. 

We hope that we have given convincing evidence that issues involving
functoriality, naturality and equivariance are likely to occur
frequently in our subject, and that it will pay for us get in the
habit of dealing with them as we go.  We also hope that we have made a
few other points along the way: our view of induced $C^*$-algebras as
an obstruction to imprimitivity, our heavy use of linking-algebra
techniques to identify imprimitivity bimodules, and the seemingly deep
and strange relations between induction and duality, should all have
applications elsewhere.
For instance, the strong connection between these ideas and equivariant
$KK$-theory is well-documented in \cite{BS-CH} and \cite{kas}. Several
of the ideas are also present in the approach of \cite{CE-TE} and
\cite{CE-PP} towards a Mackey machine for the Baum-Connes conjecture, 
and are applied to the Connes-Kasparov conjecture in \cite{CEN-CK}.

%
%

\chapter{Right-Hilbert Bimodules}
\label{hilbert-chap}

In this chapter we gather together the basic theory of right-Hilbert
bimodules.
We start with the basic definitions
and some important notation which shall be used
throughout this work.

\section{Right-Hilbert bimodules and partial imprimitivity bimodules}

Let $B$ be a
\cstar algebra. Recall that a \emph{Hilbert $B$-module} is a vector
space $X$ which is a right $B$-module equipped with a positive
definite $B$-valued sesquilinear form $\<\cdot,\cdot\>_B$
satisfying
\begin{equation}\label{Hb-mod-eq}
\<x,y \d b\>_B = \<x,y\>_Bb \midtext{and}
\<x,y\>_B^* = \<y,x\>_B
\righttext{for all}x,y\in X,b\in B,
\end{equation}
and which is complete in the norm $\|x\| =
\|\<x,x\>_B\|^{1/2}$.
Our primary reference for Hilbert modules
is \cite{RW}, and a secondary reference is \cite{lan:hilbert}.
Some notational conventions: 
we often omit the dot ($\cdot$) when writing module actions;
and in general, if $(u,v) \mapsto uv \: U \times V \to W$ is a
pairing among
vector spaces, then for $P \subseteq U$ and $Q \subseteq
V$ we write $PQ$ to mean the \emph{linear span} of the set 
$\{uv \mid u \in P, v \in Q\}$.

\begin{defn}\label{rHb-defn}
Let $A$ and $B$ be $C^*$-algebras.
A \emph{right-Hilbert $A-B$ bimodule} is a
Hilbert $B$-module $X$
which is also a nondegenerate left $A$-module (\ie, $AX=X$) satisfying
\begin{subequations}
\begin{align}
\label{rHb-eq1}
a \d (x \d b) &= (a \d x) \d b
\righttext{and}
\\
\label{rHb-eq2}
\<a \d x,y\>_B &= \<x,a^* \d y\>_B
\end{align}
\end{subequations}
for all $a \in A$, $x,y \in X$, and $b \in B$.  We write ${}_AX_B$ to
indicate all the data, and we call $X$ \emph{full} if it is full as a
Hilbert $B$-module, \ie, $\overline{\<X,X\>}_B = B$.
In general, if $X$ is not full,
we shall
write  $B_X$ for the closed ideal $\overline{\lk X,X\rk}_B\subseteq
B$, and we call $B_X$ the {\em range} of the inner product
on $X$.
\end{defn}

\begin{rem}
(1) 
In recent years, objects very similar to right-Hilbert bimodules
have been introduced into the literature: for example, the
$A-B$ correspondences of \cite{ms}.
In many cases
(as in \cite{ms}), the left module action is permitted to be
degenerate; we require it
to be nondegenerate so that we can extend it to the
multiplier algebra $M(A)$ (see below).

(2) Note that if $X$ is an
$A-B$ correspondence, then $\overline{AX}$ is a closed $A-B$ sub-bimodule
of $X$, and therefore becomes a right-Hilbert $A-B$ bimodule.
In fact we have $\overline{AX}=AX=\{ax\mid a\in A, x\in X\}$, since
it follows from Cohen's factorization theorem
that $\overline{AX}=A\overline{AX}\subseteq AX$
(we refer to  \cite[Proposition
2.33]{RW} for a statement and an easy proof of Cohen's factorization
theorem in the case where $A$ is a $C^*$-algebra).
More generally, a similar application of Cohen's theorem implies
that for any $C^*$-subalgebras $C$ and $D$ of $A$ and $B$,
respectively, we have $\overline{CX}=\{cx\mid c\in C, x\in X\}$
and $\overline{XD}=\{xd\mid x\in X, d\in D\}$.
\end{rem}

\begin{ex}\label{ex-standard}
If $B$ is a $C^*$-algebra, then $B$
becomes a full
right-Hilbert $B-B$ bimodule in a natural way by putting
\[
a \cdot b \cdot c = abc \midtext{and} \<a,b\>_B = a^*b
\righttext{for}a,b,c\in B.
\]
If $\phi \: A \to M(B)$ is a
nondegenerate $C^*$-algebra homomorphism, then $B$ becomes a
full right-Hilbert $A-B$ bimodule with left action given by
\[
a \cdot b = \phi(a)b.
\]
More generally, if
$\phi\:A\to M(B)$ is an arbitrary (possibly degenerate)
$*$-homomorphism, then $B$ becomes an $A-B$ correspondence
and, therefore, $X=\varphi(A)B$ is a right-Hilbert $A-B$ bimodule.
We call a right-Hilbert bimodule ${}_AX_B$ arising in this way
\emph{standard}. If $\varphi\:A\to M(B)$ is nondegenerate,
\ie, if $X=\varphi(A)B=B$, then we say that
$_AB_B$ is a {\em nondegenerate} standard right-Hilbert
  bimodule.
\end{ex}

\begin{rem}\label{rem-full}
    It is clear that a nondegenerate standard right-Hilbert bimodule
    $_AB_B$ is full. The converse is not true
    in general. To see an example let $B=M_2(\CC)$ and let
    $\varphi\:\CC\to M_2(\CC); \varphi(\lambda)=
    \left(\begin{smallmatrix} \lambda& 0\\ 0&0\end{smallmatrix}\right)$.
    Then $M_2(\CC)\varphi(\CC)M_2(\CC)=M_2(\CC)$ and
    $X=\varphi(\CC)M_2(\CC)\cong \CC^2$ is a full right-Hilbert
    $\CC-M_2(\CC)$ bimodule, but $\phi$ is degenerate.
\end{rem}

If $X$ and $Y$ are Hilbert $B$-modules, $\c L_B(X,Y)$ denotes the
set of maps $T\:X\to Y$
which are adjointable in the sense that there
exists $T^*\:Y\to X$ such that
\[
\<Tx,y\>_B = \<x,T^*y\>_B
\righttext{for all}x\in X,y\in Y.
\]
Such $T$ are automatically bounded and $B$-linear \cite[Lemma
2.18]{RW}.
The notation is shortened to $\c L(X,Y)$ if $B$ is understood, and $\c
L_B(X)$ (or just $\c L(X)$) if $X=Y$. 
In the latter case
$\c L(X)$ is a $C^*$-algebra with the operator norm $\|T\|=
\sup\{\|Tx\|\mid\|x\|\le 1\}$ \cite[Proposition~2.1]{RW}.

Now, if ${}_AX_B$ is a right-Hilbert $A-B$ bimodule then for each $a
\in A$ the map $x \mapsto a \d x$ is adjointable (consequently the
associativity condition \eqref{rHb-eq1} is redundant), so we get a
homomorphism $\kappa \: A \to \c L_B(X)$
such that
\[
\kappa(a)x = a \d x,
\]
and which is nondegenerate in the sense that $\kappa(A)X=X$.
Conversely, every right-Hilbert bimodule arises in this way: If $X$ is
a Hilbert $B$-module and $\kappa \: A \to \c L(X)$ is a nondegenerate
homomorphism, then $X$ becomes a right-Hilbert $A-B$ bimodule via
\[
a \d x = \kappa(a)x.
\]
Thus a right-Hilbert $A-B$ bimodule is nothing more nor less than a Hilbert
$B$-module $X$ together with a nondegenerate homomorphism $A\to \c
L(X)$.

If $X$ and $Y$ are Hilbert $B$-modules, $\c K(X,Y)$ denotes the 
\emph{compact operators} from $X$ to $Y$: by definition, it is 
the closed
span in $\c L(X,Y)$ of the maps $z\mapsto y\<x,z\>_B$ for $x\in X$
and $y\in Y$. $\c K(X) = \c K(X,X)$ is a closed ideal in $\c L(X)$, and in fact
$\c L(X)\cong M(\c K(X))$ \cite[Corollary~2.54]{RW}.
In particular, if $X$ is a Hilbert $B$-module, then the formula
\begin{equation}\label{eq-KXinnerproduct}
     _{\K(X)}\lk x,y\rk z\deq x\lk y,z\rk_B
     \end{equation}
defines a full $\K(X)$-valued inner product on $X$, which
gives $X$ the structure of a {\em left} Hilbert $\K(X)$-module.
Then $B$ acts via adjointable operators on the right of
$_{\K(X)}X$, and $B_X=\overline{\lk X,X\rk}_B$ identifies with the
compact operators of the left Hilbert $\K(X)$-module
$_{\K(X)}X$.

We pause to clear up an apparent ambiguity: Suppose ${}_AX_B$ is a
right-Hilbert bimodule, and let $\kappa\:A\to\c L_B(X)$ be the
associated homomorphism.  By definition $\kappa$ is nondegenerate in the
sense that $\kappa(A)X=X$.  Is $\kappa$ still nondegenerate when viewed as
a homomorphism of $A$ into the multiplier algebra $M(\c K(X))$?
Yes, because $\kappa(A)\c K(X)$ contains all the ``rank-one'' maps ${}_{\c
K(X)}\<x,y\>$: just factor $x=\kappa(a)z$ for some $a\in A$ and $z\in X$,
and then
\[
{}_{\c K(X)}\<x,y\>
= {}_{\c K(X)}\<\kappa(a)z,y\>
= \kappa(a){}_{\c K(X)}\<z,y\>.
\]
Thus $\kappa$ extends uniquely to a strictly continuous homomorphism
$\bar\kappa\:M(A)\to M(\c K(X))$.  Since the strict topology on $M(\c
K(X)) = \c L(X)$ is stronger than the strong operator topology
\cite[Proposition~C.7]{RW}, $\bar\kappa$ is also strict-strong operator
continuous into $\c L(X)$.  In particular, $X$ is a left
$M(A)$-module, in fact a right-Hilbert $M(A)-M(B)$ bimodule.

In this work, imprimitivity bimodules will play a very important
r\^ole, since they will represent the
equivalences in our categories (see \chapref{categories-chap} below).
But we shall also need the
more general notion of partial imprimitivity
bimodules:

\begin{defn}\label{defn-partial}
     Suppose that $A$ and $B$ are $C^*$-algebras.
     A {\em partial $A-B$ imprimitivity bimodule} is
     a complex vector space $X$ which is
     a right Hilbert $B$-module and a left Hilbert $A$-module
     such that
     $$ a\d(x\d b)=(a\d x)\d b
\quad\text{and}\quad {_A\lk x,y\rk}\d z=x\d\lk y,z\rk_B$$
     for all $a\in A$, $b\in B$, and $x,y,z\in X$.
     If ${}_A\overline{\lk X, X\rk} =A$, we
     say that $X$ is a {\em right-partial
     imprimitivity bimodule}. 
If $\overline{\lk X, X\rk}_B=B$, $X$ is a \emph{left-partial
imprimitivity bimodule}.  
If both $_A\lk\cdot, \cdot\rk$
     and $\lk\cdot,\cdot\rk_B$ are full, then
     $X$ is called an \emph{$A-B$ imprimitivity bimodule}.
     $A$ and $B$ are called {\em Morita equivalent} if there exists
     at least one $A-B$ imprimitivity bimodule.
\end{defn}

Imprimitivity bimodules for $C^*$-algebras were 
introduced by Rieffel in \cite{RieIR}. As a general reference
we recommend \cite{RW}.
Our partial imprimitivity bimodules are often
called Hilbert bimodules in the literature.
However, this terminology would be much too close to
``right-Hilbert bimodules'', so we decided to introduce new terminology
to prevent confusion.

The reason for choosing the notion of partial imprimitivity
bimodules comes from the trivial fact that if
$X$ is a partial $A-B$ imprimitivity bimodule, then
$X$ becomes an $A_X-B_X$ imprimitivity bimodule for
the ideals $A_X={}_A\overline{\lk X,X\rk}$ and 
$B_X=\overline{\lk X,X\rk}_B$
of $A$ and $B$, respectively.
In the special case where $A=C_0(V)$ and $B=C_0(W)$ are commutative,
this just means that we have a homeomorphism between the open sets
$V_X$ and $W_X$ corresponding to the ideals $A_X$ and $B_X$,
\ie, a partial homeomorphism between $V$ and $W$
(see \cite[Corollary 3.33]{RW}).

\begin{ex}\label{ex-partialimp}
     (1) Suppose that $X$ is a right-Hilbert $A-B$ bimodule.
     Let $_{\K(X)}\lk\cdot, \cdot\rk$ be the left $\K(X)$-valued
     inner product on $X$ as defined in \eqref{eq-KXinnerproduct}.
     Then $X$ becomes a right-partial $\K(X)-B$ imprimitivity
     bimodule.
     In particular, if $X$ is a full Hilbert $B$-module, then
     $X$ is a $\K(X)-B$ imprimitivity bimodule.

     (2) More generally, suppose that $X$ is a right-Hilbert $A-B$ bimodule
     such that the corresponding homomorphism $\kappa\:A\to \L_B(X)$
     is injective and  $\kappa(A)\supseteq \K(X)$.
     Define a left $A$-valued inner product on $X$ by
     $$_A\lk x,y\rk =\kappa^{-1}({_{\K(X)}\lk x,y\rk})$$
     for all $x,y \in X$. Then $X$ becomes a partial $A-B$
     imprimitivity bimodule. In particular, if $X$ is a full
     right-Hilbert $A-B$ bimodule, then $X$ is an $A-B$ imprimitivity
     bimodule if and only if $\kappa\:A\to \L(X)$ is injective
     and $\kappa(A)=\K(X)$.

     (3) Suppose that $X$ is a partial $A-B$ imprimitivity bimodule.
     Let $\widetilde{X}=\{\tilde{x}\mid x\in X\}$ with left $B$-action and
     right $A$-action defined by
     $$b\d\tilde{x}\deq\widetilde{x\d b^*}
\quad\text{and}\quad
     \tilde{x}\d a\deq\widetilde{a^*\d x}$$
     and $B$- and $A$-valued inner products given by
     $$_B\lk \tilde{x}, \tilde{y}\rk=\lk x,y\rk_B
     \quad\text{and}
     \lk \tilde{x},\tilde{y}\rk_A={_A\lk x,y\rk},$$
     for $x,y\in X$.
     Then $\tilde{X}$ is a partial $B-A$ imprimitivity bimodule,
     called the {\em conjugate} (also called the \emph{reverse}, 
or \emph{dual} in the literature) of $X$.
     If $X$ is a right-partial imprimitivity bimodule, then
     $\widetilde{X}$ is a left-partial imprimitivity bimodule and \emph{vice
     versa}.
     If $X$ is an $A-B$ imprimitivity bimodule, then
     $\widetilde{X}$ is a $B-A$ imprimitivity bimodule.

     (4) If $X$ is a partial $A-B$ imprimitivity bimodule, then
     via the extension of the actions to the multiplier algebras,
     $X$ becomes a partial $M(A)-M(B)$ imprimitivity bimodule.
\end{ex}

Another important example is given by the multiplier bimodule $M(X)$
of a right-partial imprimitivity bimodule $X$, which we shall define
in \secref{sec-multbimod} below.
We now recall the notion of induced ideals. If not explicitly
stated otherwise, by an ideal of an algebra we always mean
a two-sided ideal.

\begin{defn}\label{def-induced}
     Suppose that $X$ is a right-Hilbert $A-B$ bimodule
     (or just an $A-B$ correspondence). Let $I$ be a closed
     ideal in $B$. Then call $XI$ the \emph{closed
     $A-B$ submodule of $X$ corresponding to $I$}.
     The closed ideal
     $X$-$\ind I\deq\{a\in A\mid aX\subseteq XI\}$
     is called the {\em ideal of $A$ induced from $I$}.
\end{defn}

If $X$ is actually an $A-B$ imprimitivity bimodule then we
have the following basic result. For a proof we refer to
\cite[Theorem 3.22, Lemma 3.23 and Proposition 3.24 ]{RW}.

\begin{prop}[{Rieffel Correspondence}]\label{prop-rieffel}
     Suppose that $X$ is an $A-B$ imprimitivity bimodule. Then
     there exist inclusion-preserving one-to-one correspondences
     between
     \begin{enumerate}
	\item the closed ideals of $B$,
	\item the closed $A-B$ submodules of $X$, and
	\item the closed ideals of $A$.
     \end{enumerate}
    If $I$ is a closed ideal of $B$, then
     $XI$ is the corresponding closed $A-B$ submodule of
     $X$ and $X$-$\ind I$ is the corresponding closed
     ideal of $A$. Moreover, we have
     the identities $X$-$\ind I={}_A\overline{\lk XI,XI\rk}$,
     $(X$-$\ind I)X=XI$, and $\overline{\lk XI, XI\rk}_B=I$.
\end{prop}

\begin{rem}\label{rem-quotientimp}
     If $X$ is an $A-B$ imprimitivity bimodule, then
     it follows from the Rieffel correspondence
     that if $I$ is a closed ideal of $B$, then
     $XI$ becomes a $X$-$\ind I-I$ imprimitivity bimodule
     via restricting all actions and inner products to
     the subspaces. With a little bit more work one can
     also show that $X/XI$ can be made into
     an $A/(X$-$\ind I)-B/I$ imprimitivity bimodule
     with respect to the obvious actions and inner products
     (see \cite[Proposition 3.25]{RW} for more details).
\end{rem}

\section{Multiplier bimodules and homomorphisms}\label{sec-multbimod}

A theory of multipliers of imprimitivity
bimodules was developed in \cite{er:mult}.
We cannot adopt that theory directly, as our
right-Hilbert bimodules are decidedly left-challenged; we need something
more like the one-sided theory of \cite{BS-CH} and \cite{ng:module}.
Our gadgets and objectives are not exactly the same as in these
references; for this reason and for completeness we give a self-contained
treatment.

If $X$ is a Hilbert $B$-module, then $\c L(X)$ is a $C^*$-algebra,
hence is a right-Hilbert $\c L(X) - \c L(X)$ bimodule.  This extends
to adjointable maps between different modules:

\begin{prop}\label{LXY-mod-prop}
Let $X$ and $Y$ be Hilbert $B$-modules. Then $\L(X,Y)$ is a
right-Hilbert $\L(Y)-\L(X)$
bimodule with operations
\[
  P \d T =  P \circ T,
\quad
T \d  Q = T \circ  Q,
\midtext{and}
\<T,S\>_{\L(X)} = T^*\circ S
\]
for $ P \in \L(Y)$, $T,S \in \L(X,Y)$, and $ Q \in \L(X)$.
Moreover, $\L(X,Y)$ becomes a partial $\L(Y)-\L(X)$
imprimitivity bimodule if we define the $\L(Y)$-valued inner
product on $\L(X,Y)$ by
${_{\L(Y)}\lk T,S\rk}=T\circ S^*$.

\end{prop}

\begin{proof} Of course $PT\in \L(X,Y)$, with $(PT)^*
= T^*P^*$, and similarly for the other compositions.
Composition is associative,
so $\L(X,Y)$ is an $\L(Y)-\L(X)$ bimodule, and easy calculations
verify the sesquilinearity of $\<\cdot,\cdot\>_{\L(X)}$ and the
consistency relations \eqref{Hb-mod-eq} and \eqref{rHb-eq2} (recall
that \eqref{rHb-eq1} is redundant).

For positivity, note that for each $x \in X$,
\[
\<x,T^*Tx\>_B = \<Tx,Tx\>_B \ge 0,
\]
so that $T^*T \ge 0$ in $\L(X)$ by
\cite[Lemma~2.28]{RW}. The above
equality also implies that $\|Tx\|^2 = 0$ for all $x \in X$
whenever $T^*T = 0$, which gives definiteness of the inner
product.

Since $\L(Y)$ is unital, it certainly acts nondegenerately on
$\L(X,Y)$, so it
only remains to check that $\L(X,Y)$ is complete in the norm
\[
\|T\| = \|\<T,T\>_{\L(X)}\|^{1/2} =
\|T^*T\|^{1/2}_{\L(X)}.
\]
First we claim that this norm coincides with the operator norm on the
Banach space $\B(X,Y)$ of all bounded linear maps of $X$ into
$Y$. Since $T^*T \ge 0$ in $\L(X)$, we can find $ Q \in
\L(X)$ with $ Q^* Q = T^*T$. Then
\begin{align*}
\|T\|^2_{\B(X,Y)}
&= \sup_{\|x\| \le 1}\{\|Tx\|^2_Y\}
= \sup_{\|x\| \le 1}\{\|\<Tx,Tx\>_B\|\}
= \sup_{\|x\| \le 1}\{\|\<x,T^*Tx\>_B\|\}
\\&= \sup_{\|x\| \le 1}\{\|\<x, Q^* Qx\>_B\|\}
= \sup_{\|x\| \le 1}\{\| Qx\|^2_X\}
= \| Q\|^2_{\L(X)}
\\&= \| Q^* Q\|_{\L(X)}
= \|T^*T\|_{\L(X)}
= \|T\|^2.
\end{align*}
Next, the adjoint map
$T\mapsto T^*$ is isometric from $\L(X,Y)$ into $\L(Y,X)$:
\begin{align*}
\|T^*\|
&= \|T^*\|_{\B(Y,X)}
= \sup_{\|y\| \le 1}\{\|T^*y\|_X\}
= \sup_{\|x\|,\|y\| \le 1}\{\|\<x,T^*y\>_B\|\}
\\&= \sup_{\|x\|,\|y\| \le 1}\{\|\<Tx,y\>_B\|\}
= \sup_{\|x\| \le 1}\{\|Tx\|_Y\}
= \|T\|_{\B(X,Y)}
= \|T\|.
\end{align*}
Now a standard argument shows that $\L(X,Y)$ is closed in
$\B(X,Y)$, and therefore complete: suppose that $\{T_i\}$ is Cauchy
in $\L(X,Y)$, so $T_i \to T$ in $\B(X,Y)$ and $T_i^* \to S$ in
$\B(Y,X)$. Then $S$ is an adjoint for $T$:
\begin{equation*}
\<Tx,y\>_B
= \lim_i \<T_ix,y\>_B
= \lim_i \<x,T_i^*y\>_B
= \<x,Sy\>_B.
\end{equation*}
The last assertion now follows from the complete symmetry of the
situation.
\end{proof}

\begin{rem}\label{KXY-ib-rem}
Note that the $\L(Y)$- and $\L(X)$-valued inner products on $\L(X,Y)$
are in general not full, so $\L(X,Y)$ is in general not an
$\L(Y)-\L(X)$ imprimitivity bimodule. To see an example,
let $X=\H$ be an infinite dimensional Hilbert space, viewed as
a Hilbert $\CC$-module, and let $Y=\CC$. Then one can check that
$\L(\H, \CC)\cong \H$
and $_{\L(\H)}\lk \H,\H\rk=\K(\H)$. But $\L(\H)$ coincides with
the much larger algebra $\B(\H)$.
\end{rem}

\begin{rem}\label{LB=MB-rem}
In what follows, we are primarily interested in $\L_B(B,X)$, where $X$
is a Hilbert $B$-module, and $B$ is viewed as a Hilbert $B$-module
in the usual way. In this case $B \cong \c K_B(B)$ and (hence) $M(B)
\cong \c L_B(B)$.
\end{rem}

\begin{cor}\label{LBX-mod-cor}
Let $X$ be a right-Hilbert $A-B$ bimodule. Then
$\L_B(B,X)$ is a right-Hilbert $M(A)-M(B)$ bimodule
with operations
\begin{equation}
\label{LBX-mod-eq}
(n\d T)b = n\d (Tb),
\quad
(T\d m)b = T(mb),
\midtext{and}
\<T,S\>_{M(B)} = T^*S
\end{equation}
for $n\in M(A)$, $T,S\in \L_B(B,X)$, $b\in B$, and $m\in M(B)$.
\end{cor}

\begin{proof}
By \propref{LXY-mod-prop} and \remref{LB=MB-rem}, $\L(B,X)$ is a
right-Hilbert $\L(X)-M(B)$ bimodule, so we only need to produce the
compatible $M(A)$-action.  But this is easy: letting $\kappa\:A\to \c
L(X)$ be the associated nondegenerate homomorphism, we have the
canonical extension $\bar\kappa\:M(A)\to \c L(X)$, so $M(A)$ acts on
the left of $\c L(B,X)$ by $m\d T = \bar\kappa(m)T$.
Unraveling the identifications, one sees that this left action is
given by the first formula in \eqref{LBX-mod-eq}.
\end{proof}

\begin{defn}\label{MX-defn}
We call $\L_B(B,X)$ the \emph{multiplier bimodule} of the right-Hilbert
$A-B$ bimodule $X$, and we denote it by $M(X)$.
\end{defn}

\begin{rem}\label{rem-imp-multiplier}
     (1) If $X$ is a right-partial $A-B$ imprimitivity bimodule,
     then the left action of $M(A)$ on $X$ identifies
     $M(A)$ with $M(\K(X))\cong \L_B(X)$. Thus it follows
     from \propref{LXY-mod-prop} that $M(X)$ carries a
     left $M(A)$-valued inner product which makes
     $M(X)$ into a partial $M(A)-M(B)$ imprimitivity bimodule.

     (2) It follows from \cite[Proposition 1.3]{er:mult} that
     if $_AX_B$ is an $A-B$ imprimitivity bimodule, then
     the above-defined multiplier bimodule $M(X)$ is
     canonically isomorphic
     to the multiplier bimodule defined in \cite[Definition
     1.1]{er:mult}.

\end{rem}

\begin{defn}\label{rHb-hom-defn}
Let ${}_AX_B$ and ${}_CY_D$ be right-Hilbert bimodules, let $\phi \: A
\to M(C)$ and $\psi \: B \to M(D)$ be homomorphisms, and let $\Phi \:
X \to M(Y)$ be a linear map. We call $\Phi$ a \emph{$\phi-\psi$
compatible right-Hilbert bimodule homomorphism} if
\begin{enumerate}
\item $\Phi(a \d x) = \phi(a) \d \Phi(x)$,
\item $\Phi(x \d b) = \Phi(x) \d \psi(b)$, and
\item $\psi(\<x,z\>_B) = \<\Phi(x),\Phi(z)\>_{M(D)}$
for all $a \in A$, $x,z \in X$ and $b \in B$.
\end{enumerate}

We call $\phi$ and $\psi$
the \emph{coefficient maps} of $\Phi$; we write
${}_\phi\Phi_\psi \: {}_AX_B \to M({}_CY_D)$ to indicate all the data.

We say that $\Phi$ is \emph{nondegenerate} if $\overline{\Phi(X)
D}=Y$ (\ie, the closed linear span of $\{\Phi(x)d\mid x\in X, d\in D\}$
equals $Y$)
and both $\phi$ and $\psi$ are nondegenerate.

If $\phi$ and $\psi$ are isomorphisms of $A$ and $B$ onto $C$ and
$D$, respectively, and $\Phi$ is a bijection onto $Y$, we call $\Phi$
a \emph{right-Hilbert bimodule isomorphism} of $X$ onto $Y$.

If $A=C$, $B=D$, and both $\phi$ and $\psi$ are identity maps, we
call $\Phi$ a \emph{right-Hilbert $A-B$ bimodule homomorphism}. If in
addition the map $\Phi$ is bijective we call it a
\emph{right-Hilbert $A-B$ bimodule isomorphism}.
\end{defn}

\begin{rem}\label{rHb-hom-rem}
(1)
Note that a right-Hilbert bimodule homomorphism $\Phi\:X\to M(Y)$
is always norm-decreasing, since (with the above notation)
for each $x\in X$ we have
\[
\|\Phi(x)\|^2 = \|\<\Phi(x),\Phi(x)\>_{M(D)}\| =
\|\psi(\<x,x\>_B)\| \leq \|\<x,x\>_B\| = \|x\|^2.
\]
Of course, this also shows that a 
right-Hilbert bimodule isomorphism
is necessarily isometric.

(2)
Condition (ii) is actually implied by condition (iii) and
the linearity of $\psi$:
\begin{align*}
&\|\Phi(xb)-\Phi(x)\psi(b)\|^2
\\&\quad
= \bigl\|\<\Phi(xb),\Phi(xb)\>_{M(D)}
- \<\Phi(xb),\Phi(x)\psi(b)\>_{M(D)}
\\&\quad
\hphantom{=\biggl\|}\quad
+ \<\Phi(x)\psi(b),\Phi(x)\psi(b)\>_{M(D)}
-\<\Phi(x)\psi(b),\Phi(xb)\>_{M(D)}\bigr\|
\\&\quad=
\bigl\|\<\Phi(xb),\Phi(xb)\>_{M(D)}
- \<\Phi(xb),\Phi(x)\>_{M(D)}\psi(b)
\\&\quad
\hphantom{=\biggl\|}\quad
+ \psi(b)^*\<\Phi(x),\Phi(x)\>_{M(D)}\psi(b)
- \psi(b)^*\<\Phi(x),\Phi(xb)\>_{M(D)}\|
\\&\quad=
\Bigl\|\psi\Bigl(\<xb,xb\>_{M(B)}
- \<xb,x\>_{M(B)}b + b^*\<x,x\>_{M(B)}b
- b^*\<x,xb\>_{M(B)}\Bigr)\Bigr\|
\\&\quad
= 0.
\end{align*}
In fact, a similar computation shows that if $\Phi\:X\to M(Y)$ is only a
\emph{map} satisfying (i) and (iii) then $\Phi$ is automatically linear,
hence is a right-Hilbert bimodule homomorphism (see, for example,
\cite[Lemma 2.5]{SieME}).
\end{rem}

If $_AX_B$ and $_CY_D$ are partial imprimitivity bimodules, a
nondegenerate right-Hilbert  bimodule homomorphism automatically
preserves this extra structure (as well as possible).
Recall that in this situation
the action of $C$ on $Y$ determines a homomorphism
$\kappa_C\:M(C)\to \L_D(Y)$. If
$Y$ is a right-partial $C-D$ imprimitivity bimodule
we have $M(C)\cong \L_D(Y)$ via $\kappa_C$, and we used this
identification to define an $M(C)$-valued inner
product on $M(Y)$.

\begin{lem}\label{lem-partialhom}
     Suppose that $_AX_B$ is a partial imprimitivity bimodule,
     $_CY_D$ is a right-Hilbert bimodule
     and  $_\phi\Phi_\psi\:{X}\to{M(Y)}$ is a right-Hilbert
     bimodule homomorphism such that $\overline{\Phi(X)D}=Y$. Then
     $$_{\L(Y)}\lk \Phi(x), \Phi(z)\rk=\kappa_C\circ\phi({_A\lk
     x,z\rk}),$$
     for all $x,z\in X$.
     If, in addition, $_CY_D$ is a right-partial
     imprimitivity bimodule, then $_{\phi}\Phi_{\psi}\:X\to M(Y)$
     is a {\em partial imprimitivity bimodule homomorphism}, that is,
     $\Phi$ is a right-Hilbert bimodule homomorphism with the
     additional property that
     $$_{M(C)}\lk \Phi(x),\Phi(z)\rk= \phi({_A\lk x,z\rk})$$
     for all $x,z\in X$.
\end{lem}
\begin{proof}
  The second assertion follows directly from the first assertion and
  the definition of $_{M(C)}\lk \cdot,\cdot\rk$. For the first
  assertion we have to check that
  $$\phi({_A\lk x,z\rk})y={_{\L(Y)}\lk\Phi(x),\Phi(z)\rk}y$$
  for all $y\in Y$.
  Since $Y=\overline{\Phi(X)D}$, we may assume without loss of generality that
  $y=\Phi(w)d$ for some $w\in X$, $d\in D$, and then
  \begin{align*}
  \phi({}_A\<x,z\>)\d\Phi(w)d
&= \Phi({}_A\<x,z\>\d w)d
= \Phi(x\d\<z,w\>_B)d\\
&= \Phi(x)\d\<\Phi(z),\Phi(w)\>_{M(D)}d
=\Phi(x)\circ\Phi(z)^*\circ \Phi(w)\d d\\
&= {}_{\L(Y)}\<\Phi(x),\Phi(z)\> \d \Phi(w) \d d.
\end{align*}
\end{proof}

\begin{rem}\label{rHb-hom-rem1}
(1) If $_AX_B$ and $_CY_D$ are imprimitivity bimodules, it follows
     from the above lemma that a nondegenerate
     right-Hilbert bimodule homomorphism $_\phi\Phi_\psi\:X\to M(Y)$
     is automatically a nondegenerate imprimitivity bimodule
     homomorphism in the sense of \cite[Definition~1.8]{er:mult}.

Curiously, the above statement has a converse: if $\Phi$ is a
nondegenerate imprimitivity bimodule homomorphism,
then $\Phi$ is automatically a
nondegenerate right-Hilbert bimodule homomorphism,
since then $\overline{\Phi(X)D=Y}$ by \cite[Lemma 5.1]{kqr:resind}.
However, we must be careful in applying this
remark: if $X$ and $Y$ are imprimitivity bimodules, but we only know
that $\Phi\:X\to M(Y)$ is a right-Hilbert bimodule homomorphism, then
nondegeneracy of the coefficient homomorphisms does \emph{not} imply
nondegeneracy of $\Phi$ --- for example, if $\c H$ is a Hilbert space
of dimension greater than one, and if $\xi\in\c H$, we get a
right-Hilbert bimodule homomorphism ${}_\phi\Phi_\psi\:{}_\b C\b C_\b
C\to{}_{\c K(\c H)}\c H_\b C$ by defining
\[
\Phi(c)=c\xi,
\quad
\phi(c)=c1_\c H,
\midtext{and}
\psi(c)=c.
\]
Both $\phi$ and $\psi$ are nondegenerate, but $\Phi$ is degenerate.
Of course, $\Phi$ is \emph{not} an imprimitivity bimodule
homomorphism. For nondegeneracy of $\Phi$ to follow from that of the
coefficient maps, we must know $\Phi$ preserves \emph{both} inner
products. As a reward for checking both inner products, we can avoid
checking the left module action, as a computation similar to that in
item (2) of \remref{rHb-hom-rem} shows.

(2)
Several times we will need the following adaptation to our context of
\cite[Example~1.10]{er:mult}: if ${}_\phi\Phi_\psi\:{}_AX_B\to
M({}_CY_D)$ is a nondegenerate right-Hilbert bimodule homomorphism,
then there exists a unique map $\mu\:\c K_B(X)\to M(\c K_D(Y))$ such
that ${}_\mu\Phi_\psi$ is a nondegenerate
right-Hilbert bimodule homomorphism from ${}_{\c K(X)}X_B$ to
$M({}_{\c K(Y)}Y_D)$. Note that it follows from \lemref{lem-partialhom}
above that ${}_\mu\Phi_\psi$ is then automatically a right-Hilbert
bimodule homomorphism, so that $\mu$ is determined by
the equation
$$\mu({}_{\K(X)}\lk x,y\rk)=
{}_{M(\K_D(Y))}\lk \Phi(x), \Phi(y)\rk
\righttext{for all} x,y\in X.$$
\end{rem}

It is  useful to note the following relation between
right-partial imprimitivity bimodule homomorphisms
and induced ideals.

\begin{lem}\label{lem-induced-hom}
     Let $_{\phi}\Phi_{\psi}\:{}_AX_B\to M(_CY_D)$ be a right-Hilbert
     bimodule homomorphism. Then $\ker\Phi=X\cdot\ker\psi$.
     In particular, $\Phi$ is injective if $\psi$ is.
    If, in addition, $_AX_B$ and $_CY_D$ are
    right-partial imprimitivity
    bimodules, we also have $\ker\phi=X$-$\ind(\ker\psi)$.
    Thus, if $_AX_B$ is an imprimitivity bimodule, then
     $\ker\phi, \ker\Phi$, and $\ker\psi$ correspond to each
     other via the Rieffel correspondence.
\end{lem}
\begin{proof} Since $\Phi$ is a bimodule map,
it follows that $X\d\ker\psi\subseteq \ker\Phi$.
Conversely,
if $Y\deq\ker\Phi$, then $Y$ is a closed $A-B$ submodule
of $X$ such that $\lk Y, Y\rk_B\subseteq \ker\psi$.
It follows that $Y=Y\d\overline{\lk Y,Y\rk}_B\subseteq X\d\ker\psi$.

Suppose now that $A= \c K(X)$ and $C=\K(Y)$. As above we first see that
$(\ker\phi)X\subseteq \ker\Phi$, so we only have to check
that $aX\subseteq \ker\Phi$ implies $a\in \ker\phi$ for all
$a\in \c K(X)$.
But, using \lemref{lem-partialhom}, $aX\subseteq \ker\Phi$ implies  that
$$\phi(a \c K(X))=\phi\big(a({}_{\c K(X)}\overline{\lk X,X\rk})\big)=
\phi({}_{\c K(X)}\overline{\lk aX,X\rk})=
{}_{M(\K(Y))}\overline{\lk \Phi(aX),\Phi(X)\rk}=\{0\},$$
and hence $\phi(a)=0$.
The final assertions now follow from
Proposition \ref{prop-rieffel}.
\end{proof}

\begin{rem}\label{rem-isoimp}
One useful application of the above lemma and the Rieffel
correspondence is the following result:
If ${}_\phi\Phi_\psi\:{}_AX_B\to {}_CY_D$ is a homomorphism between
imprimitivity bimodules, and if $\phi$ and $\psi$ are isomorphisms,
then so is $\Phi$. The injectivity of $\Phi$ follows from
the injectivity of $\psi$,
and the surjectivity follows
from the Rieffel correspondence and
the fact that $\Phi(X)$ is a closed
$C-D$ submodule of $Y$ such that $\overline{\lk \Phi(X),\Phi(X)\rk}_D=
\psi(B)=D$.
\end{rem}

In \propref{tnsr-hom-prop} below, and elsewhere, we need to construct
a right-Hilbert bimodule homomorphism by extending from a
pre-right-Hilbert bimodule.  We pause to make this precise, and
streamline the process with an elementary lemma.

\begin{defn}
Let $A_0$ and $B_0$ be dense $^*$-subalgebras of $C^*$-algebras $A$
and $B$, respectively, and let $X_0$ be an $A_0-B_0$ bimodule.  We say
$X_0$ is a \emph{pre-right-Hilbert $A_0-B_0$ bimodule} if it has a
$B_0$-valued pre-inner product (with positivity interpreted in $B$),
so that \eqref{Hb-mod-eq} holds, and also
\[
\<a \d x,a \d x\>_{B_0} \le \|a\|^2\<x,x\>_{B_0}
\righttext{for all} a \in A_0, x \in X_0.
\]
\end{defn}

As with pre-Hilbert modules (\cf\ \cite[Lemma
2.16]{RW}),
the completion $X$ of $X_0$ becomes a
right-Hilbert $A-B$ bimodule by taking limits of the operations.

\begin{lem}
\label{pre-rHb-lem}
Let $X_0$ be a pre-right-Hilbert $A_0-B_0$ bimodule, with completion
${}_AX_B$.  Suppose we have homomorphisms $\phi\:A\to M(C)$
and $\psi\:B\to M(D)$, and a map $\Phi_0\:X_0\to M(Y)$ such that:
\begin{enumerate}
\item $\Phi_0(a \d x) = \phi(a) \d \Phi_0(x)$ and
\item $\<\Phi_0(x),\Phi_0(y)\>_{M(D)} = \psi(\<x,y\>_{B_0})$
\end{enumerate}
for all $a \in A_0$ and $x,y \in X_0$.  Then $\Phi$ extends uniquely
to a right-Hilbert bimodule homomorphism $\Phi \: X \to M(Y)$.
\end{lem}

\begin{proof}
Similarly to \remref{rHb-hom-rem}, the inner product property (ii)
implies first of all that $\Phi_0$ is linear, then that it is bounded,
hence extends uniquely to a bounded linear map $\Phi \: X \to M(Y)$.
The compatibility conditions (i) and (iii) of \defnref{rHb-hom-defn}
follow from continuity and density.  (Recall from \remref{rHb-hom-rem}
that \defnref{rHb-hom-defn} (ii) is redundant.)
\end{proof}

Occasionally it is convenient to take advantage of certain
labor-saving devices when we are dealing with
partial imprimitivity bimodules.
In the language we have already introduced, we observe that a
pre-imprimitivity bimodule can be defined as a pre-right-Hilbert
bimodule ${}_{A_0}(X_0)_{B_0}$ which also has a left $A_0$-valued
pre-inner product (with positivity interpreted in the completion $A$)
such that for all $x,y,z\in X_0$ and $b\in B_0$ we have
\[
{}_{A_0}\<xb,xb\>\le\|b\|^2{}_{A_0}\<x,x\>
\midtext{and}
{}_{A_0}\<x,y\>z=x\<y,z\>_{B_0}.
\]

\begin{lem}
\label{pre-imp-lem}
Let $X_0$ be a pre-imprimitivity $A_0-B_0$ bimodule, with completion
${}_AX_B$.  Suppose we have homomorphisms $\phi\:A\to M(C)$
and $\psi\:B\to M(D)$, and a map $\Phi_0\:X_0\to M(Y)$ such that:
\begin{enumerate}
\item
${}_{M(C)}\<\Phi_0(x),\Phi_0(y)\>=\phi({}_{A_0}\<x,y\>)$
and
\item $\<\Phi_0(x),\Phi_0(y)\>_{M(D)}=\psi(\<x,y\>_{B_0})$
\end{enumerate}
for all $x,y \in X_0$.  Then $\Phi$ extends uniquely
to an imprimitivity bimodule homomorphism $\Phi\:X\to M(Y)$.

Moreover, if $\phi$ and $\psi$ are nondegenerate, then so is $\Phi$.
\end{lem}

\begin{proof}
Similarly to the proof of the preceding lemma, the inner product
properties (i)--(ii) imply first that $\Phi_0$ is linear, then that it
is bounded, hence extends uniquely to a bounded linear map $\Phi\:X\to
M(Y)$ which preserves both inner products by continuity and density;
\remref{rHb-hom-rem} assures us that this suffices to give the unique
extension $\Phi$.  We have already mentioned in \remref{rHb-hom-rem1}
that nondegeneracy of $\phi$ and $\psi$ guarantee nondegeneracy of
$\Phi$.
\end{proof}

\begin{defn}
\label{def-strict}
Let ${}_AX_B$ be a right-Hilbert bimodule.  The \emph{strict topology}
on $M(X)$ is that generated by the seminorms $m \mapsto \|Tm\|,\|mb\|$ for
$T \in \c K_B(X)$ and $b \in B$.
\end{defn}

\begin{rem}\label{rem-strict}
     (1) Thus the strict topology on $M({}_AX_B)$ has nothing to do with
the left $A$-module action and only depends on the right-partial
imprimitivity bimodule structure of $_{\K(X)}X_B$.

(2) What \cite{jt} and
\cite{ng:module} call the strict topology on $M(X)$ is really the
$*$-strong topology.  This is easily seen to be weaker than the
strict topology.%
\footnote{Is it always the same?  Presumably not, but we do
not know a counterexample.}   
Indeed, if $m_i\to m$ strictly in $M(X)$,
then  $m_ib\to mb$ for all $b\in B$ by definition.
Thus we only have to check that $m_i^*x\to m^*x$ for all
$x\in X$, which follows from factoring $x=cy$ for some
$c\in \c K(X)$, $y\in X$ and computing
$$m_i^*(cy)=\<m_i,cy\>_{M(B)}=\<c^*m_i,y\>_{M(B)}\to\<c^*m,y\>_{M(B)}
=m(cy).
$$
\end{rem}

The above definition is modeled after
\cite[Definition~1.5]{er:mult},
where it is discussed for imprimitivity
bimodules. If $X$ is full,
the strict topology on $M(X)$ only depends upon the
$\c K(X)-B$
imprimitivity bimodule structure on $X$ and all results from
\cite{er:mult} are available. However, we need to generalize several
results of \cite{er:mult} to
more general right-Hilbert bimodules.

Note first that we have a canonical
right-Hilbert $\c L(X)-M(B)$ bimodule embedding
of $X$ into $M(X)$ given by
$x \mapsto m_x$, where $m_xb = x \d b$ and $m_x^*(y)
= \<x,y\>_B$. It follows directly from the definitions
that this embedding preserves the left and right
actions and the $M(B)$-valued inner product.
Moreover, since the left action of $M(A)$ on $X$ and $M(X)$
is given via the same homomorphism $\kappa\:M(A)\to \c L(X)$,
we see that $x\to m_x$ is also a right-Hilbert $M(A)-M(B)$
bimodule homomorphism.

\begin{prop}[{\cf\ \cite[Proposition 1.6]{er:mult}}]
     \label{prop-strictlydense}
     Let $X$ be a right-Hilbert $A-B$ bimodule and let us
     view $X$ as an $\c L(X)-M(B)$ sub-bimodule of $M(X)$
     via the above embedding.
     Then\textup:
     \begin{enumerate}
	\item
     $\c K(X)\cdot M(X)\subseteq X$ and $M(X)\cdot B\subseteq X$.
     \item  $M(X)$ is the strict completion of $X$.
     \item The pairings
     $$\c L(X)\times M(X)\to M(X),\quad M(A)\times M(X)\to
     M(X)\;\;\text{and}$$
     $$M(X)\times M(B)\to M(X),\quad
     \<\cdot,\cdot\>_{M(B)}\:M(X)\times M(X)\to M(B)$$
     are separately strictly continuous, where we identify $\c L(X)$
     with $M(\c K(X))$.
     \end{enumerate}
\end{prop}
\begin{proof} Let $m\in M(X)$.
     Factoring $b\in B$ as $b=cd$ for $c,d\in B$, it follows that
     $mb=(mc)(d)\in X$. On the other side,
     we have $_{\K(X)}\lk x,y\rk m=x\lk y, m\rk_{M(B)}\in X$
     for all $x,y\in X$, $m\in M(X)$.  This proves (i).

     For the proof of (ii), let $(u_i)_{i\in I}$ be an approximate
     unit in $B$. Then it follows from (i) that $mu_i$ converges
     strictly to $m$. Thus $X$ is strictly dense
     in $M(X)$. Conversely, assume that $(m_i)_{i\in I}$
     is a strict Cauchy net in $M(X)$. Then it follows from
     Remark \ref{rem-strict} that $(m_i)_{i\in I}$ is also
     a $*$-strong Cauchy net. This implies that we can define
     $m\in M(X)$ by $m(b)=\lim_{i\in I} m_i(b)$,
     $m^*(x)=\lim_{i\in I}m_i^*(x)$.
     Now let $c\in \c K(X)$. Since $(cm_i)_{i\in I}$ is a Cauchy net
     in $X$ by assumption, it follows that
     $cm_i$ converges to some $y\in X$, and since
     $yb=\lim_{i\in I}(cm_i)b=\lim_{i\in I}c(m_ib)=c(mb)=(cm)b$
     for all $b\in B$, it follows that $y=cm$. Thus $M(X)$
     is complete in the strict topology.

     The proof of (iii) now follows from an application of
     Remark \ref{rem-strict} and the norm-continuity of the
     pairings. For example, if $m_i\to m$ strictly in $M(X)$
     and $n\in M(X)$, then it follows for all $b\in B$ that
     $$
     \<m_i,n\>_{M(B)}b=m_i^*(n(b))\to m^*(n(b))=\<m,n\>_{M(B)}b,
     $$ by $*$-strong convergence of $(m_i)_{i\in I}$.
     Similarly,
     $$b\<m_i,n\>_{M(B)}=\<m_i(b^*),n\>_{M(B)}
     \to \<m(b^*),n\>_{M(B)}=b\<m,n\>_{M(B)}.$$
     It follows that $\<m_i,n\>_{M(B)}\to \<m,n\>_{M(B)}$ strictly.
     The other  assertions follow from similar arguments.
\end{proof}

In fact, $M(X)$ is maximal with respect to the above properties:

\begin{prop}\label{MX-char-prop}
Let $X$ be a right-Hilbert $A-B$ bimodule.  Suppose $M$ is a
right-Hilbert $M(A)-M(B)$ bimodule containing $X$ as a right-Hilbert
$M(A)-M(B)$ sub-bimodule \textup(extending the operations on $X$ in
the usual way\textup) such that
\[
M\d B\subseteq X.
\]
Then $M$ embeds as an  $M(A)-M(B)$ sub-bimodule of $M(X)$.
\end{prop}

\begin{proof}
For $m\in M$ define $\Phi(m)\:B\to X$ by
\begin{equation}
\label{MX-char-eq1}
\Phi(m)b=mb.
\end{equation}
Then $\Phi(m)$ is adjointable, with
\begin{equation}
\label{MX-char-eq2}
\Phi(m)^*x=\<m,x\>_{M(B)}.
\end{equation}
To see this, note
that $\<m,x\>_{M(B)}\in B$, since we may factor $x=yb$ for some
$y\in X$ and $b\in B$, and then
\[
\<m,x\>_{M(B)}=\<m,yb\>_{M(B)}=\<m,y\>_{M(B)}b\in B.
\]
Checking the
adjoint property, for $b\in B$ and $x\in X$ we have
\[
\bigl\<b,\<m,x\>_{M(B)}\bigr\>_B
= b^*\<m,x\>_{M(B)}
= \<mb,x\>_{M(B)}
= \<mb,x\>_B
= \<\Phi(m)b,x\>_B.
\]
Thus we get a map $\Phi\:M\to M(X)$.

For each $b\in B$,
$k\in M(A)$, and $m,n\in M$ we have
\[
\Phi(km)b
= kmb
= k\Phi(m)b
\]
and
\[
\<\Phi(m),\Phi(n)\>_{M(B)}b
= \Phi(m)^*\Phi(n)b
= \Phi(m)^*nb
= \<m,nb\>_{M(B)}
= \<m,n\>_{M(B)}b.
\]
Thus $\Phi\:M\to M(X)$ is a right-Hilbert $M(A)-M(B)$ bimodule
homomorphism which is necessarily isometric; of course $\Phi$
restricts to the usual embedding of $X$ in $M(X)$.

For uniqueness, 
suppose that $\Psi\:M\to M(X)$ is another right-Hilbert
$M(A)-M(B)$-bimodule homomorphism extending the embedding of $X$.
Then, for all $b\in B$,
$(\Phi(m)-\Psi(m))(b)=mb-\Psi(mb)=0$ since $mb\in X$.
Thus $\Psi=\Phi$.
\end{proof}

\begin{rem}\label{two-mult-rem}
When $X$ is an $A-B$ imprimitivity bimodule, we have two ostensibly
different $M(A)-M(B)$ bimodules: $M({}_AX) \deq \L_A(A,X)$ and
$M(X_B)\deq
\L_B(B,X)$. In fact these turn out to be naturally isomorphic. The
link is
the multiplier bimodule $M({}_AX_B)$ studied in \cite{er:mult}, which
consists of
pairs $(m_A,m_B)$ in which $m_A \: A \to X$ is $A$-linear, $m_B \:
B \to X$
is $B$-linear, and
\[
m_A(a) \d b = a \d m_B(b)
\righttext{for all} \ a \in A, b \in B.
\]
The maps $(m_A,m_B)\mapsto m_A$ and $(m_A,m_B)\mapsto m_B$ give
isomorphisms of $M({}_AX_B)$ onto $M({}_AX)$ and $M(X_B)$,
respectively \cite[Proposition~1.3]{er:mult}. This is not
obvious: there is
\emph{a priori} no assertion of adjointability in the definition of
$M({}_AX_B)$, and as a result the abstract characterization of
$M({}_AX_B)$ in \cite[Proposition~1.2]{er:mult}\ looks rather different
from its one-sided analogue in \propref{MX-char-prop}.
Unfortunately, a similar result does not hold for arbitrary
partial imprimitivity bimodules, so that we definitely have a
one-sided theory of multiplier bimodules for
partial imprimitivity bimodules.
\end{rem}

We will now show that nondegeneracy
allows us to extend right-Hilbert bimodule homomorphisms to the
multiplier bimodules.

\begin{thm}
\label{rHb-hom-ext-prop}
Let ${}_\phi\Phi_\psi\:{}_AX_B\to M({}_CY_D)$ be a nondegenerate
right-Hilbert bimodule homomorphism.  Then $\Phi$ extends uniquely to
a right-Hilbert bimodule homomorphism
$${_{\bar\phi}\bar\Phi_{\bar\psi}}\:
{}_{M(A)}M(X)_{M(B)}\to {}_{M(C)}M(Y)_{M(D)},$$
where $\bar\phi\:M(A)\to M(C)$ and $\bar\psi\:M(B)\to
M(D)$ are the unique extensions of $\phi$ and $\psi$.
Moreover, $\bar{\Phi}$ is strictly continuous.
\end{thm}

\begin{proof}
We first aim to show that $\Phi$ is continuous from the relative
strict topology of $X$ to the strict topology of $M(Y)$.  By
\remref{rHb-hom-rem1} there exists a map $\mu\:\c K(X)\to M(\c K(Y))$
such that ${}_\mu\Phi_\psi\:{}_{\c K(X)}X_B \to M({}_{\c K(Y)}Y_D)$ is
a nondegenerate right-Hilbert bimodule homomorphism.  Let $x_i\to 0$
strictly in $X$.  Then for all $d\in D$ we can factor $d=\psi(b)d'$
for some $b\in B$ and $d'\in D$, and then
\[
\|\Phi(x_i)d\|
= \|\Phi(x_i)\psi(b)d'\|
= \|\Phi(x_i \d b)d'\|
\le \|\Phi(x_i \d b)\|\|d'\|
\le \|x_i \d b\|\|d'\|
\to 0.
\]
Similarly, for all $T \in \c K(Y)$ we can factor $T=T'\mu(S)$ for
some $T' \in \c K(Y)$ and $S \in \c K(X)$, and then
\[
\|T\Phi(x_i)\|
= \|T'\mu(S)\Phi(x_i)\|
= \|T'\Phi(Sx_i)\|
\le \|T'\|\|\Phi(Sx_i)\|
\le \|T'\|\|Sx_i\|
\to 0.
\]
Thus $\Phi(x_i) \to 0$ strictly in $M(Y)$, as desired.

Therefore $\Phi$ certainly has a unique strictly continuous extension
$\bar\Phi \: M(X) \to M(Y)$.
Since $\bar\phi$ and $\bar\psi$ are also strictly continuous,
it follows from the strict continuity of all
bimodule operations that $\bar\Phi$ is a $\bar\mu-\bar\psi$ compatible
right-Hilbert $\c K(X)-B$ bimodule homomorphism.

Finally, for the uniqueness, suppose $\Psi \: M(X) \to M(Y)$ is any
right-Hilbert bimodule homomorphism extending $\Phi$. Since the
extensions of $\phi$ and $\psi$ to the multiplier algebras are unique,
$\Psi$ must be $\bar\phi-\bar\psi$ compatible. For all $m \in M(X)$ and
$d \in D$ we can factor $d=\psi(b)d'$ for some $b \in B$ and $d' \in D$,
and then
\[
\Psi(m)d
= \Psi(m)\psi(b)d'
= \Psi(mb)d'
= \Phi(mb)d'
= \bar\Phi(m)d.
\]
Therefore $\Psi(m)=\bar\Phi(m)$.
\end{proof}

The above result allows us
to compose nondegenerate
right-Hilbert bimodule homomorphisms: if $\Phi \: {}_AX_B \to
M({}_CY_D)$ and $\Psi \: {}_CY_D \to M({}_EZ_F)$ are nondegenerate
homomorphisms, then we have a nondegenerate composition $\Psi \circ
\Phi \: X \to M(Z)$.

\section{Tensor products}
\label{sec-tensor}

Let ${}_AX_B$ and ${}_BY_C$ be right-Hilbert bimodules.  Then the
algebraic tensor product $X\odot Y$ becomes a pre-right-Hilbert $A-C$
bimodule with operations
\begin{gather*}
a \d (x\otimes y) = a \d x\otimes y, \quad
(x\otimes y) \d b = x\otimes y \d b,
\midtext{and}
\\\<x\otimes y,z\otimes w\>_C = \<y,\<x,z\>_B \d w\>_C
\end{gather*}
\cite[Proposition 3.16]{RW}. The completion is a
right-Hilbert $A-C$ bimodule $X\otimes_B Y$. Note that if $X$ and $Y$
are full, then so is $X\otimes_BY$, since
\[
\overline{\<X\otimes_BY,X\otimes_BY\>}_C
= \overline{\<Y,\<X,X\>_B \d Y\>}_C
= \overline{\<Y,B \d Y\>}_C
= \overline{\<Y,Y\>}_C = C.
\]

\begin{defn}\label{rHb-tnsr-defn}
We call $X\otimes_B Y$ the \emph{balanced tensor product} of
${}_AX_B$ and ${}_BY_C$.
\end{defn}

Given two 
right-Hilbert bimodule homomorphisms $\Phi
\: {}_AX_B \to M({}_DZ_E)$ and $\Psi \: {}_BY_C \to M({}_EW_F)$, we
would like to form a tensor product
homomorphism
$$\Phi\otimes_B\Psi \: X\otimes_BY \to M(Z\otimes_EW).$$
There are two obstructions to this utopia: first, understandably we
have to require that the coefficient maps at $B$ coincide.  More
significantly, the obvious map
$x\otimes_By\mapsto\Phi(x)\otimes_{M(E)}\Psi(y)$ will be into
$M(Z)\otimes_{M(E)}M(W)$.  Fortunately, using the abstract
characterization of \propref{MX-char-prop}, we can show that
$M(Z)\otimes_{M(E)}M(W)$ always embeds in $M(Z\otimes_E W)$:

\begin{lem}\label{M-tnsr-emb-lem}
Let ${}_DZ_E$ and ${}_EW_F$ be right-Hilbert bimodules. There exists
an isometric right-Hilbert $M(D)-M(F)$ bimodule homomorphism
$\Upsilon \: M(Z)\otimes_{M(E)}M(W) \to M(Z\otimes_EW)$ such that
\begin{equation}\label{M-tnsr-emb-eq}
\Upsilon(m)f = m \d f \midtext{and}
\Upsilon(m)^*x = \<m,x\>_{M(F)}
\end{equation}
for each $m \in M(Z)\otimes_{M(E)}M(W)$, $f \in F$, and
$x \in Z\otimes_EW$.
\end{lem}

\begin{proof}
Note that $M(Z)\otimes_{M(E)}M(W)$ is a right-Hilbert $M(D)-M(F)$
bimodule which contains a copy of $Z\otimes_EW = Z\otimes_{M(E)}W$.
Moreover,
\begin{align*}
M(Z)\otimes_{M(E)}M(W) \d F
&= M(Z)\otimes_{M(E)}W = M(Z)\otimes_{M(E)}E \d W
\\&= M(Z) \d E\otimes_{M(E)}W = Z\otimes_{M(E)}W
= Z \otimes_E W.
\end{align*}
Hence \propref{MX-char-prop}\ provides an isometric
  right-Hilbert
$M(D)-M(F)$ bimodule homomorphism
$\Upsilon \: M(Z)\otimes_{M(E)}M(W) \to M(Z\otimes_EW)$
which is the identity on $Z\otimes W$, and which by
\eqref{MX-char-eq1}\ and \eqref{MX-char-eq2} in
the proof of \propref{MX-char-prop}\ satisfies \eqref{M-tnsr-emb-eq}.
\end{proof}

\begin{rem}
In the future we will suppress the map $\Upsilon$, using the above
lemma to identify $M(Z)\otimes_{M(E)}M(W)$
with its image in
$M(Z\otimes_EW)$. Note that for $m\in M(Z)$ and $n\in M(W)$ the
adjointable map $m\otimes n\: F \to Z\otimes_EW$ is given by
$(m\otimes n)f=m\otimes n\d f$.
\end{rem}

\begin{prop}\label{tnsr-hom-prop}
Let ${}_\phi\Phi_\psi\:{}_AX_B\to M({}_DZ_E)$ and
${}_\psi\Psi_\theta\:{}_BY_C\to M({}_EW_F)$ be right-Hilbert bimodule
homomorphisms.  Then there exists a $\phi-\theta$ compatible
right-Hilbert bimodule homomorphism
$\Phi\otimes_B\Psi\:{}_A(X\otimes_BY)_C\to M({}_D(Z\otimes_EW)_F)$
such that
\begin{equation}\label{tnsr-hom-eq}
(\Phi\otimes_B\Psi)(x\otimes y) = \Phi(x)\otimes \Psi(y).
\end{equation}
If $\Phi$ and $\Psi$ are nondegenerate, then $\Phi \otimes_B \Psi$ is
too.
\end{prop}

\begin{proof}
The pairing $(x,y)\mapsto\Phi(x)\otimes\Psi(y)$ is bilinear, so determines
a unique linear map \mbox{$\Phi\odot\Psi\:X\odot Y\to M(Z\otimes_EW)$}.
We verify conditions (i)--(ii) of \lemref{pre-rHb-lem}, and it suffices
to do this for elementary tensors: for all $a\in A$, $x,z\in X$, and
$y,w\in Y$ we have
\begin{align*}
(\Phi\odot\Psi)(a\d(x\otimes y))
&= (\Phi\odot\Psi)(a\d x\otimes y)
= \Phi(a\d x)\otimes\Psi(y)
\\&= \phi(a)\d\Psi(x)\otimes\Psi(y)
= \phi(a)\d(\Phi(x)\otimes\Psi(y))
\end{align*}
and
\begin{align*}
&\bigl\<(\Phi\odot\Psi)(x\otimes y),(\Phi\odot\Psi)(z\otimes
w)\bigr\>_{M(F)}
= \bigl\<\Phi(x)\otimes\Psi(y),\Phi(z)\otimes\Psi(w)\bigr\>_{M(F)}
\\&\quad=
\bigl\<\Psi(y),\<\Phi(x),\Phi(x)\>_{M(E)}\d\Psi(w)\bigr\>_{M(F)}
= \bigl\<\Psi(y),\psi(\<x,z\>_B)\d\Psi(w)\bigr\>_{M(F)}
\\&\quad= \bigl\<\Psi(y),\Psi(\<x,z\>_B\d w)\bigr\>_{M(F)}
= \theta\bigl(\<y,\<x,z\>_B\d w\>_C\bigr)
\\&\quad= \theta\bigl(\<x\otimes y,z\otimes w\>_C\bigr).
\end{align*}

Now suppose $\Phi$ and $\Psi$ are nondegenerate. Then we have
\begin{align*}
&\overline{(\Phi\otimes_B\Psi)(X\otimes_B Y)\d F}
= \overline{\bigl(\Phi(X)\odot\Psi(Y)\bigr)F}
= \overline{\Phi(X)\odot\Psi(Y)\d F}
\\&\quad= \overline{\Phi(X)\odot W}
= \overline{\Phi(X)\odot E\d W}
= \overline{\Phi(X)\d E\odot W}
= Z \otimes_E W.
\end{align*}
Since $\phi$ and $\theta$ are nondegenerate by assumption, this shows
that $\Phi \otimes_B \Psi$ is nondegenerate.
\end{proof}

We will need the exterior tensor product as well.  Let ${}_AX_B$ and
${}_CY_D$ be right-Hilbert bimodules.  Then the algebraic tensor
product $X\odot Y$ becomes a pre-right-Hilbert $(A\odot C)-(B\odot D)$
bimodule with operations
\begin{gather*}
(a\otimes c)\d(x\otimes y) = a\d x\otimes c\d y,
\quad
(x\otimes y)\d(b\otimes d) = x\d b\otimes y\d d,
\midtext{and}
\\\<x\otimes y,z\otimes w\>_{B\otimes D} = \<x,z\>_B\otimes\<y,w\>_D
\end{gather*}
\cite[Proposition 3.36]{RW}. The completion is a right-Hilbert $(A\otimes
C)-(B\otimes D)$ bimodule $X\otimes Y$, where ``$\otimes$'' always denotes
the {\em minimal} tensor product for $C^*$-algebras.
If $X$ and $Y$ are full, then
so is $X\otimes Y$, since
\[
\overline{\<X\otimes Y,X\otimes Y\>}_{B\otimes D} =
\overline{\<X,X\>_B\otimes\<Y,Y\>_D} = B\otimes D.
\]

\begin{defn}\label{xtnsr-defn}
We call $X\otimes Y$ the \emph{exterior tensor
product} of ${}_AX_B$ and ${}_CY_D$.
\end{defn}

As for balanced tensor products, we want to ``exterior tensor'' 
bimodule
homomorphisms, and again we first need to know that things go
into the right place:

\begin{lem}\label{M-xtnsr-lem}
Let ${}_EZ_F$ and ${}_GW_H$ be right-Hilbert bimodules. There exists
an isometric right-Hilbert bimodule homomorphism $\Xi \: M(Z)\otimes
M(W) \to M(Z\otimes W)$ with coefficient maps the canonical
homomorphisms $M(E)\otimes M(G) \to M(E\otimes G)$ and $M(F)\otimes
M(H) \to M(F\otimes H)$ such that
\begin{equation}\label{M-xtnsr-emb-eq}
\Xi(m)b = m \d b \midtext{and}
\Xi(m)^*x = \<m,x\>_{M(F)\otimes M(H)}
\end{equation}
for each $m \in M(Z)\otimes M(W)$, $b \in F\otimes H$, and $x \in
Z\otimes W$.
\end{lem}

\begin{proof}
Note that $M(Z)\otimes M(W)$ contains $Z\otimes W$ in the obvious way
as a right-Hilbert $(M(E)\otimes M(G))-(M(F)\otimes M(H))$
sub-bimodule, that $E\otimes G\subseteq M(E)\otimes M(G)$ and
$F\otimes H\subseteq M(F)\otimes M(H)$ are essential closed ideals,
and that
\[
(M(Z)\otimes M(W)) \d (F\otimes H)
= M(Z)\d F\otimes M(W)\d H
\subseteq Z\otimes W.
\]
Hence \propref{MX-char-prop}\ provides an isometric
right-Hilbert bimodule
homomorphism $\Xi \: M(Z)\otimes M(W) \to M(Z\otimes W)$ as
desired.  By \eqref{MX-char-eq1}\ and \eqref{MX-char-eq2}\ in the
proof of \propref{MX-char-prop}, $\Xi$ satisfies
\eqref{M-xtnsr-emb-eq}.
\end{proof}

\begin{rem}
As for balanced tensor products, we will suppress the map $\Xi$,
using the above lemma to identify $M(Z) \otimes M(W)$ with its image
in $M(Z\otimes W)$. Note that for $m\in M(Z)$ and $n\in M(W)$ the
adjointable map $m\otimes n\:F\otimes H\to Z\otimes W$ is given on
elementary tensors by $(m\otimes n)(f\otimes h) = m\d f\otimes
n\d h$.
\end{rem}

\begin{prop}\label{xtnsr-hom-prop}
Let ${}_\phi\Phi_\psi\:{}_AX_B\to M({}_EZ_F)$ and
${}_\theta\Psi_\rho\:{}_CY_D\to M({}_GW_H)$ be right-Hilbert bimodule
homomorphisms.  There exists a right-Hilbert bimodule homomorphism
$$\Phi\otimes\Psi\:{}_{A\otimes C}(X\otimes Y)_{B\otimes D}\to
M({}_{E\otimes G}(Z\otimes W)_{F\otimes H}),$$
with coefficient maps
the usual homomorphisms $\phi\otimes\theta\:A\otimes C\to M(E\otimes
G)$ and $\psi\otimes\rho\:B\otimes D\to M(F\otimes H)$, such that
\begin{equation}\label{xtnsr-hom-eq}
(\Phi\otimes\Psi)(x\otimes y) = \Phi(x)\otimes \Psi(y).
\end{equation}
If $\Phi$ and $\Psi$ are nondegenerate, then $\Phi\otimes\Psi$ is too.
\end{prop}

\begin{proof}
As in the proof of \propref{tnsr-hom-prop}, we clearly have a unique
linear map $\Phi\odot\Psi\:X\odot Y\to M(Z)\otimes M(W)$, so for the first
part it suffices to check the left module actions and the inner products
on the generators: for $a\in A$, $c\in C$, $x,z\in X$, and $y,w\in Y$
we have
\begin{align*}
&(\Phi\otimes\Psi)\bigl((a\otimes c)\d(x\otimes y)\bigr)
= (\Phi\otimes\Psi)(a\d x\otimes c\d y\bigr)
\\&\quad= \Phi(a\d x)\otimes\Psi(c \d y)
= \phi(a)\d\Phi(x)\otimes\theta(c)\d\Psi(y)
\\&\quad= \bigl(\phi(a)\otimes\theta(c)\bigr)\d
\bigl(\Phi(x)\otimes\Psi(y)\bigr)
= (\phi\otimes\theta)(a\otimes c)\d(\Phi\otimes\Psi)(x\otimes y)
\end{align*}
and
\begin{align*}
&\bigl\<(\Phi\otimes\Psi)(x\otimes y),
(\Phi\otimes\Psi)(z\otimes w)\bigr\>_{M(F\otimes H)}
= \bigl\<\Phi(x)\otimes\Psi(y),
\Phi(z)\otimes\Psi(w)\bigr\>_{M(F\otimes H)}
\\&\quad= \<\Phi(x),\Phi(z)\>_{M(F)}\otimes
\<\Psi(y),\Psi(w)\>_{M(H)}
= \psi(\<x,z\>_B)\otimes\rho(\<y,w\>_D)
\\&\quad= (\psi\otimes\rho)\bigl(\<x,z\>_B\otimes\<y,w\>_D\bigr)
= (\psi\otimes\rho)\bigl(\<x\otimes y,z\otimes w\>_{B\otimes D}\bigr).
\end{align*}

Now suppose $\Phi$ and $\Psi$ are nondegenerate. Then
\begin{align*}
\overline{(\Phi\otimes\Psi)(X\otimes Y)\d(F\otimes H)}
&= \overline{(\Phi(X)\otimes\Psi(Y))\d(F\otimes H)}
\\&= \overline{\Phi(X)\d F\otimes\Psi(Y)\d H}
= Z\otimes W,
\end{align*}
so $\Phi\otimes\Psi$ is nondegenerate as well.
\end{proof}

\section{The $C$-multiplier bimodule $M_C(X\otimes C)$}
\label{sec-Cmult}

In this section we introduce the notion of
the $C$-multiplier bimodule of the exterior tensor
product $X\otimes C$ of a right-Hilbert bimodule $X$
with a $C^*$-algebra $C$. The definition is very similar to the
definition of the $C$-multiplier algebra $M_C(A\otimes C)$
as given in \appxref{coactions-chap}, and we shall obtain 
properties for $C$-multiplier bimodules similar to those we
obtained in \propref{prop-Cstrict} and \propref{prop-Cmultextend}
for $C$-multiplier bimodules. We start with:

\begin{defn}\label{defn-Cmultmod}
     Suppose that $X$ is a right-Hilbert $A-B$ bimodule and let
     $C$ be a $C^*$-algebra. The {\em $C$-multiplier bimodule}
     $M_C(X\otimes C)$ of $X\otimes C$ is defined as the set
     $$M_C(X\otimes C)=\{m\in M(X\otimes C)\mid(1\otimes C)m\cup
     m(1\otimes C)\subseteq X\otimes C\}.$$
     The {\em $C$-strict topology} on $M_C(X\otimes C)$ is the
     locally convex topology generated by the seminorms
     $m\mapsto \|m(1\otimes c)\|$ 
and $m\mapsto\|(1\otimes c)m\|$, $c\in C$.
\end{defn}

If $C=C^*(G)$ for some locally compact group $G$, we shall
simply write $M_G(X\otimes C^*(G))$ for $M_{C^*(G)}(X\otimes C^*(G))$
and call it the {\em $G$-multiplier
bimodule} of $X\otimes C^*(G)$. The following result
is the bimodule analogue of \propref{prop-Cstrict}.

\begin{lem}\label{lem-Cstrictcont}
     Suppose that $A$, $B$, and $C$ are $C^*$-algebras and that
     $X$ is a right-Hilbert $A-B$ bimodule.
     Then\textup:
     \begin{enumerate}
	\item $M_C(X\otimes C)$ is a closed $A\otimes C- B\otimes C$
     sub-bimodule of $M(X\otimes C)$.
     \item $M_C(X\otimes C)$ is a right-Hilbert $M_C(A\otimes
     C)-M_C(B\otimes C)$ bimodule with respect to the
     bimodule operations on $M(X\otimes C)$ restricted to
     $M_C(X\otimes C)$.
     \item The $C$-strict topology on  $M_C(X\otimes C)$is stronger than the
    strict topology inherited from the  full multiplier
     bimodule $M(X\otimes C)$.
     \item The right-Hilbert $M_C(A\otimes
     C)-M_C(B\otimes C)$ bimodule operations on $M_C(X\otimes C)$
     are separately $C$-strictly continuous.
     \item $M_C(A\otimes C)$ and $M_C(B\otimes C)$ are the $C$-strict
     completions of $A\otimes C$ and $B\otimes C$, respectively,
     and $M_C(X\otimes C)$ is the $C$-strict completion of $X\otimes C$.
     \end{enumerate}
\end{lem}
\begin{proof}
     We omit the straightforward verifications of (i) and (ii).
     For (iii), suppose that $(m_i)_{i\in I}$ is
     a net in $M_C(X\otimes C)$ which converges $C$-strictly to some
     $m\in M_C(X\otimes C)$. Let $z\in B\otimes C$. Factoring
     $z=(1\otimes c)y$ for some $c\in C$, $y\in B\otimes C$ gives
     $$\lim_i m_iz=\lim_i(m_i(1\otimes c))y
     =m(1\otimes c)y=mz.$$
     A similar argument shows that $zm_i\to zm$ for all $z\in
     \K(X\otimes C)\cong \K(X)\otimes C$.
     Thus $m_i\to m$ in the strict topology of $M(X\otimes C)$.

     Part~(iv) follows directly from (iii) and Proposition \ref{prop-strictlydense}.
     For the proof of (v), let $(u_i)_{i\in I}$ be a bounded
     approximate unit for $C$. Then, if $m\in M_C(X\otimes C)$, we
     see that $m(1\otimes u_i)$ is a net in $X\otimes C$
     which converges $C$-strictly to $m$, so $X\otimes C$ is
     $C$-strictly dense in $M_C(X\otimes C)$. Conversely, if
     $(m_i)_{i\in I}$ is a $C$-strict Cauchy net in $M_C(X\otimes C)$,
     it follows from (iii) that it is also a strict Cauchy net
     in $M(X\otimes C)$. It follows then from Proposition 
\ref{prop-strictlydense}
     that there exists an $m\in M(X\otimes C)$ such that $m_i\to m$
     strictly.

     We claim that $m_i\to m$ $C$-strictly and that $m\in
     M_C(X\otimes C)$.
     For this we let $\eps>0$ and $c\in C$. We show that there exists an
     $i_0\in I$ such that $\|(m_i-m)(1\otimes c)z\|\leq\eps$
     for all $z\in B\otimes C$ with $\|z\|\leq 1$ 
and $i\geq i_0$.
     For this let $i_0\in I$ such that $\|(m_i-m_j)(1\otimes c)\|\leq\eps$
     for all $i,j\geq i_0$. Since
     $m_j(1\otimes c)z\to m(1\otimes c)z$ in norm, it follows that
     $$\|(m_i-m)(1\otimes c)z\|=\lim_j\|(m_i-m_j)(1\otimes c)z\|\leq
     \eps$$
     for all $i\geq i_0$. Thus $m_i(1\otimes c)\to m(1\otimes c)$ in
     norm,
     and a similar argument shows that $(1\otimes c)m_i\to (1\otimes
     c)m$ in norm for all $c\in C$. Thus $m_i\to m$ $C$-strictly.
     Finally, since $m_i(1\otimes c), (1\otimes c)m_i\in X\otimes C$,
     it follows from the fact that $X\otimes C$ is norm-closed
     in $M(X\otimes C)$ that $(1\otimes c)m, m(1\otimes c)\in X\otimes C$
     for all $c\in C$, hence $m\in M_C(X\otimes C)$.
\end{proof}

\begin{rem}\label{rem-multiply}
     It is useful to observe that we always have
     $$(1\otimes M(C))\cdot M_C(X\otimes C)\cup
     M_C(X\otimes C)\cdot (1\otimes M(C))
     \subseteq M_C(X\otimes C).$$
     We omit the straightforward argument.
\end{rem}

\begin{prop}\label{prop-Cmult}
     Suppose that $A$, $B$, $C$, $D$, $E$, $F$ are $C^*$-algebras.
     Let $X$ be a right-Hilbert $A-B$ bimodule and let
     $Y$ be a right-Hilbert $E-F$ bimodule. Assume further
     that $_{\varphi}\Phi_{\psi}\:{_A}X_B\to M({_E}Y_F)$ is a
right-Hilbert bimodule homomorphism,
     and that $\Psi\:C\to M(D)$ is a nondegenerate $*$-homomorphism.
     Let $\Phi\otimes\Psi\:X\otimes C\to M(Y\otimes D)$ denote
     the tensor product homomorphism of 
Proposition~\textup{\ref{xtnsr-hom-prop}}.
     Then\textup:
     \begin{enumerate}
     \item There exists a unique bimodule homomorphism
     $$\overline{\Phi\otimes \Psi}\:{_{M_C(A\otimes C)}}M_C(X\otimes
     C)_{M_C(B\otimes C)}\to {_{M(E\otimes D)}}M(Y\otimes
     D)_{M(F\otimes D)}$$
     which extends $\Phi\otimes \Psi\:X\otimes C\to M(Y\otimes D)$.
     \item The map $\overline{\Phi\otimes \Psi}$ of 
     \textup{(i)} is
     $C$-strict to strict continuous.
     \item If $\Phi(X)\subseteq Y$, then
     $\overline{\Phi\otimes \Psi}\big((M_C(X\otimes C)\big)\subseteq
     M_D(Y\otimes D)$ and $\overline{\Phi\otimes\Psi}$ is
     $C$-strict to $D$-strict continuous.
     \end{enumerate}
\end{prop}

\begin{rem}
      In the rest of the paper we shall usually write
     $\Phi\otimes\Psi$ also for its extension
     $\overline{\Phi\otimes\Psi}$ to $M_C(X\otimes C)$
     if confusion seems unlikely.
     Note that it already follows from \propref{prop-Cmultextend}
     that the coefficient maps $\phi\otimes\Psi\:A\otimes C
     \to M(E\otimes D)$ and $\psi\otimes \Psi\: B\otimes C\to
     M(F\otimes D)$ have unique extensions to
     $M_C(A\otimes C)$ and $M_C(B\otimes C)$, respectively,
     with similar properties as stated for $M_C(X\otimes C)$ above.
\end{rem}

\begin{proof}[Proof of Proposition \ref{prop-Cmult}]
     The proof follows closely the ideas of the
     proof of Theorem \ref{rHb-hom-ext-prop}, and it is also very
     similar to the proof of \propref{prop-Cmultextend}.
     We first remark that the map $\Phi\otimes\Psi\:X\otimes C\to
     M(Y\otimes D)$ is $C$-strict to strict continuous.
     Indeed, if $(x_i)_{i\in I}$ is a net in $X\otimes C$
     which converges $C$-strictly to $x\in X\otimes C$,
     and if $z\in Y\otimes D$, we can factor $z=(1\otimes\Psi(c))y$
     for some $c\in C$
     and $y\in Y\otimes D$ (since $\Psi$ is nondegenerate) to conclude that
     $$\Phi\otimes\Psi(x_i)z=\Phi\otimes \Psi(x_i(1\otimes c))y
     \to \Phi\otimes\Psi(x(1\otimes c))y=\Phi\otimes \Psi(x)z,$$
     where convergence is in norm. A similar argument shows that
     $k(\Phi\otimes \Psi(x_i))\to k(\Phi\otimes \Psi(x))$ in norm
     for all $k\in \K(Y\otimes D)\cong\K(Y)\otimes D$.
     Thus there exist unique $C$-strict to strict
     continuous linear extensions
     $\overline{\Phi\otimes\Psi}\:M_C(X\otimes C)\to M(Y\otimes D)$.

     Since all bimodule operations on $_{M_C(A\otimes C)}M_C(X\otimes
     C)_{M_C(B\otimes C)}$ are
     separately $C$-strictly continuous, and since all bimodule
     operations on $M(Y\otimes D)$ are separately strictly continuous,
     we conclude that the extensions
     $$\overline{\varphi\otimes\Psi}\:M_C(A\otimes C)\to M(E\otimes D)
     \quad\text{and}\quad
     \overline{\psi\otimes\Psi}\:M_C(B\otimes C)\to M(F\otimes D)$$
     are $*$-homomorphisms, and that
     $\overline{\Phi\otimes\Psi}$ is $\overline{\varphi\otimes\Psi}-
     \overline{\psi\otimes\Psi}$ compatible.
     A straightforward argument, similar to the one presented
     in the proof of Theorem \ref{rHb-hom-ext-prop}
     and \propref{prop-Cmultextend}, shows that
     ${_{\overline{\varphi\otimes\Psi}}}\overline{\Phi
     \otimes\Psi}_{\overline{\psi\otimes\Psi}}$ is indeed the
     only bimodule extension of $\Phi\otimes \Psi$ to
     $M_C(X\otimes C)$.  This proves~(i) and~(ii).

     For the proof of (iii) assume that $\Phi(X)\subseteq Y$.
     Then we can copy the arguments given in the proof
     of part (iv) of \propref{prop-Cmultextend} to conclude
     that $\Phi\otimes \Psi( X\otimes C)\subseteq M_D(Y\otimes D)$
     and that $\Phi\otimes\Psi\:X\otimes C\to M_D(Y\otimes D)$
     is $C$-strict to $D$-strict continuous.
     It then follows that
     $\Phi\otimes\Psi$ extends to a $C$-strict to $D$-strict 
continuous bimodule
     homomorphism of $M_C(X\otimes C)$ into $M_D(Y\otimes D)$,
     and by the uniqueness clause in (i), this must coincide
     with $\overline{\Phi\otimes \Psi}$. 
\end{proof}

It is an important observation that the $C$-multiplier bimodule
$M_C(X\otimes C)$ only depends on $X$ and not on the
coefficients.

\begin{lem}\label{lem-coefficient}
     Let $_AX_B$ be a right-Hilbert $A-B$ bimodule and let
     $C$ be a $C^*$-algebra. Let $\Phi\:{_A}X_B\to {_{\K(X)}}X_B$
     be the identity map. Then
     $$\overline{\Phi\otimes \id_C}\:M_C({_A}X_B\otimes C)
     \to M_C({_{\K(X)}}X_B\otimes C)$$
     is isometric and surjective.
     Also, if $B_X=\overline{\lk X,X\rk}_B$ and 
     $\Psi\:{_A}X_{B_X}\to {_A}X_{B}$ denotes the identity map,
     then
     $$\overline{\Psi\otimes \id_C}\:M_C({_A}X_{B_X}\otimes C)
     \to M_C({_A}X_B\otimes C)$$
     is isometric and surjective. In particular, the range of the
     $M_C(B\otimes C)$-valued inner product on
     $M_C(X\otimes C)$ lies in  $M_C(B_X\otimes C)$.
\end{lem}
\begin{proof} The result follows directly from the fact that
     $\Phi\otimes \id_C\: {_A}X_B\otimes C\to {_{\K(X)}}X_B\otimes C$
     and $\Psi\otimes \id_C\:{_A}X_{B_X}\otimes C\to {_A}X_B\otimes C$
     are linear homeomorphisms with respect to the $C$-strict
     topologies.
\end{proof}

\begin{prop}\label{prop-iso}
     Let $\Phi\:{_A}X_B\to M({_{E}}Y_F)$ and $\Psi\:C\to M(D)$ be as
     in Proposition~\textup{\ref{prop-Cmult}}. If $\Phi$ and $\Psi$ are
     isometric, then so is
     $\overline{\Phi\otimes\Psi}\:M_C(X\otimes C)\to M(Y\otimes D)$.
\end{prop}
\begin{proof}
    Let $B_X=\overline{\lk X,X\rk}_B$ and consider the composition
    $$
    \begin{CD}
    {_A}X_{B_X}\otimes C @>>> {_A}X_B\otimes C @>\Phi\otimes \Psi >>
    M(Y\otimes D).
    \end{CD}
    $$
    This composition extends to the composition
    $$
    \begin{CD}
    M_C({_A}X_{B_X}\otimes C) @>\cong >> M_C({_A}X_B\otimes C)
    @>\overline{\Phi\otimes \Psi} >>
    M(Y\otimes D).
    \end{CD}
    $$
    Thus in order to show that $\overline{\Phi\otimes \Psi}$ is
    isometric, it is enough to see that the above composition
    of maps is isometric. But by uniqueness of extensions,
    this composition equals the extension of
    $\Phi_1\otimes \Psi\:{_A}X_{B_X}\otimes C
    \to M(Y\otimes D)$, where $\Phi_1=\Phi$ with right coefficient
    map restricted to the ideal $B_X$ of $B$.
    Since $\Phi\:X\to M(Y)$ is isometric, it follows that
    the right coefficient map $\psi\:B\to M(F)$ restricts to an
    isometric $*$-homomorphism $\psi_1\:B_X\to M(F)$.
    It follows then from \propref{prop-Cmultextend}
    that the right coefficient map
    $\overline{\psi_1\otimes\Psi}\:M_C(B_X\otimes C)\to M(Y\otimes D)$
    of $\overline{\Phi_1\otimes \Psi}$
    is isometric. Thus
    $\overline{\Phi_1\otimes \Psi}$ is isometric, too.
\end{proof}

In the special case where $\Psi=\id_C$ and $\Phi(X)\subseteq Y$
it will be necessary to identify the
elements in the image of
$\overline{\Phi\otimes \id}_C\:M_C(X\otimes C)\to M_C(Y\otimes C)$.
This can be done as follows:

\begin{lem}\label{lem-Cmultidentify}
     Suppose that $\Phi\:{_AX_B}\to M({_EY_F})$ is an isometric
     right-Hilbert
     bimodule homomorphism with $\Phi(X)\subseteq Y$.
     Then the isometry
     $\overline{\Phi\otimes \id}_C\:M_C(X\otimes C)\to M_C(Y\otimes C)$
     has image
     $$M\deq\{m\in M_C(Y\otimes C)\mid (1\otimes C)m\cup m(1\otimes C)\subseteq
     \Phi\otimes \id_C(X\otimes C)\}.$$
\end{lem}
\begin{proof}
   It is clear that   $\overline{\Phi\otimes \id}_C\big(M_C(X\otimes
   C)\big)\subseteq M$. So let $m\in M$. Let
   $(c_i)_{i\in I}$ be a bounded approximate identity of $C$, and
   define $(z_i)_{i\in I}\subseteq X\otimes C$ by
   $z_i=(\Phi\otimes \id_C)^{-1}(m(1\otimes c_i))$.
   Since $m(1\otimes c_i)$ is a $C$-strict Cauchy net
   in $Y\otimes C$ which lies in the image of $\Phi\otimes \id_C$,
   it follows that $(z_i)_{i\in I}$ is a $C$-strict Cauchy net
   in $X\otimes C$. Since $M_C(X\otimes C)$ is the $C$-strict
   completion of $X\otimes C$, we find an $n\in M_C(X\otimes C)$
   with $z_i\to n$ $C$-strictly.
   Since $\overline{\Phi\otimes\id}_C\:M_C(X\otimes C)\to M_C(Y\otimes
   C)$ is $C$-strict to $C$-strict continuous, it follows that
   $$\overline{\Phi\otimes\id}_C(n)=
   \lim_i\overline{\Phi\otimes\id}_C(z_i)=\lim_i m(1\otimes c_i)=m.$$
\end{proof}

As a first application of this lemma we get

\begin{cor}\label{MCideal}
     Suppose that $I$ is a closed ideal of $B$, and let us
     view $M_C(I\otimes C)$ as a subalgebra of $M_C(B\otimes C)$
     via the unique extension of the embedding $I\otimes C\to
     B\otimes C$ as given by \propref{prop-Cmultextend}. Then
     $M_C(I\otimes C)$ is an ideal in $M_C(B\otimes C)$.
\end{cor}
\begin{proof}
Let $m\in M_C(I\otimes C)$. Then
$$m(b\otimes c)=m(1\otimes c)(b\otimes 1)\in (I\otimes C)(B\otimes 1)
=I\otimes C$$
for all elementary tensors $b\otimes c\in B\otimes C$, from which
it follows that $m(B\otimes C)\subseteq I\otimes C$.
Similarly, we get $(B\otimes C)m\subseteq I\otimes C$.

Thus, if
$n\in M_C(B\otimes C)$, it follows that
$$mn(1\otimes C)\subseteq m(B\otimes C)\subseteq I\otimes C,$$
and we also have $(1\otimes C)mn\subseteq (I\otimes C)n\in I\otimes C$,
since $I\otimes C$ is an ideal in $M(B\otimes C)$.
Thus, it follows from \lemref{lem-Cmultidentify}
that $mn\in M_C(I\otimes C)$. For symmetric reasons,
this also implies that $nm\in M_C(I\otimes C)$.
\end{proof}

\section{Linking algebras}\label{sec-link}

We collect here a few facts we will need concerning linking algebras of
partial imprimitivity bimodules; the primary references are \cite{bgr},
\cite{lan:hilbert}, and \cite{RW}, where the linking algebras
were studied for the special case of imprimitivity bimodules.
Linking algebras of Hilbert $C^*$-bimodules appear in~\cite{BMS-QM}. 

The \emph{linking algebra}
of a partial $A-B$ imprimitivity bimodule $X$ is
\[
L(X)\deq\left\{\left.
\begin{pmatrix}
a&x
\\
\widetilde z&b\end{pmatrix}\right| a \in A,\ b \in B,\ x,z \in X\right\},
\]
with the usual linear structure and $*$-algebra operations
\begin{equation}
\label{mat1}
\begin{pmatrix}
a&x
\\
\widetilde z&b
\end{pmatrix}
\begin{pmatrix}
a'&x'
\\
\widetilde{z'} & b'
\end{pmatrix}
= \begin{pmatrix}
aa'+{}_A\langle x,z'\rangle & a \d x' + x \d b'
\\
\widetilde z \d a' + b\d\widetilde{z'} &
\langle z,x'\rangle_B + bb'
\end{pmatrix}
\end{equation}
and
\begin{equation}
\label{mat2}
\begin{pmatrix}
a&x
\\
\widetilde z & b
\end{pmatrix}^*
= \begin{pmatrix}
a^* & z
\\
\widetilde x & b^*
\end{pmatrix}.
\end{equation}
$L(X)$ acts by adjointable operators on the Hilbert
$B$-module $X\oplus B$ via
\[
\begin{pmatrix}
a&x
\\
\widetilde z &b\end{pmatrix}\begin{pmatrix}
y
\\
c\end{pmatrix}=\begin{pmatrix}
a \d y+x \d c
\\
\langle z,y\rangle_B+bc\end{pmatrix},
\]
giving a homomorphism of $L(X)$ into the $C^*$-algebra ${\mathcal
L}_B(X\oplus B)$.  By considering the action on vectors of the form
$(z,0)$ or $(0,b)$, we can see that this homomorphism is injective
on $X$, $\rev X$ and $B$. Similarly, we can define a right-action
on the left Hilbert $A$-module $A\oplus X$, giving a $*$-homomorphism
of $L(X)$ into $\L_A(A\oplus X)$ which is injective on
$A$, $X$, and $\rev X$. Thus, defining a norm on $L(X)$ by
the maximum of the respective operator norms, we
obtain a complete $C^*$-norm on $L(X)$
(see \cite{bgr} and \cite[pages~50--51]{RW} ---
where the constructions have been done for imprimitivity bimodules,
in which case $L(X)$ embeds injectively into $\L_B(X\otimes B)$).
Since any
$*$-algebra has at most one complete $C^*$-norm, this unambiguously
makes $L(X)$ into a $C^*$-algebra.  The linking algebra of ${}_AX_B$
contains copies of $A$, $B$, $X$ and the conjugate
module $\widetilde{X}$.  We can recover these copies by
observing that the matrices
\[
p=p_{L(X)}\deq\begin{pmatrix}
1_{M(A)}&0
\\
0&0\end{pmatrix}\ \mbox{ and }\
q=q_{L(X)}\deq\begin{pmatrix}
0&0
\\
0&1_{M(B)}\end{pmatrix}
\]
define double centralizers of $L(X)$, so that $p,q \in M(L(X))$, and
noting that, for example,
\[
A= \begin{pmatrix}
A&0
\\
0&0\end{pmatrix}=pL(X)p;
\]
more formally, the inclusion $a\mapsto
\left(\begin{smallmatrix} a&0
\\
0&0\end{smallmatrix}\right)$ is an isomorphism of $A$ onto the corner
$pL(X)p$.

The above construction has an obvious inverse:

\begin{prop}\label{prop-partiallink}
     Suppose that $L$ is a $C^*$-algebra and $p,q\in M(L)$ are
     complementary projections \textup(\ie, $p+q=1$\textup). Then
     $pLq$ is a partial $pLp-qLq$ imprimitivity bimodule
     with operations given by multiplication and involution
     on $L$ in the canonical way \textup(\eg, the
     $qLq$-valued inner product on $pLq$ is given by
     $\lk a,b\rk_{qLq}= a^*b$\textup).
     Moreover, $pLq$ is a right-partial $pLp-qLq$ imprimitivity bimodule
     if and only if $q$ is full \textup(\ie, $\overline{LqL}=L$\textup),
and $pLq$ is a left-partial $pLp-qLq$ imprimitivity bimodule if and
only if $p$ is full. 
   In particular, $pLq$ is a $pLp-qLq$ imprimitivity
   bimodule if and only if both $p$ and $q$ are full.
\end{prop}
\begin{proof} It is straightforward to check that
     $pLq$ is a partial $pLp-qLq$ imprimitivity bimodule
     with respect to the canonical operations.
     If $q$ is full, then we get
     $${_{pLp}}\overline{\lk pLq, pLq\rk}= \overline{pLq(pLq)^*}
     =\overline{pLqLp}=pLp,$$
     which implies that $pLq$ is a right-partial
     $pLp-qLq$ imprimitivity bimodule.
     For the converse, we first observe that for any projection
     $q\in M(L)$ one has the equations
     ${LqLq}=Lq$ and (hence) $qLqL=qL$. This follows from the fact
     that $Lq$ is a Hilbert $qLq$-module, which implies
     that $qLq$ acts nondegenerately on $Lq$ by multiplication.
     Thus, if $\overline{pLqLp}={_{pLp}}\overline{\lk pLq, pLq\rk}
     =pLp$, then a short computation reveals that
     $$\overline{LqL}=\overline{pLqLp}+pLqLq+qLqLp+qLqLq=
     pLp+pLq+qLp+qLq=L,$$
     so $q$ is a full projection. Applying the same arguments
     to $p$ completes the proof.
\end{proof}

In practice it
is crucial to be able to recognize
under what conditions a given algebra $L$ with
given complementary projections $p,q\in M(L)$
is isomorphic to the linking algebra of a
given partial imprimitivity bimodule:

\begin{prop}\label{recogniselink}
Let ${}_AX_B$ be a partial imprimitivity bimodule.  Suppose that
$(L,p,q,{}_\phi\Phi_\psi)$ consists of a $C^*$-algebra $L$,
complementary projections $p,q$ in $M(L)$, and a partial imprimitivity
bimodule isomorphism ${}_\phi\Phi_\psi\:{}_AX_B\to
{}_{pLp}(pLq)_{qLq}$ \textup(with the obvious meaning\textup).  Then
\[
\theta\begin{pmatrix}
a&x
\\
\widetilde z &b\end{pmatrix}
= \phi(a)+\Phi(x)+\Phi(z)^*+\psi(b)
\]
defines an isomorphism $\theta\:L(X)\iso L$.
\end{prop}

\begin{proof}
Routine calculations show that $\theta$ is a homomorphism.  It is
injective because each of the components is; it is surjective because
every $d \in L$ can be decomposed as $d=pdp+pdq+qdp+qdq$, and because
the maps $\phi$, $\Phi$, and $\psi$ are surjective.
\end{proof}

\begin{rem}\label{L-rem}
We shall often write $L=L(X)$ to summarize an application of
this proposition; this means that we think there are obvious
candidates for the projections $p,q$ and the maps $\phi$, $\Phi$,
and $\psi$, and that these candidates satisfy the hypotheses of
Proposition~\ref{recogniselink}.  Thus, for example, if $C$ is any
other $C^*$-algebra, we write $L(X)\otimes C=L(X\otimes C)$ to mean that
the quadruple
$$\big(L(X)\otimes C,p_{L(X)}\otimes 1_{M(C)},q_{L(X)}\otimes
1_{M(C)},{}_\phi\Phi_\psi\otimes{\id}\big)
$$
satisfies the hypotheses of
Proposition~\ref{recogniselink} for the external tensor product $X\otimes
C$, which is a partial $(A\otimes C)-(B\otimes C)$ imprimitivity bimodule.
This leads to statements like $A\otimes C=p(L(X)\otimes C)p$, by which
we mean that the obvious map $\phi\otimes {\id}$ is an isomorphism of
$A\otimes C$ onto the given corner.
\end{rem}

For any right-partial $A-B$ imprimitivity bimodule $X$,
we saw in \remref{rem-imp-multiplier} that $M(X)$ is a
partial
$M(A)-M(B)$ imprimitivity bimodule, so we can also
form the linking algebra $L(M(X))=\left(\begin{smallmatrix}M(A)&M(X)
\\
\widetilde{M( X)}&M(B)\end{smallmatrix}\right)$.
It follows from
the properties of $M(X)$ as listed in
Proposition \ref{prop-strictlydense}
that $L(X)$ is a closed ideal in $L(M(X))$.
Thus, by the universal properties of
the multiplier algebras, there exists a unique
algebra homomorphism $\Phi\:L(M(X))\to M(L(X))$
which extends the identity on $L(X)$.
In fact, this map is always an isomorphism:

\begin{prop}
\label{multsoflink}
Suppose that $X$ is a right-partial $A-B$ imprimitivity bimodule.
Then the canonical homomorphism $\Phi\:L(M(X))\to M(L(X))$
is an isomorphism of $C^*$-algebras.
Moreover, $\Phi$ is a homeomorphism with respect to
the topology on
$L(M(X))=\left(\begin{smallmatrix}M(A)&M(X)
\\
M(\widetilde X)&M(B)\end{smallmatrix}\right)$ which is the product
of the strict topologies on the corners and
the strict topology on $M(L(X))$.
\end{prop}
\begin{proof}
Using the separate strict continuity of all
bimodule pairings of $_AX_B$
(compare with Proposition \ref{prop-strictlydense})
it is straightforward to check that the product of the
strict topologies on the corners of $L(X)$, viewed as
subspaces of the corners of $L(M(X))$, coincides with
the strict topology of $L(X)$, viewed as a subalgebra of $M(L(X))$.
It follows that the identity map on $L(X)$ extends to
a product-strict to strict
linear homeomorphism $\bar{\Phi}$ between
the corresponding completions $L(M(X))$ (see
Proposition \ref{prop-strictlydense})
and $M(L(X))$. Since all algebra operations are separately
continuous
with respect to these topologies (which in the case of
$L(M(X))$  also follows from Proposition \ref{prop-strictlydense}),
it follows that $\bar\Phi$ is an algebra homomorphism
which extends the identity on $L(X)$. Thus $\bar\Phi=\Phi$
by the uniqueness of $\Phi$.
\end{proof}

Just as we are going to simplify notation by writing $L=L(X)$ when we
believe it is obvious what the projections $p,q$ and the bimodule
isomorphism ${}_\phi\Phi_\psi$ are meant to be, we shall write things
like $qM(L)q=D$ when we think there is an obvious candidate for $q$
and for the embedding of $D$ in $M(L)$.  Thus, for example, the
identification $L(X)\otimes C=L(X\otimes C)$ induces an identification
$M(B\otimes C)=qM(L(X\otimes C))q$ which means, strictly speaking,
that the isomorphism $\psi\otimes\id$ of $B\otimes C$ onto the corner
$(q\otimes 1)(L(X\otimes C))(q\otimes 1)$ induces, via
Proposition~\ref{multsoflink}, an isomorphism of $M(B\otimes C)$ onto
$(q\otimes 1)M(L(X\otimes C))(q\otimes 1)$.

We close this section with a slight extension of
  \cite[Remark (2) on p. 307]{er:mult}:

\begin{lem}\label{link-nodeg}
     Assume that $_AX_B$  and
     $_CY_D$ are right-partial imprimitivity bimodules.
     Let $_\phi\Phi_\psi\:X\to M(Y)$ be a partial
     imprimitivity bimodule
     homomorphism. Then \textup(using the identification
     $M(L(Y))=L(M(Y))$ given in Proposition~\textup{\ref{multsoflink}}\textup)
the formula
     $$\Psi\left(\smtx{a&x\\ \tilde{y}& b}\right)=
     \mtx{\phi(a)&\Phi(x)\\ \widetilde{\Phi(y)}&\psi(b)}$$
defines a 
     $*$-homomorphism $\Psi\:L(X)\to M(L(Y))$, 
and $\Psi$ is nondegenerate if $\Phi$
     is nondegenerate. The
     extension $\bar{\Psi}\:M(L(X))\to M(L(Y))$ is given
     by the extensions of the corner maps.

     Conversely, assume that $\Psi\:L(X)\to M(L(Y))$ is a
     nondegenerate $*$-homomorphism
     such that $\Psi(p_Xl)=p_Y\Psi(l)$ and $\Psi(q_Xl)=q_Y\Psi(l)$, where
     $p_X,q_X$ and $p_Y, q_Y$ denote the respective corner projections.
     Then $\Psi$ determines a
     partial imprimitivity bimodule homomorphism
     $_{\phi}\Phi_\psi\:X\to M(Y)$
     by defining
     $$\phi(p_Xlp_X)=p_Y\Psi(l)p_Y,\quad
     \Phi(p_Xlq_X)=p_Y\Psi(l)q_Y,\quad\text{and}\quad
     \psi(q_Xlq_X)=q_Y\Psi(l)q_Y$$
     for $l\in L(X)$, and $\Phi$ is nondegenerate if $\Psi$
     is nondegenerate.
\end{lem}
\begin{proof} It follows immediately from the algebraic properties
     of the linking algebras that the above procedure gives
     a correspondence between partial imprimitivity bimodule
     homomorphisms of $X$ into $M(Y)$ and $*$-homomorphisms
     of $L(X)$ into $M(L(Y))$. So we only have to check
     that this correspondence preserves nondegeneracy.

     If $_{\phi}\Phi_{\psi}$ is nondegenerate, then
     $$\Psi(L(X))L(Y)=\mtx{\phi(A)C+{_A\lk \Phi(X),Y\rk}
     &\Phi(A)Y+\Phi(X)D\\
     \widetilde{\Phi(X)}C+\psi(B)\widetilde{Y}&\lk\Phi(X),Y\rk_B+\psi(B)D}.
     $$
     Since $\phi$ and $\psi$ are nondegenerate, it follows that
     $\phi(A)C=C$, $\phi(A)Y=Y$, $Y\psi(B)=Y$, and $\psi(B)D=D$,
     which implies that $\Psi(L(X))L(Y)=L(Y)$.

     Now assume conversely that $\Psi(L(X))\to M(L(Y))$ is nondegenerate.
     Using the general equation $pL(X)pL(X)=pL(X)$
     (see the proof of \lemref{prop-partiallink})
     we get:
     \begin{align*}
	\phi(A)C&=p\Psi(L(X))pL(Y)p=p\Psi(L(X))p\Psi(L(X))L(Y)p\\
	&=\Psi(pL(X)pL(X))L(Y)p=\Psi(pL(X))L(Y)p=p\Psi(L(X))L(Y)p\\
	&=pL(Y)p=C.
     \end{align*}
     (Here we've omitted the subscripts on the $p$'s and $q$'s.)
     Hence $\phi$ is nondegenerate, and a similar computation shows that
     $\psi$ is nondegenerate. Finally, using the fullness of $q_X$
     (since the left inner product on $X$ is full),
     we get
     \begin{align*}
	\overline{\Phi(X)D}&=\overline{p\Psi(L(X))qL(Y)q}=
	\overline{p\Psi(L(X))q\Psi(L(X))L(Y)q}\\
	&=\overline{p\Psi(L(X)qL(X))L(Y)q}=p\Psi(L(X))L(Y)q=Y.
     \end{align*}
     This finishes the proof.
\end{proof}

\begin{rem}\label{remdegen}
     (1) Note that if $_AX_B$ and $_CY_D$ are right-partial imprimitivity
     bimodules, then it follows from \lemref{lem-partialhom} that
     the nondegenerate right-Hilbert bimodule homomorphisms
     of $X$ into $M(Y)$ are automatically nondegenerate
     partial imprimitivity bimodule homomorphisms (and \emph{vice versa}).
     Thus the above lemma gives a one-to-one correspondence
     between the nondegenerate right-Hilbert bimodule homomorphisms
     of $X$ into $M(Y)$ and the nondegenerate $*$-homomorphisms of
$L(X)$ into $M(L(Y))$.

     (2) The proof of the lemma implies the interesting observation
     that a partial imprimitivity bimodule homomorphism
     $_{\phi}\Phi_{\psi}\:{_AX_B}\to M({_CY_D})$
     between right-partial imprimitivity bimodules is nondegenerate
     if (and only if) the coefficient homomorphisms
     $\phi\:A\to M(C)$ and $\psi\:B\to M(D)$ are nondegenerate,
     since this was all we needed for the proof of the nondegeneracy
     of the corresponding homomorphism of the linking algebra.
     For imprimitivity bimodules this was already observed
     in \remref{rHb-hom-rem1}. As remarked there, for this to be true
     it is necessary that $_{\phi}\Phi_{\psi}$ preserve
     the left inner products.
\end{rem}

%
%

\chapter{The Categories}
\label{categories-chap}

In this chapter we show that there exists a category $\c C$ in which
the objects are $C^*$-algebras, and the morphisms from $A$ to $B$ are
the isomorphism classes of right-Hilbert $A-B$ bimodules.
For any
locally compact group $G$, there are also equivariant categories $\c
A(G)$, $\c C(G)$, and $\c A \c C(G)$ which combine, respectively,
actions of $G$, coactions of $G$, and both, with the structure of $\c
C$. Note that our category differs from the one
considered in \cite{ekqr:green} by allowing
isomorphism classes of non-full
right-Hilbert modules to be morphisms in $\C$.

\section{${C}^*$-Algebras}

Let $X$ and $Y$ be right-Hilbert $A-B$ bimodules.
Recall from the preceding chapter
that a \emph{right-Hilbert $A-B$ bimodule isomorphism} of $X$ onto $Y$
is a bijective right-Hilbert bimodule homomorphism $\Phi\:X\to Y$ whose
coefficient maps are $\id_A$ and $\id_B$.  From now on we abuse the
terminology and simply say $X$ and $Y$ are \emph{isomorphic} if there
exists a right-Hilbert $A-B$ bimodule homomorphism of $X$ onto $Y$.
It is not hard to check that this notion of isomorphism is an equivalence
relation on the class of right-Hilbert $A-B$ bimodules.

\begin{rem}
When we began writing this, the idea of a category in
which the morphisms come from bimodules seemed new.  In the intervening
years, however, this category and close relatives of it have
been independently discovered by several others.
(See, for instance, \cite{LanBO, LanQR} and 
\cite{sch, sch:cro, sch:dil}.)
It also turns out that this category has been at least 
implicitly in the air for quite some time; for example, \cite{bgr} 
briefly mentions the category of $C^*$-algebras and (isomorphism 
classes of) imprimitivity bimodules.
\end{rem}

\begin{thm}
\label{C-cat-thm}
There is a category $\c C$ in which the objects are $C^*$-algebras,
and in which the morphisms from $A$ to $B$ are the isomorphism classes
of right-Hilbert $A-B$ bimodules.
The composition of $[X]\:A\to
B$ with $[Y]\:B\to C$ is the isomorphism class of the balanced tensor
product $X\otimes_BY$\textup; the identity morphism on $A$ is the isomorphism
class of the standard right-Hilbert bimodule ${}_AA_A$.
\end{thm}

\begin{proof}
We first note that the composition of morphisms is well-defined:
suppose $[X] = [X']\:A\to B$ and $[Y] = [Y']\:B\to C$, so that
we have right-Hilbert $A-B$ bimodule isomorphisms
\[
\Phi\:X\to X'
\midtext{and}
\Psi\:Y\to Y'.
\]
Then the tensor product homomorphism $\Phi\otimes_B \Psi$ maps $X
\otimes_B Y$ into $X'\otimes_B Y'$ and has inverse $\Phi^{-1}
\otimes_B \Psi^{-1}$, so gives an
isomorphism $X\otimes_B Y\cong X'\otimes_B Y'$.

Next we establish that composition of morphisms in $\c C$ is
associative; by the above, it suffices to show that $X\otimes_B (Y
\otimes_C Z)$ and $(X\otimes_B Y)\otimes_C Z$ are isomorphic for any
right-Hilbert bimodules ${}_AX_B$, ${}_BY_C$, and ${}_CZ_D$. But
straightforward calculations show that the usual linear isomorphism of
$X \odot (Y \odot Z)$ onto $(X \odot Y) \odot Z$ respects the module
actions and right inner products, so extends to the desired
isomorphism.

Finally, note that $A\otimes_A X \cong X$ and $Y\otimes_A A \cong Y$
for any right-Hilbert bimodules ${}_AX_B$ and ${}_BY_A$. Hence the
identity morphism from $A$ to $A$ is given by the isomorphism class of
the standard bimodule ${}_AA_A$.
\end{proof}

\begin{rem}
In any category, $\mor(A,B)$ is required to be a \emph{set}
(not merely a class)
for each
pair of objects $A$ and $B$.  This will fail in $\C$ unless we limit
the size of
the bimodules involved.  We can do this by considering only
$C^*$-algebras and Hilbert modules with dense subsets whose
cardinalities do not exceed a fixed large cardinal.  For example, we
could consider only separable $C^*$-algebras and bimodules.
Alternatively, for each $A,B$ we could restrict attention to
right-Hilbert $A-B$ bimodules with cardinality dominated by the
larger of the cardinalities of $A$ and $B$ (which would accommodate all
the bimodules which occur in the usual imprimitivity theorems, for
example).  In practice, these issues should never present a
real problem, and we shall ignore them.
\end{rem}

Recall that in any category, a morphism $f\:A\to B$ is called an
\emph{equivalence} if there exists $g\:B\to A$ such that $f\circ g=\id_B$ and
$g\circ f=\id_A$.  The equivalences in $\c C$ are exactly the
(isomorphism classes of) imprimitivity bimodules.  In one direction,
if $X$ is an $A-B$ imprimitivity bimodule, then the isomorphism class
of the conjugate module $\rev X$ is in $\mor(B,A)$, and satisfies
$[X]\circ[\rev X]=[B]$ and $[\rev X]\circ[X]=[A]$.  To see the
converse requires substantially more work.  Schweizer
independently discovered the following result in 
\cite[Proposition 2.3]{sch}; 
our proof is considerably different.%
\footnote{It is also different from the representation-theoretic proof
used in \cite{ekqr:green}, where this theorem was proved for full
right-Hilbert bimodules.}

\begin{lem}
\label{C-ibs-lem}
Let ${}_AX_B$ and ${}_BY_A$ be right-Hilbert bimodules such that
\[
X\otimes_B Y\cong A
\midtext{and}
Y\otimes_A X\cong B.
\]
Then $X$ is an $A-B$ imprimitivity bimodule, and $Y\cong\rev{X}$.
\end{lem}

\begin{proof}
We first note that
$Y\otimes_A X\cong B$ implies that $X$ is full, since
it follows from the definition of the
inner product on $Y\otimes_AX$ that its range is always
contained in the range of the inner product on  $X$.
A similar argument shows that $Y$ is full, too.
It therefore suffices to show that the canonical homomorphism $\phi\:A\to \c
L_B(X)$ is an isomorphism of $A$ onto $\c K_B(X)$, for this shows $X$
is an $A-B$ imprimitivity bimodule, and then the second statement of
the lemma will follow from uniqueness of inverses in a category.  We
first show that $\phi$ is faithful.  Since $X\otimes_BY \cong A$, we
know that the canonical homomorphism of $A$ into $\c L_A(X\otimes_BY)$
is faithful (because $A\to\c L_A(A)$ is).  This implies $\phi$ is
faithful, for if $a\d x = 0$ for all $x\in X$, then $a\d (x\otimes y)
= 0$ for all $x\in X, y\in Y$, so $a\d z = 0$ for all $z\in
X\otimes_BY$, hence $a = 0$.

To see that $\phi(A)=\c K_B(X)$, note that the latter coincides
with $\overline{XX^*}$, using the canonical identification of $X$ with
$\c K_B(B,X)$.  Our strategy is to show that there is a right-Hilbert
bimodule homomorphism
\[
\Phi\:{}_BY_A\to \c L_B(X,B)
= M\bigl({}_B(\c K_B(X,B))_{\c K(X)}\bigr)
\]
such that the right coefficient map of $\Phi$ is $\phi$, and
$\Phi(Y)^*=X$. This will do the job, since we will then have
\[
\overline{XX^*}=\overline{\Phi(Y)^*\Phi(Y)}=\phi(A).
\]
Let $\Psi\:{}_B(Y\otimes_AX)_B\iso {}_BB_B$, and define
$\Phi\:Y\to \c L_B(X,B)$ by
\[
\Phi(y)x=\Psi(y\otimes x).
\]
We have
\[
\Phi(b\d y)x
= \Psi(b\d y\otimes x)
= b\Psi(y\otimes x)
= b\Phi(y) x
\righttext{for all }b\in B,\ x\in X,\ y\in Y,
\]
and
\begin{align*}
\bigl\<\<\Phi(y_1),\Phi(y_2)\>_{\c K(X)}\,x_1, x_2\bigr\>_B
&= \bigl\<\Phi(y_1)^*\Phi(y_2)x_1,x_2\bigr\>_B
= \bigl\<\Phi(y_2)x_1,\Phi(y_1)x_2\bigr\>_B
\\&= \bigl\<\Psi(y_2\otimes x_1),\Psi(y_1\otimes x_2)\bigr\>_B
= \<y_2\otimes x_1,y_1\otimes x_2\>_B
\\&= \bigl\<\<y_1,y_2\>_A\d\,x_1,x_2\bigr\>_B
= \bigl<\phi\bigl(\<y_1,y_2\>_A\bigr)x_1,x_2\bigr\>_B,
\end{align*}
for all $x_i\in X$ and $y_i\in Y$,
which implies that $\Phi$ is an $\id_B-\phi$ compatible
right-Hilbert bimodule homomorphism.

Finally,
\begin{align*}
X
&= A\d X
= \overline{\Phi(Y)^* \Phi(Y) X}
= \overline{\Phi(Y)^* \Psi(Y\otimes_A X)}
= \overline{\Phi(Y)^* B}
= \Phi(B\d Y)^*
= \Phi(Y)^*,
\end{align*}
where $B$ is identified with $\c K_B(B_B)$ where appropriate.
\end{proof}

\section{Group actions}

\begin{defn}
\label{rHb-act-defn}
Let $G$ be a locally compact group, let $\alpha$ and $\beta$ be actions
of $G$ on $C^*$-algebras $A$ and $B$, and let $X$ be a right-Hilbert $A -
B$ bimodule. An \emph{$\alpha-\beta$ compatible right-Hilbert bimodule
action} of $G$ on $X$ is a homomorphism $\gamma$ of $G$ into the group
of invertible linear maps on $X$ such that
\begin{enumerate}
\item
$\gamma_s(a\d x) = \alpha_s(a)\d \gamma_s(x)$

\item
$\gamma_s(x\d b) = \gamma_s(x)\d \beta_s(b)$

\item
$\<\gamma_s(x),\gamma_s(y)\>_B =
\beta_s(\<x,y\>_B)$
\end{enumerate}
for each $s\in G$, $a\in A$, $x,y\in X$, and $b\in B$; and such that
each map $s\mapsto \gamma_s(x)$ is continuous from $G$ into $X$.
We call $\alpha$ and $\beta$ the \emph{coefficient actions} of $\gamma$.
\end{defn}

\begin{rem}
\label{rHb-act-rem}
(1)
Note that 
each $\gamma_s$ is in particular a
right-Hilbert bimodule homomorphism of $X$ onto itself with coefficient
maps $\alpha_s$ and $\beta_s$; thus condition (ii) is implied by condition
(iii) and the linearity of $\beta_s$ by \remref{rHb-hom-rem}.

(2)
If $X$ is an $A-B$ imprimitivity bimodule, and if $\gamma$ is an
$\alpha-\beta$ compatible right-Hilbert bimodule action on $X$, then
$\gamma$ is automatically an imprimitivity bimodule action in the sense
of \cite{com}: calculating as in the proof of
\lemref{lem-partialhom}, for each $s
\in G$ and $x,y,z\in X$ we have
\begin{align*}
\alpha_s({}_A\<x,y\>)\d \gamma_s(z)
&= \gamma_s({}_A\<x,y\>\d z)
\\&= \gamma_s(x\d \<y,z\>_B)
\\&= \gamma_s(x)\d \<\gamma_s(y), \gamma_s(z)\>_B
\\&= {}_A\<\gamma_s(x), \gamma_s(y)\>\d \gamma_s(z),
\end{align*}
which shows that $\alpha_s({}_A\<x,y\>) =
{}_A\<\gamma_s(x), \gamma_s(y)\>$.

By \remref{rHb-hom-rem1}, we have a sort of converse: when $X$ is an
imprimitivity bimodule we can replace (i) in \defnref{rHb-act-defn}
by
\begin{enumerate}
\item[(i${}'$)]
${}_{A}\<\gamma_s(x),\gamma_s(y)\>=\alpha_s({}_A\<x,y\>)$.
\end{enumerate}

(3) We should point out that group actions on Hilbert bimodules as introduced
above are well known in the literature. In particular they
play an important r\^ole in the construction of
Kasparov's equivariant $KK$-Theory for $C^*$-algebras (see
also \cite{kas}).
\end{rem}

Given right-Hilbert bimodule actions
${}_{(A,\alpha)}(X,\gamma)_{(B,\beta)}$ and
${}_{(B,\beta)}(Y,\rho)_{(C,\epsilon)}$, it is easy to check that the
automorphisms of $X\otimes_B Y$ defined by
\[
(\gamma\otimes_B \rho)_s = \gamma_s\otimes_B \rho_s
\]
for each $s\in G$ give rise to an $\alpha-\epsilon$ compatible
action $\gamma\otimes_B \rho$ of $G$ on $X\otimes_B Y$.

\begin{defn}
\label{AG-isom-defn}
Let $\gamma$ and $\rho$ be $\alpha-\beta$ compatible actions
of $G$ on right-Hilbert $A-B$ bimodules $X$ and $Y$. An
isomorphism $\Phi$ of $X$ onto $Y$ is
\emph{$\gamma-\rho$ equivariant}, or \emph{intertwines}
$\gamma$ and $\rho$, if
\[
\Phi\circ \gamma_s = \rho_s\circ \Phi
\]
for all $s\in G$. We say $\gamma$ and $\rho$ are
\emph{isomorphic}, or $X$ and $Y$ are \emph{equivariantly isomorphic},
if such a $\Phi$ exists.
\end{defn}

It is straightforward to check that
this notion of isomorphism is an equivalence relation on
the class of right-Hilbert $A-B$ bimodules with $\alpha-\beta$
compatible actions of $G$.

\begin{thm}
\label{AG-cat-thm}
Let $G$ be a locally compact group. There is a category $\c A(G)$ in
which the objects are $C^*$-algebras with actions of $G$, and in which
the morphisms from $(A,\alpha)$ to $(B,\beta)$ are the equivariant
isomorphism classes of 
right-Hilbert $A-B$ bimodules with
$\alpha-\beta$ compatible actions of $G$. The composition of
$[X,\gamma]\:(A,\alpha)\to (B,\beta)$ with $[Y,\rho] \:
(B,\beta)\to (C,\nu)$,
is the isomorphism class of the tensor product action $(X
\otimes_B Y,\gamma\otimes_B \rho)$\textup; the identity morphism on
$(A,\alpha)$ is the isomorphism class of the $\alpha-\alpha$
compatible right-Hilbert bimodule action $(A,\alpha)$ itself.
\end{thm}

\begin{proof}
We first note that the composition of morphisms is well-defined:
suppose $[X,\gamma] = [X',\gamma']\:(A,\alpha)\to (B,\beta)$
and $[Y,\rho] = [Y',\rho']\:(B,\beta)\to (C,\nu)$,
so that we have equivariant right-Hilbert bimodule
isomorphisms
\[
\Phi\:(X,\gamma)\to (X',\gamma')
\midtext{and}
\Psi\:(Y,\rho)\to (Y',\rho').
\]
Then straightforward calculations show that the isomorphism
$\Phi\otimes_B \Psi\:X\otimes_B Y\to X' \otimes_B Y'$ from the proof
of \thmref{C-cat-thm} satisfies
\[
(\Phi\otimes_B \Psi)\circ (\gamma\otimes_B \rho)_s =
(\gamma'\otimes_B \rho')_s\circ (\Phi\otimes_B \Psi)
\]
for all $s\in G$, hence gives an isomorphism between the $\alpha
- \nu$ compatible actions $\gamma\otimes_B \rho$ and
$\gamma'\otimes_B \rho'$.

Next we establish that composition of morphisms in $\c A(G)$ is
associative; by the above, it suffices to show that the actions
$\gamma\otimes_B (\rho\otimes_C \sigma)$ and $(\gamma\otimes_B \rho)
\otimes_C \sigma$ are isomorphic for any right-Hilbert bimodule
actions $({}_AX_B,\gamma)$, $({}_BY_C,\rho)$, and $({}_CZ_D,\sigma)$.
Again, straightforward calculations show that the isomorphism
$X\otimes_B (Y\otimes_C Z) \cong (X\otimes_B Y)\otimes_C Z$ from the
proof of \thmref{C-cat-thm}\ intertwines the actions as desired.

Finally, note that for any action $({}_AX_B,\gamma)$, with left
coefficient action $\alpha$, the canonical isomorphism $A\otimes_A X
\cong X$ intertwines the actions $\alpha\otimes_A \gamma$ and
$\gamma$.  Similarly, for any action $({}_BY_A,\rho)$ with right
coefficient action $\alpha$, the canonical isomorphism $Y\otimes_A A
\cong Y$ intertwines $\rho\otimes_A \alpha$ and $\rho$.  Hence the
identity morphism on $(A,\alpha)$ is the isomorphism class of the
$\alpha - \alpha$ compatible right-Hilbert bimodule action
$(A,\alpha)$.
\end{proof}

\begin{rem}
\label{AG-ibs-rem}
The equivalences in $\c A(G)$ are exactly the (isomorphism
classes of) imprimitivity bimodule actions of $G$ in the
sense of \cite{com} or \cite{CMW-CP}. In one direction, if
${}_{(A,\alpha)}(X,\gamma)_{(B,\beta)}$ is an imprimitivity bimodule
action, then the canonical isomorphisms $X\otimes_B \rev X \cong A$ and
$\rev X\otimes_A X \cong B$ are easily seen to be $\gamma\otimes_B
\rev \gamma-\alpha$ and $\rev \gamma\otimes_A \gamma-\beta$
equivariant, respectively. In the other direction, if $({}_AX_B,\gamma)$
and $({}_BY_A,\rho)$ are $\alpha-\beta$ and $\beta-\alpha$ compatible
right-Hilbert bimodule actions with
\[
\gamma\otimes_B \rho \cong \alpha
\midtext{and}
\rho\otimes_A \gamma \cong \beta,
\]
then $X$ is in particular an $A-B$ imprimitivity bimodule by
\lemref{C-ibs-lem}, so $\gamma$ and $\rho$ are imprimitivity
bimodule actions by \remref{rHb-act-rem}.
\end{rem}

\section{Group coactions}

In this section we are going to construct a
category $\c C(G)$, where all actions and morphisms
are equipped with coactions of the group $G$.
The necessary background on coactions of groups on
$C^*$-algebras is given in Appendix~\ref{coactions-chap}.

For any right-Hilbert $A-B$ bimodule $X$ and any locally compact group
$G$, we may take the exterior tensor product of $X$ with the
right-Hilbert $C^*(G)-C^*(G)$ module $C^*(G)$, as in
\defnref{xtnsr-defn}, and get a right-Hilbert $(A\otimes
C^*(G))-(B\otimes C^*(G))$ bimodule $X\otimes C^*(G)$.

\begin{defn}
\label{rHb-co-defn}
Let $G$ be a locally compact group, let $\delta$ and $\epsilon$ be
coactions of $G$ on $C^*$-algebras $A$ and $B$, and let $X$ be a
right-Hilbert $A-B$ bimodule.  A \emph{$\delta-\epsilon$
compatible right-Hilbert bimodule coaction} of $G$ on $X$ is a
nondegenerate%
\footnote{Notice that we have incorporated
nondegeneracy of $\zeta$ \emph{as a bimodule homomorphism} into our
definition, whereas Ng's
definition \cite[Definition~2.10]{ng:module} does not.}
right-Hilbert bimodule homomorphism
$\zeta\:X\to M(X \otimes C^*(G))$,%
with coefficient maps
$\delta$ and $\epsilon$ (called the \emph{coefficient coactions} of
$\zeta$), such that
\begin{enumerate}
     \item $(1_{M(A)}\otimes C^*(G))\zeta(X)\subseteq X\otimes
     C^*(G)$, and
\item $(\zeta\otimes \id_G)\circ \zeta = (\id_X\otimes\delta_G)
\circ \zeta$ (the \emph{coaction identity}).
\end{enumerate}
Moreover, $\zeta$ is called {\em nondegenerate} if $\delta$ and
$\epsilon$ are nondegenerate coactions and
\[
\overline{(1\otimes C^*(G))\zeta(X)}=X\otimes C^*(G).
\]
\end{defn}

\begin{rem}
\label{rHb-co-rem}
(1)
A condition (as in \cite{er:mult})
\begin{enumerate}
\item
[(iii)]
$\zeta(X)(1_{M(B)}\otimes C^*(G))\subseteq X\otimes C^*(G)$,
\end{enumerate}
symmetric to (i), would be redundant: it follows from the analogous
property for the $C^*$-coaction $\epsilon$, because
\begin{align*}
\zeta(X)(1\otimes C^*(G))
&=\zeta(XB)(1\otimes C^*(G))
=\zeta(X)\epsilon(B)(1\otimes C^*(G))
\\&\subseteq\zeta(X)(B\otimes C^*(G))
\subseteq X\otimes C^*(G).
\end{align*}
It follows from this that condition (i) of \defnref{rHb-co-defn}
is
equivalent to the requirement that $\zeta(X)\subseteq M_G(X\otimes C^*(G))$
(see \defnref{defn-Cmultmod}).

(2)
If $\zeta$ is
nondegenerate, then
it automatically
satisfies
\[
\overline{\zeta(X)\bigl(1\otimes C^*(G)\bigr)}
= X\otimes C^*(G).
\]
This is easily checked by using the right module homomorphism property and
nondegeneracy of $\epsilon$.

(3)
If $X$ is an $A-B$ imprimitivity bimodule, and if $\zeta$ is a $\delta -
\epsilon$ compatible right-Hilbert bimodule coaction on $X$, then $\zeta$
is automatically an imprimitivity bimodule coaction in the sense of
\cite[Definition~3.1]{er:mult} (\cf\ \cite[2.2]{BS-CH},\cite[2.15]{BuiME},
and \cite[Definition~3.3]{ng:module}).  We only need to check that
\[
\delta({}_A\<x,y\>) = {}_{M(A\otimes C^*(G))}\<\zeta(x), \zeta(y)\>
\]
for all $x,y\in X$; but this is immediate from \lemref{lem-partialhom}
(in fact, for this to be true we only need to assume that
$X$ is a right-partial imprimitivity bimodule).
Moreover, if $X$ is an imprimitivity bimodule, then
condition~(i) in \defnref{rHb-co-defn} is redundant, 
nondegeneracy of $\zeta$ as a bimodule homomorphism follows
automatically from nondegeneracy of the coefficient homomorphisms, and
nondegeneracy of $\zeta$ as a coaction follows automatically from 
nondegeneracy of the coefficient coactions.
\end{rem}

The following lemma will be fundamental in all computations with
bimodule coactions.

\begin{lem}
\label{Theta-lem}
Let ${}_AX_B$ and ${}_BY_C$ be
right-Hilbert bimodules, and let $D$ be a $C^*$-algebra.
There exists a right-Hilbert $A\otimes D - C\otimes D$ bimodule
isomorphism%
\[
\Theta\:(X\otimes D)\otimes_{B\otimes D}(Y\otimes
D)\to (X\otimes_B Y)\otimes D
\]
such that
\begin{equation}
\label{Theta-eq}
\Theta \bigl((x\otimes d)\otimes_{B\otimes D} (y\otimes e)\bigr)
= (x\otimes_B y)\otimes de
\end{equation}
for $x\in X$, $y\in Y$, and $d,e\in D$.%
\footnote{More generally, one can prove that $({}_AX_B\otimes {}_DZ_E)
\otimes_{B\otimes E}({}_BY_C\otimes {}_EW_F)\cong (X\otimes_B Y)
\otimes(Z\otimes_E W)$ as right-Hilbert $A\otimes D-C\otimes F$
bimodules.}

Moreover, the unique extension of $\Theta$ to the multiplier bimodules
satisfies\footnote{It is fairly obvious that the following identity
can be generalized to allow $x$ and $y$ to be multipliers as well, but
we only need the specific fact we record here.}
\[
\Theta \bigl((x\otimes m)\otimes(y\otimes n)\bigr)
= (x\otimes y)\otimes mn
\righttext{for} x\in X, y\in Y, m, n\in M(D).
\]
\end{lem}

\begin{proof}
Equation \eqref{Theta-eq} clearly determines a map $\Theta \:
(X \odot D) \odot (Y \odot D)\to (X\otimes_B Y)\otimes D$.  The
following computation implies that $\Theta$ preserves inner products:
\begin{align*}
\bigl\<
\Theta \bigl((x\otimes d)\otimes(y\otimes e)\bigr),&
\Theta \bigl((z\otimes f)\otimes(w\otimes g)\bigr)
\bigr\>_{C\otimes D}
\\&\quad= \bigl\< (x\otimes y)\otimes de,
(z\otimes w)\otimes fg\bigr\>_{C\otimes D}
\\&\quad= \< x\otimes y, z\otimes w \>_C
\otimes (de)^*fg
\\&\quad= \bigl\< y, \<x,z\>_B\d w\bigr\>_C
\otimes e^*d^*fg
\\&\quad= \bigl\< y\otimes e, \<x,z\>_B
\d w\otimes d^*fg\bigr\>_{C\otimes D}
\\&\quad= \bigl\< y\otimes e, (\<x,z\>_B\otimes (d^*f))
\d (w\otimes g)\bigr\>_{C\otimes D}
\\&\quad= \bigl\< y\otimes e,
\< x\otimes d, z\otimes f\>_{B\otimes D}
\d (w\otimes g)\bigr\>_{C\otimes D}
\\&\quad= \bigl\< (x\otimes d)\otimes(y\otimes e),
(z\otimes f)\otimes (w\otimes g)\bigr\>_{C\otimes D}.
\end{align*}
Since $\Theta$ clearly has dense range in $(X\otimes_B Y)\otimes D$,
it therefore extends to an isometry, which we continue to denote by
$\Theta$, of $(X\otimes D)\otimes_{B\otimes D}(Y\otimes D)$ onto
$(X\otimes_B Y)\otimes D$.  Straightforward calculations with
elementary tensors verify that $\Theta$ intertwines the left actions,
so that $\Theta$ is indeed a right-Hilbert 
$A\otimes D - C\otimes D$ bimodule isomorphism.

For the other part, take $c\in C$ and $d\in D$ and compute:
\begin{align*}
\Theta \bigl((x\otimes m)\otimes(y\otimes n)\bigr) \d (c\otimes d)
&= \Theta \bigl((x\otimes m)\otimes(y\otimes n) \d (c\otimes d)\bigr)\\
&= \Theta \bigl((x\otimes m)\otimes(y\d c\otimes nd)\bigr)\\
&= \Theta \bigl((x\otimes m)\otimes(y\d c\otimes ef)\bigr)
\righttext{(for some $e,f\in D$)}\\
&= \Theta \bigl((x\otimes me)\otimes(y\d c\otimes f)\bigr)\\
&= (x\otimes y\d c)\otimes mef\\
&= (x\otimes y)\d (c\otimes mnd)\\
&= \bigl((x\otimes y)\otimes mn\bigr)\d (c\otimes d).
\end{align*}
\end{proof}

The following construction of the balanced tensor product of coactions
should be compared to \cite[Proposition 2.10]{BS-CH}.

\begin{prop}
\label{tnsr-co-prop}
If ${}_{(A,\delta)}(X,\zeta)_{(B,\epsilon)}$ and
${}_{(B,\epsilon)}(Y,\eta)_{(C,\vartheta)}$ are right-Hilbert bimodule
coactions of $G$, then
\[
\zeta \cotimes_B \eta
\deq  \Theta\circ (\zeta\otimes_B \eta)
\]
defines a $\delta-\vartheta$ compatible coaction $\zeta\cotimes_B\eta$
of $G$ on $X \otimes_BY$,
where $\Theta$ is the isomorphism of Lemma~\textup{\ref{Theta-lem}}.
Moreover, if $\zeta$ and $\eta$ are nondegenerate, then so is
$\zeta\cotimes_B\eta$.
\end{prop}

\begin{proof}
By \propref{tnsr-hom-prop} and \lemref{Theta-lem}, $\Theta$ and $\zeta
\otimes_B \eta$ are nondegenerate right-Hilbert bimodule
homomorphisms, so their composition (which clearly has the desired
coefficient maps) is too.

To show that $(1_{M(A)}\otimes C^*(G))\d (\zeta
\cotimes_B \eta)(X\otimes_B Y) \subseteq (X\otimes_B Y)\otimes
C^*(G)$, we first claim that
\begin{equation}
\label{tnsr-co-eq2}
(1_{M(A)}\otimes c)\d
(\zeta(x)\otimes \eta(y))\in (X\otimes
C^*(G))\otimes_{B\otimes C^*(G)} (Y\otimes C^*(G))
\end{equation}
for all $c\in C^*(G)$, $x\in X$, and $y\in Y$.  For, since $\zeta$ is
a coaction on $X$, $(1\otimes c)\d \zeta(x)\in X\otimes C^*(G)$, so we
may write $(1\otimes c)\d \zeta(x) = z\d (1 \otimes d)$ for some $z\in
X\otimes C^*(G)$ and $d\in C^*(G)$.  Since $\eta$ is a coaction on $Y$
we have $(1\otimes d)\d \eta(y)\in Y \otimes C^*(G)$.  Thus
\begin{align*}
(1\otimes c)\d (\zeta(x)\otimes \eta(y))
&= (1\otimes c)\d \zeta(x)\otimes \eta(y)
= z\d (1\otimes d)\otimes \eta(y)
\\&= z\otimes (1\otimes d)\d \eta(y)
\in (X\otimes C^*(G))\otimes_{B\otimes C^*(G)} (Y\otimes C^*(G)),
\end{align*}
which gives \eqref{tnsr-co-eq2}.

Next, note that
\begin{equation}
\label{tnsr-co-eq3}
(1_{M(A)}\otimes c)\d \Theta(w)
= \Theta((1_{M(A)}\otimes c)\d w)
\end{equation}
for all $c\in C^*(G)$ and $w\in M((X\otimes C^*(G))\otimes_{B
\otimes C^*(G)}(Y\otimes C^*(G)))$, because the unique extension
(\propref{rHb-hom-ext-prop}) of $\Theta$ to the multiplier bimodule
has left coefficient map $\id_{M(A)}$.

Now combine \eqref{tnsr-co-eq2} and \eqref{tnsr-co-eq3}\ to get
\begin{align*}
(1_{M(A)}\otimes c)\d (\zeta \cotimes_B \eta)(x\otimes y)
&= (1\otimes c)\d \Theta\bigl(\zeta(x)\otimes \eta(y)\bigr)
\\&= \Theta \bigl((1\otimes c)\d (\zeta(x)\otimes \eta(y))\bigr)
\\&\in\Theta \bigl((X\otimes C^*(G))
\otimes_{B\otimes C^*(G)} (Y\otimes C^*(G))\bigr)
\\&= (X\otimes_B Y)\otimes C^*(G),
\end{align*}
so that $(1_{M(A)}\otimes C^*(G))\d (\zeta \cotimes_B \eta)(X
\otimes_B Y) \subseteq (X\otimes_B Y)\otimes C^*(G)$ by density and
continuity. 

If $\zeta$ and $\eta$ are nondegenerate, the coefficient
homomorphisms of $\zeta\cotimes_B\eta$ are certainly nondegenerate,
and using the above considerations we can compute
\begin{align*}
     \overline{(1\otimes C^*(G))(\zeta\cotimes_B\eta)(X\otimes_BY)}
&= \Theta\big(\overline{(1\otimes C^*(G)(\zeta(X)\otimes_{M(B\otimes
     C^*(G))}\eta(Y))}\big)\\
     &=\Theta\big(\overline{(X\otimes C^*(G))\otimes_{M(B\otimes
     C^*(G))}\eta(Y)}\big)\\
     &=\Theta\big(\overline{(X\otimes C^*(G))\otimes_{B\otimes
     C^*(G)}(1\otimes C^*(G))\eta(Y)}\big)\\
     &=\Theta\big((X\otimes C^*(G))\otimes_{B\otimes C^*(G)}(Y\otimes
     C^*(G))\big)\\
     &=(X\otimes_BY)\otimes C^*(G).
\end{align*}

It only remains to check the coaction identity:
\begin{align*}
\bigl((\zeta \cotimes_B \eta)\otimes \id\bigr)
\circ (\zeta \cotimes_B \eta)
&= \bigl(\Theta\circ (\zeta\otimes_B \eta)\otimes \id\bigr)
\circ \Theta
\circ (\zeta\otimes_B \eta)
\\&= (\Theta\otimes \id)
\circ \bigl(
(\zeta\otimes_B \eta)\otimes \id
\bigr)
\circ \Theta
\circ (\zeta\otimes_B \eta)
\\&\overset{(1)}{=} (\Theta\otimes \id)
\circ \Theta
\circ \bigl(
(\zeta\otimes \id)\otimes_{B\otimes C^*(G)}
(\eta\otimes \id)
\bigr)
\circ (\zeta\otimes_B \eta)
\\&\overset{(2)}{=} (\Theta\otimes \id)
\circ \Theta
\circ \bigl(
(\id\otimes\delta_G)\otimes_{B\otimes C^*(G)}
(\id\otimes\delta_G)
\bigr)
\circ (\zeta\otimes_B \eta)
\\&\overset{(3)}{=} (\id\otimes_G\delta)\circ \Theta
\circ (\zeta\otimes_B \eta)
\\&= (\id\otimes_G\delta)\circ (\zeta \cotimes_B \eta),
\end{align*}
where the equality at~(1) is justified by computing, for $x\in X$, $y
\in Y$, and $c,d\in C^*(G)$,
\begin{align*}
\bigl((\zeta\otimes_B \eta)\otimes \id\bigr)
\circ \Theta
\bigl((x\otimes c)\otimes (y\otimes d)\bigr)
&= \bigl((\zeta\otimes_B \eta)\otimes \id\bigr)
\bigl((x\otimes y)\otimes cd\bigr)
\\&= (\zeta\otimes_B \eta)(x\otimes y)\otimes cd
\\&= (\zeta(x)\otimes \eta(y))\otimes cd
\\&= \Theta \bigl(
(\zeta(x)\otimes c)\otimes (\eta(y)\otimes d)\bigr)
\\&= \Theta \bigl(
(\zeta\otimes \id)(x\otimes c)
\otimes
(\eta\otimes \id)(y\otimes d)\bigr)
\\&= \Theta\circ
\bigl((\zeta\otimes \id)\otimes_{B\otimes C^*(G)}
(\eta\otimes \id)\bigr)
\bigl((x\otimes c)\otimes (y\otimes d)\bigr),
\end{align*}
and then appealing to linearity, strict density, and strict
continuity; the equality at~(2) is justified by
\begin{align*}
\bigl(
(\zeta\otimes \id)\otimes_{B\otimes C^*(G)}
(\eta\otimes \id)
\bigr)
\circ (\zeta\otimes_B \eta)
&= (\zeta\otimes \id)\circ \zeta
\otimes_B
(\eta\otimes \id)\circ \eta
\\&= (\id\otimes\delta_G)\circ \zeta
\otimes_B
(\id\otimes\delta_G)\circ \eta
\\&= \bigl(
(\id\otimes\delta_G)\otimes_{B\otimes C^*(G)}
(\id\otimes\delta_G)
\bigr)
\circ (\zeta\otimes_B \eta);
\end{align*}
and finally the equality at~(3) is justified by computing, for $x\in
X$, $y\in Y$, and $s,t\in G$ (viewed as elements of $M(C^*(G))$
via the canonical embedding $G\to UM(C^*(G))$),
\begin{align*}
(\Theta\otimes \id)
\circ \Theta
\circ \bigl(
(\id\otimes\delta_G)\otimes_{B\otimes C^*(G)}
(\id\otimes\delta_G)
\bigr)&
\bigl((x\otimes s)\otimes (y\otimes t)\bigr)\\
&=(\Theta\otimes \id)
\circ \Theta
\bigl(
(\id\otimes\delta_G)(x\otimes s)
\otimes (\id\otimes\delta_G)(y\otimes t)
\bigr)
\\&= (\Theta\otimes \id)
\circ \Theta
\bigl(
(x\otimes s\otimes s)
\otimes (y\otimes t\otimes t)
\bigr)
\\&= (\Theta\otimes \id)
\bigl( \bigl(
(x\otimes s)\otimes (y\otimes t)\bigr)\otimes st \bigr)
\\&= (x\otimes y)\otimes st\otimes st
\\&= (\id\otimes\delta_G)
\bigl((x\otimes y)\otimes st\bigr)
\\&= (\id\otimes\delta_G)
\circ \Theta
\bigl((x\otimes s)\otimes (y\otimes t)\bigr),
\end{align*}
and then appealing to linearity, strict density, and strict
continuity.
\end{proof}

\begin{defn}
\label{CG-isom-defn}
Let $\zeta$ and $\eta$ be $\delta-\epsilon$ compatible coactions on
right-Hilbert $A-B$ bimodules $X$ and $Y$.  An isomorphism $\Phi$ of
$X$ onto $Y$ is \emph{$\zeta-\eta$ equivariant}, or \emph{intertwines}
$\zeta$ and $\eta$, if
\[
\eta\circ \Phi = (\Phi\otimes \id)\circ \zeta.
\]
We say $\zeta$ and $\eta$ are \emph{isomorphic}, or $X$ and $Y$ are
\emph{equivariantly isomorphic}, if such a $\Phi$ exists.
\end{defn}

It is straightforward to check that 
this notion of isomorphism is an equivalence relation on the
class of right-Hilbert $A-B$ bimodules with $\delta-\epsilon$
compatible coactions of $G$.

\begin{thm}
\label{CG-cat-thm}
Let $G$ be a locally compact group.  There is a category $\c C(G)$ in
which the objects are $C^*$-algebras with nondegenerate normal
coactions of $G$ \textup(see Definition~\textup{\ref{def-normal}} for
the meaning of normal coaction\textup), 
and in which the morphisms from $(A,\delta)$ to
$(B,\epsilon)$ are the equivariant isomorphism classes of
nondegenerate
right-Hilbert $A-B$ bimodules with $\delta-\epsilon$ compatible
coactions of $G$.  The composition of $[X,\zeta]\:(A,\delta)\to
(B,\epsilon)$ with $[Y,\eta]\:(B,\epsilon)\to (C,\vartheta)$ is the
isomorphism class of the tensor product coaction $\zeta \cotimes_B
\eta$ on $X\otimes_B Y$\textup; the identity morphism on $(A,\delta)$ is the
isomorphism class of the $\delta -\delta$ compatible right-Hilbert
bimodule coaction $(A,\delta)$ itself.
\end{thm}

\begin{proof}
We first note that the composition of morphisms is well-defined:
suppose $[X,\zeta] = [X',\zeta']\:(A,\delta)\to (B,\epsilon)$ and
$[Y,\eta] = [Y',\eta']\:(B,\epsilon)\to (C,\vartheta)$, so that
we have equivariant right-Hilbert bimodule isomorphisms
\[
\Phi\:(X,\zeta)\to (X',\zeta')
\midtext{and}
\Psi\:(Y,\eta)\to (Y',\eta').
\]
Then the isomorphism $\Phi\otimes_B
\Psi\:X\otimes_B Y\to X'\otimes_B Y'$ from the proof of
\thmref{C-cat-thm} satisfies
\begin{align*}
(\zeta' \cotimes_B \eta')\circ (\Phi\otimes_B \Psi)
&= \Theta\circ
(\zeta'\otimes_B \eta')\circ (\Phi\otimes_B \Psi)
\\&= \Theta\circ
\bigl(
(\zeta'\circ \Phi)\otimes_B (\eta'\circ \Psi)
\bigr)
\\&= \Theta\circ
\bigl(
(\Phi\otimes \id)\circ \zeta\otimes_B
(\Psi\otimes \id)\circ \eta
\bigr)
\\&= \Theta\circ
\bigl(
(\Phi\otimes \id)\otimes_{B\otimes B^*(G)} (\Psi\otimes \id)
\bigr)
\circ (\zeta\otimes_B \eta)
\\&= \bigl((\Phi\otimes_B \Psi)\otimes \id\bigr)
\circ \Theta
\circ (\zeta\otimes_B \eta)
\\&= \bigl((\Phi\otimes_B \Psi)\otimes \id\bigr)
\circ (\zeta \cotimes_B \eta).
\end{align*}

Next we establish that composition of morphisms in $\c C(G)$ is
associative; by the above it suffices to show that the coactions
$\zeta \cotimes_B (\eta \cotimes_C \tau)$ and $(\zeta \cotimes_B \eta)
\cotimes_C \tau$ are isomorphic for any tensorable right-Hilbert
bimodule coactions $({}_AX_B,\zeta)$, $({}_BY_C,\eta)$, and
$({}_CZ_D,\tau)$.  Let $\Phi\:X\otimes_B (Y\otimes_C Z)\to (X
\otimes_B Y)\otimes_C Z$ be the isomorphism from the proof of
\thmref{C-cat-thm}; we need to show that
\[
\bigl((\zeta \cotimes_B \eta) \cotimes_C \tau\bigr)
\circ \Phi
= (\Phi\otimes \id)
\circ \bigl(\zeta \cotimes_B (\eta \cotimes_C \tau)\bigr).
\]
We check this for an elementary tensor $x\otimes (y\otimes z)$:
\begin{align*}
\bigl((\zeta \cotimes_B \eta) \cotimes_C \tau\bigr)
\circ \Phi
\bigl(x\otimes (y\otimes z)\bigr)
&= \Theta \Bigl(
\Theta \bigl(\zeta(x)\otimes \eta(y)\bigr)
\otimes \tau(z)
\Bigr)
\\&\overset{(*)}{=} (\Phi\otimes \id)\circ \Theta
\Bigl(
\zeta(x)\otimes
\Theta \bigl(\eta(y)\otimes \tau(z)\bigr)
\Bigr)
\\&= (\Phi\otimes \id)\circ
\bigl(\zeta \cotimes_B (\eta \cotimes_C \tau)\bigr)
\bigl(x\otimes (y\otimes z)\bigr),
\end{align*}
where the equality at $(*)$ is justified by replacing $\zeta(x)$ by an
elementary tensor $x'\otimes c\in X\otimes C^*(G)$, replacing
$\eta(y)$ by $y'\otimes d$ and $\tau(z)$ by $z'\otimes e$,
and then appealing to linearity, strict density, and strict continuity:
\begin{align*}
\Theta \bigl(
\Theta \bigl((x'\otimes c)\otimes (y'\otimes d)\bigr)
\otimes (z'\otimes e)
\bigr)
&= \Theta \bigl(
\bigl((x'\otimes y')\otimes cd\bigr)
\otimes (z'\otimes e)
\bigr)
\\&= \bigl((x'\otimes y')\otimes z'\bigr)
\otimes cde
\\&= \Phi \bigl(x'\otimes (y'\otimes z')\bigr)
\otimes cde
\\&= (\Phi\otimes \id)
\bigl(
\bigl(x'\otimes (y'\otimes z')\bigr)
\otimes cde
\bigr)
\\&= (\Phi\otimes \id)
\circ \Theta
\bigl(
(x'\otimes c)\otimes
\bigl((y'\otimes z')\otimes de)\bigr)
\bigr)
\\&= (\Phi\otimes \id)
\circ \Theta
\bigl(
(x'\otimes c)\otimes
\Theta \bigl((y'\otimes d)\otimes (z'\otimes e)\bigr)
\bigr).
\end{align*}

Finally, note that for any coaction
$({}_AX_B,\zeta)$ of $G$, with left coefficient coaction $\delta$, the
canonical isomorphism $\Phi\:A\otimes_A X\to X$ intertwines the
coactions $\delta \cotimes_A \zeta$ and
$\zeta$.  To see this, take $a\in A$ and $x\in X$, and compute:
\begin{equation*}
\zeta\circ \Phi(a\otimes x)
= \zeta(a\d x)
=\delta(a)\d \zeta(x)
\overset{(*)}{=} (\Phi\otimes \id)\circ \Theta
\bigl(\delta(a)\otimes \zeta(x)\bigr)
= (\Phi\otimes \id)\circ (\delta \cotimes_A \zeta)(a\otimes x),
\end{equation*}
where the equality at $(*)$ is justified by replacing $\delta(a)$ by an
elementary tensor $b\otimes c\in A\otimes C^*(G)$, replacing
$\zeta(x)$ by $y\otimes d$, and then appealing to linearity, strict
density, and strict continuity:
\begin{align*}
(b\otimes c)\d (y\otimes d)
&= b\d y\otimes cd
= \Phi(b\otimes y)\otimes cd
= (\Phi\otimes \id)
\bigl((b\otimes y)\otimes cd\bigr)
\\&= (\Phi\otimes \id)\circ \Theta
\bigl((b\otimes c)\otimes (y\otimes d)\bigr).
\end{align*}
This gives $\zeta\circ \Phi = (\Phi\otimes \id)\circ(\delta
\cotimes_A \zeta)$.  A similar argument shows that for any coaction
$({}_BX_A,\zeta)$ with right coefficient coaction $\delta$, the
canonical isomorphism $X\otimes A \xrightarrow{\cong} X$
intertwines $\zeta \cotimes_A\delta$ and $\zeta$; hence the identity
morphism on $(A,\delta)$ is the isomorphism class of the $\delta -
\delta$ compatible right-Hilbert bimodule coaction $(A,\delta)$.
\end{proof}

\begin{rem}
Notice that we could consider a larger category where the
$C^*$-coactions on the objects are not required to be nondegenerate
and normal, in which our category $\c C(G)$ would sit as a 
subcategory.  For our purposes, we require the specific category $\c
C(G)$, so that the Mansfield imprimitivity theorem applies without
any further hypotheses (see \appxref{imprim-chap}).
\end{rem}

\begin{rem}
\label{CG-ibs-rem}
The equivalences in $\c C(G)$ are exactly the (isomorphism classes of)
imprimitivity bimodule coactions.  In one direction, suppose
$({}_AX_B,\zeta)$ is an imprimitivity bimodule coaction of $G$, with
left coefficient coaction $\delta$.  Then the map
\[
\rev x\mapsto \rev{\zeta(x)}
\: \rev X\to M \bigl( (X\otimes C^*(G))\sptilde\bigr)
\]
is a nondegenerate homomorphism, and it is routine to check
that the map $\Psi\:(X\otimes C^*(G))\sptilde\to \rev X\otimes
C^*(G)$ defined by
\[
\Psi(\rev{x\otimes c}) = \rev x\otimes c^*
\]
is an
isomorphism; hence we can define a nondegenerate homomorphism
$\rev \zeta\:\rev X\to M \bigl(\rev X\otimes C^*(G)\bigr)$
by
\[
\rev \zeta(\rev x) = \Psi \bigl(\rev{\zeta(x)}\bigr).
\]
Routine computations show that $\rev \zeta$ is in fact a bimodule
coaction, which is nondegenerate if $\zeta$ is nondegenerate.

We claim that the canonical isomorphism $\Phi\:X\otimes_B \rev X
\to A$ intertwines $\zeta \cotimes_B \rev \zeta$ and $\delta$.  To see
this, take $x,y\in X$, and compute:
\begin{align*}
&(\Phi\otimes \id)\circ (\zeta \cotimes_B \rev \zeta)
(x\otimes \rev y)
= (\Phi\otimes \id)\circ \Theta
\Bigl( \zeta(x)\otimes
\Psi \bigl(\rev{\zeta(y)}\bigr) \Bigr)
\\&\quad\overset{(*)}{=} {}_{A\otimes C^*(G)}\< \zeta(x), \zeta(y) \>
=\delta \bigl({}_A\<x,y\>\bigr)
=\delta\circ \Phi(x\otimes \rev y),
\end{align*}
where the equality at $(*)$ is justified by replacing $\zeta(x)$ by an
elementary tensor $z\otimes c\in X\otimes C^*(G)$, replacing
$\zeta(y)$ by $w\otimes d$, and then appealing to linearity, strict
density, and strict continuity:
\begin{align*}
&(\Phi\otimes \id)\circ \Theta
\Bigl( (z\otimes c)\otimes
\Psi \bigl(\rev{w\otimes d}\bigr) \Bigr)
= (\Phi\otimes \id)\circ \Theta
\bigl((z\otimes c)\otimes (\rev w\otimes d^*)\bigr)
\\&\quad= (\Phi\otimes \id)
\bigl((z\otimes \rev w)\otimes cd^*\bigr)
= \Phi(z\otimes \rev w)\otimes cd^*
\\&\quad= {}_A\<z,w\>\otimes cd^*
= {}_{A\otimes C^*(G)}\<z\otimes c, w\otimes d\>.
\end{align*}
A similar argument shows that the canonical homomorphism $\Psi\:\rev
X\otimes_A X\to B$ intertwines $\rev \zeta \cotimes_A \zeta$ and
$\epsilon$.

In the other direction, if $({}_AX_B,\zeta)$ and $({}_BY_A,\eta)$ are
$\delta-\epsilon$ and $\epsilon -\delta$ compatible right-Hilbert
bimodule coactions of $G$ with
\[
\zeta \cotimes_B \eta \cong\delta
\midtext{and}
\eta \cotimes_A \zeta \cong \epsilon,
\]
then in particular $X$ is an $A-B$ imprimitivity bimodule by
\lemref{C-ibs-lem}, so $\zeta$ and $\eta$ are imprimitivity bimodule
coactions by \remref{rHb-co-rem}.
\end{rem}

\section{Actions and coactions}

\begin{defn}
\label{ACG-isom-defn}
Let $\gamma$ and $\rho$ be $\alpha-\beta$ compatible actions, and let
$\zeta$ and $\eta$ be $\delta-\epsilon$ compatible coactions, of $G$
on right-Hilbert $A-B$ bimodules $X$ and $Y$.  We say the triples
$(X,\gamma,\zeta)$ and $(Y,\rho,\eta)$ are \emph{isomorphic}, or $X$
and $Y$ are \emph{equivariantly isomorphic}, if there exists a
right-Hilbert $A-B$ bimodule isomorphism of $X$ onto $Y$ which is both
$\gamma-\rho$ equivariant and $\zeta-\eta$ equivariant.
\end{defn}

Needless to say, this notion of
isomorphism is an equivalence relation on the class of right-Hilbert
$A-B$ bimodules with $\alpha-\beta$ compatible actions and $\delta
- \epsilon$ compatible coactions of $G$.

\begin{thm}
\label{ACG-cat-thm}
Let $G$ be a locally compact group.  There is a category $\c A \c
C(G)$ in which the objects are $C^*$-algebras with actions and
nondegenerate normal coactions of $G$, and in which the morphisms from
$(A,\alpha,\delta)$ to $(B,\beta,\epsilon)$ are the equivariant
isomorphism classes of
right-Hilbert $A-B$ bimodules with
$\alpha-\beta$ compatible actions and
nondegenerate
$\delta-\epsilon$ compatible
coactions of $G$.  The composition of $[X,\gamma,\zeta] \:
(A,\alpha,\delta)\to (B,\beta,\epsilon)$ with $[Y,\rho,\eta] \:
(B,\beta,\epsilon)\to (C,\nu,\vartheta)$ is the isomorphism class of
$(X\otimes_B Y,\gamma\otimes_B \rho,\zeta \cotimes_B \eta)$\textup; the
identity morphism on $(A,\alpha,\delta)$ is the isomorphism class of
$(A,\alpha,\delta)$ itself.
\end{thm}

\begin{proof}
Since the
proofs of Theorems~\ref{AG-cat-thm} and~\ref{CG-cat-thm} both use the
same isomorphisms $X\otimes_B Y \cong Z\otimes_B W$, $X\otimes_B (Y
\otimes_C Z)\cong (X\otimes_B Y)\otimes_C Z$, $A\otimes_A X \cong
X$, and $Y\otimes_A A \cong Y$, these proofs combine to show that
composition of morphisms is well-defined and associative, and that
there are identity morphisms, in $\c A \c C(G)$.
\end{proof}

\begin{rem}
It might seem that this theorem could be proven by
identifying $\c A \c C(G)$ with the ``diagonal'' subcategory of $\c
A(G)\times \c C(G)$ consisting of the objects
$((A,\alpha),(B,\delta))$ which satisfy $A = B$ and the morphisms
$([X,\gamma],[Y,\zeta])$ which satisfy $X \cong Y$.  However, this
approach would fail, because the map $[X,\gamma,\zeta]\mapsto
([X,\gamma],[X,\zeta])$ may not be injective: on the right side
there can be different action-equivariant and coaction-equivariant
isomorphisms, while on the left side there has to be one
isomorphism which is both action- and coaction-equivariant.
\end{rem}

\section{Actions and coactions on linking algebras}
\label{sec-actcoactlink}

One of the basic tools in this work is to exploit the relation
between actions and coactions on (right-) partial imprimitivity
bimodules and actions and coactions on the corresponding linking
algebras. Recall from \secref{sec-link} in \chapref{hilbert-chap}
that for any partial $A-B$ imprimitivity bimodule $X$, we can form
the linking algebra
$$L(X)=\mtx{A&X\\ \rev{X}&B}$$
with multiplication and involution as given in \eqref{mat1}.
As usual, we denote the corner projections
$\smtx{1&0\\ 0&0}$ and $\smtx{0&0\\ 0&1}$
by $p$ and $q$, respectively.
For actions, we get:

\begin{lem}\label{lem-actlink}
     Suppose that $X$ is a partial $A-B$ imprimitivity bimodule.
     Let $\alpha$ and $\beta$ be actions of $G$ on $A$ and $B$,
     respectively, and let $\gamma$ be an $\alpha-\beta$ compatible
     right-Hilbert bimodule action of $G$ on $X$. Then there exists
     an action $\nu\:G\to \Aut L(X)$ given by
     $$\nu_s\left(\smtx{a&x\\ \tilde{y}&b}\right)\deq 
     \mtx{\alpha_s(a)&\gamma_s(x)\\ \widetilde{\gamma_s(y)}&
     \beta_s(b)}.$$
     Conversely, if $\nu\:G\to \Aut L(X)$ is an action such that
     $\nu_s(p)=p$ and \textup(hence\textup) $\nu_s(q)=q$, we obtain actions
     $\alpha$, $\beta$, and $\gamma$ on $A$, $B$, and $X$, respectively,
     such that
     $$\alpha_s(plp)=p\nu_s(l)p,\quad \beta_s(qlq)=q\nu_s(l)q,\quad
     \text{and}\quad \gamma_s(plq)=p\nu_s(l)q,$$
     for all $s\in G$ and $l\in L(X)$.
\end{lem}
\begin{proof} Note first that it follows from \lemref{lem-partialhom}
     that a right-Hilbert bimodule action on $_AX_B$ is automatically
     a partial imprimitivity bimodule action (and \emph{vice versa}). 
Using this,
     the proof follows directly from the algebraic properties
     in $L(X)$, so we omit further details.
\end{proof}

We now consider coactions. For imprimitivity bimodules, the following
result can be found in \cite{er:stab}:

\begin{lem}\label{lem-coactlink}
     Suppose that $_AX_B$ is a right-partial imprimitivity bimodule
     and assume that $\zeta$ is coaction of $G$ on $X$ with
     coefficient coactions $\delta$ and $\epsilon$ of $G$ on
     $A$ and $B$, respectively.
     Then there is a unique coaction $\nu\:L(X)\to M(L(X)\otimes C^*(G))$
     such that \textup(after identifying $L(X)\otimes C^*(G)$ with
     $L(X\otimes C^*(G))$ as in Remark~\textup{\ref{L-rem}}\textup)
     $$\nu\left(\smtx{a&x\\ \tilde{y}&b}\right)=\mtx{\delta(a)&\zeta(x)\\
     \widetilde{\zeta(y)}&\epsilon(b)}$$
for all $a\in A$, $x,y\in X$, and $b\in B$. 
Moreover,
     $\nu$ will be nondegenerate if and only if
     $\zeta$ is, and $\nu$ will be normal if and only if $\epsilon$ is.

     Conversely, if $\nu\:L(X)\to M(L(X)\otimes C^*(G))$ is a coaction
     such that $\nu(p)=p\otimes 1$ and $\nu(q)=q\otimes 1$, then $\nu$
     compresses to give coactions $\delta$ and $\epsilon$ on the
     corners $A$ and $B$ and a $\delta-\epsilon$ compatible coaction
     $\zeta$ on $X$\textup; these coactions will be
     nondegenerate if and only if $\nu$ is, and $\delta$ and $\epsilon$
will be normal if and only if $\nu$ is.
\end{lem}
\begin{proof} Let $\zeta$ be as in the first part of the
     lemma. Then it follows from
     \lemref{link-nodeg} (see also \remref{remdegen})
     that $\nu\:L(X)\to M(L(X\otimes
     C^*(G))=M(L(X)\otimes C^*(G)))$ is a nondegenerate homomorphism.
     A straightforward computation then shows that it is a coaction.
     Since
     $$(1\otimes C^*(G))\nu(L(X))=
     \mtx{(1\otimes C^*(G))\delta(A)&(1\otimes C^*(G))\zeta(X)\\
     (1\otimes C^*(G))\widetilde{\zeta(X)}& (1\otimes C^*(G))\epsilon(B)}
     $$
     it follows that $\nu$ is nondegenerate (as a coaction) if and only
     $_{\delta}\zeta_{\epsilon}$ is nondegenerate.
     Similarly, since
     $$(\id_{L(X)}\otimes \lambda)\circ \nu=
     \mtx{(\id_A\otimes\lambda)\circ \delta&
     (\id_X\otimes\lambda)\circ \zeta\\
     (\id_X\otimes \lambda)\circ\zeta\!\!\widetilde{\ \ }&
     (\id_B\otimes \lambda)\circ \epsilon},$$
     and since, by the Rieffel correspondence,
     all corner maps are injective if and only if
     the lower-right-hand-corner map is injective, it follows
     that $\nu$ is normal if and only if $\epsilon$ is normal.

     Conversely, if $\nu\:L(X)\to M(L(X)\otimes C^*(G))$ is a coaction,
     it follows from \lemref{link-nodeg} that it compresses to
     a nondegenerate right-Hilbert bimodule homomorphism
     $_{\delta}\zeta_{\epsilon}$,
     and it is then easy to check that $\zeta$ is a $\delta-\epsilon$
     compatible right-Hilbert bimodule homomorphism. The final
     assertion then follows from the first part of the proof.
\end{proof}

\section{Standard factorization of morphisms}

In this section we want to show that every morphism in our categories
$\c C$, $\c A(G)$, $\c C(G)$, and $\c A\c C(G)$
can be factored as a product of a
nondegenerate standard morphism
and a right-partial equivalence (see the definitions below).
We start by defining what we mean by a
standard morphism.  This is easy in the category $\c C$:
recall from Example \ref{ex-standard}
that a standard right-Hilbert $A-B$ bimodule is one of the form
$X=\varphi(A)B$, where $\varphi\:A\to M(B)$ is a $*$-homomorphism.
We shall call the corresponding morphism
$[X]\:A\to B$ in our category $\c C$ the {\em standard morphism
associated to $\varphi$}. The following observation allows us to extend
this notion to the category $\c A(G)$.

\begin{lem}\label{lem-equiv}
     Suppose that $(A,\alpha)$ and $(B,\beta)$ are actions and
     assume that $\varphi\:A\to M(B)$ is a \textup(possibly degenerate\textup)
     $\alpha-\beta$ equivariant $*$-homomorphism.
     Let $X=\varphi(A)B\subseteq B$, and
     for each $s\in G$ let $\beta_X(s)$ denote the restriction
     of $\beta_s$ to $X$. Then $\beta_X$ is an $\alpha-\beta$
     compatible action of $G$ on $X$.
\end{lem}

We omit the straightforward proof. The coaction
analogue of the above result is a bit more complicated.

\begin{lem}\label{lem-coactstandard}
     Suppose that $(A,\delta)$ and $(B,\epsilon)$ are nondegenerate 
coactions of $G$,
     and assume that $\varphi\:A\to M(B)$ is a \textup(possibly degenerate\textup)
     $\delta-\epsilon$ equivariant $*$-homomorphism. Then the restriction
     of $\epsilon$ to $X=\varphi(A)B\subseteq B$ determines
     a $\delta-\epsilon$ compatible coaction $\epsilon_X$ on the standard
     $A-B$ bimodule $X$. Moreover, if $\delta$ and $\epsilon$ are
     nondegenerate, then so is $\epsilon_X$.
\end{lem}
\begin{proof}
     Let $\iota\:X\hookrightarrow B$ denote the inclusion map. 
     Then $\iota$ is
     a right-Hilbert bimodule homomorphism from $_AX_B$ to $M({_BB_B})$.
     Thus it follows from Proposition \ref{prop-iso} that the inclusion
     $\iota\otimes \id_G\:X\otimes C^*(G)\to M_G(B\otimes C^*(G))$
     extends to an inclusion of the $G$-multiplier bimodules
     $M_G(X\otimes C^*(G))\subseteq
     M_G(B\otimes C^*(G))$.
     We show that $\epsilon(X)\subseteq M_G(X\otimes C^*(G))$.
     Using Lemma \ref{lem-Cmultidentify}, it is enough to show that
     $$(1\otimes C^*(G))\epsilon(X)\cup \epsilon(X)(1\otimes C^*(G))\subseteq
     X\otimes C^*(G).$$
     But this follows from
     \begin{align*}
	{\epsilon(X)(1\otimes C^*(G))}&=
	{\epsilon(\varphi(A)B)(1\otimes C^*(G))}\\
	&={(\varphi\otimes\id_G)\big(\delta(A)\big)\big(\epsilon(B)(1\otimes
	C^*(G))\big)}\\
	&\subseteq
	{(\varphi\otimes\id_G)\big(\delta(A)\big)\big(B\otimes C^*(G)\big)}\\
	&={(\varphi\otimes\id_G)\big(\delta(A)\big)\big((1\otimes C^*(G))
	(B\otimes C^*(G))\big)}\\
	&= {(\varphi\otimes\id_G)
	\big(\delta(A)(1\otimes C^*(G))\big)\big(
	B\otimes C^*(G)\big)}\\
	&\subseteq {\varphi\otimes \id_G\big(A\otimes 
C^*(G)\big)(B\otimes C^*(G))}\\
	&=X\otimes C^*(G),
     \end{align*}
     and a similar computation which shows that
     $(1\otimes C^*(G))\epsilon(X)\subseteq X\otimes C^*(G)$.
     Note that if $\delta$ and $\epsilon$ are both nondegenerate, then
     taking closures of all terms in the above computation allows us to
     replace all inclusions by equal signs.

     It now makes sense to define the restriction
     $\epsilon_X\:X\to M_G(X\otimes C^*(G))$
      of $\epsilon$ to $X$.
     Then $\epsilon_X$ is a right-Hilbert
     bimodule homomorphism with coefficient maps $\delta$ and
     $\epsilon$. To see that it is nondegenerate (as a bimodule
     homomorphism), we compute as above
     \begin{align*}
	{\epsilon_X(X)(B\otimes C^*(G))}&=
	{(\phi\otimes\id_G)\circ\delta(A)(B\otimes C^*(G))}\\
	&\supseteq (\phi\otimes\id_G)\big(\delta(A)(A\otimes
	C^*(G))\big)(B\otimes C^*(G))\\
	&=(\phi\otimes\id_G)\big(A\otimes C^*(G)\big)(B\otimes C^*(G))\\
	&=X\otimes C^*(G).
     \end{align*}
     The coaction identity for $\epsilon_X$ follows directly from the
     coaction identity for $\epsilon$.
     If $\delta$ and $\epsilon$ are nondegenerate, then we already saw
     above that
     $$\overline{(1\otimes C^*(G))\epsilon_X(X)}=\overline{(1\otimes 
C^*(G))\epsilon(X)}
     =X\otimes C^*(G),$$
     so $\epsilon_X$ is a nondegenerate coaction.
     \end{proof}

\begin{defn}\label{def-standard}
     Suppose that $\varphi\:A\to M(B)$ is a $*$-homomorphism and
     let $X=\varphi(A)B$ be the standard right-Hilbert $A-B$ bimodule
     associated to $\varphi$.
\begin{enumerate}
\item
     If $\varphi\:A\to M(B)$ is $\alpha-\beta$ equivariant for the
     actions $\alpha$ and $\beta$ of $G$ on $A$ and $B$, respectively,
     and if $\beta_X$ is the restriction of $\beta$ to $X$ as in
     Lemma \ref{lem-equiv},
     then we say that $(X,\beta_X)$ is the {\em standard action}
     associated to $\varphi$.

\item
     If $\varphi\:A\to M(B)$ is $\delta-\epsilon$
     equivariant  for the nondegenerate
     coactions $\delta$ and $\epsilon$ of $G$ on $A$ and $B$,
     respectively, and if
     $\epsilon_X\:X\to M_G(X\otimes C^*(G))$ is as in Lemma
     \ref{lem-coactstandard}, then
      $(X,\epsilon_X)$ is called the
     {\em standard coaction} associated
     to $\varphi$.

\item
     If $\varphi\:A\to M(B)$ is nondegenerate, then we say that
     the above-defined standard morphisms, actions, or coactions, are
     {\em nondegenerate}.
     (Note that in this case we have $X=B$.)
\end{enumerate}
\end{defn}

We now introduce the notion of partial equivalences.

\begin{defn}\label{defn-partialequiv}
     Suppose that $_AX_B$ is right-partial imprimitivity bimodule.
     Then $[X]\:A\to B$ is called a {\em right-partial equivalence}
     between $A$ and $B$ in the category $\C$.
     Similarly, morphisms $[(X,\gamma)]$, $[(X,\zeta)]$ and $[(X,
     \gamma,\zeta)]$  in the categories
     $\A(G), \C(G)$ and $\A\C(G)$, respectively,
     are called {\em right-partial equivalences}, if the underlying
     module $X$ is a right-partial imprimitivity bimodule.
\end{defn}

Suppose now that $_AX_B$ is any right-Hilbert $A-B$ bimodule.  If we
put $K=\c K(X)$ and let $\kappa\:A\to M(K)$ be the associated
nondegenerate homomorphism, we obtain a
nondegenerate standard bimodule $_AK_K$, and $_KX_B$
becomes a right-partial $K-B$ imprimitivity bimodule.
It is trivial to check that the map
$$k\otimes x\mapsto k\d  x$$
extends to a right-Hilbert $A-B$ bimodule isomorphism
$$_A(K\otimes_KX)_B\cong {_AX_B},$$
so we see that any morphism in the category $\c
C$ can be factored as the composition of a nondegenerate
standard morphism, and a right-partial equivalence.
It is easy
to extend this observation to morphisms in the category $\c A(G)$:

\begin{prop}
\label{decomaction}
Suppose that ${}_{(A,\alpha)}(X,\gamma)_{(B,\beta)}$ is an action,
and let $K=\c K(X)$.
Then\textup:
\begin{enumerate}
\item
There exists a unique action $\mu$ on $K$ such that $\gamma$ is
$\mu-\beta$ compatible.

\item
The canonical nondegenerate homomorphism $\kappa\:A\to M(K)$ is
$\alpha-\mu$ equivariant\textup; thus $(_AK_K,\mu)$ is a
nondegenerate standard action with left coefficient action $\alpha$.

\item
The map $k\otimes x\mapsto k\d x$ implements an equivariant
isomorphism between $(_A(K\otimes_KX)_B,
\mu\otimes_K\gamma)$
and $(_AX_B,\gamma)$.
\end{enumerate}
In particular, every morphism in $\c A(G)$ can be factored as the
composition of a nondegenerate standard morphism,
and a right partial equivalence.
\end{prop}

\begin{proof}
Since by definition $\gamma_s(x)\d\<\gamma_s(y),\gamma_s(z)
\>_B=\gamma_s(x\d\<y,z\>_B)$ for all $s\in G$ and $x,y,z\in
X$ it follows from \cite[page~292]{com} that
$\mu_s({}_K\<x,y\>)\deq {}_K\<\gamma_s(x),\gamma_s(y)\>$
determines a unique strongly
continuous action $\mu$ of $G$ on $K$ such that $\gamma$ is
$\mu-\beta$ compatible.  This gives (i).

For each $a\in A$, $s\in G$, and $x\in X$ we have
\begin{align*}
\phi(\alpha_s(a))\gamma_s(x)
&= \alpha_s(a)\d\gamma_s(x)
= \gamma_s(a\d x)
= \gamma_s(\phi(a)x)=\mu_s(\phi(a))\gamma_s(x).
\end{align*}
Thus $\phi\:A\to M(K)$ is $\alpha-\mu$ equivariant and so
$(_AK_K,\mu)$ is a standard action.

For (iii) just check that $k\otimes x\mapsto k\d x$ is
$(\mu\otimes_K\gamma)-\gamma$ equivariant,
which is trivial.
\end{proof}

We have to work a bit more to get the analogous result for morphisms
in $\c C(G)$.
Recall from \remref{rHb-hom-rem1} that if
${}_\phi\Phi_\psi\:{}_AX_B\to M({}_CY_D)$
is any nondegenerate right-Hilbert
bimodule homomorphism then there exists a map $\mu\:\c K_B(X)\to M(\c
K_D(Y))$ such that ${}_{\mu}\Phi_\psi$ is a nondegenerate
partial imprimitivity bimodule homomorphism.  The
following easy lemma is quite fundamental.

\begin{lem}
\label{lemfund}
With the above notation\textup:
\begin{enumerate}
\item
If $\kappa_A\:A\to M(\c K_B(X))$ and $\kappa_C\:C\to M(\c K_D(Y))$
are the canonical
nondegenerate homomorphisms associated with the left module
actions
on $X$ and $Y$, then the following diagram commutes\textup:
\[
\xymatrix{
{A}
\ar[r]^-{\kappa_A}
\ar[d]_{\phi}
&{M(\c K_B(X))}
\ar[d]^{\mu}
\\
{M(C)}
\ar[r]_-{\kappa_C}
&{M(\c K_D(Y)).}
}
\]

\item
$\ker\mu$ is the ideal induced from $\ker\psi$ via $X$.
\end{enumerate}
\end{lem}

\begin{proof}
For (i) let $a\in A$, $x\in X$, and $d \in D$, and compute
\begin{align*}
\mu(\kappa_A(a))\d\Phi(x)d
&= \Phi(\kappa_A(a)x)d
= \Phi(a\d x)d
= \phi(a)\d\Phi(x) d
= \kappa_C(\phi(a))\Phi(x)d.
\end{align*}
Since $\overline{\Phi(X)\d D}=Y$ by the nondegeneracy of $\Phi$,
it follows that $\mu(\kappa_A(a))$ and $\kappa_C(\phi(a))$ are
identical in $\L(Y)$.   The assertion (ii) follows from
Lemma \ref{lem-induced-hom}.
\end{proof}

\begin{lem}
\label{KK-tnsr-lem}
Let ${}_AX_B$ and ${}_CY_D$ be right-Hilbert bimodules. The equation
\[
\Phi(k\otimes l)(x\otimes y) = k(x)\otimes l(y)
\]
determines a $C^*$-algebra isomorphism $\Phi$ of $\c K_B(X)\otimes \c
K_D(Y)$ onto $\c K_{B\otimes D}(X\otimes Y)$ such that
\[
\Phi\circ(\kappa_A\otimes\kappa_C) = \kappa_{A\otimes C},
\]
where $\kappa_A\:A\to M(\c K_B(X))$, $\kappa_C\:C\to M(\c K_D(Y))$,
and
$\kappa_{A\otimes C}\:A\otimes C\to M(\c K_{B\otimes D}(X\otimes Y))$
are the nondegenerate homomorphisms associated with $X$,
$Y$, and $X\otimes Y$.
\end{lem}

\begin{proof}
The first part of the lemma is \cite[pages~35--37]{lan:hilbert}. For the
second part, simply compute:
\begin{align*}
&\Phi\circ(\kappa_A\otimes\kappa_C)(a\otimes c)(x\otimes y)
= \Phi(\kappa_A(a)\otimes\kappa_C(c))(x\otimes y)
\\&\quad= \kappa_A(a) x\otimes\kappa_C(c) y
= a\d x\otimes c\d y
= (a\otimes c)\d (x\otimes y)
= \kappa_{A\otimes C}(a\otimes c)(x\otimes y).
\end{align*}
\end{proof}

\begin{prop}
\label{decomcoaction}
Suppose that ${}_{(A,\delta)}(X,\zeta)_{(B,\epsilon)}$ is a
coaction, and let $K=\c K_B(X)$.  Then\textup:
\begin{enumerate}
\item
There exists a unique coaction $\mu$ on $K$ such that
$\zeta$ is $\mu-\epsilon$ compatible.

\item
The canonical nondegenerate homomorphism $\kappa_A\:A\to M(K)$ is
$\delta-\mu$ equivariant; thus $(_AK_K,\mu)$ is a
nondegenerate standard coaction
with left coefficient coaction $\delta$.

\item
The map $k\otimes x\mapsto k\d x$ induces an equivariant
isomorphism between $(_A(K\otimes_K X)_B,\mu \cotimes_K \zeta)$
and $(_AX_B,\zeta)$.
\end{enumerate}
Moreover, $\mu$ is nondegenerate if $\zeta$ is, and $\mu$ is normal if
$\epsilon$ is. 
\end{prop}

\begin{proof}
Since $\zeta\:X\to M(X\otimes C^*(G))$ is  a nondegenerate bimodule
homomorphism, it follows from \remref{rHb-hom-rem1}
that there exists a
map $\nu\:K\to M(\c K_{B\otimes C^*(G)}(X\otimes C^*(G)))$
such that ${}_\nu\zeta_\epsilon$ is a nondegenerate
right-Hilbert bimodule homomorphism.  We compose $\nu$ with the
isomorphism $\Phi^{-1}\:M(\c K_{B\otimes C^*(G)}(X\otimes C^*(G)))
\to M(K\otimes C^*(G))$ provided by \lemref{KK-tnsr-lem} to get a
nondegenerate homomorphism $\mu\:K\to M(K\otimes C^*(G))$
such that ${}_{\mu}\zeta_\epsilon$
is a nondegenerate right-Hilbert
bimodule homomorphism.  Since $\epsilon$ is injective and since,
by Lemma \ref{KK-tnsr-lem},
$\ker\mu=\ker\nu$ is the ideal induced from $\ker\epsilon$ via
$X$, it follows that $\mu$ is injective, too.  Thus to show that
$\mu$ is a coaction on $K$ it only remains to show that
\begin{enumerate}
\item[(a)]
$(1\otimes z)\mu(k)\in K\otimes C^*(G)$
for all $z\in C^*(G)$, $k\in K$, and

\item[(b)]
$(\mu\otimes \id_G)\circ\mu
= (\id_K\otimes\delta_G)\circ\mu$.
\end{enumerate}
But (a) follows from the simple calculation
\begin{align*}
(1\otimes z)\mu({}_K\<x,y\>)
&= (1\otimes z){}_{M(K\otimes
C^*(G))}\<\zeta(x), \zeta(y)\>
= {}_{M(K\otimes C^*(G))}\<(1\otimes z)\d \zeta(x), \zeta(y)\>,
\end{align*}
which is in $K\otimes C^*(G)$ since $(1\otimes z)\d \zeta(x)\in
X\otimes C^*(G)$ and
\[
{}_{M(K\otimes C^*(G))}\bigl\<X\otimes C^*(G),
M(X\otimes C^*(G))\bigr\> \subseteq
K\otimes C^*(G).
\]
In order to prove (b) we just compute
\begin{align*}
(\mu\otimes \id_G)\circ\mu({}_K\<x,y\>)
&= (\mu\otimes \id_G)
\bigl({}_{M(K\otimes C^*(G))}\<\zeta(x), \zeta(y)\>\bigr)
\\&\quad= {}_{M(K\otimes C^*(G)\otimes C^*(G))}\bigl\<
(\zeta\otimes \id_G)\circ \zeta(x),
(\zeta\otimes \id_G)\circ \zeta(y)\bigr\>
\\&\quad= {}_{M(K\otimes C^*(G)\otimes C^*(G))}\bigl\<
(\id_X\otimes\delta_G)\circ \zeta(x),
(\id_X\otimes\delta_G)\circ \zeta(y)\bigr\>
\\&\quad= (\id_K\otimes\delta_G)\circ\mu({}_K\<x,y\>)
\end{align*}
for all $x,y\in X$.  Then (b) follows from the density of
${}_K\<X,X\>$ in $K$.  Since the left coefficient map
on $\c K(X)$ of any right-Hilbert bimodule
homomorphism $\Phi\:{_{\K(X)}}X_B\to M({_CY_D})$ is uniquely determined
by $\Phi$ (since it is compatible with the $\c K(X)$-valued
inner product on $X$), this proves~(i).

It follows from \lemref{lemfund} that
\[
\mu\circ \kappa_A = \Phi^{-1}\circ \kappa_{A\otimes
C^*(G)}\circ\delta,
\]
where $\kappa_{A\otimes C^*(G)}\:A\otimes C^*(G)\to M(\c K_{B
\otimes C^*(G)}(X\otimes C^*(G)))$ is the canonical nondegenerate
homomorphism. Thus the second part of \lemref{KK-tnsr-lem} gives
\[
\mu\circ \kappa_A = (\kappa_A\otimes \id)\circ\delta,
\]
which establishes (ii).

For (iii), we need to show that
\[
\zeta\circ \Phi
= (\Phi\otimes \id)\circ (\mu \cotimes_K \zeta),
\]
where $\Phi\:K\otimes_K X\to X$ is the canonical isomorphism.  But
this was verified in the proof 
of \thmref{CG-cat-thm}, where we showed that 
there are identities in $\c C(G)$.

Now assume that $\zeta$ (and therefore also $\delta$ and
$\epsilon$) is a nondegenerate coaction. Then
\begin{align*}
    \overline{(1\otimes C^*(G))\mu(K)}&=
    \overline{(1\otimes C^*(G)){}_{M(K\otimes C^*(G))}\lk \zeta(X),
    \zeta(X)\rk}\\
    &={}_{K\otimes C^*(G)}\overline{\lk (1\otimes
    C^*(G))\zeta(X),\zeta(X)\rk}\\
    &={}_{K\otimes C^*(G)}\overline{\lk X\otimes C^*(G),\zeta(X)\rk}\\
    &={}_{K\otimes C^*(G)}\overline{\lk (X\otimes C^*(G))(1\otimes
    C^*(G)),\zeta(X)\rk}\\
    &={}_{K\otimes C^*(G)}\overline{\lk X\otimes C^*(G),\zeta(X)(1\otimes
    C^*(G))\rk}\\
    &={}_{K\otimes C^*(G)}\overline{\lk X\otimes C^*(G),X\otimes C^*(G)\rk}\\
    &=\K(X)\otimes C^*(G).
\end{align*}
Finally, if $\epsilon$ is normal, the normality of  $\mu$
follows from the fact that
$$_{(\id_K\otimes\lambda)\circ\mu}\big((\id_X\otimes\lambda)
     \circ\zeta\big)_{(\id_B\otimes\lambda)\circ\epsilon}\:{_KX_B}\to
     M\big({_{K\otimes\K(L^2(G))}(X\otimes\K(L^2(G)))_{B\otimes
     \K(L^2(G))}}\big)$$
     is a right-partial imprimitivity bimodule homomorphism. Thus
     injectivity of the right coefficient map implies injectivity of
     the left coefficient map.
\end{proof}

\begin{rem}\label{rem-decomactioncoaction}
     Combining \propref{decomaction} with \propref{decomcoaction}
     gives a similar factorization result for morphisms in
     $\A\C(G)$. We omit the obvious details.
\end{rem}

In some situations, it becomes necessary to further factor
the right-partial equivalence part of the morphisms in
$\C$, $\A(G)$, $\C(G)$, or $\A\C(G)$, into an equivalence
and a standard morphism coming from an inclusion of
an ideal. To be more precise, if $_AX_B$ is a right-Hilbert
$A-B$  bimodule, then $X$ can be regarded as an
$_AX_{B_X}$ bimodule in the canonical way,
where $B_X=\overline{\lk X, X\rk}_B$.  It is then clear
that the map of $_A(X\otimes_{B_X}B_X)_B$ into ${_AX_B}$
determined by $x\otimes b\mapsto x\d b$
is a right-Hilbert bimodule isomorphism.
If $_AX_B$ is a right-partial $A-B$ imprimitivity bimodule, then
$_AX_{B_X}$ is an imprimitivity bimodule, and the
above isomorphism
gives a factorization
$$[{_AX_B}]=[{}_{B_X}B_B]\circ[_AX_{B_X}]$$
of $[{_AX_B}]$ as the product of the equivalence $[_AX_{B_X}]$
and the (probably degenerate) standard morphism $[{}_{B_X}B_B]$ associated
to the inclusion of the ideal $B_X$ into $B$.

More generally, if $({_AX_B},\gamma)$ is a right-Hilbert
bimodule action, then the above isomorphism is easily seen to be
$\gamma-\gamma\otimes_{B_X}\beta_X$ equivariant, where
$\beta_X$ denotes the restriction of $\beta$ to $B_X$.
Thus we get a similar factorization for
right-partial imprimitivity bimodule actions.
As usual, the case of coactions is a bit more complicated,
so we do it in a lemma. Interestingly, it seems to be necessary to
assume nondegeneracy of all coactions to do this step.

\begin{lem}\label{lem-idealcoact}
   Assume that ${}_{(A,\delta)}(X,\zeta)_{(B,\epsilon)}$ is a
   nondegenerate
coaction and let $B_X=\overline{\lk X,X\rk}_B$. Then\textup:
\begin{enumerate}
     \item There exists a unique nondegenerate coaction $\epsilon_X\:B_X\to
     M(B_X\otimes C^*(G))$
     such that $({_A}X_{B_X},\zeta)$
     becomes a nondegenerate $\delta-\epsilon_X$ compatible coaction 
on $_AX_{B_X}$.
     \item The inclusion $B_X\hookrightarrow B$ is $\epsilon_X-\epsilon$
     equivariant, and
     $${}_{(A,\delta)}(X,\zeta)_{(B,\epsilon)}\cong
     {}_{(A,\delta)}(X\otimes_{B_X}B_X,
     \zeta\cotimes_{B_X}\epsilon_X)_{(B,\epsilon)}$$
     via the canonical isomorphism $_A(X\otimes_{B_X}B_X)_B\cong {_AX_B}$.
     \item If $\epsilon$ is normal
     then so is $\epsilon_X$.
\end{enumerate}
\end{lem}
\begin{proof}
     Of course, we want to define $\epsilon_X$ as the restriction
     of $\epsilon$ to $B_X$. For this to make sense, we have to
     show, similarly to the proof of Lemma \ref{lem-coactstandard}, that
     $\epsilon(B_X)\subseteq M_G(B_X\otimes C^*(G))$, where we view
     $M_G(B_X\otimes C^*(G))$ as an ideal of $M_G(B\otimes C^*(G))$
     via the canonical inclusion (see  \corref{MCideal}).
     But this follows from
     $$\epsilon(B_X)=\epsilon(\lk X, X\rk_B)\subseteq
     \lk M_G(X\otimes C^*(G)), M_G(X\otimes C^*(G))\rk_{M_G(B\otimes
     C^*(G))},$$
     and the fact that the range of the $M_G(B\otimes C^*(G))$-valued
     inner product on $M_G(X\otimes C^*(G))$ lies in $M_G(B_X\otimes
     C^*(G))$ (see \lemref{lem-coefficient}).
     The coaction identity of $\epsilon_X$ follows directly from the
     coaction identity for $\epsilon$. Thus, to see that $\epsilon_X$
     is a nondegenerate coaction, it only remains to check that
     $\overline{\epsilon_X(B_X)(1\otimes C^*(G))}=B_X\otimes C^*(G)$. For this
     we compute
     \begin{align*}
	\overline{\epsilon_X(B_X)(1\otimes C^*(G))}&=
	\overline{\lk\zeta(X),\zeta(X)\rk_{M(B\otimes C^*(G))}(1\otimes
	C^*(G))}\\
	&=\overline{\lk\zeta(X),\zeta(X)(1\otimes
	C^*(G))\rk}_{B\otimes C^*(G)}\\
         &=\overline{\lk\zeta(X),X\otimes C^*(G))\rk}_{B\otimes C^*(G)}\\
	&=\overline{\lk(1\otimes C^*(G))\zeta(X),X\otimes 
C^*(G))\rk}_{B\otimes C^*(G)}\\
	&=\overline{\lk X\otimes C^*(G),X\otimes C^*(G))\rk}_{B\otimes 
C^*(G)}\\
	&=B_X\otimes C^*(G).
     \end{align*}

     By \lemref{lem-coefficient} we can canonically identify
     $M_G({_AX_B}\otimes C^*(G))$ with $M_G({_AX_{B_X}}\otimes C^*(G))$,
     which now shows that $\zeta$ can be viewed as a
     $\delta-\epsilon_X$ compatible coaction on $_AX_{B_X}$.
     This proves (i).

     It is clear from the construction of $\epsilon_X$ that
     the inclusion $B_X\hookrightarrow B$ is $\epsilon_X-\epsilon$ equivariant.
     For~(ii), it only remains to show that the canonical isomorphism
     $X\otimes_{B_X}B_X\cong X$ is
     $\zeta\cotimes_{B_X}\epsilon_X-\zeta$ equivariant.
     But this follows from the last paragraph in the proof of Theorem
     \ref{CG-cat-thm}.

Finally, using the identifications made in \corref{MCideal}, 
$(\id_{B_X}\otimes\lambda)\circ \epsilon_X$ is precisely the
restriction to $B_X$ of $(\id_B\otimes\lambda)\circ\epsilon$. 
Thus, if $\epsilon$ is normal 
(see \remref{rem-normal}), then $\epsilon_X$ is too.

     Finally, if $\epsilon$ is normal, and if 
$\iota\:B_X\hookrightarrow B$ is the
     inclusion map, it follows from \propref{prop-Cmultextend} that
     $\iota\otimes \id_{\K}$, with $\K=\K(L^2(G))$, extends to an 
injective inclusion of
     $M_{\K}(B_X\otimes \K)$ into $M_{\K}(B\otimes\K)$, and then we get
     $$(\iota\otimes\id_{\K})\circ (\id_{B_X}\otimes\lambda)\circ\epsilon_X
     =(\id_B\otimes\lambda)\circ\epsilon\circ\iota,$$
where $\lambda\colon C^*(G)\to\B(L^2(G))=M(\K)$ is the regular
representation.
     If $\epsilon$ is normal, the right side of this equality is
     injective, which then implies that
     $(\id_{B_X}\otimes\lambda)\circ\epsilon_X$ is injective, too.
     But 
(see \propref{prop-ind} and \remref{rem-normal})
this implies that $\epsilon_X$ is normal.
\end{proof}

\section{Morphisms and induced representations}

In this section we  give an outline of
the relationship between our categories $\C$, $\A(G)$, $\C(G)$,
and $\A\C(G)$, and (covariant) representations on Hilbert space.
As always in this paper, we assume that the inner product
on a Hilbert space $\H$ is linear in the second and conjugate
linear in the first variable.
Thus, if $A$ is a $C^*$-algebra, it follows from this convention that
a nondegenerate representation $\pi\:A\to \B(\H)$ gives
$\H$ the structure of a right-Hilbert $A-\CC$ bimodule.
Moreover, two such modules $\H_1$ and $\H_2$
are isomorphic as right-Hilbert $A-\CC$ bimodules if and only
if there exists a unitary $U\:\H_1\to \H_2$ which
intertwines the left $A$-actions, \ie, if and only if the
corresponding nondegenerate representations of
$A$ on $\H_1$ and $\H_2$ are unitarily equivalent.
It follows that the unitary equivalence classes of
nondegenerate Hilbert space representations of $A$ are
precisely the morphisms from $A$ to $\CC$ in the
category $\C$. Note also that $\H=\{0\}$ is the only
right-Hilbert $A-\CC$ bimodule which is not full.

\begin{defn}\label{defn-rep}
     Let $A$ be a $C^*$-algebra.
     We denote by $\Rep(A)$ the class $\Mor(A,\CC)$
     of morphisms from $A$ to $\CC$ in the category $\C$.
     Similarly, if $(A,\alpha)$ is an action of a locally compact
     group $G$ on $A$, then we put
     $$\Rep(A,\alpha)=\Mor\big((A,\alpha), (\CC,\id)\big)$$
     in $\A(G)$, and if
     $(A,\delta)$ is a coaction of $G$, we put
     $$\Rep(A,\delta)=\Mor\big((A,\delta),(\CC,\id\otimes 1)\big)$$
     in $\C(G)$. (Here $\id\otimes 1\:\CC\to M(\CC\otimes C^*(G))$
     denotes the trivial coaction of $G$ on $\CC$.)
     We similarly define $\Rep(A,\alpha,\delta)$
     for our mixed category $\A\C(G)$.
\end{defn}

We already saw above that $\Rep(A)$ coincides with the class of
all unitary equivalence classes of nondegenerate
Hilbert-space representations
of the $C^*$-algebra $A$ (including the zero-representation
$\H=\{0\}$). We shall later check that
$\Rep(A,\alpha)$ and $\Rep(A,\delta)$ are also
precisely the unitary equivalence classes of nondegenerate
covariant representations of $(A,\alpha)$ and $(A,\delta)$,
respectively (also including the zero-representation).

But before we do this, we note that the definition of
$\Rep(A)$ as $\Mor(A,\CC)$ in the category $\C$ directly gives
us a procedure to use any right-Hilbert $B-A$ bimodule
$_BX_A$ to induce representations from $A$ to $B$:
we simply compose the morphisms in $\Rep(A)$
with the morphism $[X]$ to obtain the map
$$\Rep(A)\to \Rep(B)\colon  [\H]\mapsto [\H]\circ [X]=[X\otimes_A\H],$$
which is a bijection if $[X]$ is an equivalence in $\C$,
\ie, if $X$ is a $B-A$ imprimitivity bimodule.
Exactly the same arguments work in the equivariant cases:
  if $(X,\gamma)$ is a $\beta-\alpha$
compatible right-Hilbert bimodule action, we get an induction map
$$\Rep(A,\alpha)\to \Rep(B,\beta)\colon  [(\H,U)]\mapsto
[(\H,U)]\circ[(X,\gamma)],$$
where composition is in $\A(G)$,
and if $(X,\zeta)$ is an $\epsilon-\delta$ compatible
right-Hilbert bimodule coaction we get a
map
$$\Rep(A,\delta)\to \Rep(B,\epsilon)\colon  [(\H,\eta)]\mapsto
[(\H,\eta)]\circ[(X,\zeta)],$$
where composition is in $\C(G)$; the maps are bijections
if the respective morphisms are equivalences in the appropriate
categories. 
Of course, similar
observations can be made for $\A\C(G)$.

We now turn to identifying the spaces
$\Rep(A,\alpha)$ and $\Rep(A,\delta)$ with
the equivalence classes of the covariant representations
of $(A,\alpha)$ and $(A,\delta)$, respectively.
We start with the easy case of an action:

\begin{lem}\label{lem-repaction}
     Assume that  $(\pi, U)$ is a covariant representation
     of $(A,\alpha)$ on a Hilbert space $\H_{\pi}$.
     Then $(\H_{\pi}, U)$ is an $\alpha-\id$ compatible
     right-Hilbert bimodule action, and the assignment
     $$[(\pi,U)]\to [(\H_{\pi},U)]$$
     is a one-to-one correspondence between the unitary equivalence
     classes of nondegenerate covariant representations
     of $(A,\alpha)$ \textup(including the zero-representation\textup)
     and $\Rep(A,\alpha)$.
\end{lem}
\begin{proof}  Note that in our definition of covariant
representations $(\pi,U)$, we already assume that
$\pi\:A\to \B(\H_{\pi})$ is nondegenerate (see
Section~\ref{sec-actions} of Appendix~\ref{coactions-chap}). Thus
it follows from the definition of covariance and
    the definition of right-Hilbert bimodule
    actions (see Definition \ref{rHb-act-defn})
     that $(\pi,U)\mapsto (\H_{\pi},U)$
     is actually a one-to-one correspondence between the
     nondegenerate covariant Hilbert-space representations of $(A,\alpha)$
     and the $\alpha-\id$ compatible
     right-Hilbert $A-\CC$ bimodule actions.
     It is then straightforward to check that
     two covariant representations are unitarily equivalent if and
     only if the corresponding bimodule actions are isomorphic.
\end{proof}

As usual, the coaction case requires a bit more work.
Recall from Definition \ref{def-covariant} that a
covariant representation of a coaction $(A,\delta)$
on a Hilbert space $\H_{\pi}$
consists of a pair $(\pi,\mu)$ such that
$\pi\:A\to \B(\H_{\pi})$ and $\mu\:C_0(G)\to \B(\H_{\pi})$
are nondegenerate representations satisfying
$$(\pi\otimes\id_G)\circ \delta(a)=(\mu\otimes\id_G)(w_G)(\pi(a)\otimes 1)
(\mu\otimes\id_G)(w_G)^*,$$
where $w_G\in C^b(G, M^{\beta}(C^*(G)))=M(C_0(G)\otimes C^*(G))$
denotes the canonical embedding 
$s\mapsto u(s)$ of
$G$ into $UM(C^*(G))$.

\begin{lem}\label{lem-repcoaction}
     Assume that $(\pi,\mu)$ is a nondegenerate covariant
     representation of the coaction $(A,\delta)$.
     Then the map
     $\zeta_{\mu}\:\H_{\pi}\to M(\H_{\pi}\otimes C^*(G))$ defined by
\[
     \zeta_{\mu}(x)=(\mu\otimes \id_G)(w_G)(x\otimes 1)
\]
     is a $\delta-\id\otimes 1$ compatible coaction
     on $\H_{\pi}$.
\end{lem}
\begin{proof}
     We first check that $\zeta_{\mu}$ is indeed a
     right-Hilbert bimodule coaction in the sense of Definition
     \ref{rHb-co-defn}. To see that
     $\zeta_{\mu}(\H_{\pi})\subseteq M(\H_{\pi}\otimes C^*(G))$,
     note that $(\mu\otimes \id_G)(w_G)\in
     M(\K(\H_{\pi})\otimes C^*(G))=
     \L_{\CC\otimes C^*(G)}(\H_{\pi}\otimes C^*(G))$.
     Moreover, if $x\in \H_{\pi}$,
     then $x\otimes 1\in M(\H_{\pi}\otimes C^*(G))=
     \L_{\CC\otimes C^*(G)}(\CC\otimes C^*(G),\H_{\pi}\otimes C^*(G))$.
     Thus it is clear that
     $$\zeta_{\mu}(x)=(\mu\otimes \id_G)(w_G)(x\otimes 1)\in
     \L_{\CC\otimes C^*(G)}(\CC\otimes C^*(G), \H_{\pi}\otimes C^*(G))=
     M(\H_{\pi}\otimes C^*(G)).$$

     Next we check that
     $\zeta_{\mu}\:\H_{\pi}\to M(\H_{\pi}\otimes C^*(G))$
     is a bimodule homomorphism with
     coefficient maps $\delta$ and $\id\otimes 1$.
     Indeed, if we write $a\cdot x$ for $\pi(a)x$,
     the covariance condition for $(\delta,\mu)$
     implies that
     \begin{align*}
	\zeta_{\mu}(a\cdot x)&=(\mu\otimes \id_G)(w_G)(a\cdot x\otimes 1)
     =(\mu\otimes \id_G)(w_G)(\pi(a)\otimes 1)(x\otimes 1)\\
     &=(\pi\otimes\id_G)\circ \delta(a)\big((\mu\otimes
     \id_G)(w_G)(x\otimes 1)\big)
     =\delta(a)\cdot\zeta_{\mu}(x)
     \end{align*}
     for all $a\in A$ and $x\in \H_{\pi}$. On the other side, 
     viewing bimodule elements as adjointable operators, we
     have
     \begin{align*}
     \lk \zeta_{\mu}(x),\zeta_{\mu}(y)\rk_{\CC\otimes C^*(G)}&=
     (x\otimes 1)^*(\mu\otimes\id_G)(w_G)^*(\mu\otimes\id_G)(w_G)(y\otimes
     1)\\
     &=x^*y\otimes 1=\lk x,y\rk_{\CC}\otimes 1
     =(\id\otimes 1)(\lk x,y\rk_{\CC})
     \end{align*}
     for all $x,y\in \H_{\pi}$. 

Since $(\mu\otimes\id_G)(w_G)$
     is actually a unitary operator on $\H_{\pi}\otimes C^*(G)$,
     we get
     \begin{align*}
	\overline{\zeta_{\mu}(\H_{\pi})(\CC\otimes C^*(G))}&
     =(\mu\otimes \id_G)(w_G)\overline{\big((H_{\pi}\otimes 1)(1\otimes
     C^*(G))\big)}\\
     &=(\mu\otimes\id_G)(w_G)(\H_{\pi}\otimes C^*(G))\\
     &=\H_{\pi}\otimes C^*(G),
     \end{align*}
     which proves that $\zeta_{\mu}$ is nondegenerate.

     We now check that
     $(1\otimes C^*(G))\zeta_{\mu}(\H_{\pi})\subseteq \H_{\pi}\otimes
     C^*(G)$. If $z\in C_c(G)\subseteq C^*(G)$, then
     $(1\otimes z)(w_G)\in M(C_0(G)\otimes C^*(G))$ is given by
     the continuous function $s\mapsto zu(s)$ of $G$ into $C^*(G)$.
     Since $u(s)$ is the canonical image of $s$
     in $UM(C^*(G))$, it follows from the formulas given in
     Remark \ref{rem-universal} (where $u(s)$ is denoted $i_G(s)$)
     that $zu(s)\in C_c(G)\subseteq C^*(G)$ and
     $\big(zu(s)\big)(t)=\Delta(s^{-1})z(ts^{-1})$ for $t\in G$.
     It follows that the function
     $s\mapsto zu(s)$ lies in $C_c(G, C^*(G))$ for all $z\in
     C_c(G)\subseteq C^*(G)$; hence we see that
     $(1\otimes C^*(G))w_G\subseteq C_0(G)\otimes C^*(G)$.
     From this we obtain
     \begin{align*}
     (1\otimes C^*(G))\zeta_{\mu}(\H_{\pi})&=
     (1\otimes C^*(G))(\mu\otimes \id_G)(w_G)(\H_{\pi}\otimes 1)\\
     &=(\mu\otimes\id_G)\big((1\otimes C^*(G))w_G\big)(\H_{\pi}\otimes
     1)\\
     &\subseteq (\mu\otimes\id_G)(C_0(G)\otimes C^*(G))(\H_{\pi}\otimes
     1)\\
     &= (\mu\otimes\id_G)(C_0(G)\otimes C^*(G))(1\otimes
     C^*(G))(\H_{\pi}\otimes 1)\\
     &\subseteq (\mu\otimes\id_G)(C_0(G)\otimes
     C^*(G))(\H_{\pi}\otimes C^*(G))\\
     &\subseteq \H_{\pi}\otimes C^*(G).
     \end{align*}

     So it only remains to check the coaction identity for
     $\zeta_{\mu}$.  For this recall from Proposition \ref{slice-unitary}
     that if $w\deq (\mu\otimes \id_G)(w_G)$, then we have
     $w_{12}w_{13}=(\mu\otimes\delta_G)(w_G)$,
     with notation as in that proposition.
     Using this equation we compute
     \begin{align*}
     (\zeta_{\mu}\otimes \id_G)\circ \zeta_{\mu}(x)&=
     (\zeta_{\mu}\otimes \id_G)\big(w(x\otimes 1)\big)\\
     &=w_{12}w_{13}(x\otimes 1\otimes 1)\\
     &=(\mu\otimes\delta_G)(w_G)(x\otimes 1\otimes 1)\\
     &=(\id_{\H_{\pi}}\otimes\delta_G)\big((\mu\otimes
     \id_G)(w_G)(x\otimes 1)\big)\\
     &=(\id_{\H_{\pi}}\otimes\delta_G)\circ \zeta_{\mu}(x)
     \end{align*}
     for all $x\in X$. This completes the proof of the lemma.
\end{proof}

Thus, every covariant representation $(\pi,\mu)$
of $(A,\delta)$
determines a $\delta-\id\otimes 1$ equivariant coaction
on $\H_{\pi}$. Conversely, we get:

\begin{lem}\label{lem-conversecoact}
     Suppose that $(\H, \zeta)$ is a $\delta-\id\otimes 1$
     right-Hilbert $A-\CC$ bimodule coaction.
     Let $\pi\:A\to \B(\H)$ be the nondegenerate representation
     coming from the action of $A$ on $\H$.
     Then there exists a unique nondegenerate representation
     $\mu\colon C_0(G)\to \B(\H)$ such that $(\pi,\mu)$ is
     a covariant representation of $(A,\delta)$ on $\H$
     and such that $\zeta=\zeta_{\mu}$ with $\zeta_{\mu}$
     as in Lemma~\textup{\ref{lem-repcoaction}}.
\end{lem}
\begin{proof} We may assume without loss of generality that
     $\H\neq\{0\}$.
     We first show that there exists a unitary $w\in U(\H\otimes
     C^*(G))$ such that
     $\zeta(x)=w(x\otimes 1)$ for all $x\in \H$.
     For this we first define a map
     $w\:\H\odot C^*(G)\to \H\otimes C^*(G)$ by
     $$w(x\otimes u)=\zeta(x)(1\otimes u).$$
     Note that $\zeta(x)(1\otimes u)\in \H\otimes C^*(G)$,
     since $\zeta(\H)\subseteq M_G(\H\otimes C^*(G))$.
     Moreover,
     $w$ preserves the $\CC\otimes C^*(G)$-valued
     inner products: for $x,y\in \H$ and $u,v\in C^*(G)$
     we have
     \begin{align*}
	\lk w(x\otimes u), w(y\otimes v)\rk_{\CC\otimes C^*(G)}&=
	\lk \zeta(x)(1\otimes u), \zeta(y)(1\otimes v)\rk_{\CC\otimes
	C^*(G)}\\
	&=(1\otimes u^*)\lk \zeta(x),\zeta(y)\rk_{M(\CC\otimes C^*(G))}
	(1\otimes v)\\
	&{=}(1\otimes u^*)\big((\id_{\CC}\otimes 1)(\lk x,y\rk_{\CC})\big)
	(1\otimes v)\\
	&{=}(1\otimes u^*)\big(\lk x,y\rk_{\CC}\otimes 1\big)
	(1\otimes v)\\
	&=\lk x,y\rk_{\CC}\otimes u^*v\\
&=\lk x\otimes u, y\otimes v\rk_{\CC\otimes C^*(G)}.
     \end{align*}
     It follows that $w$ extends to an isometry
     $w\:\H\otimes C^*(G)\to \H\otimes C^*(G)$, which is surjective since
     $w(\H\otimes C^*(G))\supseteq\zeta(\H)(1\otimes C^*(G))
     =\zeta(\H)(\CC\otimes C^*(G))=\H\otimes C^*(G)$
     by nondegeneracy of $\zeta$.

    Now, if $(c_i)_i$ is a bounded approximate identity of $C^*(G)$, we
    get
    $$\zeta(x)=\lim_i\zeta(x)(1\otimes c_i)=\lim_iw(x\otimes
    c_i)=w(x\otimes 1)$$
    (note that both $\zeta(x)$ and $w(x\otimes 1)$ are in
    $M_G(\H\otimes C^*(G))$ --- compare with the proof of
    Lemma~\ref{lem-Cstrictcont}).
    Thus we have $\zeta(x)=w(x\otimes 1)$ for all $x\in \H$.
    Using the coaction identity for $\zeta$, it follows
    that
    \begin{align*}
        &w_{12}w_{13}(x\otimes 1\otimes 1)
      =(\zeta\otimes \id_G)\circ \zeta(x)\\
      &=(\id_{\H}\otimes\delta_G)\circ \zeta(x)
     =(\id_{\K(\H)}\otimes\delta_G)(w)(x\otimes 1\otimes 1)
     \end{align*}
     for all $x\in \H$, from which we get the equation
     $w_{12}w_{13}=(\id_{\K(\H)}\otimes\delta_G)(w)$.
     It follows then from Remark \ref{rem-slice} that
     there exists a unique nondegenerate homomorphism $\mu\:C_0(G)\to
     \B(\H_{\pi})$
     such that $w=(\mu\otimes\id_G)(w_G)$. It clearly
     follows from our constructions
     that $\mu$ is then uniquely determined by the
     property that $\zeta=\zeta_{\mu}$.

     It only remains to show that
     $(\pi,\mu)$ is a covariant homomorphism
     of $(A,\delta)$, where $\pi\:A\to \B(\H)$ comes from the left
     action of $A$ on $X$. For this we just compute
     \begin{multline*}
	\big((\pi\otimes \id_G)\circ \delta(a)\big)(w(x\otimes 1))=
     \big((\pi\otimes \id_G)\circ \delta(a)\big)\zeta(x)\\
     =\zeta(\pi(a)x)=w(\pi(a)x\otimes 1)
     \big(w(\pi(a)\otimes 1)w^*\big)w(x\otimes 1),
     \end{multline*}
     from which it follows that $(\pi\otimes \id_G)\circ \delta(a)
     =w(\pi(a)\otimes 1)w^*$ for all $a\in A$.
     Since $w=(\mu\otimes\id_G)(w_G)$, this just means
     that $(\pi,\mu)$ is covariant.
\end{proof}

We now combine the above results to get:

\begin{prop}\label{prop-repcoaction}
     Let $(A,\delta)$ be a coaction. For a covariant
     representation $(\pi,\mu)$ of $(A,\delta)$ let
     $(\H_{\pi}, \zeta_{\mu})$ be the corresponding
     $\delta-\id\otimes 1$ compatible
     right-Hilbert $A-\CC$ bimodule coaction as in
     Lemma~\textup{\ref{lem-repcoaction}}.
Then
     $$[(\pi,\mu)]\mapsto [(\H_{\pi},\zeta_{\mu})]$$
     is a one-to-one correspondence between the unitary
     equivalence classes of covariant Hilbert-space representations
     of $(A,\delta)$ and $\Rep(A,\delta)$.
\end{prop}
\begin{proof}
     Having Lemmas \ref{lem-repcoaction} and \ref{lem-conversecoact}
     at hand, we only have to show that
     two representations $(\pi,\mu)$ and $(\rho,\nu)$
     of $(A,\delta)$ are unitarily equivalent if and only if
     $(\H_{\pi},\zeta_{\mu})$ and $(\H_{\rho},\zeta_{\nu})$
     are isomorphic right-Hilbert bimodule coactions.
     So let $U\:\H_{\pi}\to \H_{\rho}$ be a unitary such that
     $U\pi(a)=\rho(a)U$ for all $a\in A$
     and $U\mu(f)=\nu(f)U$ for all $f\in C_0(G)$.
     It then follows that
     $$(U\otimes
     1)(\mu\otimes\id_G)(w_G)=(\nu\otimes\id_G)(w_G)(U\otimes 1).$$
     This implies that
     \begin{multline*}
	\zeta_{\nu}(Ux)=(\nu\otimes \id_G)(w_G)(Ux\otimes 1)
	=(\nu\otimes \id_G)(w_G)(U\otimes 1)(x\otimes 1)\\
	=(U\otimes 1)(\mu\otimes\id_G)(w_G)(x\otimes 1)
	=(U\otimes 1) \zeta_{\mu}(x),
     \end{multline*}
      which means that $U$ is a $\zeta_{\mu}-\zeta_{\nu}$
     compatible isomorphism between $\H_{\pi}$ and $\H_{\rho}$.

     Conversely, if $V\:\H_{\pi}\to\H_{\rho}$ is such a compatible
     isomorphism, then we have $V\pi(a)=\rho(a)V$ for all
     $a\in A$, and a short computation as above shows that
     $$(V\otimes
     1)(\mu\otimes\id_G)(w_G)=(\nu\otimes\id_G)(w_G)(V\otimes 1).$$
     Using slice maps $S_f$ for $f$ in the Fourier algebra $A(G)$, 
it follows from
     Proposition \ref{slice-unitary} that
     $$V\mu(f)=S_f\big((V\otimes 1)(\mu\otimes\id_G)(w_G)\big)
     =S_f\big((\nu\otimes\id_G)(w_G)(V\otimes 1)\big)=
     \nu(f)V.$$
     Since $A(G)$ is norm-dense in $C_0(G)$, we see
     that $V$ is a unitary which intertwines $(\pi,\mu)$
     and $(\rho,\nu)$.
\end{proof}

We close this chapter with a short remark on the additive structure
of our categories. Recall that if
$X$ and $Y$ are two right-Hilbert $A-B$ bimodules, then
we can equip the  direct sum
$X\oplus Y$ with the structure of a right-Hilbert $A-B$
bimodule by defining
\begin{align*}
     &\;\;\lk x_1+y_1, x_2+y_2\rk_B\deq \lk x_1, x_2\rk_B+\lk y_1,y_2\rk_B,\\
     &a\d(x+y)=a\d x+a\d y,\quad\text{and}\quad(x+y)\d b=x\d b+y\d b,
\end{align*}
for $x,x_1,x_2\in X, y, y_1,y_2\in Y$, $a\in A$ and $b\in B$.
It is not hard to check that
$[X]+[Y]=[X\oplus Y]$
defines an
additive structure on the morphisms $\Mor(A,B)$ in our category
$\C$, so that
the zero-morphism (\ie, $X=\{0\}$) serves as
an additive identity. It is also straightforward to check that
the distributive laws
$$ [Z]\circ ([X]+[Y])=[Z]\circ [X]+[Z]\circ [Y]
\quad\text{and}\quad
([X]+[Y])\circ [Z]= [X]\circ [Z]+[Y]\circ [Z]$$
hold in $\C$.
Of course we can define similar additive structures on
the categories $\A(G)$, $\C(G)$, and $\A\C(G)$ --- we omit further
details on this. One can even define infinite direct sums
of right-Hilbert bimodules,
and then we get a notion of
infinite sums of morphisms from $A$ to $B$:
if $\{X_i\}_{i\in I}$
is a family of right-Hilbert $A-B$ bimodules, then, as usual,
one defines $\oplus_{i\in I}X_i$ as the completion of
the vector space
$$\{(x_i)_{i\in I}\mid x_i\neq 0 \;\text{for only
finitely many $i\in I$}\}$$
with respect to the $B$-valued inner product
$\lk(x_i), (x'_i)\rk_B=\sum_{i\in I}\lk x_i, x'_i\rk_B$ and
the obvious left and right actions of $A$ and $B$.
We then put $\sum_{i\in I}[X_i]\deq [\oplus_{i\in I}X_i]$.

Note that these additive structures on $\C$, $\A(G)$, $\C(G)$,
and $\A\C(G)$, are something which is not available in the
category of $C^*$-algebras with $*$-homomorphisms as
morphisms, since in general the sum of
two $*$-homomorphisms $\phi_1,\phi_2\:A\to B$
is not a $*$-homomorphism any more.

%
%

\chapter{The Functors}
\label{functors-chap}

In this chapter we show that many fundamental $C^*$-algebraic
constructions---including restricting, inflating, decomposing, and
taking crossed products by actions and coactions---have
right-Hilbert bimodule counterparts, and that these constructions give
functors among the various categories $\c A(G)$, $\c C(G)$, and $\c A
\c C(G)$.

\section{Crossed products}
\label{xpr}

\subsection{Actions}

In \thmref{xpr-fun} below we will define a functor from $\c A(G)$ to
$\c C(G)$ with object map
$(A,\alpha)\mapsto (A\times_r G,\hat\alpha)$.
Because the objects in $\c A(G)$ are actions on $C^*$-algebras, rather
than isomorphism classes of such, we need to choose a single version
of the reduced crossed product and stick with it.  We use
\defnref{def-red crossed}: the action $\alpha$ can
also be regarded as a nondegenerate homomorphism $\alpha\:A\to
C_b(G,A)\subseteq M(A\otimes C_0(G))$ via $\alpha(a)(s) =
\alpha_{s^{-1}}(a)$, and we define
\[
A\times_{\alpha,r} G
= (i^r_A\times i^r_G)(A\times_\alpha G)
\subseteq M(A\otimes \c K(L^2(G))),
\]
where $(i^r_A,i^r_G)$ is the covariant homomorphism
$\bigl((\id_A\otimes M)\circ\alpha,1\otimes\lambda\bigr)$.

The dual coaction of $G$ on $A\times_{\alpha,r} G$ is defined 
(see \exref{ex-dualcoaction} in~\appxref{coactions-chap})
by
\[
\hat\alpha(i^r_A(a) i^r_G(f))
= (i^r_A(a)\otimes 1) (i^r_G\otimes u)(f)
\righttext{for} a\in A, f\in C_c(G),
\]
where $u\:G\to UM(C^*(G))$ is the canonical inclusion.
Recall that we require all coactions in $\c
C(G)$ to be normal and nondegenerate; every dual coaction
$\hat\alpha$ has these properties (see \exref{ex-dualcoaction} and
\propref{prop-dualnormal}).

We find it convenient to do many of our calculations
with vector-valued $C_c$-functions.  To facilitate this, by a slight
abuse of notation we \emph{identify} $C_c(G,A)$ with its image
$(i^r_A\times i^r_G)(C_c(G,A))$ in $A\times_{\alpha,r} G$.  The
calculations will sometimes involve integrals of functions with values
in locally convex spaces, rather than just $C^*$-algebras.  Because we
(almost always) only need to integrate continuous functions of compact
support, the standard theorems about vector-valued integration apply,
and we have chosen to use this theory to avoid having to insert
$C^*$-algebra elements and linear functionals at every turn; the
technical foundations we require are laid out in
\appxref{indlim-chap}.  As an example of how such integrals arise,
notice that the embedding of $C_c(G,A)$ in $A\times_{\alpha,r} G$ is
described in terms of the canonical embeddings $(i^r_A,i^r_G)$ by the
strictly convergent integral
\[
f=\int_G i^r_A(f(s)) i^r_G(s)\,ds.
\]
As a point of notation, we usually identify $G$ with its canonical image
in $UM(C^*(G))$, so that, for example, we can say that the dual coaction
of $G$ on $A \times_{\alpha,r} G$ is given on $C_c(G,A)$ by the strictly
convergent integral
\[
\hat\alpha(f)=\int_G i^r_A(f(s)) i^r_G(s)\otimes s\,ds.
\]

We now show how to regard $\hat\alpha(f)$ as a $C_c$-function.
To prepare for this, we require some technical flexibility involving
tensor products.  If $B$ is a $C^*$-algebra, we (often without
comment) make extensive use of the canonical embedding of the
algebraic tensor product $C_c(G,A) \odot B$ into $C_c(G,A\otimes B)$
given by
$(f\otimes b)(s)=f(s)\otimes b$.
(We occasionally extend this convention to other function spaces, such as
$C_b$, the bounded continuous functions.)
Note that this
embedding extends to the canonical isomorphism
of $(A\times_{\alpha,r} G)\otimes B$ onto 
$(A\otimes B)\times_{\alpha\otimes \id,r} G$
from \lemref{lem-mintensor}.
We employ this trick many times throughout this work when
dealing with coactions in terms of $C_c$-functions. The idea is to treat
the extra copy of $C^*(G)$ as a freely moving object.

This allows us to embed $C_c(G, A\otimes B)$ in $(A\times_{\alpha,r}
G)\otimes B$, and $C_c(G, M^\beta(A\otimes B))$ in $M((A
\times_{\alpha,r} G)\otimes B)$; for the latter we use the embedding
of $C_c(G, M^{\beta}(A))$ into $M(A\times_{\alpha,r} G)$ provided
by \corref{corconv}. (Recall that we write $M^{\beta}(A)$
for $M(A)$ equipped with the strict topology.)
Note that for $f\in C_c(G,A)$ and
$b\in B$, we have
\[
\begin{split}
f\otimes b
&
=\int_G i_A^r(f(s))i_G^r(s)\,ds\otimes b
=\int_G i_A^r(f(s))i_G^r(s)\otimes b\,ds
\\&
=\int_G \bigl(i_A^r(f(s))\otimes b\bigr)(i_G^r(s)\otimes 1)\,ds
=\int_G (i_A^r\otimes\id)\bigl((f\otimes b)(s)\bigr)
(i_G^r(s)\otimes 1)\,ds,
\end{split}
\]
so that, by linearity, density and continuity, if $f\in
C_c(G,M^\beta(A\otimes B))$ then as an element of $M((A
\times_{\alpha,r} G)\otimes B)$ we have (by the same abuse of
notation as above)
\[
f=\int_G (i_B^r\otimes \id)(f(s)) (i_G^r(s)\otimes 1)\,ds.
\]

\begin{lem}
\label{alpha hat}
Let $(A,G,\alpha)$ be an action.  For each $f\in C_c(G,A)$, we have
\[
\hat\alpha(f)\in \mathbox{C_c(G,M^\beta(A\otimes C^*(G)))}\subseteq
M\big((A\times_r G)\otimes C^*(G)\big),
\]
with 
$\hat\alpha(f)(s)=f(s)\otimes s$ for $s\in G$.
\end{lem}

\begin{proof}
Compute:
\[
\begin{split}
\hat\alpha(f)
&
=\hat\alpha\left(\int_G i_A^r(f(s)) i_G^r(s)\,ds \right)
=\int_G i_A^r(f(s)) i_G^r(s)\otimes s\,ds
\\&
=\int_G (i_A^r\otimes \id)(f(s)\otimes s)
(i_G^r(s)\otimes 1)\,ds.
\end{split}
\]
Thus, $\hat\alpha(f)$ agrees with the element $s\mapsto f(s)\otimes
s$ of $C_c(G,M^\beta(A\otimes C^*(G)))$.
\end{proof}

Functoriality requires that we also give a construction of crossed
products of right-Hilbert bimodules.  Such constructions have appeared
in several places in the literature (\eg, see
\cite{CMW-CP, com, kas, BS-CH, ekqr:green}
for several constructions of full and reduced
crossed products by Hilbert modules).
In order to have a construction which is best suited
for our needs, and for completeness, we shall give our own construction
of reduced crossed products by
a right-Hilbert bimodule actions. Of course we closely follow
the ideas presented in the literature cited above.

Let
$_{(A,\alpha)}(X,\gamma)_{(B,\beta)}$ be a right-Hilbert bimodule
action. Consider the dense subalgebras
$C_c(G,A)$ and $C_c(G,B)$ of $A\times_{\alpha,r}G$ and
$B\times_{\beta,r}G$.
We want to equip $C_c(G,X)$ with a
pre-right-Hilbert $C_c(G,A)-C_c(G,B)$ bimodule structure with operations
\begin{equation}
\begin{split}
\label{X x G}
f\d x(s)
&=\int_G f(t)\d\gamma_t(x(t^{-1}s))\, dt,
\\
x\d g(s)
&=\int_G x(t)\d\beta_t(g(t^{-1}s))\, dt,\quad\text{and}
\\
\<x,y\>_{B\times_\beta G}(s)
&=\int_G\beta_{t^{-1}}\bigl(\<x(t),y(ts)\>_B\bigr)\,dt
\end{split}
\end{equation}
for $f\in C_c(G,A)$, $x,y\in C_c(G,X)$, and $g\in C_c(G,B)$.

\begin{prop}\label{X xr G}
     With the above operations, $C_c(G,X)$ completes
     to give a right-Hilbert $A\times_{\alpha,r}G-B\times_{\beta,r}G$
     bimodule.
Moreover, if $B_X=\overline{\lk X,X\rk}_B$, then
     $\overline{\lk X\times_{\gamma,r}G,
     X\times_{\gamma,r}G\rk}_{B\times_{\beta,r}G}=B_X\times_{\beta,r}G$.
\end{prop}

We do the proof in two steps. In the first step we
consider the right-Hilbert bimodule
action $({_KX_B}, {_\mu\gamma_{\beta}})$, with $K=\K(X)$,
as provided by \propref{decomaction}, and use
a linking algebra argument to see that the proposition is true
in this special situation. We then observe that the
homomorphism $\kappa\times_rG\:A\times_{\alpha,r}G\to
M(K\times_{\mu,r}G)$ coming from the
$\alpha-\mu$
equivariant homomorphism $\kappa\:A\to M(K)\cong \L(X)$
induces a left action of $A\times_{\alpha,r}G$ on $X\times_{\gamma,r}G$
which coincides on $C_c$-functions with the actions
as given in \eqref{X x G}.

If $({_KX_B}, {_\mu\gamma_\beta})$
is a right-Hilbert imprimitivity bimodule action
     of $G$ on the partial imprimitivity
     bimodule $_KX_B$, then we obtain a strongly continuous
     action $\nu\:G\to \Aut L(X)$ by defining
\begin{equation}\label{eq-linkactionformula}
     \nu_s\left(\mtx{k&x\\ \tilde{y}& b}\right)=
\mtx{\mu_s(k) &\gamma_s(x)\\ \widetilde{\gamma_s(y)}& \beta_s(b)}
\end{equation}
(see \lemref{lem-actlink}).
The convolution algebra $C_c(G, L(X))$ has a canonical
decomposition as two-by-two matrices
\begin{equation}\label{eq-linkingCc}
C_c(G, L(X))=\mtx{C_c(G,K) & C_c(G,X)\\ C_c(G,\rev{X})& C_c(G,B)},
\end{equation}
and it is straightforward to check that the pairings among
the corners $C_c(G,K)$, $C_c(G,X)$ and $C_c(G, B)$
given by convolution on $C_c(G, L(X))$ are precisely the ones
given by~\eqref{X x G}.

Now let $p=\smtx{1 &0\\ 0&0}$ and $q=\smtx{0 &0\\ 0&1}$
    denote the corner projections in $M(L(X))$, and let
     $i_L^r(p)$ and $i_L^r(q)$ denote the images of $p$ and $q$
     in $M(L(X)\times_{\nu,r}G)$.
     It follows from \lemref{lem-equivdegenerate} that
     the canonical inclusions of $C_c(G, K)$ and $C_c(G,B)$
     into $C_c(G, L(X))$ provided by \eqref{eq-linkingCc}
     extend to isomorphisms
     $$K\times_{\mu,r}G\cong i_L^r(p)\big(L(X)\times_{\nu,r}G\big)i_L^r(p)
     \quad\text{and}\quad
     B\times_{\beta,r}G\cong i_L^r(q)\big(L(X)\times_{\nu,r}G\big)i_L^r(q),
     $$
     and we have
     $$\overline{C_c(G,X)}=i_L^r(p)\big(L(X)\times_{\nu,r}G\big)i_L^r(q),$$
     where we identify $C_c(G,X)$ with the
     upper right corner of $C_c(G, L(X))$ as in \eqref{eq-linkingCc}.
     In particular, it follows from \lemref{prop-partiallink}
     that
     $C_c(G,X)$ completes to give a
     partial $K\times_{\mu,r}G-B\times_{\beta,r}G$
     imprimitivity bimodule $X\times_{\gamma, r}G$, with actions
     and inner products given by \eqref{X x G}.
     Using \lemref{recogniselink} we also see that
     $L(X)\times_{\nu,r}G$ is then canonically isomorphic to
     $L(X\times_{\gamma,r}G)$.

We now gather the above observations:

\begin{lem}\label{lem-linkingaction}
     Suppose that $({_KX_B},{_\mu\gamma_\beta})$ is a
     right-Hilbert bimodule action of $G$ on the
     partial imprimitivity bimodule $_KX_B$. Then
     $C_c(G,X)$, equipped with the $C_c(G,K)-C_c(G,B)$
     pre-right-Hilbert bimodule structure of \eqref{X x G},
     completes to a
     partial $K\times_{\mu,r}G-B\times_{\beta,r}G$
     imprimitivity bimodule such that the identification
     \eqref{eq-linkingCc} extends to a canonical isomorphism
     $$L(X)\times_{\nu,r}G\cong L(X\times_{\gamma,r}G),$$
     with $\nu$ as in \eqref{eq-linkactionformula}.
     Moreover, the ranges of the left- and right-valued
     inner products on $X\times_{\gamma,r}G$ are given
     by the ideals
     $K_X\times_{\mu,r}G$ and $B_X\times_{\beta,r}G$,
     where $K_X$ and $B_X$ denote the ranges of the inner
     products on $X$.
     In particular,
     $X\times_{\gamma,r}G$ is a
     $K\times_{\mu,r}G-B\times_{\beta,r}G$ imprimitivity bimodule
     if and only if $X$ is a $K-B$ imprimitivity bimodule.
\end{lem}
\begin{proof} Everything but the assertion on the ranges
     of the inner products
     follows from the above considerations.

     For the ranges we first consider the case where $K_X=K$
     and $B_X=B$, \ie, where $X$ is a $K-B$ imprimitivity bimodule.
     It follows then from
     \lemref{prop-partiallink} that $p$ and $q$
     are full projections in $L(X)$,
     which implies that
     $i_L^r(q)$ and $i_L^r(p)$ are full projections
     in $M(L(X)\times_{\nu,r}G)$.
     To see this, compute
     \begin{align*}
	L(X)\times_{\nu,r}G&=\overline{i_G^r(C^*(G))i_L^r(L(X))i_G^r(C^*(G))}\\
	&=\overline{i_G^r(C^*(G))i_L^r(L(X)qL(X))i_G^r(C^*(G))}\\
	&=\overline{\big(i_G^r(C^*(G))i_L^r(L(X))\big)i_L^r(q)
	\big(i_L^r(L(X))i_G^r(C^*(G))\big)}\\
	&=\overline{\big(L(X)\times_{\nu,r}G\big)i_L^r(q)
	\big(L(X)\times_{\nu,r}G\big)}.
     \end{align*}
     Thus it follows that $X\times_{\gamma,r}G$ is a
     $K\times_{\mu,r}G-B\times_{\beta,r}G$ imprimitivity bimodule.

     If $_KX_B$ is a partial imprimitivity bimodule,
     then $\mtx{K_X&X\\ \rev{X}& B_X}$, 
which we denote by $L^{\ip}(X)$,  embeds
     as an ideal in $L(X)$ and therefore
     $L^{\ip}(X)\times_{\nu,r}G=\mtx{K_X\times_{\mu,r}G&
     X\times_{\gamma,r}G\\ (X\times_{\gamma,r}G)\!\!\widetilde{\ \ }&
     B_X\times_{\beta,r}G}$ canonically
     embeds as an ideal in $L(X\times_{\gamma,r}G)$.
     Since this embedding is the identity on $X\times_{\gamma,r}G$,
     and since $_{K_X}X_{B_X}$ is an imprimitivity bimodule,
     it follows that the ranges of the inner products
     on $X\times_{\gamma,r}G$ are precisely
     $K_X\times_{\mu,r}G$ and $B_X\times_{\beta,r}G$, respectively.
\end{proof}

\begin{proof}[Proof of \propref{X xr G}]
     Let $K=\K(X)$. Then it follows from \propref{decomaction}
     that there exists an action $\mu$ of $G$ on $K$
     such that $({_KX_B}, {_\mu\gamma_\beta})$ is a right-Hilbert
     bimodule action on the right-partial $K-A$ imprimitivity bimodule
     $_KX_B$, and such that the homomorphism
     $\kappa\:A\to M(K)\cong \L_B(X)$ coming from the left
     action of $A$ on $X$ is $\alpha-\mu$ equivariant.
     \lemref{lem-linkingaction} implies
     that $C_c(G,X)$ completes to a right-partial
     $K\times_{\mu,r}G-B\times_{\beta,r}G$
     imprimitivity bimodule, and it follows from
     \lemref{lem-equivdegenerate}
     that there exists a homomorphism
     $\kappa\times_rG\:A\times_{\alpha,r}G\to M(K\times_{\nu,r}G)
     \cong \L_{B\times_rG}(X\times_{\gamma,r}G)$
     given on $C_c(G,A)$ by $f\mapsto \kappa\circ f$.
     Using this, it is straightforward to check that the resulting left
     action of $A\times_{\alpha,r}G$ on $X\times_{\gamma,r}G$
     coincides on $C_c(G,A)$ with the action given by
     \eqref{X x G}.
\end{proof}

As with reduced $C^*$-crossed products, we identify $C_c(G,X)$
with its image in $X\times_{\gamma,r} G$.

To complete our functor, we need a dual coaction on $X
\times_{\gamma,r} G$, and for this we will
view $C_c(G, M^\beta(X\otimes C^*(G)))$ as a subspace
of $M((X\times_{\gamma,r} G)\otimes C^*(G))$.
The following general lemma justifies this embedding,
since for any $C^*$-algebra $C$, it
follows from \lemref{lemmor} that we have an embedding
$C_c(G,M^\beta(X\otimes C)) \hookrightarrow M((X\otimes C)
\times_{\gamma\otimes \id, r} G)$.

\begin{lem}
\label{X-mult}
For any action
$_{(A,\alpha)}(X,\gamma)_{(B,\beta)}$ and for any $C^*$-algebra $C$
the canonical embedding of $C_c(G,X) \odot C$ in $C_c(G,X\otimes C)$
extends to an isomorphism
\begin{multline*}
\Phi\:
_{(A\times_{\alpha,r} G)\otimes C}
((X\times_{\gamma,r} G)\otimes C)
_{(B\times_{\beta,r} G)\otimes C}
\\
\iso
_{(A\otimes C)\times_{\alpha\otimes \id,r} G}
\bigl((X\otimes C)\times_{\gamma\otimes \id,r} G\bigr)
_{(B\otimes C)\times_{\beta\otimes \id,r} G}.
\end{multline*}
\end{lem}

\begin{proof}
We know that the embedding
of $C_c(G,A) \odot C$ into $C_c(G,A\otimes C)$
determines an isomorphism of 
$(A\times_{\gamma,r} G)\otimes C$ with 
$(A\otimes C)\times_{\gamma\otimes \id,r} G$,
and similarly for $B$.  A routine calculation
now shows that the embedding of
$C_c(G,X) \odot C \subseteq (X\times_{\gamma,r} G) \otimes C$ 
into 
$C_c(G, X\otimes C) \subseteq (X\otimes C)\times_{\gamma\otimes \id,r} G$
is compatible with these coefficient maps and all
bimodule operations.  Since the image of $C_c(G,X)\otimes C$ is
inductive limit dense in $C_c(G,X\otimes C)$, the result follows.
\end{proof}

\begin{prop}
\label{gamma hat}
If $\gamma$ is an $\alpha-\beta$ compatible action of $G$ on a
right-Hilbert bimodule ${}_AX_B$, then
there is a unique nondegenerate $\hat\alpha -
\hat\beta$ compatible coaction $\hat\gamma$ of
$G$ on $X\times_{\gamma,r} G$ such that for 
each $x\in C_c(G,X)$ we have
\[
\hat\gamma(x)\in C_c(G,M^\beta(X\otimes C^*(G)))
\subseteq M((X\times_{\gamma,r} G)\otimes C^*(G)),
\]
with
$\hat\gamma(x)(s)=x(s)\otimes s$ for $s\in G$.
\end{prop}

\begin{proof}
As for the construction of the crossed product $X\times_{\gamma,r}G$,
we use a linking algebra argument. For this let $K=\K(X)$,
let $\mu\:G\to \Aut K$ denote the action of $G$ on $K$ provided
by \propref{decomaction}, and let $\nu\:G\to \Aut L(X)$ denote
the corresponding action on the linking algebra.
It then follows from \lemref{lem-linkingaction} that
$X\times_{\gamma,r}G$ is a right-partial
$K\times_{\mu,r}G-B\times_{\beta,r}G$ imprimitivity bimodule
and $L(X\times_{\gamma,r}G)\cong L(X)\times_{\nu,r}G$.
By \lemref{alpha hat}, the dual coaction $\hat{\nu}$
maps $C_c(G,L(X))\subseteq L(X)\times_{\nu,r}G$ into
$C_c(G, M^{\beta}(L(X)\otimes C^*(G)))\subseteq
M((L(X)\times_{\nu,r}G)\otimes C^*(G))$ according to the rule
$\hat{\nu}(f)(s)=f(s)\otimes s$
for $s\in G$.
We make the canonical identifications
$$C_c(G, L(X))=\mtx{C_c(G,A)&C_c(G,X)\\ C_c(G,X)\!\!\widetilde{\ \ }
&C_c(G,B)}$$
and
$$C_c\big(G, M^{\beta}(L(X\otimes C^*(G)))\big)=
\mtx{C_c(G, M^{\beta}(A\otimes C^*(G)))& C_c(G, M^{\beta}(X\otimes
C^*(G)))\\
C_c(G, M^{\beta}(X\otimes C^*(G)))\!\!\widetilde{\ \ \ }&
C_c(G, M^{\beta}(B\otimes C^*(G)))},$$
where the latter identification is allowed since the
strict topology on $M(L(X)\otimes C^*(G))$ coincides
with the product of the strict topologies of the corners
by \propref{multsoflink}.
With this identification, the formula for $\hat{\nu}$ becomes
$$\hat{\nu}\left(\mtx{g&x\\ \tilde{y}&l}\right)(s)=
\mtx{g(s)\otimes s& x(s)\otimes s\\ \widetilde{y(s)\otimes s}& 
l(s)\otimes s},$$
for $g\in C_c(G,A)$, $x,y\in C_c(G,X)$, and $l\in C_c(G,B)$.
In particular, it follows from this that $\hat{\nu}(p)=p\otimes 1$
and $\hat{\nu}(q)=q\otimes 1$, where $p,q\in M(L(X\times_{\gamma,r}G))
=M(L(K)\times_{\nu,r}G)$ denote
the usual corner projections.
Thus it follows from \lemref{lem-coactlink} that $\hat\nu$
compresses to a $\hat\mu-\hat\beta$ compatible coaction $\hat\gamma$
on the upper right corner $X\times_{\gamma,r}G$ which
then has to be given by the formula
$$\hat\gamma(x)(s)=x(s)\otimes s$$
for $x\in C_c(G,X)$.
Since $\hat\nu$ is nondegenerate by \exref{ex-reduceddual}, it follows
from \lemref{lem-coactlink} that $\hat\gamma$ is nondegenerate, too.

We now complete the proof by observing that it follows directly from
the formulas for $\hat\alpha$ and $\hat\mu$ that
the homomorphism $\kappa\times_rG\:A\times_{\alpha,r}G\to
M(K\times_{\mu,r}G)$ is $\hat\alpha-\hat\beta$ equivariant
(where $\kappa\:A\to M(K)$ is determined by the left $A$-action on $X$).
This shows that $\hat\gamma$ is 
a nondegenerate $\hat\alpha-\hat\beta$ compatible coaction on
$X\times_{\gamma,r}G$.
\end{proof}

For later use it is important to see that the
crossed product by a standard morphism $({_AX_B},\gamma)$
associated to an $\alpha-\beta$ equivariant homomorphism
$\varphi\:A\to M(B)$ is again standard.

\begin{prop}\label{prop-standard}
     Suppose that $\phi\:A\to M(B)$ is a \textup(possibly degenerate\textup)
     $\alpha-\beta$ equivariant homomorphism, and
     let $({_AX_B},\gamma)$ be the standard
     right-Hilbert $A-B$ bimodule action associated to
     $\phi$ \textup(see Definition~\textup{\ref{def-standard}}\textup).
     Then $_{A\times_rG}(X\times_{\gamma,r}G)_{B\times_rG}$
     is canonically isomorphic to the standard
     right-Hilbert $A\times_{\alpha,r}G-B\times_{\beta,r}G$
     bimodule
     $Y=\phi\times_rG(A\times_{\alpha,r}G)(B\times_{\beta,r}G)$,
     and this isomorphism is equivariant with
     respect to $\hat{\gamma}$ and the restriction
     $\hat{\beta}_Y$ of $\hat{\beta}$ to $Y$.
\end{prop}
\begin{proof}
     Since $_AX_B=\phi(A)B$, we see that
     $$\phi\times_rG(C_c(G,A))C_c(G,B)\supseteq
     (C_c(G)\odot \varphi(A))(C_c(G)\odot B)=C_c(G)\odot X$$
     is inductive limit dense in $C_c(G,X)$.
     The result then follows from the fact that
     the above inclusion of $\phi\times_rG(C_c(G,A))C_c(G,B)$
     preserves the left and right $C_c(G,A)$- and $C_c(G,B)$-actions
     and the $C_c(G,B)$-valued inner products.
     The $\hat{\gamma}-\hat{\beta}_Y$ equivariance
     follows directly from the formulas for the respective
     dual coactions as given in \lemref{alpha hat} and
     \propref{gamma hat}.
\end{proof}

\begin{thm}
\label{xpr-fun}
The object map $(A,\alpha)\mapsto(A\times_{\alpha,r} G,\hat\alpha)$ 
and the
morphism map $[_AX_B,\gamma]\mapsto [_{A\times_{r} G}(X\times_{\gamma,r}
G)_{B\times_r G},\hat\gamma]$ define a functor from $\c A(G)$ to $\c
C(G)$.
\end{thm}

\begin{proof}
We first show that the map on morphisms is well-defined.  Suppose $\Phi\:
X\to Y$ is an isomorphism of right-Hilbert $A-B$ bimodules which is
equivariant for $\alpha -\beta$ compatible actions $\gamma$ and $\rho$
of $G$.  Then $(\Phi\times_rG)(x)(s) \deq  \Phi(x(s))$ is easily seen to
give a bijective map $\Phi\times_rG\:C_c(G,X)\to C_c(G,Y)$ which
respects the pre-right-Hilbert bimodule structures given by \eqeqref{X
x G} and hence extends to a right-Hilbert $(A \times_{\alpha,r}
G)-(B\times_{\beta,r} G)$ bimodule isomorphism of
$X\times_{\gamma,r} G$ onto $Y\times_{\rho,r} G$.  It suffices to
check $\hat\gamma -\hat\rho$ equivariance of $\Phi\times_rG$ for $x
\in C_c(G,X)$:
\begin{multline*}
((\Phi\times_rG)\otimes \id)\circ\hat\gamma(x)(s)
= (\Phi\otimes \id)(\hat\gamma(x)(s))
= \Phi(x(s))\otimes s\\
= (\Phi\times_rG)(x)(s)\otimes s
=\hat\rho\circ (\Phi\times_rG)(x)(s).
\end{multline*}

If $(A,G,\alpha)$ is an
action then it follows from \propref{prop-standard} that
the bimodule crossed product
$_AA_A\times_r G$ is equivariantly isomorphic to
the standard right-Hilbert bimodule $_{A
\times_r G}(A\times_r G)_{A\times_r G}$,
which shows that the crossed-product functor preserves
identities.

It only remains to see that the assignment $[X,\gamma]\mapsto [X
\times_{\gamma,r} G]$ respects composition of morphisms;
\footnote{This is also proven (using non-categorical language) 
in \cite[Lemma 3.10]{kas}; we include our own proof for completeness.} 
that is, if
$(_AX_B,\gamma)$ is $\alpha -\beta$ compatible and $(_BY_C,\rho)$ is
$\beta-\mu$ compatible, we need to find a right-Hilbert $(A
\times_{\alpha,r} G)-(C\times_{\mu,r} G)$ bimodule isomorphism
\[
\Upsilon\:(X\times_{\gamma,r} G)\otimes_{B\times_r G}
(Y\times_{\rho,r} G) \iso
(X\otimes_B Y)\times_{\gamma\otimes \rho,r} G
\]
which is
$(\hat\gamma\cotimes_{B\times_rG}\hat\rho)-\what{\gamma\otimes\rho}$
equivariant.  The rule
\[
\Upsilon(x\otimes y)(s)
=\int_G x(t)\otimes\rho_t(y(t^{-1}s))\,dt
\]
defines a map from $C_c(G,X) \odot C_c(G,Y)$ to $C_c(G,X\otimes_B
Y)$ which is easily checked to preserve the pre-right-Hilbert bimodule
structures.  In order to see that $\Upsilon$ extends to an isomorphism
of the completions, we need only verify that $\Upsilon$ has dense range
for the inductive limit topology.  For this, let $x_1\in X$ and $f
\in C_c(G,B)$, and define $x\in C_c(G,X)$ by $x(s)=x_1\d f(s)$.
Then for $y\in C_c(G,Y)$ we have
\begin{align*}
\Upsilon(x\otimes y)(s)
&
= x_1\otimes\int_G f(t)\d \rho_t(y(t^{-1}s))\,dt
= x_1\otimes (f\d y)(s).
\end{align*}
Now, we can approximate $y$ by $f\d y$ in the inductive limit
topology, and taking $y$ of the form $y(s)=y_1 g(s)$ for $y_1\in Y$
and $g\in C_c(G)$ we can thus approximate the function $s\mapsto
(x_1\otimes y_1) g(s)$.  But such functions have
inductive limit dense span in $C_c(G,X\otimes_B Y)$.

For the equivariance, it suffices to show that for $x\in C_c(G,X)
\subseteq X\times_{\gamma,r} G$ and $y\in C_c(G,Y)\subseteq Y
\times_{\rho,r} G$ we have
\[
\what{\gamma\otimes \rho}\circ \Upsilon(x\otimes y)
= (\Upsilon\otimes \id)\circ
(\hat\gamma\cotimes\hat\rho)(x\otimes y)
\]
as elements of
$C_c\bigl(G,M^\beta\bigl((X\otimes_B Y)\otimes C^*(G)\bigr)\bigr)
\subseteq
M\bigl(\bigl((X\otimes_B Y)\times_{\gamma\otimes \rho,r} G\bigr)
\otimes C^*(G)\bigr)$.
First note that for $x'\in C_c(G,X)$, $y'\in C_c(G,Y)$, and $u,v,s
\in G$ we have
\begin{align*}
&(\Upsilon\otimes \id)\circ \Theta
  \bigl((x'\otimes u)\otimes (y'\otimes v)\bigr)(s)
= (\Upsilon\otimes \id)
  \bigl((x'\otimes y')\otimes uv\bigr)(s)\\
&\qquad= \bigl(\Upsilon(x'\otimes y')\otimes uv\bigr)(s)
= \Upsilon(x'\otimes y')(s)\otimes uv\\
&\qquad= \int_G x'(t)\otimes \rho_t(y'(t^{-1}s))\,dt\otimes uv\\
&\qquad= \int_G\bigl(x'(t)\otimes \rho_t(y'(t^{-1}s))\bigr)\otimes uv\,dt\\
&\qquad= \int_G\Theta\bigl((x'(t)\otimes u)\otimes 
  (\rho_t\otimes \id)(y'(t^{-1}s)\otimes v)\bigr)\,dt\\
&\qquad= \int_G\Theta\Bigl((x'\otimes u)(t)\otimes
  (\rho_t\otimes \id)\bigl((y'\otimes v)(t^{-1}s)\bigr)\Bigr)\,dt.
\end{align*}
By linearity, density, and continuity we conclude that for $x\in
C_c(G,X)$, $y\in C_c(G,Y)$, and $s\in G$ we have
\begin{align*}
&(\Upsilon\otimes \id)\circ (\hat\gamma\cotimes\hat\rho)
(x\otimes y)(s)
= (\Upsilon\otimes \id)\circ \Theta
\bigl(\hat\gamma(x)\otimes\hat\rho(y)\bigr)(s)
\\
&\qquad
=\int_G \Theta\Bigl(
\hat\gamma(x)(t)\otimes
(\rho_t\otimes \id)\bigl(\hat\rho(y)(t^{-1}s)\bigr)
\Bigr)\,dt
\\
&\qquad
=\int_G \Theta\Bigl(
(x(t)\otimes t)\otimes
(\rho_t(y(t^{-1}s))\otimes t^{-1}s\bigr)
\Bigr)\,dt
\\
&\qquad
=\int_G \bigl(x(t)\otimes \rho_t(y(t^{-1}s))\bigr)
\otimes s\,dt
=\int_G x(t)\otimes \rho_t(y(t^{-1}s))\,dt\otimes s
\\
&\qquad
= \Upsilon(x\otimes y)(s)\otimes s
= \what{\gamma\otimes \rho}\circ
\Upsilon(x\otimes y)(s).
\end{align*}
\end{proof}

\begin{rem}\label{rem-decompaseactioncross}
    For later use, it is worthwhile to observe that the proof
    of \propref{X xr G} implies that if $_AX_B={_A(K\otimes_KX)}_B$
    is the (equivariant) standard factorization of $X$, as
    described in \propref{decomaction}, then
    we get the standard factorization
    $$_{A\times_{\alpha,r}G}(X\times_{\gamma,r}G)_{B\times_{\beta,r}G}
    \cong {_{A\times_{\alpha,r}G}
    \big((K\times_{\mu,r}G)\otimes_{K\times_{\nu,r}G}
    (X\times_{\gamma,r}G)\big)_{B\times_{\beta,r}G}},$$
    which is clearly equivariant with respect to the dual coactions.
    It follows from \propref{prop-standard} below  that
    $_{A\times_{\alpha,r}G}
    (K\times_{\mu,r}G)_{K\times_{\mu,r}G}$ is equivariantly isomorphic
    to $_AK_K\times_{\mu,r}G$, so that the action-crossed-product
    functor of \thmref{xpr-fun} preserves standard factorizations of
    morphisms.
\end{rem}

\subsection{Coactions}

In \thmref{co-xpr-fun} below we define a functor from $\c C(G)$
to $\c A(G)$ with object map
\[
(A, \delta)\mapsto(A\times_\delta G,\hat\delta).
\]
As for actions, we first need to choose once and for all a single
version of the crossed product, and we use \defnref{def-crossed}: the
coaction $\delta$ is a nondegenerate homomorphism into $M(A\otimes
C^*(G))$, and composing with the nondegenerate homomorphism
$\id\otimes\lambda\:A\otimes C^*(G)\to M(A\otimes \c K(L^2(G)))$ gives
a nondegenerate homomorphism
\[
j_A=(\id\otimes\lambda)\circ \delta\:A\to M(A\otimes \c
K(L^2(G))).
\]
Writing $j_G=1\otimes M$, we get a covariant homomorphism
$(j_A,j_G)$ of $(A,G,\delta)$ to $M(A\otimes \c K(L^2(G)))$.  (We
sometimes write $j^A_G$ for $j_G$ if there is danger of confusion.)
Then the crossed product is defined as
\[
A\times_\delta G=\overline{ j_A(A) j_G(C_0(G)) }
\subseteq M(A\otimes \c K(L^2(G))).
\]

The standard theory of coactions (see \appxref{coactions-chap}) tells
us that $A\times_\delta G$ is a $C^*$-algebra and that the triple $(A
\times_\delta G, j_A, j_G)$ has the universal property that for every
covariant homomorphism $(\pi,\mu)$ of $(A,C_0(G))$ into a multiplier
algebra $M(B)$ there exists a unique nondegenerate homomorphism $\pi
\times \mu\:A\times_\delta G\to M(B)$ such that $(\pi\times
\mu)(j_A(a) j_G(f))=\pi(a) \mu(f)$ for all $a\in A$ and $f\in
C_0(G)$.

The \emph{dual action} $\hat\delta$ of $G$ on $A\times_{\delta} G$
is defined (see \defnref{def-dualact}) by $\hat\delta_s=j_A\times
(j_G\circ \rho_s)$, where $\rho_s$ denotes right translation on
$C_0(G)$.

If $\phi\:A\to M(B)$ is a homomorphism which is
equivariant for coactions $\delta$ and $\epsilon$ of $G$ on $A$ and
$B$, there is an associated homomorphism
\[
\phi\times G=(j_B\circ \phi)\times j_G
\]
of $A\times_\delta G$ into $M(B\times_\epsilon G)$,
which is nondegenerate if $\phi$ is nondegenerate (see
\lemref{lem-equivcrossed}).

Let $\zeta\:X\to M(X\otimes C^*(G))$ be a $\delta-\epsilon$
compatible nondegenerate coaction on a right-Hilbert $A-B$ bimodule $X$.
We will
define the crossed product $X\times_\zeta G$ following \cite[Theorem
3.2]{er:mult} (except here we use full coactions and avoid
representing things on Hilbert space).  Put
\[
j_X = (\id\otimes\lambda)\circ\zeta.
\]

\begin{prop}\label{prop-coactmodcross}
Assume that $(_AX_B,\zeta)$ is a nondegenerate
right-Hilbert bimodule coaction. With the above notation, the subspace
\[
X\times_\zeta G \deq  \overline{ j_X(X)\d j^B_G(C_0(G)) }
\]
of $M(X\otimes \c K(L^2(G)))$ is a right-Hilbert $(A\times_\delta G) -
(B\times_\epsilon G)$ bimodule.  Indeed, $X\times_\zeta G$ is closed
under the bimodule actions of $A\times_\delta G$ and $B\times_\epsilon
G$, and the $M(B\otimes \c K(L^2(G)))$-valued inner product on 
$X\times_\zeta G$
actually takes values in $B\times_\epsilon G$.
Moreover, if 
$B_X=\overline{\lk X,X\rk}_B$, then
$\overline{\lk X\times_{\zeta}G, X\times_{\zeta}G\rk}_{B\times_{\epsilon}G}
=B_X\times_{\epsilon_X}G$.
\end{prop}

As for actions, we first give the proof for the case of
a nondegenerate right-partial imprimitivity bimodule coaction
$({_KX_B}, \zeta)$,
using the linking algebra approach. So let
$L(X)=\smtx{K& X\\ \rev{X} & B}$ denote the linking algebra
of $X$. If we identify $L(X)\otimes C^*(G)$ with
$L(X\otimes C^*(G))$ as described in \secref{sec-link}
of~\chapref{hilbert-chap}, we
can use \propref{multsoflink} to get a canonical identification
$M(L(X)\otimes C^*(G))\cong L(M(X\otimes C^*(G)))$.
Thus, by \lemref{lem-coactlink} we obtain a coaction
$\nu$ of $G$ on $L(X)$ determined by
\begin{equation}\label{eq-LXcoaction}
\nu\left(\mtx{ k& x\\
  \tilde{y}& b}\right)
=\mtx{\mu(k) &\zeta(x)\\ 
  \widetilde{\zeta(y)}& \epsilon(b)}.
\end{equation}
If we similarly identify $M(L(X)\otimes \K(L^2(G)))$ with
$L(M(X\otimes \K(L^2(G))))$, we see that the
maps $(j_L, j_G^L)\:(L(X), C_0(G))\to M(L(X)\otimes \K(L^2(G)))$
can be written as
$$j_L=(\id_L\otimes \lambda)\circ \nu=
\mtx{(\id_K\otimes\lambda)\circ \mu& (\id_X\otimes\lambda)\circ \zeta\\
(\id_{\rev{X}}\otimes\lambda)\circ \tilde{\zeta}&
(\id_B\otimes\lambda)\circ \epsilon}
=\mtx{j_K& j_X\\ \widetilde{j_X}& j_B}$$
and
$$j_G^L=1_{M(L)}\otimes M=\mtx{1_{M(K)}\otimes M&0\\0&
1_{M(B)}\otimes M}=\mtx{j_G^K& 0\\ 0& j_G^B}.$$
Thus, for the crossed product $L(X)\times_{\nu}G$, we can write
$$L(X)\times_{\nu}G=\overline{j_L(L(X))j_G^L(C_0(G))}
=\mtx{\overline{j_K(K)j_G^K(C_0(G))}& \overline{j_X(X)j_G^B(C_0(G))}\\
\overline{\widetilde{j_X(X)}j_G^K(C_0(G))}&\overline{j_B(B)j_G^B(C_0(G))}}.
$$
Note that if  $p=\smtx{1&0\\ 0&0}$
and $q=\smtx{0&0\\ 0&1}$ denote the corner projections in $L(X)$,
then $j_L(p)$ and $j_L(q)$ are the corner projections
of $L(X)\times_{\nu}G\cong L(X\times_{\zeta}G)$.
We now gather these results into:

\begin{lem}\label{lem-linkcoactcross}
     Suppose that $({_KX_B}, {_\mu\zeta_\epsilon})$ is a
     nondegenerate right-Hilbert bimodule
     coaction of $G$ on the right-partial imprimitivity bimodule
     $_KX_B$. Then
     $$X\times_{\zeta}G= \overline{j_X(X)j_G^B(C_0(G))}=
     \overline{j_G^K(C_0(G))j_X(X)}$$
     is a right-partial $K\times_{\mu}G-B\times_{\epsilon}G$
     imprimitivity bimodule with respect to the bimodule operations
     as described in \propref{prop-coactmodcross}.

     If $_KX_B$ is an imprimitivity bimodule, then so is
     $_{K\times_{\mu}G}(X\times_{\zeta}G)_{B\times_{\epsilon}G}$ --- in
     general, the range of the $B\times_{\epsilon}G$-valued inner
     product on $X\times_{\zeta}G$ is $B_X\times_{\epsilon_X}G$,
     where $B_X$ denotes the range of the $B$-valued
     inner product on $X$ and $\epsilon_X$ denotes the restriction
     of $\epsilon$ to $B_X$ 
\textup(see Lemma~\textup{\ref{lem-idealcoact}}\textup).

     Finally, if
     $\nu=\smtx{\mu &\zeta\\ \tilde{\zeta}& \epsilon}$ is the
     corresponding coaction on the linking algebra $L(X)$,
     then $L(X)\times_{\nu}G$ is canonically isomorphic to
     $L(X\times_{\zeta}G)$.
\end{lem}
\begin{proof} All assertions follow directly from the
above considerations, except the fullness of the left
inner product on $X\times_{\zeta}G$ and the identification
of the range of the right inner product with $B_X\times_{\epsilon_X}G$.
(The identity $X\times_{\zeta}G=\overline{j_G^K(C_0(G))j_X(X)}$
follows from identifying $\overline{j_G^K(C_0(G))j_X(X)}$
as the upper right corner of $\overline{j_G^L(C_0(G))j_L(L(X))}$.)

The fullness of the left inner product follows from the fact that
if $q\in M(L(X))$ is a full projection, then $j_L(q)$ is a
full projection in $L(X)\times_\nu G$, which follows from
\begin{multline*}
     \overline{(L(X)\times_{\nu}G) j_L(q)(L(X)\times_{\nu}G)}
=\overline{j_G^L(C_0(G))j_L(L(X)qL(X))j_G^L(C_0(G))}\\
=\overline{j_G^L(C_0(G))j_L(L(X))j_G^L(C_0(G))}=L(X)\times_{\nu}G.
\end{multline*}
(Compare with the proof of \lemref{lem-linkingaction}.)
If $B_X=B$, the same argument applies to give fullness
of the right inner product on $X\times_{\zeta}G$.
Thus, in this case $X\times_{\zeta}G$ is an imprimitivity bimodule.

To get the general formula for the right inner product we
do the same trick as in the action case.
Let $L^{\ip}(X)=\smtx{K&X\\ \rev{X}&B_X}$. We know from
\lemref{lem-idealcoact} that $\zeta$ is also $\mu-\epsilon_X$
compatible, so we can view $L^{\ip}(X)$ as a
coaction invariant ideal in $L(X)$, and obtain a
corresponding embedding of
$$L^{\ip}(X)\times_{\nu} G=\mtx{K\times_{\mu}G&X\times_{\zeta}G\\
(X\times_{\zeta}G)\!\!\widetilde{\ \ }& B_X\times_{\epsilon_X}G}$$
into $L(X)\times_{\nu}G= L(X\times_{\zeta}G)$, which is the
identity on the upper right corner $X\times_{\zeta}G$
and which compresses to the inclusion
$B_X\times_{\epsilon_X}G\subseteq B\times_{\epsilon}G$ in the
lower right corner.
\end{proof}

\begin{proof}[Proof of \propref{prop-coactmodcross}]
     Let $K=\K(X)$. By \propref{decomcoaction} we see that
     there is a unique nondegenerate
     coaction $\mu\:K\to M(K\otimes C^*(G))$
     such that the homomorphism $\kappa\:A\to M(K)$
     coming from the left action of $A$ on $X$ is
     $\delta-\mu$ equivariant and such that
     $({_KX_B},\zeta)$ is $\mu-\epsilon$ compatible.
     By \lemref{lem-linkcoactcross} we get an identification
     $K\times_{\mu}G\cong\K(X\times_{\zeta}G)$.
     Thus  the
     homomorphism
     $$\kappa\times G\: A\times_{\delta}G\to M(K\times_{\mu}G)
     \cong\L_{B\times_{\epsilon}G}(X\times_{\zeta}G)$$
     determines a left action of $A\times_{\delta}G$
     on $X\times_{\zeta}G$. Since
     $$\kappa\times G=(j_K\circ \kappa)\times j_G^K
     =(\kappa\otimes \id_\K)\circ (j_A\times j_G^A)$$
     (see the proof of \lemref{lem-equivcrossed}),
     it follows that this left action coincides with
     the action described in the proposition.
\end{proof}

We now introduce the dual action of $G$ on $X\times_{\zeta}G$.

\begin{defn}
Given a $\delta-\epsilon$ compatible nondegenerate coaction $\zeta$ of $G$ on a
right-Hilbert $A-B$ bimodule $X$, the \emph{dual action} $\hat\zeta$
of $G$ on the crossed product $X\times_\zeta G$ is defined on 
generators by
\[
\hat\zeta_s(j_X(x)\d j_G(f))=j_X(x)\d j_G(\rho_s(f)).
\]
(Recall that $\rho_s$ denotes right translation on $C_0(G)$ by $s$.)
\end{defn}

It is straightforward to check that the above formula determines a
$\hat\delta -\hat\epsilon$ compatible action on the right-Hilbert
$(A\times_\delta G)-(B\times_\epsilon G)$ bimodule $X\times_\zeta G$.

For later use it is necessary to know that
coaction-crossed products by standard right-Hilbert modules are
standard.

\begin{lem}
\label{co-std-xpr}
Suppose that $\phi\: A\to M(B)$ is a \textup(possibly degenerate\textup)
homomorphism which is
equivariant for nondegenerate coactions $\delta$ and $\epsilon$ of $G$,
and let $({_AX_B}, \epsilon_X)$ be the standard coaction
associated to $\phi$ 
\textup(see Definition~\textup{\ref{def-standard}}\textup).
Then 
$\big({_{A\times_{\delta}G}(X\times_{\epsilon_X}G)_{B\times_{\epsilon}G}},\;
\hat{\epsilon}_X\big)$ is equivariantly isomorphic to the standard
action
$\big({_{A\times_{\delta}G}Z_{B\times_{\epsilon}G}},\hat{\epsilon}_Z\big)$,
$Z\deq \phi\times G(A\times_{\delta}G)
(B\times_{\epsilon}G)$,
associated to $\phi\times G\:A\times_{\delta}G\to
M(B\times_{\epsilon}G)$.
\end{lem}

\begin{proof}
     Since $\epsilon_X$ is defined as the restriction of
     $\epsilon$ to $X=\phi(A)B$, it is clear that the
     inclusion $\iota\:{_AX_B}\hookrightarrow {_BB_B}$ is 
$\epsilon_X-\epsilon$
     equivariant.
     It follows from \lemref{lem-Cmultidentify} that
     $\iota\otimes\id_{\K}\:M_{\K(L^2(G))}(X\otimes \K(L^2(G)))
\hookrightarrow
     M_{\K(L^2(G))}(B\otimes \K(L^2(G)))$ is an isometric inclusion.
     Thus, for $a\in A$ and $b\in B$, we can compute
     \begin{align*}
	(\iota\otimes \id_{\K})\circ j_X(\phi(a)b)
	=(\id_B\otimes\lambda)\circ
	(\iota \otimes \id_G)\circ \epsilon_X(\phi(a)b)
	=(\id_B\otimes\lambda)\circ \epsilon(\phi(a)b)
	=j_B(\phi(a)b).
     \end{align*}
     Restricting to
     $X\times_{\epsilon_X}G=\overline{j_G(C_0(G))j_X(\phi(A)B)j_G(C_0(G))}$,
     it follows that
     \begin{align*}
     \iota\otimes\id_\K(X\times_{\epsilon_X}G)&=
     \overline{j_G(C_0(G))j_B(\phi(A))j_B(B)j_G(C_0(G))}\\
     &=\overline{\phi\times G(j_G(C_0(G))j_A(A))(j_B(B)j_G(C_0(G))}\\
     &=\phi\times G(A\times_{\delta}G)(B\times_{\epsilon}G).
     \end{align*}
     It is clear that this isomorphism is
     $\widehat{\epsilon_X}-\hat{\epsilon}_Z$ equivariant and preserves
     all bimodule operations.
\end{proof}

\begin{thm}
\label{co-xpr-fun}
The object map $(A, \delta)\mapsto (A\times_\delta G,\hat\delta)$
and the morphism map $[_AX_B, \zeta]\mapsto [_{A\times G}(X
\times_\zeta G)_{B\times G},\hat\zeta]$ define a functor from $\c
C(G)$ to $\c A(G)$.
\end{thm}

\begin{proof}
To see that the morphism map is well-defined, let
\[
\Phi\:_{(A,\delta)}(X,\zeta)_{(B,\epsilon)}
\iso _{(A,\delta)}(Y,\eta)_{(B,\epsilon)}
\]
be an equivariant isomorphism.  Then the
isomorphism $\Phi\otimes \id\:X\otimes \c K(L^2(G))\to Y\otimes
\c K(L^2(G))$ takes $j_X$ onto $j_Y$ (after extending to the
multiplier bimodule), hence maps $X\times_\zeta G$ onto $Y
\times_\eta G$.  Obviously, $\Phi\otimes \id$ preserves the
right-Hilbert bimodule operations and is equivariant for the dual
actions.  Thus $[X\times_\zeta G,\hat\zeta]=[Y\times_\eta
G,\hat\eta]$.

It is a very special case of \lemref{co-std-xpr} that 
the morphism map preserves identities.

For compositions, if $(_AX_B,\zeta)$ is $\delta-\epsilon$ compatible
and $(_BY_C,\eta)$ is $\epsilon-\vartheta$ compatible, we need to
find a right-Hilbert $(A\times_{\delta} G)-(C\times_{\vartheta} G)$
bimodule isomorphism
\[
\Theta\:(X\times_{\zeta} G)\otimes_{B\times_\epsilon G}
(Y\times_{\eta} G) \iso
(X\otimes_B Y)\times_{\zeta\cotimes \eta} G
\]
which is $(\hat\zeta\otimes\hat\eta)-\what{\zeta\cotimes \eta}$
equivariant.  The notation $\Theta$ is suggestive: by
\lemref{Theta-lem} we already have a right-Hilbert $(A\otimes \c
K(L^2(G)))-(C\otimes \c K(L^2(G)))$ bimodule isomorphism
\[
\Theta\:
(X\otimes \c K(L^2(G)))\otimes_{B\otimes \c K(L^2(G))}
(Y\otimes \c K(L^2(G)))
\iso (X\otimes_B Y)\otimes \c K(L^2(G)),
\]
and we will show that it takes 
$(X\times_\zeta G)\otimes_{B\times_\epsilon G} (Y\times_\eta G)$ onto 
$(X\otimes_B Y)\times_{\zeta\cotimes\eta}G$ 
(after extending to multiplier bimodules).  
For $f,g\in C_0(G)$, $x\in X$, and $y\in Y$ we have
\begin{multline*}
\Theta\bigl(j^A_G(f)\d j_X(x)\otimes j_Y(y)\d j^C_G(g)\bigr)
= j^A_G(f)\d \Theta\bigl((\id\otimes\lambda)\circ 
  \zeta(x)\otimes (\id\otimes\lambda)\circ \eta(y)\bigr) \d j^C_G(g)\\
= j^A_G(f)\d (\id\otimes\lambda)\circ 
  (\zeta\cotimes \eta)(x\otimes y) \d j^C_G(g)
= j^A_G(f)\d j_{X\otimes_B Y}(x\otimes y)\d j^C_G(g),
\end{multline*}
and such elements densely span $(X\otimes_B Y)\times_{\zeta\cotimes\eta}G$.  
It follows
immediately from the construction that $\Theta$ gives an
isomorphism $(X\times_{\zeta}G)\otimes_{B\times_\epsilon G} (Y\times_{\eta}G)
\cong (X\otimes_B Y)\times_{\zeta\cotimes\eta}G$, 
and a quick calculation verifies that
$\Theta$ intertwines the actions $\hat\zeta\otimes\hat\eta$ and
$\what{\zeta\cotimes \eta}$.
\end{proof}

\begin{rem}\label{rem-decompose}
     It follows from the
     proof of \propref{prop-coactmodcross} that
     the right-Hilbert bimodule
     $_{A\times_{\delta}G}(X\times_{\zeta}G)_{B\times_\epsilon G}$
     factors as the product
     $$_{A\times_{\delta}G}\big((K\times_\nu G)\otimes_{K\times_{\mu}G}
     (X\times_\zeta G)\big)_{B\times_{\epsilon}G},$$
     where the first factor is the nondegenerate
     standard  bimodule  corresponding
     to the homomorphism $\kappa\times G\:A\times_{\delta}G
     \to M(K\times_\nu G)$, and the second factor is the
     right-partial imprimitivity bimodule
     $_{K\times_{\mu}G}(X\times_{\zeta}G)_{B\times_{\epsilon}G}$.
     This decomposition is clearly equivariant with respect to the
     dual actions. It
     follows from \lemref{co-std-xpr} that
     $_{A\times_{\delta}G}(K\times_\nu G)_{K\times_{\mu}G}$
     is equivariantly isomorphic to $_AK_K\times_{\mu}G$, so we
     see that the coaction crossed-product functor
     of \thmref{co-xpr-fun} preserves standard decomposition
     of morphisms. (Compare with \remref{rem-decompaseactioncross}.)
\end{rem}

\section{Restriction and inflation}
\label{res-inf}

\subsection{Actions}

If $H$ is a closed subgroup of a locally compact group $G$, we can
restrict any action of $G$ to an action of $H$, and if $N$ is a closed
normal subgroup of $G$, we can inflate any action of $G/N$ to an action of $G$
(by composing with the quotient map).  
More generally, suppose we have a continuous homomorphism $\phi\colon
G\to F$ of locally compact groups.  
Then any action $\alpha$
of $F$ gives rise to an action $\alpha\circ \phi$ of $G$.  We want to
make this (and hence restriction and inflation)
into a functor, so we need to handle the morphisms.

\begin{lem}
Let $\gamma$ be an $\alpha -\beta$ compatible
action of $F$ on ${}_AX_B$, and let $\phi\:G\to F$ be a continuous
homomorphism.  Then $\gamma\circ \phi$ is an $\alpha\circ \phi -
\beta\circ \phi$ compatible action of $G$ on ${}_AX_B$.
\end{lem}

\begin{proof}
Just note that for each $x\in X$, 
the map $s\mapsto (\gamma\circ\phi)_s(x) = \gamma_{\phi(s)}(x)$
is continuous from $G$ into $X$.
\end{proof}

\begin{prop}
Let $G$ and $F$ be locally compact groups, and let $\phi\colon G\to F$
be a continuous homomorphism. 
Then the object map $(A,\alpha)\mapsto(A,\alpha\circ \phi)$ and the
morphism map $[_AX_B,\gamma]\mapsto[_AX_B,\gamma\circ \phi]$
define a functor from $\c A(F)$ to $\c A(G)$.
\end{prop}

\begin{proof}
It is obvious that the morphism map is well-defined on isomorphism
classes (because any $F$-equivariant isomorphism is also
$G$-equivariant) and sends identities to identities.  For
compositions, let $(_AX_B,\gamma)$ and $(_BY_C, \rho)$ be
actions of $F$.  Then just observe that
\[
(\gamma\otimes_B \rho)\circ \phi=
(\gamma\circ \phi)\otimes_B (\rho\circ \phi)
\]
on $_AX\otimes_BY_C$.
\end{proof}

Specializing to restriction and inflation, we immediately get the
following corollaries:

\begin{cor}
\label{res-fun}
For a closed subgroup $H$ of $G$, the object map $(A, G,\alpha)
\mapsto(A, H,\alpha|_H)$ and the morphism map $[_AX_B, G,\gamma]
\mapsto[_AX_B, H,\gamma|_H]$ define a functor from $\c A(G)$ to $\c
A(H)$.
\end{cor}

\begin{cor}
For a closed normal subgroup $N$ of $G$, the object map $(A, G/N,
\alpha)\mapsto (A, G,\infl\alpha)$ and the morphism map $[_AX_B,
G/N,\gamma]\mapsto[_AX_B, G,\infl\gamma]$ define a functor from
$\c A(G/N)$ to $\c A(G)$.
\end{cor}

\subsection{Coactions}

If $H$ is a closed subgroup of a locally compact group $G$, we can
inflate any coaction of $H$ to a coaction of $G$ (see
\defnref{ex-inflate}), and if $N$ is a closed normal subgroup we can
restrict any coaction of $G$ to a coaction of $G/N$ (see
\defnref{ex-restrict}).  Again, we want to make functors out of
restriction $\delta\mapsto\delta|$ and inflation
$\delta\mapsto\infl\delta$, and again there is a more general situation
which unifies both: that of a continuous homomorphism $\phi\colon G\to
F$ of locally compact groups.  We can integrate $\phi$ up to a nondegenerate
homomorphism, still denoted by $\phi$, from $C^*(G)$ to $M(C^*(F))$,
and this allows us to pass from coactions of $G$ to coactions of $F$.

\begin{lem}
\label{co-res-inf-lem}
\textup{(i)} If $\delta$ is a coaction of $G$ on a $C^*$-algebra $A$, then
$(\id\otimes\phi)\circ\delta$ is a coaction of $F$ on $A$.  Moreover,
if $\delta$ is normal or nondegenerate, then so is
$(\id\otimes\phi)\circ\delta$.

\textup{(ii)}
If $\zeta$ is a coaction of $G$ on a right-Hilbert bimodule $X$, then
$(\id\otimes \phi)\circ \zeta$ is a coaction of $F$ on $X$.
If $\zeta$ is nondegenerate, then so is $(\id\otimes \phi)\circ \zeta$.
\end{lem}

\begin{proof}
(i)
Certainly $(\id\otimes \phi)\circ \delta$ is a nondegenerate
injective homomorphism from $A$ to $M(A\otimes C^*(F))$.
Since $\delta(A)\subseteq M_G(A\otimes C^*(G))$ it follows from
\propref{prop-Cmultextend} that $(\id_A\otimes\phi)\circ
\delta(A)\subseteq M_F(A\otimes C^*(F))$.

For the coaction identity,
\begin{align*}
&(\id\otimes \delta_F)\circ(\id\otimes \phi)\circ \delta
= (\id\otimes \delta_F\circ \phi)\circ \delta
= \bigl(\id\otimes(\phi\otimes \phi)\circ \delta_G\bigr)\circ \delta
= (\id\otimes \phi\otimes \phi) \circ(\id\otimes \delta_G)\circ \delta\\
&\qquad
= (\id\otimes \phi\otimes \phi) \circ(\delta\otimes \id)\circ \delta
= \bigl((\id\otimes \phi)\circ \delta\otimes \phi\bigr)\circ \delta
= \bigl((\id\otimes \phi)\circ \delta\otimes \id\bigr)\circ
(\id\otimes \phi)\circ \delta.
\end{align*}

For the normality, it suffices to show that if $(\pi,\mu)$ is a
covariant homomorphism of $(A,G,\delta)$ then $(\pi,\mu\circ \phi^*)$
is a covariant homomorphism of $(A,F,(\id\otimes \phi)\circ \delta)$,
where $\phi^*\:C_0(F)\to M(C_0(G))=C_b(G)$ is defined by
$\phi^*(f)=f\circ \phi$. Normality of $(\id\otimes \phi)\circ \delta$
then follows from normality of $\delta$ by applying this fact
to $(j_A, j_G)$. We have
\begin{align*}
&\ad (\mu\circ \phi^*\otimes \id)(w_F)(\pi(a)\otimes 1)
= \ad (\mu\otimes \id) \bigl((\phi^*\otimes \id)(w_F)\bigr)
(\pi(a)\otimes 1)\\
&\qquad
= \ad (\mu\otimes \id) \bigl((\id\otimes \phi)(w_G)\bigr)
(\pi(a)\otimes 1)
= \ad (\id\otimes \phi) \bigl((\mu\otimes \id)(w_G)\bigr)
(\pi(a)\otimes 1)\\
&\qquad
= (\id\otimes \phi)\circ \ad (\mu\otimes \id)(w_G)
(\pi(a)\otimes 1)
= (\id\otimes \phi)\circ (\pi\otimes \id)\circ \delta(a)\\
&\qquad
= (\pi\otimes \id)\circ (\id\otimes \phi)\circ \delta(a),
\end{align*}
as desired.

For the nondegeneracy, we have
\begin{align*}
\overline{(\id\otimes \phi)\circ \delta(A)
(1\otimes C^*(F))}
&
= \overline{(\id\otimes \phi)\circ \delta(A)
\bigl(1\otimes \phi(C^*(G)) C^*(F)\bigr)}
\\
&
= \overline{(\id\otimes \phi)\circ \delta(A)
\bigl(1\otimes \phi(C^*(G))\bigr)
(1\otimes C^*(F))}
\\
&
= \overline{(\id\otimes \phi)\bigl(\delta(A) (1\otimes C^*(G))\bigr)
(1\otimes C^*(F))}
\\
&
= \overline{(\id\otimes \phi) (A\otimes C^*(G)) (1\otimes C^*(F))}
\\
&
= A\otimes \phi(C^*(G)) C^*(F)
= A\otimes C^*(F).
\end{align*}

(ii)
A similar argument handles bimodules.
\end{proof}

\begin{thm}
\label{co-res-inf-fun}
Let $F$ and $G$ be locally compact groups, and let $\phi\colon G\to F$
be a continuous homomorphism. Then
the object map $(A, \delta)\mapsto(A, (\id\otimes \phi)\circ
\delta)$ and the morphism map $[_AX_B, \zeta]\mapsto [_AX_B, (\id
\otimes \phi)\circ \zeta]$ define a functor from $\c C(G)$ to $\c
C(F)$.
\end{thm}

\begin{proof}
The morphism map is clearly well-defined and preserves identities.
For compositions, let $(_AX_B, \zeta)$ be $\delta-\epsilon$
compatible and $(_BY_C, \eta)$ be $\epsilon-\vartheta$ compatible.
We must show that the $(\id\otimes \phi)\circ \delta-(\id\otimes
\phi)\circ \vartheta$ compatible coactions $(\id\otimes \phi)
\circ(\zeta\cotimes_B \eta)$ and $((\id\otimes \phi)\circ \zeta)
\cotimes_B ((\id\otimes \phi)\circ \eta)$ on the right-Hilbert $A -
C$ bimodule $X\otimes_BY$ agree.  Because this requires mildly fussy
bimodule-multiplier calculations, we go through it in
detail.  By definition of the tensor product coactions, it suffices to
show that the triangle
\[
\xymatrix{
{X\otimes_BY}
\ar[r]^-{\zeta\otimes_B \eta}
\ar[dr]|-{(\id\otimes \phi)\circ \zeta\otimes_B
(\id\otimes \phi)\circ \eta}
&{M \bigl((X\otimes C^*(G))\otimes_{B\otimes C^*(G)}
(Y\otimes C^*(G) \bigr)}
\ar[d]^{(\id\otimes \phi)\otimes_{B\otimes C^*(G)}
(\id\otimes \phi)}
\\
&{M \bigl((X\otimes C^*(F))\otimes_{B\otimes C^*(F)}
(Y\otimes C^*(F)) \bigr)}
}
\]
and the square
\[
\xymatrix{
{(X\otimes C^*(G))\otimes_{B\otimes C^*(G)}
(Y\otimes C^*(G))}
\ar[r]^-{\Theta_G}
\ar[d]_{(\id\otimes \phi)\otimes_{B\otimes C^*(G)}(\id\otimes \phi)}
&{(X\otimes_BY)\otimes C^*(G)}
\ar[d]^{\id\otimes \phi}
\\
{M \bigl((X\otimes C^*(F))\otimes_{B\otimes C^*(F)}
(Y\otimes C^*(F)) \bigr)}
\ar[r]_-{\Theta_F}
&{M \bigl((X\otimes_BY)\otimes C^*(F) \bigr)}
}
\]
commute, where $\Theta_G$ is the canonical isomorphism of $(X\otimes
C^*(G))\otimes_{B\otimes C^*(G)}(Y\otimes C^*(G))$ onto $(X
\otimes_BY)\otimes C^*(G)$ from \lemref{Theta-lem}, and similarly
for $\Theta_F$. A trivial computation on elementary tensors shows
that the triangle commutes.

We show that the square commutes: for $x\in X$, $y\in
Y$, and $f, g\in C^*(G)$ we have
\begin{align*}
&\Theta_F\circ \bigl(
(\id\otimes \phi)\otimes_{B\otimes C^*(G)}
(\id\otimes \phi) \bigr)
\bigl((x\otimes f)\otimes(y\otimes g) \bigr)\\
&\qquad
= \Theta_F\bigl((x\otimes \phi(f))\otimes
(y\otimes \phi(g)\bigr)
= (x\otimes y)\otimes \phi(fg)\\
&\qquad
= (\id\otimes \phi) \bigl((x\otimes y)\otimes fg \bigr)
= (\id\otimes \phi)\circ \Theta_G
\bigl((x\otimes f)\otimes(y\otimes g) \bigr).
\end{align*}
\end{proof}

Specializing to restriction and inflation, we immediately get the
following corollaries:

\begin{cor}
For a closed normal subgroup $N$ of $G$, the object map $(A, G,
\delta)\mapsto(A, G/N, \delta|_{G/N})$ and the morphism map $[_AX_B, G,
\zeta]\mapsto[_AX_B, G/N, \zeta|_{G/N}]$ define a functor from $\c C(G)$
to $\c C(G/N)$.
\end{cor}

\begin{cor}
\label{co-inf-fun}
For a closed subgroup $H$ of $G$, the object map $(A, H, \delta)
\mapsto(A, G,\infl \delta)$ and the morphism map $[_AX_B, H, \zeta]
\mapsto[_AX_B, G,\infl \zeta]$ define a functor from $\c C(H)$ to $\c
C(G)$.
\end{cor}

\section{Decomposition}
\label{dec-act}

\subsection{Actions}

Let $(A, G,\alpha)$ be an action, and let $N$ be a closed normal
subgroup of $G$.  Then there is a canonical action $\tilde\alpha$ of
$G$ on the restricted crossed product $A\times_{\alpha|}N$, and this
is usually called the decomposition action because
\cite[Proposition~1]{gre:local} tells us that $A\times_\alpha G$ can
be decomposed into a twisted crossed product of $A\times_{\alpha|}N$
by $\tilde\alpha$.  We will define a corresponding decomposition
action $\alpha\dec$ on the reduced crossed product, and show that
$(A,\alpha)\mapsto (A\times_{\alpha|,r} N,\alpha\dec)$ is the
object map of a functor on $\c A(G)$.

Decomposition actions use the modular function for conjugation of $G$
on $N$:
\[
\Delta_{G,N}(s)\deq \Delta_G(s)\Delta_{G/N}(sN)^{-1},
\]
which has the property that
\[
\int_N f(n)\,dn
= \Delta_{G,N}(s)\int_N f(s^{-1}ns)\,dn
\righttext{for} f\in C_c(N),\ s\in G.
\]

\begin{lem}
Let $N$ be a closed normal subgroup of $G$.
\begin{enumerate}
\item
If $\alpha$ is an action of $G$ on a $C^*$-algebra $A$, then
there exists a unique action $\alpha\dec$ of $G$ on $A
\times_{\alpha|,r} N$ given on $C_c(N,A)$ by
\[
\alpha\dec_s(f)(n)=
\Delta_{G,N}(s)\alpha_s(f(s^{-1}ns)).
\]

\item
If $\gamma$ is an $\alpha-\beta$ compatible action of $G$ on 
a right-Hilbert bimodule ${}_AX_B$,
then there exists a unique $\alpha\dec-\beta\dec$ compatible action
$\gamma\dec$ of $G$ on $X\times_{\gamma|,r} N$ given on $C_c(N,X)$ by
\[
\gamma\dec_s(x)(n)=
\Delta_{G,N}(s)\gamma_s(x(s^{-1}ns)).
\]
\end{enumerate}
\end{lem}

\begin{proof}
For (i), since the right side of the given formula is the one
defining the decomposition action $\tilde\alpha$ on the full crossed
product $A\times_{\alpha|}N$, and since $\tilde\alpha$ leaves the
kernel of the regular homomorphism $i^r_A\times i^r_N$ invariant
\cite[Lemma 10]{gre:local}, the desired action $\alpha\dec$ on the
reduced crossed product $A\times_{\alpha|,r}N$ exists.

For (ii), the right side of the formula gives an
automorphism of the pre-right-Hilbert $C_c(N,A)-C_c(N,B)$ bimodule
$C_c(N,X)$ which is continuous for the inductive limit topology
and is compatible with the
the actions $\alpha\dec$ and $\beta\dec$ on the coefficients.
It therefore extends to an automorphism $\gamma\dec_s$
of $X\times_{\gamma,r}N$ by \lemref{pre-rHb-lem}.
\end{proof}

\begin{thm}
\label{dec-fun}
The object map $(A,\alpha)\mapsto(A\times_{\alpha|,r}N,\alpha\dec)$ and the
morphism map $[_AX_B,\gamma]\mapsto [_{A\times_r N}(X\times_{\gamma|,r}
N)_{B\times_r N},\gamma\dec]$ define a functor from $\c A(G)$ to
itself.
\end{thm}

\begin{proof}
Composing the reduced-crossed-product functor with the restriction
functor, and then forgetting about the actions on the image, we get a
functor $\c A(G)\to \c C$; we must trace the decomposition actions
through the various isomorphisms.  First we show that, given a
right-Hilbert $A-B$ bimodule isomorphism $\Phi\:X\to Y$ which is
equivariant for $\alpha -\beta$ compatible actions $\gamma$ and
$\rho$, the associated isomorphism $\Phi\times_rN\:X\times_{\gamma|,r} N
\to Y\times_{\rho|,r} N$ from \thmref{xpr-fun} is $\gamma\dec -
\rho\dec$ equivariant.  But this is easily checked for $x\in
C_c(N,X)$:
\begin{align*}
(\Phi\times_rN)\circ\gamma\dec_s(x)(n)
&
= \Phi\bigl(\gamma\dec_s(x(n))\bigr)
= \Phi\bigl(\Delta_{G,N}(s)
\gamma_s(x(s^{-1}ns))\bigr)
\\&
= \rho_s\bigl(\Delta_{G,N}(s)
\Phi(x(s^{-1}ns))\bigr)
= \rho\dec_s\circ (\Phi\times_rN)(x)(n).
\end{align*}

For the identities, let $(A,\alpha)$ be an action.  From
\thmref{xpr-fun} we know that the right-Hilbert bimodule crossed
product $_AA_A\times_r N$ is isomorphic to the right-Hilbert bimodule
$_{A\times_r N}(A\times_r N)_{A\times_r N}$.  Since the
decomposition actions coincide on $C_c$-functions, the isomorphism is
equivariant.

Turning to compositions, let $(_AX_B,\gamma)$ be
$\alpha -\beta$ compatible and $(_BY_C,\rho)$ be
$\beta-\nu$ compatible.  We need to show that the isomorphism
$\Upsilon$ of $(X\times_r N)\otimes_{B\times_r N}
(Y\times_r N)$ onto $(X\otimes_B Y)\times_r N$ 
from \thmref{xpr-fun} is $(\gamma\dec\otimes
\rho\dec)-(\gamma\otimes \rho)\dec$ equivariant.  
For $s\in G$, $x\in C_c(N,X)$, $y\in C_c(N,Y)$, and $n\in N$, we have:
\begin{align*}
&\Upsilon\circ (\gamma\dec\otimes \rho\dec)_s
(x\otimes y)(n)
= \Upsilon\bigl(\gamma\dec_s(x)\otimes \rho\dec_s(y)\bigr)(n)
\\&\qquad
=\int_N\gamma\dec_s(x)(k)\otimes
\nu_k\bigl(\rho\dec_s(y)(k^{-1}n)\bigr)\,dk
\\&\qquad
=\int_N \Delta_{G,N}(s)
\gamma_s(x(s^{-1}ks))\otimes
\nu_k\bigl(\Delta_{G,N}(s)
\rho_s(y(s^{-1}k^{-1}ns))\bigr)\,dk
\\&\qquad
= \Delta_{G,N}(s)
\int_N\gamma_s(x(k))\otimes
\nu_k\bigl(\rho_s(y(k^{-1}s^{-1}ns))\bigr)\,dk
\\&\qquad
= \Delta_{G,N}(s)
(\gamma\otimes \rho)_s
\left(\int_N x(k)\otimes
\nu_k\bigl(y(k^{-1}s^{-1}ns)\bigr)\,dk \right)
\\&\qquad
= \Delta_{G,N}(s)
(\gamma\otimes \rho)_s \bigl(\Upsilon(x\otimes y)(s^{-1}ns)\bigr)
= (\gamma\otimes \rho)\dec_s\circ \Upsilon
(x\otimes y)(n).
\end{align*}
\end{proof}

\subsection{Coactions}

Let $(A,G,\delta)$
be a nondegenerate normal coaction, and let $N$ be a closed normal
subgroup of $G$.
Recall from Lemma~\ref{decom-coact}
that the decomposition coaction $\delta\dec$ of $G$ on $A
\times_{\delta|}G/N$ is defined on the
generators by
\[
\delta\dec(j_A(a)j_{G/N}(f))
= (j_A\otimes \id)\circ \delta(a)(j_{G/N}(f)\otimes 1),
\]
and moreover $\delta\dec$ is also nondegenerate and normal.
We will define a functor from $\c C(G)$ to itself with object map
\[
(A, \delta)\mapsto (A\times_{\delta|}G/N, \delta\dec).
\]
For the morphism map, we need decomposition coactions for bimodules.

\begin{prop}\label{prop-zetadec}
If $\zeta$ is a nondegenerate
$\delta-\epsilon$ compatible coaction of $G$ on ${}_AX_B$,
then there is a unique $\delta\dec-\epsilon\dec$ compatible coaction
$\zeta\dec$ of $G$ on $X\times_{\zeta|}G/N$ which is given 
for $x\in X$ and $f\in C_0(G/N)$ by 
\begin{equation}\label{eq-zetadec}
\zeta\dec(j_X(x)\d j^B_{G/N}(f))
= (j_X\otimes \id)\circ \zeta(x)\d
(j^B_{G/N}(f)\otimes 1).
\end{equation}
\end{prop}

\begin{proof}
     Let $K=\K(X)$ and let $\mu\:K\to M(K\otimes C^*(G))$
     denote the unique coaction on $K$ which makes
     $({_KX_B}, \zeta)$ a $\mu-\epsilon$ compatible
     coaction (see \propref{decomcoaction}).
     Moreover, let $L(X)=\smtx{K&X\\ \rev{X}& B}$ denote
     the associated linking algebra and let
     $\nu|=\smtx{\mu|&\zeta|\\ \tilde\zeta|& \epsilon|}$
     be the corresponding coaction of $G/N$ on $L(X)$.
     By the proof of \lemref{lem-linkcoactcross}
     we get an identification
     $$L(X)\times_{\nu|}G/N=L(X\times_{\zeta|}G/N)=
     \mtx{K\times_{\mu|}G/N& X\times_{\zeta|}G/N \\
     (X\times_{\zeta|}G/N)\!\widetilde{\ }& B\times_{\epsilon|}G/N}$$
     which identifies the dense subspace $j_L(L(X))j_{G/N}^L(C_0(G/N))$
     of $L(X)\times_{\nu|}G/N$ with
     $$\mtx{j_K(K)j_{G/N}^K(C_0({G/N}))& j_X(X)j_{G/N}^B(C_0({G/N}))\\
     \widetilde{j_X(X)}j_{G/N}^K(C_0({G/N}))& j_B(B)j_{G/N}^B(C_0({G/N}))}.$$
     Using \propref{multsoflink} we get a canonical identification
     $$M\big((L(X)\times_{\nu|}G/N)\otimes C^*(G)\big)\cong
     \mtx{M\big((K\times_{\nu|}G/N)\otimes C^*(G)\big)&
     M\big((X\times_{\zeta|}G/N)\otimes C^*(G)\big)\\
     M\big((X\times_{\zeta|}G/N)\otimes C^*(G)\big)\!\!\widetilde{\ \ }&
     M\big((B\times_{\epsilon|}G/N)\otimes C^*(G)\big)}.
     $$
     With these identifications,
     the decomposition coaction $\nu\dec$ of $G$ on
     $L(X)\times_{\nu|}G/N$ is then given on the generators
     $$j_L\left(\smtx{k&x\\ \tilde{y}& b}\right)j_{G/N}^L(f)=
     \mtx{j_K(k)j_{G/N}^K(f)& j_X(x)j_{G/N}^B(f)\\
     \widetilde{j_X(y)}j_{G/N}^K(f)& j_B(b)j_{G/N}^B(f)}$$
     by
     \begin{align*}
&\nu\dec\left(j_L\left(\smtx{k&x\\ \tilde{y}& b}\right)j_{G/N}^L(f)\right)
     =(j_L\otimes \id)\circ \nu\left(\smtx{k&x\\ \tilde{y}& b}\right)
     (j_{G/N}^L(f)\otimes 1)\\
&\qquad
     =\mtx{(j_K\otimes \id)\circ \mu(k)& (j_X\otimes\id)\circ
     \zeta(x)\\
     (j_X\otimes\id)\circ \zeta(y)\!\!\widetilde{\ \ }&
     (j_B\otimes\id)\circ \epsilon(b)}
     \mtx{ j_{G/N}^K(f)\otimes 1 &  0\\  0&  j_{G/N}^B(f)\otimes 1}\\
&\qquad
     =\mtx{(j_K\otimes \id)\circ \mu(k)(j_{G/N}^K(f)\otimes 1)&
     (j_X\otimes\id)\circ \zeta(x)(j_{G/N}^B(f)\otimes 1)\\
     (j_X\otimes\id)\circ \zeta(y)\!\!\widetilde{\ \ }(j_{G/N}^K(f)\otimes 1)&
     (j_B\otimes\id)\circ \epsilon(b)(j_{G/N}^B(f)\otimes 1)}.
     \end{align*}
     Thus we see that $\nu\dec$ compresses on the corners
     $K\times_{\nu|}G/N$ and $B\times_{\epsilon|}G/N$ to
     the coactions $\mu\dec$ and $\beta\dec$, and it follows
     from \lemref{lem-coactlink} that
     $\nu\dec$ compresses to a nondegenerate $\mu\dec-\beta\dec$ compatible
     coaction on $X\times_{\zeta|}G/N$, which is then given on the
     generators by \eqref{eq-zetadec}.

     To see that $\zeta\dec$ is also $\delta\dec-\epsilon\dec$
     compatible with respect to the given left action of
     $A\times_{\alpha|}G/N$ on $X\times_{\zeta|}G/N$, it 
     suffices to show that the homomorphism
     $\kappa\times G/N\: A\times_{\alpha|}G/N\to M(K\times_{\mu|}G/N)$
     is $\delta\dec-\mu\dec$ equivariant, where $\kappa\:A\to
     M(K)\cong\L_B(X)$ is determined by the left $A$-action on $X$.
     Using the equations
     $$(\kappa\times G/N)\circ j_A= j_K\circ\kappa,\quad
     (\kappa\times G/N)\circ j_{G/N}^A=j_{G/N}^K,\quad
     \text{and} \quad (\kappa\otimes \id_G)\circ \delta=\mu\circ \kappa$$
     (which follow from the $\delta-\mu$ equivariance of $\kappa\:A\to
     M(K)$ and the definition of $\kappa\times G/N$), we can compute
     \begin{align*}
\big((\kappa\times G/N)\otimes\id_G\big)\circ
\delta\dec\big(j_A(a)j_{G/N}^A(f)\big)
&=\big((\kappa\times G/N)\otimes\id_G\big)\big((j_A\otimes\id_G)\circ \delta(a)
(j_{G/N}^A(f)\otimes 1)\big)\\
&=((j_K\circ \kappa)\otimes \id_G)\circ \delta(a)(j_{G/N}^K(f)\otimes
1)\\
&=(j_K\otimes \id_G)\circ (\kappa\otimes\id_G)\circ \delta(a)
(j_{G/N}^K(f)\otimes 1)\\
&=(j_K\otimes\id_G)\circ\mu(\kappa(a))(j_{G/N}^K(f)\otimes 1)\\
&=\mu\dec\big(j_K(\kappa(a))j_{G/N}^K(f)\big)\\
&=\mu\dec\circ (\kappa\times G/N)\big(j_A(a)j_{G/N}^A(f)\big)
\end{align*}
for all $a\in A$ and $f\in C_0(G/N)$.
Thus $((\kappa\times G/N)\otimes\id_G)\circ \delta\dec=
\mu\dec\circ(\kappa\times G/N)$
and the proof is complete.
\end{proof}

\begin{thm}
\label{co-dec-fun}
The object map $(A, \delta)\mapsto (A\times_{\delta|}G/N, \delta\dec)$ and
the morphism map $[_AX_B, \zeta]\mapsto [_{A\times G/N}(X\times_{\zeta|}
G/N)_{B\times G/N}, \zeta\dec]$ define a functor from $\c C(G)$ to
itself.
\end{thm}

\begin{proof}%
The proof is almost the same as that of \thmref{dec-fun}; we merely
indicate the differences.  We know already that the composition
$(A,\delta)\mapsto (A,\delta|)\mapsto (A\times_{\delta|} G/N)$ is a
functor from $\c C(G)$ to $\c C$; we must trace the decomposition
coactions through the various isomorphisms.  First, given an
isomorphism $\Phi\:_{(A,\delta)}(X,\zeta)_{(B,\epsilon)}\to
_{(A,\delta)}(Y,\eta)_{(B,\epsilon)}$, the associated isomorphism
$\Phi\otimes \id\:X\times_{\zeta|} G/N\to Y\times_{\eta|} G/N$
(from the proof of \thmref{co-xpr-fun})
is $\zeta\dec-\eta\dec$ equivariant:
\begin{align*}
&((\Phi\otimes \id)\otimes \id)\circ \zeta\dec(j_X(x)\d j_{G/N}(f))
= (\Phi\otimes \id\otimes \id)\bigl(
(j_X\otimes \id)\circ \zeta(x)\d (j_{G/N}(f)\otimes 1)\bigr)
\\&\qquad
= (j_Y\otimes \id)\circ \eta\circ \Phi(x)
\d (j_{G/N}(f)\otimes 1)
= \eta\dec\bigl(j_Y(\Phi(x))\d j_{G/N}(f)\bigr)
\\&\qquad
= \eta\dec\circ (\Phi\otimes \id) (j_X(x)\d j_{G/N}(f))
\end{align*}
for $x\in X$, $f\in C_0(G/N)$.

Next, it follows straight from the definitions that the decomposition
coactions coincide on a standard-bimodule crossed product 
$_AA_A \times_{\delta|}G/N
={}_{A\times G/N}(A\times_{\delta|}G/N)_{A\times G/N}$.

Finally, fix 
$_{(A,\delta)}(X,\zeta)_{(B,\epsilon)}$ and
$_{(B,\epsilon)}(Y,\eta)_{(C,\vartheta)}$, and
let 
\[
\Theta\colon
(X\times_{\zeta|}G/N)\otimes_{B\times G/N}(Y\times_{\eta|}G/N)
\to
(X\otimes_B Y)\times_{\zeta|\cotimes_B\eta|}G/N
\]
be the bimodule isomorphism 
from \thmref{co-xpr-fun} (see also \lemref{Theta-lem}).  
Further letting 
\[
\begin{split}
\Theta_G\colon
\bigl((X\times_{\zeta|}G/N)\otimes C^*(G)\bigr)
  \otimes_{(B\times G/N)\otimes C^*(G)}
  \bigl((Y\times_{\eta|}G/N)\otimes C^*(G)\bigr)\\
\to
\bigl((X\times_{\zeta|}G/N)\otimes_{B\times G/N}
  (Y\times_{\eta|}G/N)\bigr)\otimes C^*(G)
\end{split}
\]
be the canonical isomorphism 
which appears in the definition of 
$\zeta\dec\cotimes_{B\times G/N}\eta\dec$, we have
\begin{align*}
&(\Theta\otimes \id)\circ (\zeta\dec\cotimes_{B\times G/N}\eta\dec)
(j_{G/N}(f)\d j_X(x)\otimes j_Y(y)\d j_{G/N}(g))\\
&\qquad=
(\Theta\otimes \id)\circ \Theta_G\bigl(
(j_{G/N}(f)\otimes 1)\d (j_X\otimes \id)\circ \zeta(x)
\otimes
(j_Y\otimes \id)\circ \eta(y)\d (j_{G/N}(g)\otimes 1)\bigr)\\
&\qquad=
(\Theta\otimes \id)
\Bigl((j_{G/N}(f)\otimes 1)\d
\Theta_G\bigl(
(j_X\otimes \id)\circ \zeta(x)
\otimes
(j_Y\otimes \id)\circ \eta(y)\bigr)
\d (j_{G/N}(g)\otimes 1)\Bigr)\\
&\qquad=
(j_{G/N}(f)\otimes 1)\d
(j_{X\otimes_B Y}\otimes \id)\circ
\Theta_G(\zeta(x)\otimes \eta(y))
\d (j_{G/N}(g)\otimes 1)\\
&\qquad=
(j_{G/N}(f)\otimes 1)\d
(j_{X\otimes_B Y}\otimes \id)\circ
(\zeta\cotimes_B \eta)(x\otimes y)
\d (j_{G/N}(g)\otimes 1)\\
&\qquad=
(\zeta\cotimes_B \eta)\dec\bigl(
(j_{G/N}(f)\d j_{X\otimes_B Y}(x\otimes y)\d j_{G/N}(g)\bigr)\\
&\qquad=
(\zeta\cotimes_B \eta)\dec\circ \Theta
\bigl(j_{G/N}(f)\d j_X(x)\otimes j_Y(y)\d j_{G/N}(g)\bigr)
\end{align*}
for all $f,g\in C_0(G/N)$, $x\in X$, and $y\in Y$. 
Thus $\Theta$ 
is $(\zeta\dec\cotimes_{B\times G/N} \eta\dec)-(\zeta\cotimes_B \eta)\dec$
equivariant, which shows that the morphism map preserves compositions. 
\end{proof}

\section{Induced actions}

Let $H$ be a closed subgroup of $G$.  We will define a functor from
$\c A(H)$ to $\c A(G)$ with object map $(A,\alpha)\mapsto(\ind^G_HA,
\ind\alpha)$. (See \eqref{ind-a-eq} in \appxref{imprim-chap} for the 
definition of the induced action $\ind\alpha$.)
Recall from \cite[Corollary 3.2]{qr:induced} that if $(A,H,\alpha)$
is an action, then
\[
M(\ind_H^GA)=\{x\in M(A\otimes C_0(G))\mid
x(sh)=\alpha_{h^{-1}}(x(s))
\text{ for all }s\in G,h\in H\}.
\]
It follows that if $\phi\:A\to M(B)$ is a nondegenerate
homomorphism which is equivariant for actions of $H$, there is an
associated nondegenerate homomorphism $\ind \phi\:\ind A\to M(\ind
B)$ defined by
\[
\ind \phi(f)(s)=\phi(f(s)).
\]
For the morphism map we need to induce right-Hilbert bimodules.  Let
$_{(A,\alpha)}(X,\gamma)_{(B,\beta)}$ be a right-Hilbert bimodule
action of $H$.  We define
\[
\ind_H^G X=\left\{ x\in C_b(G, X) \mid
x(sh)=\gamma_{h^{-1}}(x(s)) \text{ and }
(sH\mapsto \|x(s)\|)\in C_0(G/H)
\right\},
\]
and for $s\in G$ we define $\ind\gamma_s\:\ind X\to\ind X$ by
\[
\ind\gamma_s(x)(t)=x(s^{-1}t).
\]

\begin{lem}
With the above notation, $\ind X$ becomes a right-Hilbert $\ind A -
\ind B$ bimodule with pointwise operations
$$
(f\d x)(s)=f(s)\d x(s),
\quad (x\d g)(s)=x(s)\d g(s),\quad
\text{and}\quad
\<x,y\>_{\ind B}(s)=\<x(s),y(s)\>_B,
$$
and $\ind\gamma$ is an $\ind\alpha-\ind\beta$
compatible action of $G$ on $\ind X$. Moreover, if
$B_X=\overline{\lk X,X\rk}_B$ is the range  of the $B$-valued inner
product on $X$, then $\overline{\lk \ind X,\ind X\rk}_{\ind B}=\ind B_X$.
\end{lem}

\begin{proof}
Routine calculations show that the above formulas satisfy the
algebraic properties of a right-Hilbert bimodule. The inner
product $\<\cdot,\cdot\>_{\ind B}$ is positive-definite because
$\<\cdot,\cdot\>_B$ is.
It is also straightforward to check that for each $s\in G$
the map $\ind\gamma_s$ is a right-Hilbert bimodule automorphism with
coefficient maps $\ind\alpha_s$ and $\ind\beta_s$, and that for each
$x\in\ind X$ we have $\ind\gamma_s(x)\to x$ uniformly as $s\to e$
in $G$.

  We next observe that there are enough elements in
     $\ind X$. For $f\in C_c(G,X)$ define
     $F(s)=\int_H\gamma_h(f(sh))\,dh$; then  $F\in \ind X$.
     Moreover, if $x\in X$ and if we choose
     $f=g\otimes x$ with $g\in C_c(G)^+$ supported in a small
     neighborhood of $e\in G$ such that $\int_H g(h)\,dh=1$,
     we see that $F(e)=\int_H \gamma_h(f(h))\,dh$ is close to $x$.
     Using translation, we see that there are enough elements in
     $\ind X$ such that all evaluations at the points $s\in G$
     have dense range in $X$.

     To see that $\overline{\<\ind X,\ind X\>}_{\ind B}=\ind B_X$
     (viewed canonically as a closed ideal of $\ind B$),
we appeal to the Lemma in \cite{ech:induced}. Let $C_0$ be the linear span
of the range of the $\ind B$-valued inner product.  To see that $C_0$
is dense in $\ind B_X$ it suffices to show that for each $s\in G$ the set
$\{g(s) \mid g\in C_0\}$ is dense in $B_X$, and that $C_0$ is closed under
multiplication by $C_c(G/H)$.
But the first follows from the density in $X$ of the images of the
evaluation maps on $\ind X$ and $\overline{\lk X, X\rk}_B=B_X$,
and the second from closure of $\ind X$ under the
right module action of $C_c(G/H)\subseteq M(\ind B)$.
\end{proof}

\begin{rem}\label{rem-linkinduced}
     If $X$ is a partial $A-B$ imprimitivity bimodule, we
     also obtain a left $\ind A$-valued inner product
     on $X$ which turns $X$ into a partial $\ind A-\ind B$
     imprimitivity bimodule and, as in the
     proof for the range of the $\ind B$-valued inner
     product on $\ind X$ given above, we see
     that ${}_{\ind A}\overline{\lk \ind X,\ind X\rk}=\ind A_X$ with
     $A_X={}_A\overline{\lk X,X\rk}$. This shows that
     $\ind X$ is a right-partial imprimitivity
     bimodule, or a left-partial imprimitivity bimodule, whenever $X$ is.

     Moreover, if $L(X)=\smtx{A&X\\ \rev{X}&B}$ is the linking algebra
     of $X$, equipped with the $H$-action induced from the
     $H$-action on $X$ (compare with \lemref{lem-linkingaction}),
     then it is straightforward to check that the canonical identification
     $C_b(G,L(X))=\smtx{C_b(G,A)&C_b(G,X)\\ C_b(G,X)\!\!\widetilde{\ \ }&
     C_b(G,B)}$ restricts to a canonical identification
     $$\ind L(X)=\mtx{\ind A&\ind X\\ \widetilde{\ind X} &\ind B}.$$
\end{rem}

We use the above observation for the proof of:

\begin{lem}
\label{Ind L}
Let $H$ be a closed subgroup of $G$, $(_AX_B, H,\gamma)$ a
right-Hilbert bimodule action, and $\c F$ a vector subspace of
$\ind X$ such that for each $s\in G$ the set $\{x(s) \mid x\in \c F \}$
is dense in $X$, and such that $\c F$ is closed under multiplication
by $C_c(G/H)$.  Then $\c F$ is dense in $\ind X$.
\end{lem}

\begin{proof}
This follows very quickly from the $C^*$-version, 
which is the Lemma in \cite{ech:induced}, 
modulo a linking-algebra connection.
Forgetting about the left $A$-module structure, we may view $X$ as a
right-partial $K-B$ imprimitivity bimodule, where $K=\c K_B(X)$.  Let $L(X)$
be the associated linking algebra. Then $\ind L(X)=L(\ind X)$, so 
by \cite{ech:induced}, the set
\[
\c F_1 \deq \begin{pmatrix}\ind K &\c F
\\
\c F^*&\ind B \end{pmatrix}
\]
is dense in $\ind L(X)$.  But then $\c F$ is dense in the
upper right corner $\ind X$ of $\ind L(X)$, and we are done.
\end{proof}

\begin{prop}
\label{ind-prop}
The object map $(A,\alpha)\mapsto (\ind A,\ind\alpha)$ and the
morphism map $[_AX_B,\gamma]\mapsto [_{\ind A}(\ind X)_{\ind B},\ind
\gamma]$ define a functor from $\c A(H)$ to $\c A(G)$.
\end{prop}

\begin{proof}
We first show that the morphism map is well-defined.  Suppose $\Phi\:
X\to Y$ is an isomorphism of right-Hilbert $A-B$ bimodules which is
equivariant for $\alpha -\beta$ compatible actions $\gamma$ and $\rho$
of $G$.  Then the map $\ind\Phi\colon \ind X\to \ind Y$
defined by $\ind\Phi(x)(s) \deq  \Phi(x(s))$ is easily seen to give a
right-Hilbert $\ind A-\ind B$ bimodule isomorphism of 
$\ind X$ onto $\ind Y$. It is trivial to see that the morphism map
preserves identities.

For compositions, let $(_AX_B,\gamma)$ be $\alpha -\beta$ compatible
and let $(_BY_C,\rho)$ be $\beta-\epsilon$ compatible.  Then $\Phi(x
\otimes y)(s) \deq  x(s)\otimes y(s)$ is easily seen to give a
right-Hilbert $\ind A-\ind C$ bimodule homomorphism $\Phi\:\ind X
\otimes_{\ind B}\ind Y\to\ind (X\otimes_B Y)$.  To see that
$\Phi$ is an isomorphism it remains to show that its range is dense.
But this follows from \lemref{Ind L}, since for each $s\in
G$ the set $\{x(s) \mid x\in \Phi(\ind X \odot\ind Y)\}$ is dense in $X
\otimes_BY$ and $\Phi(\ind X \odot\ind Y)$ is closed under
multiplication by $C_c(G/H)$.
\end{proof}


When the action of the subgroup $H$ is restricted from 
an action $\gamma$ of the larger group $G$ on ${}_AX_B$,
the induced action on $\ind X$ is isomorphic to, but not identical to,
the action $\gamma\otimes\tau$ on $X\otimes C_0(G/H)$, where $\tau$ is
the canonical action of $G$ on $C_0(G/H)$ by left translation.  So to 
deduce from \propref{ind-prop}
that this latter, more common construction is also functorial, we just
need to make sure that this isomorphism --- which is defined for $x\in
C_b(G,X)$ by $\phi(x)(sH) = \gamma_s(x(s))$ ---  and the related coefficient
isomorphisms transport the induced actions to the respective diagonal
actions on the tensor products.
Unsurprisingly, this is completely
straightforward, and has nothing particularly to do with the nature of
the translation action $\tau$.  The result could also be
proved directly via elementary means.

\begin{cor}
For any closed subgroup $H$ of $G$ the object map $(A,\alpha)
\mapsto(A\otimes C_0(G/H),\alpha\otimes \tau)$ and the morphism map
$[_AX_B,\gamma]\mapsto[_{A\otimes C_0(G/H)}(X\otimes C_0(G/H))_{B
\otimes C_0(G/H)},\gamma\otimes \tau]$ define a functor from $\c
A(G)$ to itself.
\end{cor}

\section{Combined functors}

In the next chapter we will need
to combine several functors we have already constructed,
involving the category $\c A \c C(G)$:

\begin{prop}
\label{combined-fun}
Let $N$ be closed normal subgroup of $G$.
\begin{enumerate}
\item The object map
$(A,\alpha)\mapsto
(A\times_{\alpha|,r} N,\alpha\dec,\infl \what{\alpha|_N})$
and the morphism map
$[X,\gamma]\mapsto
\mathbox{[X\times_{\gamma|,r} N,\gamma\dec,\infl \what{\gamma|_N}]}$
define a functor from $\c A(G)$ to $\c A \c C(G)$.

\item The object map
$(A, \delta)\mapsto
(A\times_{\delta|}G/N,\infl \what{\delta|_{G/N}},
\delta\dec)$
and the morphism map
$[X, \zeta]\mapsto
\mathbox{[X\times_{\zeta|}G/N,\infl \what{\zeta|_{G/N}},\zeta\dec]}$
define a functor from $\c C(G)$ to $\c A \c C(G)$.
\end{enumerate}
\end{prop}

\begin{proof}
(i) We know that $(A,\alpha)\mapsto(A\times_r N,\alpha\dec)$ is a
functor from $\c A(G)$ to $\c A(G)$.  On the other hand, $(A,\alpha)
\mapsto(A\times_r N,\infl \what{\alpha|_N})$ is a functor from $\c
A(G)$ to $\c C(G)$ because it is the composition of the restriction,
crossed-product, and inflation functors:
\begin{gather*}
(A, G,\alpha)\mapsto(A, N,\alpha|_N)\mapsto
(A\times_r N, \what{\alpha|_N})\mapsto
(A\times_r N,\infl \what{\alpha|_N}).
\end{gather*}
To see that the combined map $(A,\alpha)\mapsto(A\times_r N,
\alpha\dec,\infl \what{\alpha|_N})$ is a functor, we need (as usual)
to check that the morphism map is well-defined and preserves
identities and compositions.  But this will be easy; for example,
suppose we have an equivariant isomorphism $(X,
\gamma)\cong (Y, \rho)$.  Then the proof of \thmref{dec-fun} gives a
specific isomorphism $X\times_r N\cong Y
\times_r N$.  On the other hand, the proofs of \corref{res-fun},
\thmref{xpr-fun}, and \corref{co-inf-fun} give corresponding
isomorphisms which, when composed, produce the same isomorphism as
\thmref{dec-fun}.  Thus this common isomorphism is simultaneously
equivariant for both the decomposition action $\gamma\dec$ and the
inflated coaction $\infl\what{\gamma|_N}$.  Therefore, the morphism map
of the present theorem is well-defined.  Similar arguments show that
it preserves identities and compositions, and (ii) is proved
similarly.
\end{proof}

%
%

\chapter{The Natural Equivalences}
\label{natural-chap}

\section{Statement of the main results}

This chapter contains our main results.  The main idea is that any
reasonable imprimitivity theorem can be viewed as expressing a natural
equivalence of appropriate functors.  We shall illustrate this for two
versions of Green's Imprimitivity Theorem (one with induced algebras
and another which is the basis for the usual characterization of
induced representations of actions), and for Mansfield's Imprimitivity
Theorem (for an extensive survey of these theorems we refer to
\appxref{imprim-chap}).
More precisely, as our first main result we shall prove:

\begin{thm}
\label{ind-act natural}
Let $H$ be a closed subgroup of a locally compact group $G$.  
The Green imprimitivity bimodules $V_H^G(A,\alpha)$ implement a natural
equivalence between the functors
\[
(A,\alpha)\mapsto
((\ind_H^G A)\times_{\ind\alpha,r} G,\what{\ind\alpha})
\midtext{and}
(A,\alpha)\mapsto
(A\times_{\alpha,r} H,\infl_H^G \hat\alpha)
\]
from $\c A(H)$ to $\c C(G)$.
\end{thm}

Here $V_H^G(A,\alpha)$ denotes the completion of the pre-imprimitivity
bimodule $C_c(G,A)$ with actions and inner products described in
\eqeqref{eq-ind-imp}.  We say ``implement'' here because there
are some implicit assertions being made.  The natural equivalence must
assign to each object $(A,\alpha) \in \c A(H)$ an equivalence in the
category $\c C(G)$ from the object 
$((\ind A)\times_r G,\what{\ind\alpha})$ to the object $(A\times_r H,\infl\hat\alpha)$.
Green's theorem says that the isomorphism class $[V_H^G(A,\alpha)]$ is
an equivalence in the category $\c C$, so we need to construct a
coaction $\delta_V$ of $G$ on $V_H^G(A,\alpha)$ such that
$[V_H^G(A,\alpha),\delta_V]$ is an equivalence in $\c C(G)$.  Then to
prove that the assignment $(A,\alpha)\mapsto
[V_H^G(A,\alpha),\delta_V]$ is a natural equivalence we shall need to
show that any morphism $[{}_{(A,\alpha)}(X,\gamma)_{(B,\beta)}]$ in
$\c A(H)$ gives a commutative diagram
\[
\xymatrix
@C+30pt
{
{((\ind_H^G A)\times_r G,\what{\ind\alpha})}
\ar[r]^-{[V_H^G(A),\delta_V]}
\ar[d]_{[(\ind X)\times_r G,\what{\ind\gamma}]}
& {(A\times_r H,\infl\hat\alpha)}
\ar[d]^{[X\times_r H,\infl\hat\gamma]}
\\
{((\ind_H^G B)\times_r G,\what{\ind\beta})}
\ar[r]_-{[V_H^G(B),\delta_V]}
& {(B\times_r H,\infl\hat\beta)}
}
\]
in the category $\c C(G)$ (see \thmref{Ind-thm}).  For obvious reasons
we shall state this more simply as ``the diagram
\begin{equation}
\label{ind-diag}
\xymatrix
{
{(\ind_H^G A)\times_r G}
\ar[r]^-{V_H^G(A)}
\ar[d]_{(\ind X)\times_r G}
& {A\times_r H}
\ar[d]^{X\times_r H}
\\
{(\ind_H^G B)\times_r G}
\ar[r]_-{V_H^G(B)}
& {B\times_r H}
}
\end{equation}
commutes equivariantly for the appropriate coactions.''  We shall
construct the coaction $\delta_V$ in \secref{sec-green1} below and we
shall also give a proof of \thmref{ind-act natural} in that section.

If we start with an action $\alpha$ of $G$ on a $C^*$-algebra $A$,
then $(\ind_H^G A)\times_r G$ 
is naturally isomorphic to $C_0(G/H,A)\times_{\alpha\otimes
\tau,r} G$, and (after identifying the corresponding coaction on this
algebra) \thmref{ind-act natural} immediately gives another natural
equivalence between crossed-product functors; equivalently,
commutativity of the diagram
\[
\xymatrix{
{C_0(G/H,A)\times_r G}
\ar[r]^-{X_H^G(A)}
\ar[d]_-{C_0(G/H,X)\times_r G}
& {A\times_r H}
\ar[d]^{X\times_r H}
\\
{C_0(G/H,B)\times_r G}
\ar[r]_-{X_H^G(B)}
& {B\times_r H,}
}
\]
equivariantly for the appropriate coactions, where $X_H^G(A)$ is the
$C_0(G/H,A)\times_{\alpha\otimes\tau,r} G - A\times_{\alpha|,r} H$
imprimitivity bimodule of \thmref{thm-green-imp}.

If we require the subgroup to be \emph{normal}, we can also get
equivariance for appropriate \emph{actions}.  (In the
following discussion we call the subgroup $N$ instead of $H$ to
emphasize that it will be normal in $G$.)  In this situation there
exists a natural isomorphism between the crossed product $C_0(G/N,A)
\times_{\alpha\otimes \tau,r} G$ and the iterated crossed product $A
\times_{\alpha,r} G\times_{\hat\alpha|} G/N$, which transforms the
dual coaction $\what{\alpha\otimes \tau}$ into the decomposition
coaction $\hat\alpha\dec$ (see the proof of \lemref{left alg} below).
The $G$-action on the iterated crossed product is the action
$\infl\what{\hat\alpha|}$: we first take the dual coaction
$\hat\alpha$, then the restriction $\hat\alpha|$ to $G/N$, then the
dual action $\what{\hat\alpha|}$ of $G/N$,
and then inflate this to an action of $G$.  The counterpart on the $A
\times_r N$-side is given by the decomposition action $\alpha\dec$ as
described in \secref{functors-chap}.\ref{dec-act}.  Mainly as a
consequence of \thmref{ind-act natural} we then derive:

\begin{thm}
\label{G-act natural}
Let $N$ be a closed normal subgroup of a locally compact group $G$.
The Green imprimitivity bimodules $X_N^G(A,\alpha)$ implement a natural
equivalence between the functors
\[
(A,\alpha)\mapsto
(A\times_{\alpha,r} G\times_{\hat\alpha|,r} G/N,
\infl\what{\hat\alpha|},
\hat\alpha\dec)
\midtext{and}
(A,\alpha)\mapsto
(A\times_{\alpha,r} N,
\alpha\dec,
\infl\what{\alpha|})
\]
from $\c A \c C(G)$ to itself.
\end{thm}

Finally, in \secref{sec-mans} we shall derive a complete dual version
of~\thmref{G-act natural}.  For this let $\delta$ be a nondegenerate
normal coaction of~$G$ on the $C^*$-algebra $A$, and let $N$ be a
closed normal subgroup of~$G$.  Using the generalization of
Mansfield's imprimitivity theorem for coactions as presented in
\cite{kq:imprimitivity}, we obtain an imprimitivity bimodule
$Y_{G/N}^G(A,\delta)$ for the crossed products $A\times_{\delta} G
\times_{\hat\delta|,r} N$ and $A\times_{\delta|} G/N$ (see
\thmref{thm-mans-imp}).  Again, both crossed products carry canonical
actions and coactions of~$G$: the action on $A\times_{\delta} G
\times_{\hat\delta|,r} N$ is the decomposition action
$\hat\delta\dec$, and the action on $A\times_{\delta|} G/N$ is the
inflation $\infl\what{\delta|}$ of the dual action $\what{\delta|}$ 
of~$G/N$ on $A\times_{\delta|} G/N$.  Similarly, the coaction on $A
\times_{\delta} G\times_{\hat\delta|,r} N$ is the inflation
$\infl\what{\hat\delta|}$, and the coaction on $A\times_{\delta|}
G/N$ is the decomposition coaction $\delta\dec$ as described in
\appxref{coactions-chap} (see \lemref{decom-coact}).

\begin{thm}
\label{coact natural}
Let $N$ be a closed normal subgroup of a locally compact group $G$.
The Mansfield imprimitivity bimodules $Y_{G/N}^G(A,\delta)$ 
implement a natural
equivalence between the functors
\[
(A,\delta)\mapsto
(A\times_{\delta,r} G\times_{\hat\delta|,r} N,
\hat\delta\dec,
\infl\what{\hat\delta|})
\midtext{and}
(A,\delta)\mapsto
(A\times_{\delta,r} G/N,
\infl\what{\delta|},
\delta\dec)
\]
from $\c A \c C(G)$ to itself.
\end{thm}

The construction of the appropriate actions and coactions on the
bimodules, and the proof of \thmref{coact natural}, will be given in
\secref{sec-mans}.  The strategy of the proofs of all three theorems
discussed above is to factor the right-Hilbert bimodule $X$
which appears in Diagram~\eqref{ind-diag}
(or in the appropriate diagrams corresponding to \thmref{G-act
natural} and \thmref{coact natural}) into the tensor product of a
standard bimodule and a right-partial imprimitivity bimodule.  
Then we prove the theorem in these special cases, and put the results
together (using functoriality) 
to get the desired result for
general right-Hilbert bimodules.

\section{Some further linking algebra techniques}
\label{sec-more-link}

\begin{lem}
\label{pq-hom-lem}
Suppose $Z$ is a right-Hilbert $E - F$ bimodule. Let $p \in M(E)$
and $q \in M(F)$ be  projections.  Then,
via restriction of the actions and inner product, $pZq$ becomes a
right-Hilbert $pEp - qFq$ bimodule. If $p$ is a full projection,
then $pZq$ is a full right-Hilbert $pEp - qFq$ bimodule.

Suppose in addition that
${}_\phi\Phi_\psi\:{}_EZ_F\to M({}_RW_S)$ is a nondegenerate
right-Hilbert bimodule homomorphism and that $q\in M(F)$ is full.  Then
$\psi(q)$ is a full projection in $M(S)$,
and $\Phi$ restricts to give a nondegenerate
right-Hilbert bimodule homomorphism
\[
{}_{pEp}(pZq)_{qFq}\to
M\Bigl({}_{\phi(p) R\phi(p)}
\bigl(\phi(p) W\psi(q)\bigr)_{\psi(q) S\psi(q)}\Bigr).
\]
\end{lem}

\begin{proof}
     First note that $pEpE=pE$, which follows from the fact that
     $pE$ is a left $pEp$-Hilbert module with respect to the
     canonical module operations. Similarly, we have $FqFq=Fq$.

It is clear that $pEp\d pZq\subseteq pZq$ and
$pZq\d qFq\subseteq pZq$, so the actions restrict.  Likewise,
$\<pZq,pZq\>_F\subseteq qFq$, so the right inner product also
restricts.  Nondegeneracy of the left
$pEp$-module action follows from the nondegeneracy of the
left module action of $E$ on $Z$ and
\[
pEp\d pZq=pEpE\d Zq=pE\d Zq=pZq.
\]
If $p$ is full, \ie, if $EpE$ is dense in $E$, it follows that
\[
\<pZq,pZq\>_{qFq}
=q \<pE\d Z,pE\d Z\>_F q
=q \<Z,EpE\d Z\>_F q
\]
is dense in $qFq$. Hence $pZq$ is full.

For the other part, first note that
\[
R\phi(p) R
=R\phi(E)\phi(p)\phi(E) R
=R\phi(EpE) R
\]
is dense in $R$, because $p$ is full and $\phi$ is
nondegenerate.

Now clearly $\Phi$ maps $pEp$ into $M(\phi(p) R\phi(p))$, $pZq$ into
$M(\phi(p) W\psi(q))$, and $qFq$ into $M(\psi(q) S\psi(q))$; the
only remaining issue is the nondegeneracy of the restrictions of
$\phi$, $\Phi$, and $\psi$.  This is clear for $\phi$ and $\psi$:
for example, using $pEpE=pE$ and the nondegeneracy of $\phi$ we get
\[
\phi(pEp)R=\phi(pEp)(\phi(E)R)=\phi(pEpE)R=\phi(pE)R=\phi(p)R,
\]
which clearly implies $\phi(pEp)R\phi(p)=\phi(p)R\phi(p)$.
To prove nondegeneracy of $\Phi$ we need fullness of
$q$: we then get
\begin{align*}
&\overline{\Phi(pZq)\d\psi(q) S\psi(q) }
=\overline{\Phi(pZ\d Fq)\d\psi(F) S\psi(q) }
=\overline{\phi(p)\Phi(Z)\d\psi(FqF) S\psi(q) }\\
&\qquad
=\overline{\phi(p)\Phi(Z)\d\psi(F) S\psi(q) }
=\overline{\phi(p)\Phi(Z)\d S\psi(q) }
=\phi(p) W\psi(q).
\end{align*}
\end{proof}

\begin{rem}
In the above lemma, if $Z$ is a (partial) $E-F$ imprimitivity bimodule, then
$pZq$ becomes a (partial) $pEp-qFq$ imprimitivity bimodule---just 
argue as above
for the left inner product as well.  If, in addition,
  ${}_EZ_F$ and ${}_RW_S$
are imprimitivity bimodules, the restriction of $\Phi$ to $pZq$ is
also a nondegenerate imprimitivity bimodule homomorphism, by
\remref{rHb-hom-rem1}.
\end{rem}

The following result generalizes \cite[Lemma
4.6]{er:stab}.

\begin{lem}
\label{link-lem}
Suppose ${}_AX_B$ and ${}_CY_D$ are partial imprimitivity bimodules,
and $Z$ is a right-Hilbert $L(X) - L(Y)$ bimodule.  Let $p=\pproj$ and $q=
\qproj$ denote the canonical full projections in both $M(L(X))$ and
$M(L(Y))$.  Then\textup:
\begin{enumerate}
\item
Via restriction of the actions and inner products, $pZp$ becomes a
right-Hilbert $A - C$ bimodule, $pZq$ becomes a right-Hilbert $A - D$
bimodule, and $qZq$ becomes a right-Hilbert $B - D$ bimodule.

\item
There are natural isometric right-Hilbert $A - D$ bimodule
homomorphisms
\[
\Phi\:X\otimes_B qZq\to pZq
\midtext{and}
\Psi\:pZp\otimes_C Y\to pZq
\]
given by
$\Phi(x\otimes z)=x\d z$
{and}
$\Psi(w\otimes y)=w\d y$.

\item
Let us now denote by $q_X$ the $q$-projection in $M(L(X))$
and by $p_Y$ the $p$-projection in $M(L(Y))$.
Then, if $q_X$ is full,  the homomorphism $\Phi$ of~\textup{(ii)}
is an isomorphism. Similarly, if $p_Y$ is full, then $\Psi$ is
an isomorphism.
In particular, both maps are isomorphisms if $_AX_B$ and $_CY_D$
are imprimitivity bimodules.
\item If $Z$ is a partial $L(X)-L(Y)$ imprimitivity bimodule
and if 
${}_A{\lk pZq, pZq\rk}\subseteq{}_A\overline{\lk X,X\rk}$, 
then the map $\Phi$ of~\textup{(ii)} is an isomorphism.
Similarly, if 
$\lk pZq, pZq\rk_D\subseteq\overline{\lk Y,Y\rk}_D$,
then $\Psi$ is an isomorphism.
\end{enumerate}
\end{lem}

\begin{proof}
Since $pL(X)p=A$, $qL(X)q=B$, and so on, (i) follows directly from
\lemref{pq-hom-lem}.

In order to prove (ii), first observe that $\Phi$ clearly intertwines
the left actions.  To check that $\Phi$ preserves the $D$-valued inner
products, fix $x,x'\in X$ and $z,z'\in qZq$, and compute:
\begin{equation*}
\<x\otimes z,x'\otimes z'\>_D
=\<z,\<x,x'\>_B\d z'\>_D
=\<z,x^*x'\d z'\>_D
=\<x\d z,x'\d z'\>_D
=\<\Phi(x\otimes z),\Phi(x'\otimes z')\>_D.
\end{equation*}
A similar argument applies to $\Psi$.

Suppose now that $q_X$ is full.
Then $L(X)q_XL(X)$ is dense in $L(X)$, and using $L(X)\d Z=Z$
by nondegeneracy of the left action, we see that
\[
\Phi(X\otimes qZq)
=pL(X)q_X\d Zq
=pL(X)q_XL(X)\d Zq
\]
is dense in $pZq$.
Similarly, the fullness of $p_Y$ implies that $\Psi$ is surjective.
This proves (iii).

For the proof of (iv) assume that $Z$ is a
partial $L(X)-L(Y)$ imprimitivity bimodule. Let
$A_X={}_A\overline{\lk X, X\rk}$ and let
$L^r(X)=\smtx{A_X&X\\\rev{X}&B}$ be the linking algebra
of the right-partial imprimitivity bimodule $_{A_X}X_B$.
Then $L^r(X)$ is a closed ideal of $L(X)$, and
we can consider the
partial $L^r(X)-L(Y)$ bimodule $L^r(X)\d Z\subseteq Z$.
We claim that $pL^r(X)Zq=pZq$. To see this
we simply note that $A_XpZq=pZq$, since
$pZq$ is a left $A_X$-Hilbert module by assumption;
this implies that $p(L^r(X)\d Z)q\subseteq A_X\d pZq= pZq$.
Since $_{A_X}X_B$ is a right-partial imprimitivity bimodule,
it follows from \propref{prop-partiallink} that $q_X$ is a
full projection in $M(L^r(X))$, and thus it follows from
the proof of (iii) that $pL(X)q_XZq\supseteq pL^r(X)q_XZq=pZq$.
But this implies surjectivity of $\Phi$.
A similar argument shows that
if $\lk pZq, pZq\rk_D\subseteq \overline{\lk Y,Y\rk}_D$,
then $\Psi$ is surjective.
\end{proof}

Before continuing to develop our general techniques, let us
first explain how a result like \lemref{link-lem} can be used in the
proofs of the main theorems.  For this we go back to the situation of
\thmref{ind-act natural}.  As explained earlier, the proof of this
theorem requires us to show the equivariant commutativity of 
Diagram~\exref{ind-diag}, 
\ie, that the bimodule compositions $V_H^G(A)\otimes_{A\times_r H}
\big(X\times_r H\big)$ and $\big((\ind X)\times_r G\big)
\otimes_{(\ind B)\times_r G} V_H^G(B)$ are equivariantly isomorphic.
As mentioned at the end of the previous section, our strategy is to
factor the problem into two cases, and it is the
partial imprimitivity bimodule case where the above lemma becomes extremely
helpful.  To see this, assume now that $X$ is a partial $A-B$ imprimitivity
bimodule, and let
$V_H^G(L)$ be Green's imprimitivity bimodule for the associated
$G$-action on the linking algebra $L(X)$. By \lemref{lem-linkingaction}
and \remref{rem-linkinduced}  we have the identifications
\[
(\ind L(X))\times_r G=L((\ind X)\times_r G)
\midtext{and}
L(X)\times_r H=L(X\times_r H).
\]
If we now apply part (ii) of \lemref{link-lem} to the right-Hilbert
$L((\ind X)\times_r G) - L(X\times_r H)$ bimodule $Z=V_H^G(L)$, we
get $pZp=V_H^G(A)$ and $qZq=V_H^G(B)$ (which follows easily from
the construction of $V_H^G$ as given in \thmref{thm-green-ind}).
If we can show that in this situation
the maps $\Phi$ and $\Psi$ are both surjective,
we directly obtain the desired isomorphism
\[
V_H^G(A)\otimes_{A\times_r H}\big(X\times_r H\big) \cong
\big((\ind X)\times_r G\big)\otimes_{(\ind B)\times_r G} V_H^G(B),
\]
and it only remains to show the equivariance of the diagram with
respect to appropriate actions and coactions.
So we pause to obtain the desired surjectivity result.
It follows from item (iv) of \lemref{link-lem} together with:

\begin{lem}\label{lem-phipsisurjective}
     Suppose that $({_AX_B}, \gamma)$
     is an $\alpha-\beta$ compatible action on the
     partial $A-B$ imprimitivity bimodule $X$,
     and let $L(X)$ denote the linking algebra of $_AX_B$
     with corresponding action $\nu$.
     As above, we regard $V_H^G(L)$ as a right-Hilbert
     $L((\ind X)\times_rG)-L(X\times_rH)$ bimodule.
     Then, with the usual meanings of $p$ and $q$, we have
     $$_{(\ind A)\times_rG}\lk pV_H^G(L)q, pV_H^G(L)q\rk \subseteq
     {}_{(\ind A)\times_rG}\overline{\lk (\ind X)\times_rG, 
(\ind X)\times_rG\rk }
     \quad\text{and}$$
     $$\lk pV_H^G(L)q, pV_H^G(L)q\rk_{B\times_rH}\subseteq
     \overline{\lk X\times_rH, X\times_rH\rk}_{B\times_rH}.$$
\end{lem}
\begin{proof}
     By construction (see \thmref{thm-green-ind}),
     $V_H^G(L)$ is a completion of
     $$C_c(G, L(X))=\mtx{C_c(G,A)& C_c(G,X)\\ C_c(G,\rev{X})&C_c(G,B)}$$
     with module actions given by the formulas \eqref{eq-ind-imp}.
     It follows that $pV_H^G(L)q$ is the closure of
     the upper right corner $C_c(G, X)$.
     Using the the multiplication rule in $L(X)$ together
     with the formulas
     for the $L(X)\times_rH$-valued inner product on
     $C_c(G, L(X))$ (and the fact that
     this inner product restricted to $C_c(G,X)$
     takes values in $C_c(H,B)\subseteq
     B\times_rH$) we get for all $x,y\in C_c(G,X)$
     $$
      \lk x,y\rk_{C_c(H,B)}(h)=
      \Delta_H(h)^{-1/2}
       \int_G \lk x(t^{-1}), \gamma_h(y(t^{-1}h))\rk_B \,dt,$$
     which determines an element in $C_c(H, B_X)$, with
     $B_X=\overline{\lk X,X\rk}_B$. But this implies
     that
     $$\lk C_c(G,X), C_c(G,X)\rk_{B\times_rH}\subseteq C_c(G,B_X)
     \subseteq B_X\times_rH=\overline{\lk X\times_rH,
     X\times_rH\rk}_{B\times_rH},$$
     where the last equation follows from \lemref{lem-linkingaction}.
     A very similar argument shows that
     $$_{(\ind A)\times_rG}\lk C_c(G,X), C_c(G,X)\rk
     \subseteq C_c(G, \ind A_X)\subseteq
     {}_{(\ind A)\times_rG}\overline{\lk (\ind X)\times_rG,(\ind X)\times_rG\rk}.
     $$
\end{proof}

Of course, as the statement of \thmref{ind-act natural} indicates, we
also need to show that the above isomorphism is equivariant with
respect to the given coactions.  For this to make sense we first have
to define the appropriate coactions on Green's bimodules $V_H^G$,
which we shall do in the next section.  However, as soon as this has
been done in a coherent way (\ie, in such a way that the coaction on
$V_H^G(L)$ compresses to the corresponding coactions on the corners
$V_H^G(A)$ and $V_H^G(B)$, respectively), then this equivariance will
follow from:

\begin{lem}
\label{link-coact-lem}
Suppose ${}_AX_B$ and ${}_CY_D$ are right-partial
imprimitivity bimodules and $Z$ is
a right-Hilbert $L(X)-L(Y)$ bimodule, and let $p$ and $q$ be as in
\lemref{link-lem}.  Suppose in addition that $\tau$ is a $\mu-\nu$
compatible coaction of $G$ on $Z$, where $\mu$ and $\nu$ are coactions
of $G$ on $L(X)$ and $L(Y)$ arising from coactions $\zeta$ and $\eta$
on $X$ and $Y$, respectively \textup(as in \eqeqref{eq-LXcoaction}\textup).
Then\textup:
\begin{enumerate}
\item
$\tau$ restricts to give right-Hilbert bimodule coactions of $G$ on
$pZq$ and $qZq$ which are compatible with the appropriate
coefficient coactions.

\item
Assume further that $\tau$ restricts to a right-Hilbert bimodule coaction
on $pZp$.
Then the isometric homomorphisms $\Phi\:X\otimes_B qZq\to pZq$ and $\Psi\:pZp
\otimes_C Y\to pZq$ of \lemref{link-lem} are equivariant for the
appropriate coactions.
\end{enumerate}
\end{lem}

\begin{proof}
By definition, if $\delta$ is the left coefficient coaction of
$\zeta$, then
\[
\mu(p)
=\smtx{
\delta(1) & 0
\\
0 & 0}
=\pproj\otimes1
=p\otimes1
\]
in $M(L(X)\otimes C^*(G))$, and similarly $\mu(q)=q\otimes 1$,
$\nu(p)=p\otimes 1$, and $\nu(q)=q\otimes 1$.  Thus,
since $q$ is full (since the left inner products on
$_AX_B$ and $_CY_D$ are full by assumption), it follows from
\lemref{pq-hom-lem} that the restriction of $\tau$ to $pZq$ is a
nondegenerate right-Hilbert bimodule homomorphism into $M((p\otimes
1)(Z\otimes C^*(G))(p\otimes 1))=M(pZp\otimes C^*(G))$, and similarly
for $qZq$.  Since it is immediate that each of these
restrictions satisfies the appropriate versions of 
conditions~(i)--(ii) in \defnref{rHb-co-defn}, this establishes~(i) above.

To check the equivariance of $\Phi$ in (ii), first note that if
\[
\Theta \:
(X\otimes C^*(G))\otimes_{B\otimes C^*(G)} (qZq\otimes C^*(G))
\xrightarrow{\cong}
(X\otimes_B qZq)\otimes C^*(G)
\]
is the isomorphism from \lemref{Theta-lem}, then the extension to the
$G$-multiplier bimodules satisfies
\[
(\Phi\otimes \id)\circ\Theta(u\otimes v)=u\d v
\]
for $u \in M_G(X\otimes C^*(G)) \subseteq M_G(L(X)\otimes C^*(G))$ and
$v \in M_G(qZq\otimes C^*(G)) \subseteq M_G(Z\otimes C^*(G))$.  This can
be seen by taking elementary tensors $x\otimes c \in X \odot C^*(G)$
and $z\otimes d \in qZq \odot C^*(G)$ and computing
\[
(\Phi\otimes \id)\circ\Theta
\bigl((x\otimes c)\otimes(z\otimes d)\bigr)
=(\Phi\otimes \id)
\bigl((x\otimes z)\otimes cd\bigr)
=x\d z\otimes cd
=(x\otimes c)\d(z\otimes d),
\]
and then extending to the $G$-multiplier bimodules.

Now, by definition, $(\zeta \cotimes_B \tau)(x\otimes z)=
\Theta(\zeta(x)\otimes\tau(z))$ for $x \in X$ and $z \in qZq$.
Thus, we can compute:
\[
(\Phi\otimes \id)\circ (\zeta \cotimes_B \tau)(x\otimes z)
=(\Phi\otimes \id)\circ\Theta(\zeta(x)\otimes \tau(z))
=\zeta(x)\d \tau(z)
=\tau(x\d z)
=\tau\circ\Phi(x\otimes z).
\]
The equivariance of $\Psi$ is proved similarly.
\end{proof}

For the proofs of \thmref{G-act natural} and \thmref{coact natural} we
shall also need a version of \lemref{link-coact-lem} which handles
actions:

\begin{lem}
\label{link-act-lem}
Suppose ${}_AX_B$ and ${}_CY_D$ are partial imprimitivity bimodules and $Z$ is
a right-Hilbert $L(X) - L(Y)$ bimodule, and let $p$ and $q$ be as in
\lemref{link-lem}.  Suppose in addition that $\sigma$ is a
$\mu - \nu$ compatible action of $G$ on $Z$, where
$\mu$ and $\nu$ are actions of $G$ on $L(X)$ and
$L(Y)$ arising from actions $\gamma$ and $\rho$ on $X$ and $Y$,
respectively.  Then\textup:
\begin{enumerate}
\item
$\sigma$ restricts to give right-Hilbert bimodule actions of $G$ on
$pZp$, $pZq$, and $qZq$ which are compatible with the appropriate
coefficient actions.

\item
The isometric homomorphisms $\Phi\:X\otimes_B qZq\to pZq$ and $\Psi\:pZp
\otimes_C Y\to pZq$ of \lemref{link-lem} are equivariant for the
appropriate actions.
\end{enumerate}
\end{lem}

\begin{proof}
If $\alpha$ is the left coefficient action of $\gamma$, then for all
$s\in G$ we have
\[
\mu_s(p)=
\left(
\begin{matrix}
\alpha_s(1) & 0
\\
0 & 0
\end{matrix}
\right)
=\left(
\begin{matrix}
1 & 0
\\
0 & 0
\end{matrix}
\right)
=p,
\]
and similarly
$\mu_s(q)=q$, $\nu_s(p)=p$, and $\nu_s(q)
=q$. Thus
\[
\sigma_s(pZp)
=\mu_s(p)\sigma_s(Z)\nu_s(p)
=pZp,
\]
and similarly for $pZq$ and $qZq$, so $\sigma$ restricts appropriately
to the various corners.  It is immediate that these restrictions are
compatible with the coefficient actions as indicated, so this
establishes~(i).

For~(ii), fix $x \in X$, $z \in qZq$, and $s \in G$ and compute, for
example, that
\[
\Phi\circ (\gamma\otimes \sigma)_s(x\otimes z)
=\Phi(\gamma_s(x)\otimes \sigma_s(z))
=\gamma_s(x)\d\sigma_s(z)
=\sigma_s(x\d z)
=\sigma_s\circ\Phi(x\otimes z).
\]
Similar calculations show that $\Psi$ is appropriately equivariant;
this completes the proof.
\end{proof}

Finally, the following lemma will help to do the
nondegenerate homomorphism parts of our main theorems.

\begin{lem}
\label{link-hom-lem}
Suppose ${}_\phi\Phi_\psi\:{}_AX_B\to M({}_CY_D)$ is a nondegenerate
imprimitivity bimodule homomorphism.  Then\textup:
\begin{enumerate}
\item
The following diagram \textup(of morphisms\textup) commutes\textup:
\begin{equation}
\label{link-hom-diag}
\xymatrix{
{A}
\ar[r]^-{X}
\ar[d]_{\phi}
& {B}
\ar[d]^{\psi}
\\
{C}
\ar[r]_-{Y}
& {D.}
}
\end{equation}

\item
If $\Phi$ is equivariant for imprimitivity bimodule coactions of $G$
on $X$ and $Y$, then \diagref{link-hom-diag} commutes
equivariantly for the appropriate coactions\textup.

\item
If $\Phi$ is equivariant for imprimitivity bimodule actions of $G$ on
$X$ and $Y$, then \diagref{link-hom-diag} commutes equivariantly
for the appropriate actions.
\end{enumerate}
\end{lem}

\begin{proof}
Part (i) is exactly \cite[Lemma~5.3]{kqr:resind}; we re-prove it here for
completeness, and because it follows easily from \lemref{link-lem}.
Simply note that, by \cite[Remark (2) of Appendix]{er:mult}, $\Phi$
gives rise to a nondegenerate homomorphism
\[
\smtx{
\phi &\Phi
\\
\widetilde{\Phi} &\psi
}
\: L(X)\to M(L(Y)).
\]
This makes $L(Y)$ into a right-Hilbert $L(X) - L(Y)$ bimodule $Z$ such
that $pZp$ is $C$ with the $A - C$ bimodule structure from $\phi$
and $qZq$ is $D$ with the $B - D$ bimodule structure from $\psi$.
Hence \lemref{link-lem} (iii) provides an isomorphism
\begin{equation}
\label{link-hom-isom}
X\otimes_B D \cong C\otimes_C Y
\end{equation}
of right-Hilbert $A - D$ bimodules, which establishes the
commutativity of \diagref{link-hom-diag}.

For (ii), suppose the coactions on $X$ and $Y$ are $\zeta$ and $\eta$.
The $\zeta - \eta$ equivariance of $\Phi$ implies that $\nu$
is $\mu - \nu$ compatible, where $\mu$ and
$\nu$ are the associated coactions on the linking algebras, so
\lemref{link-coact-lem} gives the equivariance of
\diagref{link-hom-diag} for the coactions.  Similarly, (iii)
follows from \lemref{link-act-lem}.
\end{proof}

\section{Green's Theorem for induced algebras}
\label{sec-green1}

This section is devoted to the proof of \thmref{ind-act natural}.
We start with the construction of the appropriate coaction
$\delta_V$ on Green's
$((\ind_H^G A)\times_{\ind\alpha,r} G) - (A\times_{\alpha,r} H)$
imprimitivity bimodule $V_H^G(A,\alpha)$.
\thmref{thm-green-ind}
tells us that $V_H^G(A)$ (we streamline the notation when confusion is
unlikely) is a completion of the $C_c(G,\ind_H^GA) - C_c(H,A)$
pre-imprimitivity bimodule $C_c(G,A)$ with actions and inner products
given in \eqeqref{eq-ind-imp}.
We shall make
extensive use of the embedding of $C_c(G, M^{\beta}(A))$ into
$M(A\times_r G)$ as provided by \corref{corconv}.
Moreover, \propref{propgreen} shows that $C_c(G,M^\beta(A))$ also
embeds in $M(V_H^G(A))$.
As we shall make more precise in \lemref{V-mult} below,
we may similarly view $C_c(G, M^\beta(A\otimes C^*(G)))$ as a
subspace of $M(V_H^G(A)\otimes C^*(G))$, which then allows us
to state:

\begin{thm}
\label{delta V}
Let $G$ be a locally compact group, and let 
$(A,H,\alpha)$ be an action of a closed subgroup $H$ of $G$. There
is a unique imprimitivity bimodule coaction $\delta_V$ of $G$ on
$V^G_H(A,\alpha)$ implementing a Morita equivalence between the
coactions $\what{\ind\alpha}$ on $(\ind A)\times_{\ind\alpha,r} G$ and $
\infl\hat\alpha$ on $A\times_{\alpha,r} H$, such that for $x \in C_c(G,A)$ we
have
\begin{equation}
\label{delta V eq}
\delta_V(x) \in C_c(G,M^\beta(A\otimes C^*(G))) \subseteq
M\big(V_H^G(A)\otimes C^*(G)\big),
\end{equation}
with 
$\delta_V(x)(s)=x(s)\otimes s$ for $s\in G$.
\end{thm}

\begin{rem}
We should point out that a Morita equivalence for the coactions
$ \infl\hat\alpha$ and $\what{\ind\alpha}$ has been obtained
also in \cite{er:induced} as a corollary of a more general
equivariance result for the dual coactions appearing in the
symmetric imprimitivity theorem. However, the realization of
the coaction on the bimodule $V_H^G(A)$ was quite different
from that given in the above theorem;
the realization we give here
makes it much easier to obtain certain equivariance results
for our bimodules.
\end{rem}

Before we come to the proof of this theorem, we have to explain a bit
more how to work with strictly continuous compactly-supported functions
in the various multiplier algebras and bimodules.
Observe first of all that if $H$ is a closed subgroup of $G$
and $(A,H,\alpha)$ is an action, then for any
$C^*$-algebra $C$ the canonical embedding of $C_b(G,A)\odot C$ into
$C_b(G,A\otimes C)$ determines an isomorphism
$\psi$ of $(\ind A)\otimes C$ onto $\ind(A\otimes C)$
which is equivariant for the actions $(\ind\alpha)\otimes \id$ and
$\ind(\alpha\otimes \id)$.
We use the (inverse of the) composition
\[
\xymatrix{
{((\ind A)\times_r G)\otimes C}
\ar[r]^-{\phi}
& {((\ind A)\otimes C)\times_r G}
\ar[r]^{\psi\times_r G}
& {(\ind(A\otimes C))\times_r G}
}
\]
to embed $C_c(G,\ind(A\otimes C))$
(regarded as functions of two
variables as in the discussion following
\eqeqref{eq-sym-pre} in 
\appxref{imprim-chap})
into $((\ind A)\times_r G)\otimes C$, and to embed
$C_c(G,M^\beta(\ind(A\otimes C)))$ into 
$M(((\ind A)\times_r G)\otimes C)$%
\footnote{Note that for $f \in C_c(G,\ind A)$ and $c \in C$ we have
\[
\begin{split}
(\psi\times_r G)\circ\phi(f\otimes c)(s,t)
&
=\psi(\phi(f\otimes c)(s,\cdot))(t)
=\psi(f(s,\cdot)\otimes c)(t)
=f(s,t)\otimes c.
\end{split}
\]
}.
This allows us to work with the dual coaction $\what{\ind\alpha}$ in
terms of two-variable functions:

\begin{lem}
\label{lem-ind-funct}
With notation as above, 
for $f \in C_c(G,\ind A)$ we have
\[
\what{\ind\alpha}(f) \in
C_c(G,M^\beta(\ind(A\otimes C^*(G)))) \subseteq
M\big(((\ind A)\times_r G)\otimes C^*(G)\big),
\]
with $\what{\ind\alpha}(f)(s,t)=f(s,t)\otimes s$ for $s,t\in G$.
\end{lem}

\begin{proof}
By \lemref{alpha hat} we have
$\what{\ind\alpha}(f) \in C_c(G,M^\beta((\ind A)\otimes C^*(G)))$,
so via $\psi$
we can identify $\what{\ind\alpha}(f)$ with an
element of $C_c(G,M^\beta(\ind(A\otimes C^*(G))))$.
Using our
conventions and the formula from \lemref{alpha hat}, we compute:
\begin{equation*}
\what{\ind\alpha}(f)(s,t)
=\what{\ind\alpha}(f)(s)(t)
=(f(s)\otimes s)(t)
=f(s)(t)\otimes s
=f(s,t)\otimes s.
\end{equation*}
\end{proof}

The following lemma justifies \eqref{delta V eq}:
for any $C^*$-algebra $C$, 
it allows us to
embed $C_c(G,A\otimes C)$ into $V^G_H(A)\otimes C$, and
to embed
$C_c(G,M^\beta(A\otimes C))$ into $M(V^G_H(A)\otimes C)$.

\begin{lem}
\label{V-mult}
With notation as above, 
the canonical embedding of $C_c(G,A)\odot C$ in
$C_c(G,A\otimes C)$ extends to an imprimitivity bimodule isomorphism
\[
{}_{(\ind A\times_r G)\otimes C}
(V^G_H(A,\alpha)\otimes C)
_{(A\times_r H)\otimes C}
\cong
{}_{\ind(A\otimes C)\times_r G}
V^G_H(A\otimes C,\alpha\otimes\id)
_{(A\otimes C)\times_r H}.
\]
\end{lem}

\begin{proof}
We have introduced (in the above discussion and in
\chapref{functors-chap}) the canonical isomorphisms $\phi_L$ and
$\phi_R$ of the left and right coefficient algebras of
$V_H^G(A)\otimes C$ onto those of $V_H^G(A\otimes C)$.  
A routine calculation
now shows that the embedding of $C_c(G,A)\odot C \subseteq
V_H^G(A)\otimes C$ into $C_c(G, A\otimes C)\subseteq V_H^G(A\otimes C)$
respects both inner products.
By \lemref{pre-imp-lem}, this embedding extends uniquely to an
imprimitivity bimodule homomorphism of $V_H^G(A)\otimes C$ onto
$V_H^G(A\otimes C)$;  by \remref{rem-isoimp}, 
since the coefficient maps are isomorphisms, we are done.
\end{proof}

\begin{proof}[Proof of \thmref{delta V}]
The rule $\delta_V(x)(s) = x(s)\otimes s$ 
certainly defines a map from $C_c(G,A)
\subseteq V_H^G(A)$ to $C_c(G,M^\beta(A\otimes C^*(G))) \subseteq
M(V_H^G(A)\otimes C^*(G))$.
We first show that it preserves both inner products. By
\lemref{pre-imp-lem}, we will then know that $\delta_V$ extends
uniquely to a nondegenerate imprimitivity bimodule homomorphism from
$V_H^G(A)$ to $M(V_H^G(A)\otimes C^*(G))$.

Using the formulas given in \eqeqref{eq-ind-imp} (extended to the
various $C_c$-functions in the multiplier algebras and bimodules) we
compute, for $x,y \in C_c(G,A) \subseteq V_H^G(A)$ and $s,t\in G$:
\begin{align*}
&{}_{M(((\ind A)\times_r G)\otimes C^*(G))}
\<\delta_V(x),\delta_V(y)\>(s,t)
\\
&\qquad=\int_H (\alpha\otimes\id)_h
\bigl(\delta_V(x)(th)\delta_V(y)(s^{-1}th)^*\bigr)
\,dh\,\Delta(s)^{-1/2}\\
&\qquad=\int_H (\alpha_h\otimes\id)\bigl((x(th)\otimes th)
(y(s^{-1}th)^*\otimes h^{-1}t^{-1}s)\bigr)\,dh\,\Delta(s)^{-1/2}
\\
&\qquad=\int_H (\alpha_h\otimes\id)
\bigl(x(th)y(s^{-1}th)^*\otimes s\bigr)\,dh\,\Delta(s)^{-1/2}
\\
&\qquad=\int_H \alpha_h(x(th)y(s^{-1}th)^*)\,dh\,\Delta(s)^{-1/2}
\otimes s\\
&\qquad={}_{(\ind A)\times_r G}\< x,y\>(s,t)\otimes s
\\
&\qquad=\widehat{\ind\alpha}({}_{(\ind A)\times_r G}
\< x,y\>)(s,t),
\end{align*}
and for $h\in H$, 
\begin{align*}
&\<\delta_V(x),\delta_V(y)\>
_{M((A\times_r H)\otimes C^*(G))}(h)
\\
&\qquad=\int_G\delta_V(x)(s^{-1})^*
(\alpha\otimes\id)_h(\delta_V(y)(s^{-1}h))\,ds\,\Delta(h)^{-1/2}
\\
&\qquad=\int_G (x(s^{-1})^*\otimes s)
\alpha_h\otimes\id(y(s^{-1}h)\otimes s^{-1}h)
\,ds\,\Delta(h)^{-1/2}
\\
&\qquad=\int_G x(s^{-1})^*\alpha_h(y(s^{-1}h))\otimes h
\,ds\,\Delta(h)^{-1/2}
\\
&\qquad=\int_G x(s^{-1})^*
\alpha_h(y(s^{-1}h))\,ds\,\Delta(h)^{-1/2}\otimes h
\\
&\qquad=\< x,y\>_{A\times_r H}(h)\otimes h
=\hat\alpha(\< x,y\>_{A\times_r H})(h).
\end{align*}

By \remref{rHb-co-rem}, it now suffices to verify the coaction identity. 
For $s\in G$:
\[
(\delta_V\otimes \id)\circ\delta_V(x)(s)
=x(s)\otimes s\otimes s
=(\id\otimes\delta_G)\circ\delta_V(x)(s).
\]
\end{proof}

We are finally ready to prove \thmref{ind-act natural}.
As discussed earlier, the result is equivalent to:

\begin{thm}
\label{Ind-thm}
Let $G$ be a locally compact group, and let $H$ be a closed subgroup of
$G$.
If $({}_AX_B,H,\gamma)$ is a right-Hilbert bimodule action, then
the diagram
\begin{equation}
\label{Ind-thm-diag}
\xymatrix{
{(\ind_H^G A)\times_r G}
\ar[d]_{(\ind_H^G X)\times_r G}
\ar[r]^-{V_H^G(A)}
& {A\times_r H}
\ar[d]^{X\times_r H}
\\
{(\ind_H^G B)\times_r G}
\ar[r]_-{V_H^G(B)}
& {B\times_r H}
}
\end{equation}
commutes equivariantly for
the coactions $\delta_{V(A)}$, $\delta_{V(B)}$,
$\what{\ind\gamma}$, and $ \infl\hat\gamma$ of $G$.
\end{thm}

\begin{proof}
As indicated previously, 
our strategy is to factor the right-Hilbert
bimodule.
If $X$ is a partial imprimitivity bimodule, \lemref{link-coact-lem}
together with the discussion preceding that lemma show
that our
construction of the coaction $\delta_V$ on $V_H^G$
completes the proof in this case, since it follows from
the formula for
$\delta_V$ as given in \thmref{delta V} that the
associated
coaction $\delta_{V(L)}$ on the linking algebra $L=L(X)$
compresses to the coactions $\delta_{V(A)}$ and $\delta_{V(B)}$
on the corners.

For the nondegenerate homomorphism version, assume that $\phi\:A\to
M(B)$ is a nondegenerate homomorphism which is equivariant for actions
$\alpha$ and $\beta$ of $H$.  We must show that the diagram
\begin{equation}
\label{Ind-hom-diag}
\xymatrix{
{(\ind_H^G A)\times_r G}
\ar[d]_{\ind_H^G\phi\times_r G}
\ar[r]^-{V_H^G(A)}
& {A\times_r H}
\ar[d]^{\phi\times_r H}
\\
{(\ind_H^G B)\times_r G}
\ar[r]_-{V_H^G(B)}
& {B\times_r H}
}
\end{equation}
commutes equivariantly for the appropriate coactions.  For this we
show that there is a nondegenerate imprimitivity bimodule homomorphism
$\Psi \: V_H^G(A)\to M(V_H^G(B))$ which has coefficient maps
$\ind_H^G\phi \times_r G$ and $\phi\times_r H$ and also is
$\delta_{V(A)} - \delta_{V(B)}$ equivariant.  The commutativity of the
diagram will then follow from \lemref{link-hom-lem}.

For $x \in C_c(G,A) \subseteq V_H^G(A)$, define
$\Psi(x)\:G\to M(B)$ by
\begin{equation}
\label{eq-Ind-thm}
\Psi(x)(s)=\phi(x(s)).
\end{equation}
Then $\Psi$ is a map from $C_c(G,A)$ to
$C_c(G,M^\beta(B)) \subseteq M(V_H^G(B))$.
By \lemref{pre-imp-lem}, if we can show that $\Psi$ preserves both
inner products, then we will know that it extends uniquely to a
nondegenerate imprimitivity bimodule homomorphism from $V_H^G(A)$ to
$M(V_H^G(B))$.
For $x,y \in C_c(G,A)$ and $s,t\in G$ we have
\begin{align*}
{}_{M((\ind B)\times_r G)}\<\Psi(x),\Psi(y)\>(s,t)
&=\int_H \beta_h\bigl(\Psi(x)(th)\Psi(y)(s^{-1}th)^*\bigr)\,dh
\,\Delta_G(s)^{-1/2}
\\&=\int_H \beta_h\bigl(\phi(x(th))\phi(y(s^{-1}th)^*\bigr)\,dh
\,\Delta_G(s)^{-1/2}
\\&=\phi\left(\int_H \alpha_h\bigl(x(th)y(s^{-1}th)^*\bigr)\,dh
\,\Delta_G(s)^{-1/2}\right)
\\&=\phi\bigl({}_{(\ind A)\times_r G}\<x,y\>(s,t)\bigr)
\\&=(\ind_H^G\phi\times_r G)
\bigl({}_{(\ind A)\times_r G}\<x,y\>\bigr)(s,t),
\end{align*}
and for $h\in H$, 
\begin{align*}
\<\Psi(x),\Psi(y)\>_{M(B\times_r H)}(h)
&=\int_G \Psi(x)(s^{-1})^*\beta_h(\Psi(y)(s^{-1}h))\,ds
\,\Delta(h)^{-1/2}
\\ &=\int_G\phi(x(s))^*\beta_h(\phi(y(s^{-1}h)))\,ds
\,\Delta(h)^{-1/2}
\\ &=\phi\left(\int_G x(s^{-1})^*\alpha_h(y(s^{-1}h))\,ds
\,\Delta(h)^{-1/2}\right)
\\
&=\phi\bigl(\< x,y{\>}_{A\times_r H}(h)\bigr)
=(\phi\times_r H)\bigl(\< x,y{\>}_{A\times_r H}\bigr)(h).
\end{align*}

To check equivariance of the coactions, for $x \in C_c(G,A)$  and 
$s\in G$ we have:
\begin{equation*}
(\Psi\otimes \id)\circ\delta_{V(A)}(x)(s)
=(\phi\otimes \id)\bigl(\delta_{V(A)}(x)(s)\bigr)
=\phi(x(s))\otimes s
=\Psi(x)(s)\otimes s
=\delta_{V(B)}\bigl(\Psi(x)\bigr) (s).
\end{equation*}

Assume now that $_AX_B$ is an arbitrary right-Hilbert $A-B$ bimodule.
Using \propref{decomaction}, we equivariantly factor ${}_AX_B$ as the
composition of a  right-partial
imprimitivity bimodule ${}_KX_B$ and a nondegenerate
homomorphism $\phi\:A\to M(K)$.  By our previous results we now know
that for both pieces the appropriate diagrams commute.  This tells us
the upper and lower trapezoids of the diagram
\[
\xymatrix{
{(\ind A)\times_r G}
\ar[rrr]^-{V_H^G(A)}
\ar[dr]^{\ind\phi\times_r G}
\ar[dd]_{(\ind X)\times_r G}
& & &{A\times_r H}
\ar[dl]^{\phi\times_r H}
\ar[dd]^{X\times_r H}
\\
& {(\ind K)\times_r G}
\ar[r]^-{V_H^G(K)}
\ar[dl]^{(\ind X)\times_r G}
& {K\times_r H}
\ar[dr]_{X\times_r H}
\\
{(\ind B)\times_r G}
\ar[rrr]_-{V_H^G(B)}
& & &{B\times_r H}
}
\]
commute equivariantly for the appropriate coactions---that is, these
trapezoids commute in the category $\c C(G)$. The left and right
triangles commute in the category $\c C(G)$ by functoriality of $\ind$
and $\times_r$. We conclude that
the outer rectangle commutes equivariantly for the appropriate
coactions, as desired.
\end{proof}

\section{Green's Theorem for induced representations}

In this section we shall derive our second main theorem, \thmref{G-act
natural}, mainly as a consequence of \thmref{ind-act natural}.  As
mentioned in the discussion preceding
\eqeqref{eq-green-imp-thm}, if we start with an action $\alpha$ of $G$
on $A$, then there is a natural isomorphism between 
$(\ind_H^G A)\times_{\ind\alpha,r}G$ and $C_0(G/H,A)\times_{\alpha\otimes
\tau,r}G$, where $\ind\alpha$ denotes the action of $G$ induced from
the restriction of $\alpha$ to the closed subgroup $H$.  This
isomorphism is the crossed product $\phi\times_r G$ of the
isomorphism
$\phi$ of $\ind_H^G A$ onto $C_0(G/H,A)$
given by 
$\phi(f)(sH)=\alpha_s(f(s))$.
In particular, $\phi\times_rG$ is equivariant for the dual coactions
$\what{\ind\alpha}$ and $\what{\alpha\otimes\tau}$.  Together with
the map of 
$C_c(G,A)$ into itself given by 
$\Psi(x)(s)=\alpha_s(x(s))$,
we obtain an isomorphism
\[
\Psi\:
{}_{(\ind_H^G A)\times_{\ind\alpha,r}G}V_H^G(A)_{A\times_{\alpha,r}H}
\iso
{}_{C_0(G/H,A)\times_{\alpha\otimes\tau,r}G}
X_H^G(A)_{A\times_{\alpha,r}H}
\]
(see \thmref{thm-green-imp}).
This isomorphism carries the coaction $\delta_V$ on $V_H^G(A)$,
as constructed in \thmref{delta V}, to a unique coaction, 
$\delta_X$, of $G$ on $X_H^G(A)$, which makes the diagram
\[
\xymatrix{
{V^G_H(A)}
\ar[r]^-{\delta_V}
\ar[d]_{\Psi}
& {M(V^G_H(A)\otimes C^*(G))}
\ar[d]^{\Psi\otimes \id}
\\
{X^G_H(A)}
\ar[r]_-{\delta_X}
& {M(X^G_H(A)\otimes C^*(G))}
}
\]
commute.  Thus $\delta_X$ implements a Morita equivalence between
the $G$-coactions $\what{\alpha\otimes\tau}$ on $C_0(G/H,A)\times_r
G$ and $ \infl\what{\alpha|}$ on $A\times_r H$. In this special
situation, \thmref{Ind-thm} becomes:

\begin{thm}
\label{Ind-green thm}
Let $G$ be a locally compact group, and let $H$ be a closed subgroup of
$G$. 
If $({}_AX_B,G,\gamma)$ is a right-Hilbert bimodule action, then the
diagram
\begin{equation}
\label{Ind-green-thm-diag}
\xymatrix{
{C_0(G/H,A)\times_{\alpha\otimes\tau,r} G}
\ar[d]_{C_0(G/H,X)\times_{\gamma\otimes\tau,r} G}
\ar[r]^-{X_H^G(A)}
& {A\times_{\alpha,r} H}
\ar[d]^{X\times_{\gamma,r} H}
\\
{C_0(G/H,B)\times_{\beta\otimes\tau,r} G}
\ar[r]_-{X_H^G(B)}
& {B\times_{\beta,r} H}
}
\end{equation}
commutes equivariantly for the coactions $\delta_{X(A)}$,
$\delta_{X(B)}$, $\what{\gamma\otimes\tau}$, and $
\infl\what{\gamma|}$ of $G$.
\end{thm}

\begin{rem}
\label{delta X rem}
Later (in the proof of \thmref{eqvt-ekr-thm}) we will need a
$C_c$-formula for the coaction $\delta_X$ of $G$ on $X_H^G(A)$.  Using
the isomorphism $\Psi\:V\to X$, it is routine to check that in fact
$\delta_X$ is given by the same formula as $\delta_V$: for $x \in
C_c(G,A) \subseteq X_H^G(A)$ we have
\[
\delta_X(x) \in C_c\bigl(G,M^{\beta}(A\otimes C^*(G))\bigr)
\subseteq M\bigl(X_H^G(A)\otimes C^*(G)\bigr),
\]
and $\delta_X(x)(s)=x(s)\otimes s$ for $s\in G$.
\end{rem}

We now want to specialize to the case where $H$ is \emph{normal} in
$G$, and in order to make this more evident we shall from now on
write $N$ instead of $H$ for our subgroup.  In this situation the
above results may be restated in terms of duality theory and dual
coactions, as follows from:

\begin{lem}
\label{left alg}
If $(A,G,\alpha)$ is an action and $N$ is a closed normal subgroup of
$G$, then
\[
C_0(G/N,A)\times_{\alpha\otimes\tau,r} G
\cong A\times_{\alpha,r} G\times_{\hat\alpha|} G/N,
\]
equivariantly for the coactions
$\what{\alpha\otimes \tau}$ and $\hat\alpha\dec$ of $G$.
\end{lem}

\begin{proof}
\thmref{thm-isom-red} gives
an isomorphism
$\psi$ of $C_0(G/N,A)\times_{\alpha\otimes\tau,r}G$ onto 
$A\times_{\alpha,r}G \times_{\what\alpha|}G/N$
such that
\[
\psi\bigl(i^r_{C_0(G/N,A)}(a\otimes f)i^r_G(g)\bigr)
=j_{G/N}(f)j_{A\times_r G}(i^r_A(a)i^r_G(g))
\]
for $a \in A$, $f \in C_0(G/N)$, and $g \in C_c(G)$.
We only have to check that this isomorphism
is equivariant for the coactions $\what{\alpha\otimes\tau}$
and $\hat\alpha\dec$: 
\begin{align*}
&
\hat\alpha\dec\circ
\psi\bigl(i^r_{C_0(G/N,A)}(a\otimes f) i^r_G(g)\bigr)
=\hat\alpha\dec\bigl(
j_{G/N}(f) j_{A\times_r G} (i^r_A(a) i^r_G(g))\bigr) \\
&\qquad
=(j_{G/N}(f)\otimes 1)(j_{A\times_r G}\otimes \id)
\circ \hat\alpha(i^r_A(a)i^r_G(g)) \\
&\qquad
=(j_{G/N}(f)\otimes 1)(j_{A\times_r G}\otimes \id)
\Bigl(\int_G g(s)i^r_A(a)i^r_G(s)\otimes s \,ds \Bigr) \\
&\qquad
=(j_{G/N}(f)\otimes 1) \int_G g(s)j_{A\times_r G}(i^r_A(a)i^r_G(s))
\otimes s \,ds \\
&\qquad
=\int_G g(s)j_{G/N}(f)j_{A\times_r G}(i^r_A(a)i^r_G(s))\otimes s \,ds \\
&\qquad
=\int_G g(s)\psi\bigl(i^r_{C_0(G/N,A)}(a\otimes f)i^r_G(s)\bigr)
\otimes s \,ds \\
&\qquad
=\int_G g(s)(\psi\otimes \id)\bigl(
i^r_{C_0(G/N,A)}(a\otimes f)i^r_G(s)\otimes s\bigr) \,ds \\
&\qquad
=(\psi\otimes \id) \Bigl(\int_G g(s)\what{\alpha\otimes\tau}
\bigl(i^r_{C_0(G/N,A)}(a\otimes f)i^r_G(s)\bigr) \,ds \Bigr) \\
&\qquad
=(\psi\otimes \id)\circ \what{\alpha\otimes\tau}
\bigl(i^r_{C_0(G/N,A)}(a\otimes f)i^r_G(g)\bigr) .
\end{align*}
\end{proof}

Of course it follows from the above lemma that, in case of a
normal subgroup $N$ of $G$, we may replace
Diagram~\eqref{Ind-green-thm-diag}
by the diagram
\begin{equation}
\label{G-thm-diag}
\xymatrix{
{A\times_r G\times G/N}
\ar[r]^-{X_N^G(A)}
\ar[d]_{X\times_r G\times G/N}
& {A\times_r N}
\ar[d]^{X\times_r N}
\\
{B\times_r G\times G/N}
\ar[r]_-{X_N^G(B)}
& {B\times_r N,}
}
\end{equation}
which commutes equivariantly for the coactions $\hat\gamma\dec$ and
$ \infl\what{\gamma|}$ (where $\gamma$ denotes
the given action on
$_AX_B$). However, in
order to get a complete proof of
\thmref{G-act natural} we need to show that
the above diagram is also equivariant with respect to the $G$-actions
$ \infl\what{\hat\gamma|}$ on
$X\times_r G\times G/N$ and
$\gamma\dec$ on
$X\times_r N$.
For this we have to construct a $G$-action
on the
bimodule $X_N^G(A)$ which is compatible
with the actions
$ \infl\what{\hat\alpha|}$ and $\alpha\dec$.
Such an action was constructed in
\cite[Theorem~1]{EchME}, 
for actions twisted over $N$, full (twisted) crossed products, and 
with slightly different modular functions in the
formulas for the module actions and inner products.
The following lemma and its proof are an adaptation of that
construction to the present situation.

\begin{lem}
\label{alpha X}
If $(A,G,\alpha)$ is an action and $N$ is a closed normal subgroup of
$G$, there is a unique imprimitivity bimodule action $\alpha^X$
of $G$ on $X_N^G(A)$, implementing a Morita equivalence between the
actions $\infl\what{\hat\alpha|}$ on $A\times_r G\times G/N$ and
$\alpha\dec$ on $A\times_r N$, such that for $x \in C_c(G,A)$ and 
$t\in G$ we have
\begin{equation}
\label{alpha X eq}
\alpha^X_s(x)(t)=\Delta_{G,N}(s)^{1/2} x(ts).
\end{equation}
\end{lem}

Recall that the modular function $\Delta_{G,N}$ for conjugation of $G$
on $N$ is given by
\[
\Delta_{G,N}(s)=\Delta_G(s)\Delta_{G/N}(sN)^{-1}.
\]

\begin{proof}
Formula~\eqref{alpha X eq} certainly defines a homomorphism of $G$ into
the group of linear automorphisms of the vector space $C_c(G,A)$. We
need to show that for $x,y \in C_c(G,A) \subseteq X_N^G(A)$,
$s\in G$,  and $b \in
C_c(G\times G/N,A) \subseteq A\times_r G\times G/N$ we have
\begin{enumerate}
\item
$\alpha^X_s(b\d x)
=\infl\what{\hat\alpha|}_s(b)\d\alpha^X_s(x)$,

\item
$\<\alpha^X_s(x),\alpha^X_s(y) \>_{A\times_r N}
=\alpha\dec_s(\<x,y\>_{A\times_r N})$, and

\item
$s\mapsto\alpha^X_s(x)$ is continuous for the norm of $X_N^G(A)$.
\end{enumerate}
Note that \thmref{thm-isom-red} not only
allows us to identify $C_c(G\times
G/N,A)$ as a dense subalgebra of $A\times_r G\times G/N$, but also
implies that the action $\infl\what{\hat\alpha|}$ is given on
$C_c$-functions by
\begin{equation}
\label{alpha hat hat}
\infl\what{\hat\alpha|}_s(f)(t,rN)=f(t,rsN).
\end{equation}
For (i) and (ii), we compute:
\begin{align*}
\alpha^X_s(b\d x)(t)
&=\Delta_{G,N}(s)^{1/2} (b\d x)(ts)
\\&=\Delta_{G,N}(s)^{1/2}
\int_G b(r,tsN)\alpha_r(x(r^{-1}ts))
\Delta_G(r)^{1/2} \,dr
\\&=\int_G b(r,tsN)
\alpha_r\bigl(\Delta_{G,N}(s)^{1/2} x(r^{-1}ts)\bigr)
\Delta_G(r)^{1/2} \,dr
\\&=\int_G \infl\what{\hat\alpha|}_s(b)(r,tN)
\alpha_r\bigl(\alpha^X_s(x)(r^{-1}t)\bigr)
\Delta_G(r)^{1/2} \,dr
\\&=\bigl(\infl\what{\hat\alpha|}_s(b)
\cdot\alpha^X_s(x)\bigr)(t),
\end{align*}
and for $n\in N$,
\begin{align*}
\<\alpha^X_s(x),\alpha^X_s(y) \>_{A\times_r N}(n)
&=\Delta_N(n)^{-1/2}
\int_G\alpha_t\bigl(
\alpha^X_s(x)(t^{-1})^*\alpha^X_s(y)(t^{-1}n)\bigr) \,dt
\\&=\Delta_N(n)^{-1/2} \Delta_{G,N}(s)
\int_G\alpha_t\bigl(x(t^{-1}s)^* y(t^{-1}ns)\bigr) \,dt
\\&\overset{t\mapsto st}{=} \Delta_N(n)^{-1/2} \Delta_{G,N}(s)
\int_G\alpha_{st}\bigl(x(t^{-1})^* y(t^{-1}s^{-1}ns)\bigr) \,dt
\\&= \Delta_{G,N}(s)
\alpha_s \left(\Delta_N(n)^{-1/2}\int_G
\alpha_t\bigl(x(t^{-1})^* y(t^{-1}s^{-1}ns)\bigr) \,dt \right)
\\&=\Delta_{G,N}(s)
\alpha_s\bigl(\<x,y\>_{A\times_r N}(s^{-1}ns)\bigr)
\\&=\alpha\dec_s\bigl(\<x,y\>_{A\times_r N}\bigr)(n).
\end{align*}

For (iii), it suffices to observe that if $s\to e$ in $G$ then
$\alpha^X_s(x)\to x$ in the inductive limit topology of $C_c(G,A)$,
which is stronger than the norm topology from $X_N^G(A)$.
\end{proof}

We are now almost done with the proof of \thmref{G-act natural},
which is equivalent to:

\begin{thm}
\label{G-thm}
Let $G$ be a locally compact group, and let $N$ be a closed normal
subgroup of $G$. 
If $({}_AX_B,G,\gamma)$ is a right-Hilbert bimodule action,
then \diagref{G-thm-diag} commutes equivariantly for the actions
$\alpha^{X(A)}$, $\alpha^{X(B)}$, $ \infl \what{\hat\gamma|}$, and
$\gamma\dec$ of $G$,
and also for the coactions
$\delta_{X(A)}$, $\delta_{X(B)}$, $\hat\gamma\dec$, and $ \infl
\what{\gamma|}$ of $G$.
\end{thm}

\begin{proof}
We already observed that \diagref{G-thm-diag}
commutes equivariantly for the
coactions. We need to know that it is also
equivariant for the actions, and for this we must trace through the
steps of the argument of \thmref{Ind-thm}, where we factored 
${}_AX_B$ into a right-partial
imprimitivity bimodule ${}_KX_B$ and a standard bimodule ${}_AK_K$. In the
proof of the partial imprimitivity bimodule case we appealed
to Lemmas~\ref{link-lem}, \ref{lem-phipsisurjective},
and \ref{link-coact-lem} to get an isomorphism
which is equivariant for the coactions (see the discussion preceding
\lemref{lem-phipsisurjective}).
\lemref{link-act-lem} does the same job for the
actions.

Next, for the standard bimodule case, \ie, for a
nondegenerate equivariant homomorphism $\phi\:A\to M(B)$,
we appealed to \lemref{link-hom-lem} to see that
coaction-equivariance of the bimodule homomorphism $V(A)\to M(V(B))$
implied that the diagram commuted equivariantly for the coactions.
The action part of that lemma does the same job for our actions
as soon as we can show
that the bimodule homomorphism $X_N^G(A)\to M(X_N^G(B))$
is equivariant for
the appropriate actions. For this first note that the
isomorphism $V_N^G\cong X_N^G$ transforms the bimodule
homomorphism $V(A)\to M(V(B))$ of \eqeqref{eq-Ind-thm} to a
bimodule homomorphism $\Psi\:X(A)\to M(X(B))$ given by exactly the
same formula.
Thus, for $x \in C_c(G,A) \subseteq X_H^G(A)$ and $t\in G$ we
can compute:
\begin{align*}
\Psi(\alpha^{X(A)}_s(x))(t)
&=\phi(\alpha^{X(A)}_s(x)(t))
=\phi(x(ts))
=\Psi(x)(ts)
=\alpha^{X(B)}_s(\Psi(x))(t).
\end{align*}
This establishes the action-equivariance in the standard bimodule case.

Exactly as in the proof of \thmref{Ind-thm} we combine both special
cases in order to get the desired result.
\end{proof}

\section{Mansfield's Theorem}
\label{sec-mans}

In this section we are going to prove \thmref{coact natural}.  The
basic idea is the same as for the corresponding results concerning
actions as given in the previous sections: using
\propref{decomcoaction}, we equivariantly factor a given equivariant
right-Hilbert $A-B$ bimodule $X$ into a right-partial
imprimitivity bimodule and a
nondegenerate homomorphism and show that the appropriate diagrams
commute equivariantly in both special cases. 

\begin{thm}
\label{M-thm}
Suppose that $({}_AX_B,G,\zeta)$ is a right-Hilbert bimodule coaction
whose coefficient coactions are nondegenerate and normal, and that $N$
is a closed normal subgroup of $G$.  Then there is an action
$\alpha^Y$ and a coaction $\delta_Y$ of $G$ on Mansfield's bimodule
$Y_{G/N}^G$ \textup(given by the formulas \eqref{alpha Y eq} and
\eqref{delta Y eq} below\textup) such that the diagram
\begin{equation}
\label{M-thm-diag}
\xymatrix{
{A\times G\times_r N}
\ar[r]^-{Y_{G/N}^G(A)}
\ar[d]_{X\times G\times_r N}
& {A\times {G/N}}
\ar[d]^{X\times {G/N}}
\\
{B\times G\times_r N}
\ar[r]_-{Y_{G/N}^G(B)}
& {B\times {G/N}}
}
\end{equation}
commutes equivariantly for
the actions
$\alpha^{Y(A)}$, $\alpha^{Y(B)}$,
$\hat\zeta\dec$,
and
$\infl\what{\zeta|}$ of $G$,
and also for the coactions
$\delta_{Y(A)}$, $\delta_{Y(B)}$,
$\infl\what{\hat\zeta|}$,
and
$\zeta\dec$ of $G$.
\end{thm}

As already discussed, this theorem is equivalent to
\thmref{coact natural}.  We start with the construction of the action
$\alpha^Y$.  Recall from \appxref{imprim-chap} that for a
nondegenerate normal coaction $\delta$ of $G$ on a $C^*$-algebra $A$
and closed normal subgroup $N$ of $G$, the 
$(A\times_\delta G\times_{\hat\delta|,r} N)-(A\times_{\delta|} G/N)$ 
imprimitivity bimodule $Y_{G/N}^G(A)$ is obtained as
the completion of the $C_c(N,\D)-\D_N$ pre-imprimitivity bimodule $\D$
with operations given by \eqeqref{eq-actions-mans}.

\begin{prop}
\label{alpha Y}
For any nondegenerate normal coaction $(A,G,\delta)$, there is a
unique imprimitivity bimodule action $\alpha^Y$ of $G$ on
$Y_{G/N}^G(A)$, implementing a Morita equivalence between the actions
$\hat\delta\dec$ on $A\times_\delta G\times_{\hat\delta|,r} N$ 
and $ \infl\what{\delta|}$
on $A\times_{\delta|} G/N$, such that for $x \in \c D$ we have
\begin{equation}
\label{alpha Y eq}
\alpha^Y_s(x)=\Delta_{G,N}(s)^{1/2}\hat\delta_s(x).
\end{equation}
\end{prop}

\begin{proof}
Formula \eqref{alpha Y eq} certainly defines a homomorphism of $G$
into the group of linear automorphisms of the vector space $\D$.  By
\lemref{pre-imp-lem} and \remref{rHb-act-rem} it suffices to show that
for $x,y\in \D$, $s\in G$, and $g\in C_c(N,\D)$ we have
\begin{enumerate}
\item
${}_{A\times G\times_r N}\<\alpha^Y_s(x),\alpha^Y_s(y)\>
=\hat\delta\dec_s({}_{A\times G\times_r N}\<x,y\>)$,

\item
$\<\alpha^Y_s(x),\alpha^Y_s(y)\>_{A\times G/N}
=\infl\what{\delta|}_s(\< x,y\>_{A\times G/N})$, and

\item
$s\mapsto\alpha^Y_s(x)$ is continuous for the norm of $Y^G_{G/N}$.
\end{enumerate}
For (i) we compute, with $n\in N$:
\begin{align*}
{}_{A\times G\times_r N}\<\alpha^Y_s(x),\alpha^Y_s(y)\>(n)
&=\alpha^Y_s(x)\hat\delta_n(\alpha^Y_s(y)^*)\Delta(n)^{-1/2}
\\&=\Delta_{G,N}(s)\hat\delta_s(x)\hat\delta_{ns}(y^*)
\Delta(n)^{-1/2}
\\&=\Delta_{G,N}(s)\hat\delta_s\bigl(
x\hat\delta_{s^{-1}ns}(y^*)\Delta(s^{-1}ns)^{-1/2}\bigr)
\\&=\Delta_{G,N}(s)\hat\delta_s\bigl(
{}_{A\times G\times_r N}\<x,y\>(s^{-1}ns)\bigr)
\\&=\hat\delta\dec_s\bigl({}_{A\times G\times_r N}\<x,y\>\bigr)(n).
\end{align*}
In order to prove (ii), we mention first that the embedding $j_A\times
j_G|\:A\times_{\delta|}G/N\to M(A\times_{\delta}G)$ is
$\infl\what{\delta|}-\hat\delta$ equivariant.  Thus, since we identify
$A\times_{\delta|}G/N$ with its image in $M(A\times_{\delta}G)$ when
working with the dense subalgebra $\D_N$, we may compute
\begin{align*}
&\<\alpha^Y_s(x),\alpha^Y_s(y)\>_{A\times G/N}
=\int_N \hat\delta_n(\alpha^Y_s(x)^*\alpha^Y_s(y))\,dn
=\int_N \hat\delta_n\bigl(\Delta_{G,N}(s)\hat\delta_s(x^*y)\bigr)\,dn
\\&\quad\overset{n\mapsto sns^{-1}}{=}
\int_N \hat\delta_{sn}(x^*y) \,dn
=\hat\delta_s\left(\int_N \hat\delta_n(x^*y)\,dn\right)
=\infl\what{\delta|}_s(\< x,y\>_{A\times G/N}).
\end{align*}

For (iii) it suffices to observe that if $s\to e$ in $G$ then
\[
{}_{A\times G\times_r N}\< \hat\delta_s(x)-x,
\hat\delta_s(x)-x\>\to 0
\]
in the inductive limit topology of $C_c(N,A\times_\delta G)$.
\end{proof}

The $ \infl\what{\hat\delta|} -\delta\dec$
equivariant coaction $\delta_Y$ on the
bimodule $Y_{G/N}^G$ was actually constructed in
\cite[Proposition~3.2]{er:stab}. Translating the construction
given there into a more abstract setting, we get the formula
\begin{equation}
\label{delta Y eq}
\delta_Y(x)=(x\otimes 1)(j_G\otimes \id)(w_G^*)
\righttext{for} x \in \c D.
\end{equation}

\begin{proof}[Proof of \thmref{M-thm}]
We first derive the commutativity of \diagref{M-thm-diag}
in case where $_AX_B$ is a right-partial imprimitivity bimodule.
As in the action case, this follows almost automatically from
our general linking algebra techniques.
We want to apply Lemmas~\ref{link-lem}, \ref{link-act-lem}, and
\ref{link-coact-lem} to the Mansfield bimodule
$Y_{G/N}^G(L)$ of the linking algebra $L(X)$, using the identifications
\[
L(X)\times G\times_r N=L(X\times G\times_r N)
\midtext{and}
L(X)\times G/N=L(X\times G/N).
\]
Since $qY_{G/N}^G(L)q=Y_{G/N}^G(B)$ and $pY_{G/N}^G(L)p=
Y_{G/N}^G(A)$, item (ii) of \lemref{link-lem} provides
isometric $A\times G\times_rN-B\times G/N$ bimodule homomorphisms
$$\Phi\:\big(X\times G\times_r N\big)\otimes_{B\times G\times_r N}Y_{G/N}^G(B)
\to pY_{G/N}^G(L)q$$
and
$$\Psi\:Y_{G/N}^G(A)\otimes_{A\times G/N}\big(X\times G/N\big)\to
pY_{G/N}^G(L)q,$$
which, by Lemmas~\ref{link-act-lem} and~\ref{link-coact-lem},
are equivariant for the appropriate coactions. We have to show that
these maps are surjective in order to obtain an action- and 
coaction-equivariant isomorphism
\[
\big(X\times G\times_r N\big)\otimes_{B\times G\times_r N}Y_{G/N}^G(B)
\cong Y_{G/N}^G(A)\otimes_{A\times G/N}\big(X\times G/N\big).
\]
To see this we appeal to item (iv) of \lemref{link-lem} and
the following lemma, which is an analogue of
\lemref{lem-phipsisurjective}.

\pause

\begin{lem}\label{lem-phipsisurjectivemans}
The ranges of the $A\times G\times_rN$- and $B\times G/N$-
valued inner products on $pY_{G/N}^G(L)q$ lie in
range of the $A\times G\times_rN$-valued inner product
on $X\times G\times_rN$ and the range of the $B\times G/N$-valued
inner product of $X\times G/N$, respectively.
\end{lem}
\begin{proof}
  Since $X$ is a right-partial $A-B$ imprimitivity bimodule,
  it follows from \lemref{lem-linkcoactcross} and
  \lemref{lem-linkingaction} that
  $X\times G\times_rN$ is a right-partial
$A\times G\times_rN-B\times G\times_rN$ imprimitivity bimodule.
Thus the $A\times G\times_rN$-valued inner product
on $X\times G\times_rN$ is full, which clearly implies that it
contains the $A\times G\times_rN$-valued inner product
on $pY_{G/N}^G(L)q$. To see the other inclusion, let
$B_X$ be the range of the $B$-valued inner product on $X$.
It then follows from \lemref{lem-linkcoactcross} that the
$B\times G/N$-valued inner product on $X\times G/N$ has
range $B_X\times G/N$.

On the other hand, if we follow the construction of
$Y_{G/N}^G(L)$ as described in~\secref{sec-appmansfield}
of~\appxref{imprim-chap},
we see that it is the closed linear span
(with respect to the Hilbert-module norm) of certain elements
of the form
$$j_L\left(\mtx{a&x\\ \tilde{y}& b}\right)j_G(f)$$
with $f\in C_c(G)$,
from which it follows that $pY_{G/N}^G(L)q$ is generated
by certain elements
of the form $j_X(x)j_G(f)$. It follows then from the formula
for the right inner product on $Y_{G/N}^G$ as given in
\eqeqref{eq-actions-mans}, that  the
$B\times G/N$-valued inner
product
of two such elements $ j_X(x)j_G(f), j_X(y)j_G(g)\in pY_{G/N}^G(L)q$
is given by
$$\lk j_X(x)j_{G/N}(f), j_X(y)j_{G/N}(g)\rk_{B\times G/N}=
j_G(\phi(\bar{f}))j_B(\lk x,y\rk_B)j_G(\phi(g)),$$
with $\phi(h)(sN)=\int_N h(sn)\, dn$ for $h\in C_c(G)$.
Since $\lk x,y\rk_B\in B_X$, it follows that
$$\lk j_X(x)j_{G/N}(f), j_X(y)j_{G/N}(g)\rk_{B\times G/N}
\in B_X\times G/N.$$
\end{proof}

\resume{\thmref{M-thm}}
It's a bit harder to prove the commutativity of \diagref{M-thm-diag}
in the nondegenerate homomorphism case.  
Given a nondegenerate homomorphism
$\phi\:A\to M(B)$ which is equivariant for nondegenerate normal
coactions $\delta$ and $\eps$ of $G$ on $A$ and $B$, respectively, we
shall show that the restriction of $\phi\times G\:A\times_{\delta}G\to
M(B\times_{\eps}G)$ to $\D(A)\subseteq A\times_{\delta}G$ extends to a
nondegenerate imprimitivity bimodule homomorphism
\[
\Psi\:Y_{G/N}^G(A)\to M(Y_{G/N}^G(B))
\]
which has coefficient maps $\phi\times G\times_r N$ and $\phi\times
G/N$, and is also both $\delta_{Y(A)}-\delta_{Y(B)}$ and
$\alpha^{Y(A)}-\alpha^{Y(B)}$ equivariant.
The result will then follow from \lemref{link-hom-lem}.

We first claim that if $x\in \c D(A)$ and $c\in \c D_N(B)$, then
\[
(\phi\times G)(x) c \in \c D(B).
\]
By the definition of $\c D(A)$, $\c D(B)$, and $\c D_N(B)$ (see
\defnref{defn-DN}),
it suffices to show that for each pair $u,v \in A_c(G)$ and
each pair of compact subsets $E,E' \subseteq G$, there exist
$w \in A_c(G)$ and a compact set $F \subseteq G$ such that
$(\phi\times G)(x) c \in \c D_{(w,F)}(B)$
for each $x \in \c D_{(u,E)}(A)$ and $c \in \c D_{(v,E',N)}(B)$.
Since the pairing $(x,c)\mapsto (\phi\times G)(x) c$ is certainly
norm-continuous in both variables, it suffices to take
\[
x=j_A(\delta_u(a)) j_G^A(f) \in \c D_{(u,E)}(A)
\midtext{and}
c=j_B(\epsilon_v(b)) j_G^B(g) \in \c D_{(v, E', N)}(B).
\]
To verify the claim note first that by \cite[Lemma 9]{ManIR2} there
exists a compact subset $F'$ (only depending on $E$) such that
$x$ can be approximated in norm by elements of the form
$\sum_{i=1}^n j_G^A(f_i)j_A(\delta_u(a_i))$, with
$\supp f_i \subseteq F'$ for all $i$. By the norm continuity of the
pairing $(x,c)\to (\phi\times G)(x)c$
it follows from this that we may as well assume that
$x=j_G^A(f)j_A(\delta_u(a))$ with $\supp f \subseteq F'$.
Then
\[
\begin{split}
(\phi\times G)(x)c
&
=j^B_G(f)j_B(\phi(\delta_u(a)))j_B(\epsilon_v(b))j^B_G(g)
=j^B_G(f)j_B\bigl(\epsilon_u(\phi(a))\epsilon_v(b)\bigr) j^B_G(g).
\end{split}
\]
Choose $w \in A_c(G)$ which is identically $1$ on $(\supp u)(\supp v)$.
If $\{e_i\}$ is a bounded approximate identity for $B$, then
\[
\epsilon_w\bigl(\epsilon_u(e_i\phi(a))\epsilon_v(b)\bigr)
=\epsilon_u(e_i\phi(a))\epsilon_v(b)\quad\text{for all $i$},
\]
by \cite[Lemma~1~(iii)]{ManIR2}. Since
$\epsilon_u(e_i\phi(a))\to\epsilon_u(\phi(a))$
by strict continuity of $\epsilon_u$ on $M(B)$, we get
$\epsilon_w(\epsilon_u(\phi(a))\epsilon_v(b))
=\epsilon_u(\phi(a))\epsilon_v(b)$
by norm continuity of $\epsilon_w$ on $B\times_\epsilon G$.
Thus
\[
(\phi\times G)(x)c=j_G^B(f)j_B(\epsilon_w(d))j_G^B(g),
\]
with $d=\epsilon_u(\phi(a))\epsilon_v(b) \in B$.
Another use of \cite[Lemma 9]{ManIR2} reveals that there
exists a compact subset $F \subseteq G$, only depending on
$F'$ (and hence on $E$) such that
$j_G^B(f)j_B(\epsilon_w(d))j_G^B(g) \in \D(B)_{(w,F)}$.
This proves the claim. For later use note
that when $N=\{e\}$ the above
computations give
\begin{equation}
\label{eq-pairing1}
(\phi\times G)(x)y \in \D(B)
\righttext{for}
x \in \D(A),\  y \in \D(B).
\end{equation}

In the next step we show that the pairing
$(x,c)\mapsto (\phi\times G)(x)c$ is actually continuous
with respect to the norms on $\D(A)$ and $\D(B)$ inherited
from the bimodules $Y_{G/N}^G(A)$ and $Y_{G/N}^G(B)$, respectively.
For this recall from \eqeqref{eq-actions-mans} that the
$A\times_{\delta|} G/N$-valued inner product is given by
\[
\< x,y\>_{A\times G/N}=\int_N \hat\delta_n(x^*y) \,dn,
\]
where for $x$ in $\D(A)$, the integral $ \int_N \hat\delta_n(x) \,dn$
is determined by the equations
\[
\omega \left(\int_N \hat\delta_n(x) \,dn \right)
=\int_N \omega(\hat\delta_n(x)) \,dn,\quad\quad \omega \in
(A\times_{\delta}G)^*.
\]
Using this, and the identity $\hat\epsilon_n(c)=c$ for all
$c \in \D_N(B)$ and $n\in N$, we now compute
for any $\omega \in (B\times_\epsilon G)^*$
\[
\begin{split}
& \omega\bigl(\<(\phi\times G)(x)c,(\phi\times G)(x)c\>
_{B\times G/N}\bigr)
=\omega \left(\int_N
\hat\epsilon_n\bigl(c^*(\phi\times G)(x^*x)c\bigr) \,dn \right)
\\
&
\qquad
=\int_N \omega\bigr(c^*
\hat\epsilon_n((\phi\times G)(x^*x))c\bigr) \,dn
=\int_N c\d\omega\d c^*
\big((\phi\times G)(\hat\delta_n(x^*x))\big) \,dn
\\
&
\qquad
=\int_N (\phi\times G)^*(c\d\omega\d c^*)
(\hat\delta_n(x^*x)) \,dn
=(\phi\times G)^*(c\d\omega\d c^*)
\left(\int_N \hat\delta_n(x^*x) \,dn \right)
\\
&
\qquad
=\omega\bigl(c^*
(\phi\times G/N)(\< x,x\>_{A\times G/N})c\bigr) ,
\end{split}
\]
so that
\begin{equation}
\begin{split}
\label{bounded}
& \|(\phi\times G)(x)c\|^2_{Y_{G/N}^G(B)}
=\|\<(\phi\times G)(x)c,(\phi\times G)(x)c\>_{B\times G/N}\|
\\
&
\qquad
=\|c^*(\phi\times G)(\< x,x\>_{A\times G/N})c\|
\le\|\< x,x\>_{A\times G/N}\|\|c^*c\|
=\|x\|^2_{Y_{G/N}^G(A)} \|c\|^2.
\end{split}
\end{equation}

It follows that for every $x \in \D(A)$ the formula
$T(c)=(\phi\times G)(x) c$
defines a bounded linear map from $\D_N(B)$ to $\D(B)$, where the latter
is given the Hilbert module norm of $Y_{G/N}^G(B)$. We show that $T$ is
adjointable. For $y \in \D(B)$ and $c \in \D_N(B)$ we have
\[
\begin{split}
\< T(c),y\>_{B\times G/N}
&
=\int_N \hat\epsilon_n(c^*(\phi\times G)(x^*)y) \,dn
=\int_N c^*\hat\epsilon_n((\phi\times G)(x^*)y) \,dn
\\
&
=c^* \int_N \hat\epsilon_n((\phi\times G)(x^*)y) \,dn
=\left\< c, \int_N \hat\epsilon_n((\phi\times G)(x^*)y) \,dn
\right\>_{B\times G/N},
\end{split}
\]
where we use \eqref{eq-pairing1} to see that all integrals above
are well-defined.
Thus
\[
T^*(y)=\int_N \hat\epsilon_n((\phi\times G)(x^*)y) \,dn
\]
defines an adjoint $T^*$ for $T$. It follows that $T$ extends uniquely
to an adjointable linear map from $B\times_{\epsilon|} G/N$ 
to $Y_{G/N}^G(B)$, \ie, to a multiplier of the 
$(B\times_\epsilon G\times_{\hat\epsilon|,r} N)-(B\times_{\epsilon|} G/N)$ 
imprimitivity
bimodule $Y_{G/N}^G(B)$, by \cite[Proposition~1.3]{er:mult}.
We have shown that for all $x \in \c D(A)$ there exists $\Psi(x) \in
M(Y_{G/N}^G(B))$ such that
$\Psi(x) c=(\phi\times G)(x) c$ for 
$c \in \c D_N(B)$.

A computation similar to the derivation of the inequality
\eqref{bounded} shows that for $x,y \in \D(A)$ and $c,d \in \D_N(B)$
we have
\[
\<\Psi(x)c,\Psi(y)d\>_{B\times G/N}
=c^* (\phi\times G/N)(\< x,y\>_{A\times G/N}) d,
\]
so
\[
\<\Psi(x),\Psi(y)\>_{M(B\times G/N)}
=(\phi\times G/N)(\< x,y\>_{A\times G/N}).
\]

By \lemref{pre-imp-lem}, the following computation allows us to 
extend $\Psi$ uniquely to a nondegenerate imprimitivity bimodule
homomorphism, which we still denote by $\Psi$, from $Y_{G/N}^G(A)$ to
$M(Y_{G/N}^G(B))$, with coefficient maps $\phi\times G\times_r N$ and $\phi
\times G/N$: for $x,y\in\D(A)$ and $z\in\D(B)$ we have
\begin{align*}
&{}_{M(B\times G\times_rN)}\<\Psi(x),\Psi(y)\>z
=\Psi(x)\<\Psi(y),z\>_{B\times G/N}
=\Psi(x)\Psi(y)^*z
\\&\qquad=(\phi\times G)(x)\int_N
\hat\epsilon_n\bigl((\phi\times G)(y^*)z\bigr)\,dn
\\&\qquad=\int_N
(\phi\times G)(x)\hat\epsilon_n\bigl((\phi\times G)(y^*)\bigr)
\hat\epsilon_n(z)\,dn
\\&\qquad=\int_N
(\phi\times G)(x)(\phi\times G)\bigl(\hat\delta_n(y^*)\bigr)
\hat\epsilon_n(z)\,dn
\\&\qquad=\int_N
(\phi\times G)\bigl(x\hat\delta_n(y^*)\Delta(n)^{-1/2}\bigr)
\hat\epsilon_n(z)\Delta(n)^{1/2}\,dn
\\&\qquad=\int_N(\phi\times G)\bigl({}_{A\times G\times_rN}
\< x,y\>(n)\bigr)\hat\epsilon_n(z)\Delta(n)^{1/2}\,dn
\\&\qquad=\int_N(\phi\times G\times_rN)\bigl({}_{A\times G\times_rN}
\< x,y\>\bigr)(n)\hat\epsilon_n(z)\Delta(n)^{1/2}\,dn
\\&\qquad=(\phi\times G\times_rN)\bigl({}_{A\times G\times_rN}
\< x,y\>\bigr)z,
\end{align*}
so
\[
{}_{M(B\times G\times_rN)}\<\Psi(x),\Psi(y)\>
=(\phi\times G\times_rN)
\bigl({}_{A\times G\times_rN}\< x,y\>\bigr).
\]

We now check the coaction equivariance. For $x \in \D(A)$ we have
\[
\begin{split}
& (\Psi\otimes \id)\circ\delta_{Y(A)}(x)
=((\phi\times G)\otimes \id)\bigl((x\otimes 1)(j^A_G\otimes \id)
(w^*_G)\bigr)
\\
&
\qquad
=((\phi\times G)(x)\otimes 1)
\bigl((\phi\times G)\circ j^A_G\otimes \id\bigr) (w^*_G)
=((\phi\times G)(x)\otimes 1)(j^B_G\otimes \id)(w^*_G).
\end{split}
\]
Now, for any $c \in \D_N(B)$,
\[
\begin{split}
&\delta_{Y(B)}((\phi\times G)(x))\epsilon\dec(c)
=\delta_{Y(B)}((\phi\times G)(x)c)
=((\phi\times G)(x)c\otimes 1)(j^B_G\otimes \id)(w^*_G)
\\
&\qquad
=((\phi\times G)(x)\otimes 1)(c\otimes 1)(j^B_G\otimes \id)(w^*_G)
=((\phi\times G)(x)\otimes 1)(j^B_G\otimes \id)(w^*_G)\epsilon\dec(c),
\end{split}
\]
so
\[
\begin{split}
(\Psi\otimes \id)\circ\delta_{Y(A)}(x)
&
=\delta_{Y(B)}((\phi\times G)(x))
=\delta_{Y(B)}\circ \Psi(x).
\end{split}
\]
This gives the coaction equivariance.
For the action equivariance, if
$x \in \D(A)$ we have
\[
\begin{split}
\Psi\circ\alpha^{Y(A)}_s(x) &
=(\phi\times G)(\Delta_{G,N}(s)^{1/2} \hat\delta_s(x))
=\Delta_{G,N}(s)^{1/2} \hat\epsilon_s((\phi\times G)(x))
=\alpha^{Y(B)}_s\circ \Psi(x).
\end{split}
\]
As with Theorems~\ref{Ind-thm} and \ref{G-thm}, this suffices
to complete the proof.
\end{proof}

%
%

\chapter{Applications}
\label{apps-chap}

In this chapter we give some applications of the naturality theorems
from the preceding chapter to our motivating problem of
understanding the relationships between induction and duality.
First we uncover some new relationships
between Green and Mansfield induction. Important
special cases of these results say that the Green bimodules
$X_{\{e\}}^G(A)$ and Mansfield bimodules $Y_{G/G}^G(A)$ are in
duality:
\[
X_{\{e\}}^G(A)\times G\cong Y_{G/G}^G(A\times G)
\quad\text{and}\quad
X_{\{e\}}^G(A\times G)\cong Y_{G/G}^G(A)\times G.
\]
Results of this type require several applications of our main theorems,
and it is vital that we know everything is appropriately equivariant.
The same is true of our other applications to the restriction-induction
duality program of \cite{EchDI, kqr:resind, ekr}. We close with a new
application of linking algebra techniques to the symmetric
imprimitivity theorem of \cite{RaeIC}.

\section{Equivariant triangles}
\label{eqvt-tri-sec}

Our first application concerns a curious relationship between the
Green and Mansfield bimodules.
We prove two results (Theorems~\ref{eqvt-ekr-thm}
and~\ref{dual-ekr-thm}) which say,
very roughly, that Green and Mansfield imprimitivity are
``inverse to each other'', at least up to crossed product duality.

\subsection{Dual Mansfield equivariant triangle}

We begin with an action $(A,G,\alpha)$ and a closed normal subgroup
$N$ of $G$.  We consider various imprimitivity bimodules arising from
this data: first of all, there is the Green bimodule ${}_{A\times_r
G\times G/N}X_N^G(A)_{A\times_r N}$.  Temporarily 
forgetting about $N$, we also
have the Green bimodule 
${}_{A\times_r G\times G}X_e^G(A)_A$.
Recall
from \lemref{alpha X} that every Green bimodule carries an action
$\alpha^X$; in the particular case of $X_e^G$, the action $\alpha^X$
is $\hat{\hat\alpha}-\alpha$ compatible.  Restricting these actions
to $N$ and taking crossed products gives an
$(A\times_r G\times G\times_r N)-(A\times_r N)$ imprimitivity bimodule
$X_e^G(A)\times_r N$.
On the other hand,
the dual coaction $\hat\alpha$ of $G$ on $A\times_rG$
gives rise to the Mansfield bimodule
${}_{A\times_r G\times G\times_r N}Y_{G/N}^G(A\times_r G)_{A \times_r
G\times G/N}$.  The following theorem ties these bimodules together:

\begin{thm}
\label{eqvt-ekr-thm}
For any action $(A,G,\alpha)$ and any closed normal subgroup $N$ of
$G$, the diagram
\begin{equation}
\label{eqvt-ekr-diag}
\xymatrix
@C+2pc
{
{A\times_{\alpha,r}G\times_{\hat\alpha}G\times_{\hat{\hat\alpha},r}N}
\ar[r]^-{Y_{G/N}^G(A\times_rG)}
\ar[d]_{X_e^G(A)\times_rN}
&{A\times_{\alpha,r}G\times_{\hat\alpha|}G/N}
\ar[dl]^{X_N^G(A)}
\\
{A\times_{\alpha|,r}N}
}
\end{equation}
commutes equivariantly for the appropriate actions and coactions of $G$.
\end{thm}

\begin{proof}
First of all, recall that the dual coaction $\hat\alpha$ is
automatically normal and nondegenerate, so the Mansfield imprimitivity
bimodule $Y^G_{G/N}(A\times_r G)$ exists.  Commutativity as a diagram of
imprimitivity bimodules, but at the level of full crossed products
(and without the equivariance), is in \cite[Theorem 3.1]{ekr}.  We
will adapt the isomorphism of \cite{ekr} to the present context of
reduced crossed products.
\pause

Before we get into the details of this isomorphism, we should make
sure we know what the ``appropriate'' actions and coactions are:
on $X_N^G(A)$ we take the action $\alpha^X$ from \lemref{alpha X},
given on $x$ in the dense subspace $C_c(G,A)$ by
\[
\alpha^X_s(x)(t)=\Delta_{G,N}(s)^{1/2} x(ts).
\]
(Recall that $\Delta_{G,N}(s)=\Delta_G(s) \Delta_{G/N}(sN)^{-1}$.)
The coaction on $X_N^G(A)$ will be the $\delta_X$ from
\remref{delta X rem}, given on $x$ in $C_c(G,A)$ by
\[
\delta_X(x)(s)=x(s)\otimes s,
\]
where we regard
$\delta_X(x)$ as an element of
$C_c\bigl(G,M^\beta(A\otimes C^*(G))\bigr)$ 
$\subseteq M\bigl(X_N^G(A)\otimes C^*(G)\bigr)$.
The crossed product $X_e^G(A)\times_r N$ carries a decomposition
action, which for ease of writing in this proof we denote simply by
$\alpha^{X\times N}$.  Note that $C_c(N\times G,A)$ embeds in
$X_e^G(A)\times_r N$ via the chain of inclusions 
\[
C_c(N\times G,A)
\hookrightarrow C_c(N,C_c(G,A))
\hookrightarrow C_c(N,X)
\hookrightarrow X\times_r N,
\]
and $\alpha^{X\times N}$ is given on $z$ in the dense subspace 
$C_c(N \times G,A)$ by
\[
\alpha^{X\times N}_s(z)(n,t)
= \Delta_{G,N}(s) z(s^{-1}ns,ts).
\]
The coaction on $X_e^G(A)\times_r N$, which we shall simply denote
by $\delta_{X\times N}$, will be the inflation to $G$ of
the dual coaction of $N$.
However, before we can give a formula for this coaction which will
suit our purposes in this proof we must prepare some more tools
involving $C_c$-functions. 

\begin{lem}
The $(A\times_r G\times G\times_r N)-(A\times_r N)$ imprimitivity
bimodule $X^G_e\times_r N$ is the completion of the $C_c(N\times G
\times G,A)-C_c(N,A)$ pre-imprimitivity bimodule $C_c(N\times
G,A)$, with operations
\begin{align*}
f\d x(n,t)
&
=\int_N\int_G f(k,s,t)\alpha_s\bigl(x(k^{-1}n,s^{-1}tk)\bigr)
\Delta(s)^{1/2}\,ds\,dk
\\
x\d g(n,t)
&
=\int_N x(k,t)\alpha_{tk}(g(k^{-1}n))\,dk
\\
{}_{A\times_r G\times G\times_r N}\<x,y\>(n,s,t)
&
= \Delta(s)^{-1/2} \Delta_{G,N}(n)^{1/2}\int_N
x(nk,t)\alpha_s(y(k,s^{-1}tn)^*) \Delta(k)\,dk
\\
\<x,y\>_{A\times_r N}(n)
&
=\int_N\int_G\alpha_{k^{-1}s}
\bigl(x(k,s^{-1})^* y(kn,s^{-1})\bigr)\,ds\,dk,
\end{align*}
for $x,y\in C_c(N\times G,A)$, $f\in C_c(N\times G\times G,A)$,
and $g\in C_c(N,A)$.
\end{lem}

Just to be clear about how $C_c(N\times G\times G,A)$ is sitting
inside $A\times_r G\times G\times_r N$, note that $C_c(N\times G
\times G,A)$ embeds continuously in $C_c(N,C_c(G\times G,A))$ (for
the respective inductive limit topologies), hence in $A\times_r G
\times G\times_r N$ since $C_c(G\times G,A)$ embeds continuously in
$A\times_r G\times G$ via the isomorphism $A\times_r G\times G
\cong C_0(G,A)\times_r G$
from (the special case $N=\{e\}$ of)
\thmref{thm-isom-red}. Also recall that in this latter special case
we must remember that the second variable comes from
$C_0(G,A)$.

\begin{proof}
We just compute:
\begin{align*}
&
f\d x(n,t)
= f\d x(n,\cdot)(t)
\\
&\quad
= \left(\int_N f(k,\cdot,\cdot) \d
\alpha^X_k\bigl(x(k^{-1}n,\cdot)\bigr)\,dk \right)(t)
\\
&\quad
=\int_N \bigl(f(k,\cdot,\cdot) \d
\alpha^X_k(x(k^{-1}n,\cdot))\bigr)(t)\,dk
\\
&\quad
=\int_N\int_G f(k,s,t)
\alpha_s\bigl(\alpha^X_k(x(k^{-1}n,\cdot))(s^{-1}t)\bigr)
\Delta(s)^{1/2}\,ds\,dk
\\
&\quad
=\int_N\int_G f(k,s,t)
\alpha_s\bigl(x(k^{-1}n,s^{-1}tk)\bigr)
\Delta(s)^{1/2}\,ds\,dk,
\end{align*}
\begin{align*}
x\d g(n,t)
&
= x\d g(n)(t)
= \left(\int_N x(k,\cdot) \d\alpha_k(g(k^{-1}n))\,dk \right)(t)
\\
&
=\int_N \bigl(x(k,\cdot) \d
\alpha_k(g(k^{-1}n))\bigr)(t)\,dk
\\
&
=\int_N x(k,t)\alpha_t\bigl(\alpha_k(g(k^{-1}n))\bigr)\,dk
\\
&
=\int_N x(k,t)\alpha_{tk}(g(k^{-1}n))\,dk,
\end{align*}
\begin{align*}
&
{}_{A\times_r G\times G\times_r N}\<x,y\>(n,s,t)
= {}_{A\times_r G\times G\times_r N}\<x,y\>(n,\cdot,\cdot)(s,t)
\\
&\quad
= \left(\int_N
{}_{A\times_r G\times G}\bigl\< x(nk,\cdot),
\alpha^X_n(y(k,\cdot))\bigr\>
\Delta(k)\,dk \right)(s,t)
\\
&\quad
=\int_N
{}_{A\times_r G\times G}\bigl\< x(nk,\cdot),
\alpha^X_n(y(k,\cdot))\bigr\>(s,t)
\Delta(k)\,dk
\\
&\quad
=\int_N \Delta(s)^{-1/2} x(nk,t)
\alpha_s\bigl(\alpha^X_n(y(k,\cdot))(s^{-1}t)^*\bigr)
\Delta(k)\,dk
\\
&\quad
=\int_N \Delta(s)^{-1/2} x(nk,t)
\alpha_s\bigl(\Delta_{G,N}(n)^{1/2} y(k,s^{-1}tn)^*\bigr)
\Delta(k)\,dk
\\
&\quad
= \Delta(s)^{-1/2} \Delta_{G,N}(n)^{1/2}\int_N
x(nk,t)\alpha_s(y(k,s^{-1}tn)^*) \Delta(k)\,dk,
\end{align*}
and
\begin{align*}
\<x,y\>_{A\times_r N}(n)
&
=\int_N\alpha_{k^{-1}}\bigl(
\<x(k,\cdot),y(kn,\cdot)\>_A\bigr)\,dk
\\
&
=\int_N\alpha_{k^{-1}}\left(\int_G
\alpha_s\bigl(x(k,s^{-1})^* y(kn,s^{-1})\,ds \right)\,dk
\\
&
=\int_N\int_G
\alpha_{k^{-1}s}\bigl(x(k,s^{-1})^* y(kn,s^{-1})\bigr)\,ds\,dk.
\end{align*}
\end{proof}

\begin{lem}
With notation as above,
for $x\in C_c(N\times G,A)$ we have
\begin{gather*}
\delta_{X\times N}(x)\in
C_c\bigl(N\times G,M^\beta(A\otimes C^*(G))\bigr)
\subseteq M\bigl((X^G_e\times_r N)\otimes C^*(G)\bigr),
\end{gather*}
with $\delta_{X\times N}(x)(n,s)=x(n,s)\otimes n$ 
for $(n,s)\in N\times G$.
\end{lem}

\begin{proof}
Since $x\in C_c(N,X^G_e)$, by \propref{gamma hat} the dual coaction
of $N$ takes $x$ to the element of
$C_c\bigl(N,M^\beta(X^G_e\otimes C^*(N))\bigr)$
{given by} $n \mapsto x(n,\cdot)\otimes n$.
After inflating to a coaction of $G$, we get
\[
\delta_{X\times N}(x)\in C_c\bigl(N,M^\beta(X^G_e\otimes C^*(G))\bigr)
\midtext{and}
\delta_{X\times N}(x)(n)=x(n,\cdot)\otimes n.
\]
Further, since $x(n,\cdot)\in C_c(G,A)$, our usual canonical
embedding gives
\[
x(n,\cdot)\otimes n\in C_c\bigl(G,M^\beta(A\otimes C^*(N)\bigr)
\midtext{and}
\bigl(x(n,\cdot)\otimes n\bigr)(s)=x(n,s)\otimes n.
\]
Now, the map $(n,s) \mapsto x(n,s)\otimes n$ is in $C_c(N\times
G,M^\beta(A\otimes C^*(G)))$,
and we have just seen that
$\delta_{X\times N}(x)$ agrees with this map, so we are done.
\end{proof}

We must now do similar (and a little harder) work 
to obtain formulas for the action $\alpha^Y$ from \propref{alpha Y} and the
coaction $\delta_Y$ from \eqref{delta Y eq}
on the Mansfield
bimodule $Y_{G/N}^G(A\times_r G)$ at the level of $C_c$-functions.  
When $N$ is amenable,
\cite{ekr} shows that $C_c(G\times G,A)$ embeds in Mansfield's
bimodule $Y^G_{G/N}(A\times_r G)$.  The techniques there involve the
symmetric imprimitivity theorem of \cite{RaeIC}, but this can be
avoided:

\begin{lem}
The $\hat{\hat\alpha}$-invariant $*$-subalgebra $C_c(G\times G,A)$ of
$A\times_r G\times G$ is contained in $\c D(A\times_r G)$
\textup(see \appxref{imprim-chap}\textup),
and has algebraic operations
\begin{align*}
xy(s,t)
&
=\int_G x(r,t)\alpha_r(y(r^{-1}s,r^{-1}t))\,dr
\\
x^*(s,t)
&
=\alpha_s(x(s^{-1},s^{-1}t)^*) \Delta(s)^{-1}
\\
\hat{\hat\alpha}_r(x)(s,t)
&
= x(s,tr).
\end{align*}
\end{lem}

\begin{proof}
Let $x\in C_c(G\times G,A)$. Then $\supp x$ is
contained in $F\times F$ for some compact subset $F$ of $G$. Choose a
compact subset $E$ of $G$ whose interior contains $F$, and then choose
$u\in A_c(G)$ which is identically $1$ on $E$. We claim that $x$ is
$(u,E)$. Approximate $x$ in the norm of $A\times_r G\times G$ by a
finite sum
\[
\sum_k j_G(f_k)j_{A\times_r G}(i^r_A(a_k)i^r_G(c_k))
\quad\text{with }a_k\in A,f_k,c_k\in C_E(G).
\]
This sum is $(u,E)$ since $\supp f_k\subseteq E$ and
\[
i^r_A(a_k)i^r_G(c_k)=i^r_A(a_k)i^r_G(uc_k)
= \hat\alpha_u(i^r_A(a_k)i^r_G(c_k)).
\]
This shows that $C_c(G\times G,A) \subseteq \c D(A\times_r G)$.

The next two assertions 
follow quite easily from the isomorphism
$A\times_{\alpha,r} G\times_{\hat\alpha} G \cong C_0(G,A)\times_{\alpha
\otimes \tau,r} G$, which we freely abuse:
\begin{align*}
xy(s,t)
&
= xy(s)(t)
= \left(\int_G x(r)
(\alpha\otimes \tau)_r \bigl(y(r^{-1}s)\bigr)\,dr \right) (t)
\\
&
=\int_G x(r)(t)
(\alpha_r\otimes \tau_r) \bigl(y(r^{-1}s)\bigr)(t)\,dr
\\
&
=\int_G x(r,t)
\alpha_r \bigl(y(r^{-1}s)(r^{-1}t)\bigr)\,dr
\\
&
=\int_G x(r,t)
\alpha_r (y(r^{-1}s,r^{-1}t))\,dr,
\end{align*}
and
\begin{align*}
x^*(s,t)
&
= x^*(s)(t)
= (\alpha\otimes \tau)_s\bigl(x(s^{-1})^*\bigr)(t)
\Delta(s)^{-1}
\\
&
=\alpha_s\bigl(x(s^{-1})^*(s^{-1}t)\bigr)
\Delta(s)^{-1}
\\
&
=\alpha_s\bigl(x(s^{-1})(s^{-1}t)^*\bigr)
\Delta(s)^{-1}
\\
&
=\alpha_s(x(s^{-1},s^{-1}t)^*)
\Delta(s)^{-1}.
\end{align*}
The formula for $\hat{\hat\alpha}$ follows from
Equation~\eqref{alpha hat hat}.
\end{proof}

The above lemma and \propref{alpha Y} immediately give us the
following $C_c$-formula for the bimodule action $\alpha^Y$: for $x\in
C_c(G\times G,A) \subseteq \c D(A\times_r G)$ we have
\[
\alpha^Y_s(x)(t,r)=\Delta_{G,N}(s)^{1/2} x(t,rs).
\]

The following $C_c$-description of the Mansfield bimodule is essentially
\cite[Proposition 1.1]{ekr}, but we do it for reduced crossed products,
and we give a different argument.

\begin{lem}
For any action $(A,G,\alpha)$ and any closed normal subgroup $N$ of
$G$, the $(A\times_r G\times G\times_r N)-(A\times_r G\times G/N)$
imprimitivity bimodule $Y^G_{G/N}(A\times_r G)$ is the completion of
the $C_c(N\times G\times G,A)-C_c(G\times G/N,A)$
pre-imprimitivity bimodule $C_c(G\times G,A)$, with operations
\begin{align*}
f\d x(s,t)
&
=\int_N\int_G f(n,r,t)\alpha_r(x(r^{-1}s,r^{-1}tn))
\Delta(n)^{1/2}\,dr\,dn
\\
x\d g(s,t)
&
=\int_G x(r,t)\alpha_r(g(r^{-1}s,r^{-1}tN))\,dr
\\
{}_{A\times_r G\times G\times_r N}\<x,y\>(n,s,t)
&
= \Delta(n)^{-1/2}
\int_G x(r,t)\alpha_{s^{-1}}(y(sr,stn)^*) \Delta(r)\,dr
\\
\<x,y\>_{A\times_r G\times G/N}(s,tN)
&
=\int_N\int_G\alpha_{r^{-1}}\bigl(x(r,rtn)^*
y(rs,rtn)\bigr)\,dr\,dn,
\end{align*}
for $f\in C_c(N\times G\times G,A)$, $x,y\in C_c(G\times G,A)$, and
$g\in C_c(G\times G/N,A)$.
\end{lem}

\begin{proof}
Note first of all that $C_c(N\times G\times G,A)$, $C_c(G\times G,A)$,
and $C_c(G\times G/N,A)$ are dense in the respective normed spaces
$A\times_r G\times G\times_r N$, $Y^G_{G/N}(A\times_r G)$, and
$A\times_r G\times G/N$.
We have
\[
\begin{split}
f\d x(s,t)
&
=\int_N \bigl(f(n,\cdot,\cdot)\d n\d x\bigr)(s,t)\,dn
=\int_N \bigl(f(n,\cdot,\cdot) \hat{\hat\alpha}_n(x)
\Delta(n)^{1/2}\bigr)(s,t)\,dn
\\
&
=\int_N\int_G f(n,r,t)
\alpha_r\bigl(\hat{\hat\alpha}_n(x)(r^{-1}s,r^{-1}t)\bigr)
\Delta(n)^{1/2}\,dr\,dn
\\
&
=\int_N\int_G f(n,r,t)\alpha_r(x(r^{-1}s,r^{-1}tn))
\Delta(n)^{1/2}\,dr\,dn,
\end{split}
\]
\[
\begin{split}
&
{}_{A\times_r G\times G\times_r N}\<x,y\>(n,s,t)
= {}_{A\times_r G\times G\times_r N}\<x,y\>
(n,\cdot,\cdot)(s,t)
\\
&\quad
= \Delta(n)^{-1/2} (x \hat{\hat\alpha}_n(y^*))(s,t)
= \Delta(n)^{-1/2}\int_G x(r,t)
\alpha_r\bigl(\hat{\hat\alpha}_n(y^*) (r^{-1}s,r^{-1}t)\bigr)\,dr
\\
&\quad
= \Delta(n)^{-1/2}\int_G x(r,t)\alpha_r(y^*(r^{-1}s,r^{-1}tn))\,dr
\\
&\quad
= \Delta(n)^{-1/2}\int_G x(r,t)
\alpha_r\bigl(\alpha_{r^{-1}s}(y(s^{-1}r,s^{-1}rr^{-1}tn)^*)\bigr)
\Delta(r^{-1}s)^{-1}\,dr
\\
&\quad
= \Delta(n)^{-1/2}
\int_G x(r,t)\alpha_{s^{-1}}(y(sr,stn)^*) \Delta(r)\,dr,
\end{split}
\]
and
\[
\begin{split}
&
\<x,y\>_{A\times_r G\times G/N}(s,tN)
=\int_N \hat{\hat\alpha}_n(x^*y)(s,t)\,dn
=\int_N x^*y(s,tn)\,dn
\\
&\quad
=\int_N\int_G x^*(r,tn)
\alpha_r(y(r^{-1}s,r^{-1}tn))\,dr\,dn
\\
&\quad
=\int_N\int_G\alpha_r\bigl(x(r^{-1},r^{-1}tn)^*
y(r^{-1}s,r^{-1}tn)\bigr) \Delta(r)^{-1}\,dr\,dn
\\
&\quad
=\int_N\int_G\alpha_{r^{-1}}\bigl(x(r,rtn)^*
y(rs,rtn)\bigr)\,dr\,dn.
\end{split}
\]
The formula for $x\d g$ follows from membership of $g$ in $M(A
\times_r G\times G)$, since $C_c(G\times G/N,A)$ embeds continuously
in $A\times_r G\times G/N$, hence also continuously in $M^\beta(A
\times_r G\times G)$.
\end{proof}

The following lemma prepares us to deal with the coaction
$\delta_Y$ in terms of $C_c$-functions; as usual, the trick is to allow
the extra copy of $C^*(G)$ to be a freely moving object.

\begin{lem}
\label{prep-lem}
For any action $(A,G,\alpha)$, any closed normal subgroup $N$ of $G$,
and any $C^*$-algebra $C$, the canonical embedding of
$C_c(G\times G,A) \odot C$ in $C_c(G\times G,A\otimes C)$ extends
to an imprimitivity bimodule isomorphism
\begin{multline*}
\Psi \:
{}_{(A\times_r G\times G\times_r N)\otimes C}\bigl(
Y^G_{G/N}(A\times_{\alpha,r} G)\otimes C
\bigr)_{(A\times_r G\times G/N)\otimes C}
\\
\xrightarrow{\cong}
{}_{(A\otimes C)\times_r G\times G\times_r N}\bigl(
Y^G_{G/N}((A\otimes C)\times_{\alpha\otimes \id,r} G)
\bigr)_{(A\otimes C)\times_r G\times G/N}.
\end{multline*}
\end{lem}

\begin{proof}
We first must exhibit isomorphisms
\begin{align*}
\phi_L\:(A\times_r G\times G\times_r N)\otimes C
&
\iso
(A\otimes C)\times_r G\times G\times_r N
\righttext{and}
\\
\phi_R\:(A\times_r G\times G/N)\otimes C
&
\iso
(A\otimes C)\times_r G\times G/N
\end{align*}
between the left and right  coefficient algebras.  This mainly
involves a few applications of \lemref{lem-mintensor} and
\thmref{thm-isom-red}: for the right  coefficients, we have
\begin{align*}
(A\times_r G\times G/N)\otimes C
&
\cong \bigl((A\otimes C_0(G/N))\times_r G\bigr)\otimes C
\\
&
\cong \bigl((A\otimes C_0(G/N)\otimes C)\times_r G\bigr)
\\
&
\cong \bigl((A\otimes C\otimes C_0(G/N))\times_r G\bigr)
\\
&
\cong \bigl((A\otimes C)\times_r G\bigr)\times G/N,
\end{align*}
giving the isomorphism $\phi_R$. Note that $\phi_R$ is equivariant for
the actions 
$(\infl\what{\hat\alpha|})\otimes \id$ and 
$\infl (\what{\alpha\otimes \id}|)\sphat$.
For the left coefficients, we use \lemref{lem-mintensor} once
more, together with what we have already shown about $\phi_R$:
\begin{align*}
(A\times_r G\times G\times_r N)\otimes C
\cong \bigl((A\times_r G\times G)\otimes C\bigr)\times_r N
\cong \bigl((A\otimes C)\times_r G\times G\bigr)\times_r N,
\end{align*}
which gives $\phi_L$.

We are now ready to prove the required bimodule isomorphism.  By
\lemref{pre-imp-lem} and \remref{rem-isoimp}, it suffices to show
that $\Psi$ respects both inner products. For $y,z \in C_c(G\times
G,A)$ and $c,d\in C$ we have
\begin{align*}
&{}_L\<\Psi(y\otimes c),\Psi(z\otimes d)\>(n,s,t)
\\&\quad=\Delta(n)^{-1/2}\int_G
\Psi(y\otimes c)(r,t)
(\alpha\otimes\id)_{s^{-1}}\bigl(\Psi(z\otimes d)(sr,stn)^*\bigr)
\Delta(r)\,dr
\\&\quad=\Delta(n)^{-1/2}\int_G
\bigl(y(r,t)\otimes c\bigr)
\bigl(\alpha_{s^{-1}}\otimes\id\bigr)
\bigl(z(sr,stn)^*\otimes d^*\bigr)
\Delta(r)\,dr
\\&\quad=\Delta(n)^{-1/2}\int_G
y(r,t)\alpha_{s^{-1}}(z(sr,stn)^*)\Delta(r)\,dr\otimes cd^*
\\&\quad={}_L\<y,z\>(n,s,t)\otimes cd^*
=\phi_L\bigl({}_L\<y,z\>\otimes cd^*\bigr)(n,s,t)
=\phi_L\bigl({}_L\<y\otimes c,z\otimes d\>\bigr)(n,s,t)
\end{align*}
and
\begin{align*}
&\<\Psi(y\otimes c),\Psi(z\otimes d)\>_R(s,tN)
\\
&\quad=\int_N\int_G (\alpha\otimes \id)_{r^{-1}}
\bigl(\Psi(y\otimes c)(r,rtn)^*
\Psi(z\otimes d)(rs,rtn)\bigr)\,dr\,dn
\\
&\quad=\int_N\int_G (\alpha_{r^{-1}}\otimes \id)
\bigl((y(r,rtn)^*\otimes c^*) (z(rs,rtn)\otimes d)\bigr)\,dr\,dn
\\
&\quad=\int_N\int_G\alpha_{r^{-1}}\bigl(y(r,rtn)^* z(rs,rtn)\bigr)
\,dr\,dn\otimes c^*d
\\
&\quad= \<y,z\>_R(s,tN)\otimes c^*d
= \phi_R\bigl(\<y,z\>_R\otimes c^*d\bigr)(s,tN)
= \phi_R\bigl(\<y\otimes c,z\otimes d\>_R\bigr)(s,tN).
\end{align*}
\end{proof}

We use \lemref{prep-lem} to embed $C_c(G\times G,A\otimes C)$
in $Y^G_{G/N}((A\times_r G)\otimes C)$ and $C_c(G\times
G,M^\beta(A\otimes C))$ in $M(Y^G_{G/N}((A\times_r
G)\otimes C))$.
Recall from \eqref{delta Y eq} that the decomposition
coaction $\delta_Y$ of $G$ on $Y_{G/N}^G$
from \cite{er:stab} is given by
\[
\delta_Y(y)=(y\otimes 1)(j_G\otimes \id)(w^*_G)
\]
for $y\in\D$, and in particular for $y\in C_c(G\times G,A)$, and the
multiplication takes place in $M((A\times_r G\times G)\otimes C^*(G))$.

\begin{lem}
With notation as above,
for $y\in C_c(G\times G,A)$ we have
\begin{gather*}
\delta_Y(y)\in C_c(G\times G,M^\beta(A\otimes C^*(G))),
\end{gather*}
with $\delta_Y(y)(s,t)=y(s,t)\otimes t^{-1}s$ for $s,t\in G$.
\end{lem}

\begin{proof}
We first claim that for $f\in C_0(G)$ we have $y j_G(f)\in C_c(G
\times G,A)$ and
\[
(yj_G(f))(s,t)=y(s,t) f(s^{-1}t).
\]
Because of the isomorphism from \thmref{thm-isom-red}, we can
work in $C_0(G,A)\times_{\alpha\otimes \tau,r} G$ instead of $A
\times_{\alpha,r} G\times_{\hat\alpha} G$.
First, for $s\in G$ we have
\begin{align*}
\bigl(y i^r_{C_0(G,A)}(1\otimes f)\bigr)(s)
= y(s,\cdot) (\alpha\otimes \tau)_s(1\otimes f)
= y(s,\cdot) (1\otimes \tau_s(f)),
\end{align*}
and evaluating this at $t\in G$ we get
\[
\bigl(y i^r_{C_0(G,A)}(1\otimes f)\bigr)(s)(t)
= y(s,t) \tau_s(f)(t)
= y(s,t) f(s^{-1}t).
\]
This is jointly continuous in $(s,t)$ and has compact support
(because $y$ does), and this verifies the claim.

Hence, if we further have $c\in C^*(G)$ then
\[
(y\otimes 1) (j_G\otimes \id) (f\otimes c)
= y j_G(f)\otimes c\in C_c(G\times G,A\otimes C^*(G))
\]
and
\[
\begin{split}
&
\bigl((y\otimes 1) (j_G\otimes \id) (f\otimes c)\bigr)(s,t)
= (y j_G(f)\otimes c)(s,t)
= (y j_G(f))(s,t)\otimes c
\\
&\quad
= y(s,t) f(s^{-1}t)\otimes c
= y(s,t)\otimes f(s^{-1}t) c
= y(s,t)\otimes (f\otimes c)(s^{-1}t).
\end{split}
\]
The lemma follows by density and continuity, 
since $w_G\in M(C_0(G)\otimes C^*(G))$.
\end{proof}

\resume{\thmref{eqvt-ekr-thm}}
The isomorphism
$\Phi$ of $Y\otimes_{A\times_r G\times G/N} X$ onto 
$X\times_r N$
from \cite{ekr} is given on elementary tensors $y\otimes x\in
C_c(G\times G,A) \odot C_c(G,A)$ by
\[
\Phi(y\otimes x)(n,t)
= \Delta(n)^{-1/2}\int_G y(s,t)\alpha_s(x(s^{-1}tn))
\Delta(s)^{1/2}\,ds.
\]
We show that this is equivariant for the appropriate actions: for $y
\in C_c(G\times G,A)$ and $x\in X_N^G(A)$ we have
\begin{align*}
&
\Phi \circ \bigl(\alpha^Y\otimes\alpha^X\bigr)_r
(y\otimes x)(n,t)
\\
&\quad
= \Delta(n)^{-1/2}\int_G
\alpha^Y(y)(s,t)\alpha_s\bigl(\alpha^X_r(x)(s^{-1}tn)\bigr)
\Delta(s)^{1/2}\,ds
\\
&\quad
= \Delta(n)^{-1/2}\int_G
\Delta_{G,N}(r)^{1/2} y(s,tr)
\alpha_s\bigl(\Delta_{G,N}(r)^{1/2} x(s^{-1}tnr)\bigr)
\Delta(s)^{1/2}\,ds
\\
&\quad
= \Delta(n)^{-1/2} \Delta_{G,N}(r)\int_G
y(s,tr)\alpha_s\bigl(x(s^{-1}tnr)\bigr)
\Delta(s)^{1/2}\,ds
\\
&\quad
= \Delta_{G,N}(r)
\Phi(y\otimes x)(r^{-1}nr,tr)
\\
&\quad
=\alpha^{X\times N}_r \circ \Phi(y\otimes x)(n,t).
\end{align*}

We turn to the coaction equivariance, and for this we must show that
the diagram
\[
\xymatrix
@C+4pc
{
{Y\otimes_{A\times_r G\times G/N} X}
\ar[r]^-{\delta_Y\cotimes_{A\times_r G\times G/N} \delta_X}
\ar[d]_{\Phi}
&{M\bigl((Y\otimes_{A\otimes_r G\times G/N} X)
\otimes C^*(G)\bigr)}
\ar[d]^{\Phi\otimes \id}
\\
{X\times_r N}
\ar[r]_-{\delta_{X\times N}}
&{M\bigl((X\times_r N)
\otimes C^*(G)\bigr)}
}
\]
commutes. Recall from \propref{tnsr-co-prop} that
$\delta_Y\cotimes_{A\times_r G\times G/N} \delta_X
= \Theta \circ (\delta_Y\otimes_{A\times_r G\times G/N} \delta_X)$,
where
\[
\Theta\:\bigl(Y\otimes C^*(G)\bigr)
\otimes_{(A\times_r G\times G/N)\otimes C^*(G)}
\bigl(X\otimes C^*(G)\bigr)
\xrightarrow{\cong}
(Y\otimes_{A\times_r G\times G/N} X)\otimes C^*(G)
\]
is the isomorphism defined by
\[
\Theta\bigl((y\otimes c)\otimes (x\otimes d)\bigr)
= (y\otimes x)\otimes cd
\righttext{for} y\in Y,\ x\in X,\ c,d\in C^*(G).
\]
For $y\in C_c(G\times G,A)$ and $x\in C_c(G,A)$ we have
\[
\begin{split}
&
(\Phi\otimes \id) \circ
(\delta_Y\cotimes_{A\times_r G\times G/N} \delta_X)(y\otimes x)(n,t)
= (\Phi\otimes \id) \circ \Theta(\delta_Y(y)\otimes \delta_X(x))(n,t)
\\
&\quad
\overset{(*)}{=} \Delta(n)^{-1/2}\int_G \delta_Y(y)(s,t)
(\alpha_s\otimes \id)\bigl(\delta_X(x)(s^{-1}tn)\bigr)
\Delta(s)^{1/2}\,ds
\\
&\quad
= \Delta(n)^{-1/2}\int_G \bigl(y(s,t)\otimes t^{-1}s\bigr)
(\alpha_s\otimes \id)\bigl(x(s^{-1}tn)\otimes s^{-1}tn\bigr)
\Delta(s)^{1/2}\,ds
\\
&\quad
= \Delta(n)^{-1/2}\int_G y(s,t)\alpha_s(x(s^{-1}tn))
\Delta(s)^{1/2}\,ds\otimes n
\\
&\quad
= \Phi(y\otimes x)(n,t)\otimes n
= \delta_{X\times N} \circ \Phi(y\otimes x)(n,t),
\end{split}
\]
where the equality at $(*)$ is justified by replacing $\delta_Y(y)$ and
$\delta_X(x)$ by elementary tensors $u\otimes c\in C_c(G\times G,A)
\odot C^*(G)$ and $v\otimes d\in C_c(G,A) \odot C^*(G)$,
respectively, and appealing to density and continuity:
\[
\begin{split}
&(\Phi\otimes \id) \circ \Theta
\bigl((u\otimes c)\otimes (v\otimes d)\bigr)(n,t)
= (\Phi\otimes \id)\bigl((u\otimes v)\otimes cd\bigr)(n,t)
\\
&\quad
= \bigl(\Phi(u\otimes v)\otimes cd\bigr)(n,t)
= \Phi(u\otimes v)(n,t)\otimes cd
\\
&\quad
= \Delta(n)^{-1/2}\int_G u(s,t)\alpha_s(v(s^{-1}tn))
\Delta(s)^{1/2}\,ds\otimes cd
\\
&\quad
= \Delta(n)^{-1/2}\int_G (u(s,t)\otimes c)
(\alpha_s\otimes \id)\bigl(v(s^{-1}tn)\otimes d\bigr)
\Delta(s)^{1/2}\,ds
\\
&\quad
= \Delta(n)^{-1/2}\int_G (u\otimes c)(s,t)
(\alpha_s\otimes \id)\bigl((v\otimes d)(s^{-1}tn)\bigr)
\Delta(s)^{1/2}\,ds.
\end{split}
\]

\end{proof}

As an immediate and interesting corollary, taking $N=G$ we have:

\begin{cor}
\label{eqvt-ekr-cor}
Let $(A,G,\alpha)$ be an action. Then
\[
X_e^G(A)\times_{\alpha^X,r}G
\cong Y_{G/G}^G(A\times_{\alpha,r}G)
\]
as $(A\times_{\alpha,r} G\times_{\hat\alpha} G\times_{\hat{\hat\alpha},r} G)
-(A\times_{\alpha,r} G)$ imprimitivity
bimodules, equivariantly for the appropriate actions and coactions.
\end{cor}
Note that in this case the appropriate coactions are 
$(\alpha^X)\hat{}$ and $(\hat\alpha)_Y$, 
and the appropriate
actions are the trivial ones.

\subsection{Dual Green equivariant triangle}

The following result is dual to \thmref{eqvt-ekr-thm}, starting with
a coaction rather than an action.

\begin{thm}
\label{dual-ekr-thm}
For any nondegenerate normal coaction $(A,G,\delta)$
and any closed normal subgroup $N$ of $G$,
the diagram
\begin{equation}
\label{dual-ekr-diag}
\xymatrix
@L+5pt
{
&{A\times_\delta G\times_{\hat\delta,r}G\times_{\deltahathat|}G/N}
\ar[dl]_<(.6){X_N^G(A\times G)}
\ar[d]^{Y_{G/G}^G(A)\times G/N}
\\
{A\times_{\delta} G\times_{\deltahat|,r}N}
\ar[r]_-{Y_{G/N}^G(A)}
&{A\times_{\delta|} G/N}
}
\end{equation}
commutes equivariantly for the appropriate actions and coactions.
\end{thm}

\begin{proof}
First apply \thmref{eqvt-ekr-thm} to the dual action
$(A\times G,G,\hat{\delta})$ with $N=G$. This gives an
$(A\times G\times_rG\times G\times_rG)-(A\times G\times_rG)$
imprimitivity bimodule isomorphism
\[
X_e^G(A\times G)\times_rG\cong Y_{G/G}^G(A\times G\times_rG)
\]
which is equivariant for the coactions
$\what{{\hat\delta}^X}$ and $(\hat{\hat{\delta}})_Y$ (and in this
case the appropriate actions are the trivial ones). Also,
\thmref{M-thm} applied to the $A\times G\times_rG-A$ imprimitivity
bimodule $Y_{G/G}^G(A)$ with subgroup $G$ gives an imprimitivity
bimodule isomorphism
\[
Y_{G/G}^G(A\times G\times_rG)\cong Y_{G/G}^G(A)\times G\times_rG
\]
which is equivariant for the dual coactions
$(\deltahathat)_Y$ and $(\delta_Y)\hat{\hat{}}$.
Combining these, we have
\begin{equation}
\label{xpr iso}
X_e^G(A\times G)\times_rG\cong Y_{G/G}^G(A)\times G\times_rG,
\end{equation}
equivariantly for $\what{{\hat\delta}^X}$ and $(\delta_Y)\hat{\hat{}}$.
We wish to conclude that
\begin{equation}
\label{no-xpr iso}
X_e^G(A\times G)\cong Y_{G/G}^G(A)\times G,
\end{equation}
equivariantly for the actions $\hat\delta^X$ and $(\delta_Y)\hat{}$.
For this, we need
the following ``duality'' result for coactions on right-Hilbert
bimodules, which is an easy corollary of 
functoriality of the crossed products and \thmref{M-thm}.
\pause

\begin{lem}
\label{fwd-back-coact}
Suppose $Z$ and $W$ are right-Hilbert $A-B$ bimodules, $\delta$ and
$\epsilon$ are nondegenerate normal coactions of $G$ on $A$ and $B$,
and $\zeta$ and $\eta$ are $\delta-\epsilon$ compatible coactions of
$G$ on $Z$ and $W$, respectively. Then
\[
(Z,\zeta)\cong (W,\eta)
\midtext{if and only if}
(Z\times_\zeta G,\hat{\zeta}) \cong (W\times_\eta G,\hat{\eta}).
\]
\end{lem}

\begin{proof}
The forward implication follows from \thmref{co-xpr-fun}, 
so assume $(Z\times G,\what{\zeta})\cong (W\times G,\what\eta)$. 
Then \thmref{xpr-fun} gives
\begin{equation}
\label{double-dual-isom}
\Bigl(Z\times G\times_rG,\what{\what\zeta} \Bigr) \cong
\Bigl(W\times G\times_rG,\what{\what\eta} \Bigr).
\end{equation}
Now \thmref{M-thm} with subgroup $N=G$ gives diagrams
\[
\xymatrix
@C+1pc
{
{A\times G\times_rG}
\ar[r]^-{Y_{G/G}^G(A)}
\ar[d]_{Z\times G\times_rG}
&{A}
\ar[d]^{Z}
\\
{B\times G\times_rG}
\ar[r]_-{Y_{G/G}^G(B)}
&{B}
}
\midtext{and}
\xymatrix
@C+1pc
{
{A\times G\times_rG}
\ar[r]^-{Y_{G/G}^G(A)}
\ar[d]_{W\times G\times_rG}
&{A}
\ar[d]^{W}
\\
{B\times G\times_rG}
\ar[r]_-{Y_{G/G}^G(B)}
&{B}
}
\]
which commute equivariantly for the various coactions. Piecing these
together using the equivariant isomorphism \eqref{double-dual-isom}\
and canceling the equivalences $Y_{G/G}^G(A)$ and $Y_{G/G}^G(B)$ gives
the proposition.
\end{proof}

\begin{rem}
Although we won't need it, we remark that the result dual to
\propref{fwd-back-coact},
starting with actions rather than coactions, can be proved analogously
by using \thmref{G-thm} in place of \thmref{M-thm}.
\end{rem}

\resume{\thmref{dual-ekr-thm}}
As indicated immediately before the statement of
\lemref{fwd-back-coact}, the crossed product
in the isomorphism
\eqref{xpr iso} can be removed to give the isomorphism \eqref{no-xpr
iso}.
Taking crossed products by $N$, we
have
\begin{equation}
\label{times-N-isom}
X_e^G(A\times G)\times_rN\cong Y_{G/G}^G(A)\times G\times_rN,
\end{equation}
and by functoriality (\thmref{combined-fun})
this isomorphism is equivariant for the appropriate actions and
coactions.

Next, apply \thmref{eqvt-ekr-thm} to the dual action
$(A\times G,G,\hat{\delta})$ with subgroup $N$ to get a commutative
diagram
\[
\xymatrix
@C+4pc
@L+5pt
{
{A\times G\times_rG\times G\times_r N}
\ar[r]^-{Y_{G/N}^G(A\times G\times_rG)}
\ar[d]_{X_e^G(A\times G)\times_rN}
&{A\times G\times_rG\times G/N}
\ar[dl]^<(.3){X_N^G(A\times G)}
\\
{A\times G\times_rN}
}
\]
which on using \eqref{times-N-isom} gives equivariant commutativity
for actions and coactions of the upper left  triangle of the
diagram
\begin{equation}
\label{dual-ekr-pf-diag}
\xymatrix
@C+4pc
@L+3pt
{
{A\times G\times_rG\times G\times_rN}
\ar[r]^-{Y_{G/N}^G(A\times G\times_rG)}
\ar[d]_{Y_{G/G}^G(A)\times G\times_rN}
&{A\times G\times_rG\times G/N}
\ar[dl]_<(.6){X_N^G(A\times G)}
\ar[d]^{Y_{G/G}^G(A)\times G/N}
\\
{A\times G\times_rN}
\ar[r]_-{Y_{G/N}^G(A)}
&{A\times G/N}.
}
\end{equation}
Now note that the outer square in \eqref{dual-ekr-pf-diag} commutes
equivariantly for actions and coactions by applying \thmref{M-thm} to
the imprimitivity bimodule $Y_{G/G}^G(A)$; we conclude that the lower
right  triangle of \eqref{dual-ekr-pf-diag} commutes equivariantly
for the appropriate actions and coactions, which proves the theorem.
\end{proof}

An interesting corollary of the proof is dual to
\corref{eqvt-ekr-cor}:

\begin{cor}
\label{dual-ekr-cor}
Let $(A,G,\delta)$ be a nondegenerate normal coaction. Then
\[
X_e^G(A\times_\delta G)\cong Y_{G/G}^G(A)\times_{\delta_Y} G
\]
as $(A\times_\delta G\times_{\hat\delta,r} G\times_{\deltahathat} G)
-(A\times_\delta G)$ imprimitivity bimodules,
equivariantly for the appropriate actions and coactions.
\end{cor}
\noindent
In this case the appropriate actions are
$(\hat\delta)^X$ and $(\delta_Y)\hat{}$, and the appropriate
coactions are the trivial ones.

\section{Restriction and induction}
\label{res-ind-sec}

Various combinations of the authors \cite{EchDI},
\cite{ekr}, \cite{kqr:resind} have established a duality between
restriction and induction for actions and coactions.  A little more
precisely, given an action or coaction, restricting representations in
the given system is dual to inducing representations in the dual
system, and similarly with restricting and inducing reversed.  These
``dualities'' are actually expressed as commutative diagrams in the
category $\C$, and for convenience we refer to them as ``Ind-Res''
and ``Res-Ind'' diagrams, respectively.  Here we will
apply the results of the current paper to give simplified proofs of
the Ind-Res duality.  The strategy is to deduce the Res-Ind diagram
for the given system from the Ind-Res diagram 
for the dual system.

The Ind-Res diagram for actions,  
which essentially amounts to induction in stages,
is largely due to Green:

\begin{prop}[\cite{gre:local}]
\label{act-ind-res}
For any action $(A,G,\alpha)$
and any closed normal
subgroups $N \subseteq M$ of $G$, the diagram
\[
\xymatrix{
{A\times_{\alpha,r} G\times_{\what\alpha|}G/M}
\ar[r]^-{X_M^G}
\ar[d]_{A\times_{\alpha,r} G\times_{\what\alpha|}G/N}
&{A\times_{\alpha|,r} M}
\ar[d]^{X_N^M}
\\
{A\times_{\alpha,r} G\times_{\what\alpha|}G/N}
\ar[r]_-{X_N^G}
&{A\times_{\alpha|,r} N}
}
\]
commutes in the category $\C$.
\end{prop}

\begin{proof}
We must show
\[
X_M^G\otimes_{A\times M}X_N^M\cong
(A\times_r G\times G/N)\otimes_{A\times_r G\times G/N}X_N^G
\]
as right-Hilbert $(A\times_r G\times G/M)-(A\times_r N)$ bimodules.
Momentarily leaving off the second crossed products by $G/M$ and
$G/N$, this becomes
\begin{equation}
\label{act-stages}
X_M^G\otimes_{A\times M}X_N^M\cong
(A\times_r G)\otimes_{A\times_r G}X_N^G
\end{equation}
as right-Hilbert $(A\times_r G)-(A\times_r N)$ bimodules. Since
\[
(A\times_r G)\otimes_{A\times_r G}X_N^G\cong X_N^G
\]
by cancellation, the isomorphism \eqref{act-stages} is just Green's
induction in stages for actions \cite[Proposition 8]{gre:local}, and is
given on the generators by
\[
\Phi(x\otimes y)=x\d y
\righttext{for} x\in C_c(G,A),y\in C_c(M,A),
\]
where $y$ is viewed as an element of $A\times_r M$, acting on the
right of $X_M^G$. 
We just need to check that $\Phi$ preserves the left action of
$C_0(G/M)$: for $f\in C_0(G/M)$ we have
\begin{align*}
\Phi(f\d(x\otimes y))
= \Phi(f\d x\otimes y)
= (f\d x)\d y
= f\d (x\d y)
= f\d\Phi(x\otimes y).
\end{align*}
\end{proof}

We now recover the Res-Ind diagram for coactions \cite[Theorem
3.1]{kqr:resind}\footnote{Actually, the theorem we state here is slightly
weaker than in \cite{kqr:resind}, since there we get to assume only
as much ``normality'' as necessary for the bits of the diagram, namely
``Mansfield imprimitivity works for $M$''.  Here we need to require
$\delta$ itself to be normal, because of our method of proof.
Similarly for \propref{coact-ind-res} below.}:

\begin{cor}[\cite{kqr:resind}]
\label{coact-res-ind}
For any nondegenerate normal
coaction $(B,G,\delta)$ and any closed normal
subgroups $N \subseteq M$ of $G$, the diagram
\[
\xymatrix
@C+1pc
{
{B\times_\delta G\times_{\what\delta|,r} M}
\ar[r]^-{Y_{G/M}^G(B)}
\ar[d]_{X_N^M(B\times_\delta G)}
&{B\times_{\delta|} G/M}
\ar[d]^{B\times_{\delta|} G/N}
\\
{B\times_\delta G\times_{\what\delta|,r} N}
\ar[r]_-{Y_{G/N}^G(B)}
&{B\times_{\delta|} G/N}
}
\]
commutes in the category $\C$.
\end{cor}

\begin{proof}
The desired diagram is the outer rectangle of the diagram
\[
\xymatrix{
{B\times G\times_r M}
\ar[rr]^-{Y_{G/M}^G(B)}
\ar[ddd]_{X_N^M(B\times G)}
&&{B\times G/M}
\ar[ddd]^{B\times G/N}
\\
&{B\times G\times_r G\times G/M}
\ar[ul]|{X_M^G(B\times G)}
\ar[ur]|{Y_{G/G}^G(B)\times G/M}
\ar[d]_{B\times G\times_r G\times G/N}
\\
&{B\times G\times_r G\times G/N}
\ar[dl]|{X_N^G(B\times G)}
\ar[dr]|{Y_{G/G}^G(B)\times G/N}
\\
{B\times G\times_r N}
\ar[rr]_-{Y_{G/N}^G(B)}
&&{B\times G/N.}
}
\]
Since the slanted arrows are isomorphisms, it suffices to show the
inner polygons commute.

The left quadrilateral commutes by the preceding proposition.
The right quadrilateral commutes by
\cite[left face of Diagram (5.1)]{kqr:resind}. Finally, the top and
bottom triangles commute by \thmref{dual-ekr-thm}.
\end{proof}

Dually, suppose $(B,G,\delta)$ is a coaction. The following
Ind-Res diagram is \cite[Theorem 4.1]{kqr:resind}:

\begin{prop}[\cite{kqr:resind}]
\label{coact-ind-res}
For any nondegenerate normal
coaction $(B,G,\delta)$ and any closed normal
subgroups $N \subseteq M$ of $G$, the diagram
\[
\xymatrix{
{B\times_\delta G\times_{\what\delta|,r}N}
\ar[r]^-{Y_{G/N}^G}
\ar[d]_{B\times_\delta G\times_{\what\delta|,r}M}
&{B\times_{\delta|} G/N}
\ar[d]^{Y_{G/N}^{G/M}}
\\
{B\times_\delta G\times_{\what\delta|,r}M}
\ar[r]_-{Y_{G/M}^G}
&{B\times_{\delta|} G/M}
}
\]
commutes in the category $\C$.
\end{prop}

The next theorem substantially improves Corollary~3.3 of \cite{ekr}, 
which covers the case $N=\{e\}$.

\begin{thm}
\label{act-res-ind}
For any action $(A,G,\alpha)$
and any closed normal
subgroups $N \subseteq M$ of $G$, the diagram
\[
\xymatrix{
{A\times_{\alpha,r} G\times_{\hat\alpha|}G/N}
\ar[r]^-{X_{N}^G(A)}
\ar[d]_{Y_{G/M}^{G/N}(A\times_{\alpha,r} G)}
&{A\times_{\alpha|,r} N}
\ar[d]^{A\times_{\alpha|,r} M}
\\
{A\times_{\alpha,r} G\times_{\hat\alpha|} G/M}
\ar[r]_-{X_{M}^G(A)}
&{A\times_{\alpha|,r} M}
}
\]
commutes in the category $\C$.
\end{thm}

\begin{proof}
The proof is patterned after that of \corref{coact-res-ind}.  The
desired diagram is the outer rectangle of the diagram
\[
\xymatrix{
{A\times_r G\times G/N}
\ar[rr]^-{X_{N}^G(A)}
\ar[ddd]_{Y_{G/M}^{G/N}(A\times_r G)}
&&{A\times_r N}
\ar[ddd]^{A\times_r M}
\\
&{A\times_r G\times G\times_r N}
\ar[ul]^{Y_{G/N}^G(A\times_r G)}
\ar[ur]_{X_{e}^G(A)\times_r N}
\ar[d]_{A\times_r G\times G\times_r M}
\\
&{A\times_r G\times G\times_r M}
\ar[dl]^{Y_{G/M}^G(A\times_r G)}
\ar[dr]_{X_{e}^G(A)\times_r M}
\\
{A\times_r G\times G/M}
\ar[rr]_-{Y_{G/N}^G(A)}
&&{A\times_r M.}
}
\]
The left quadrilateral commutes by the preceding proposition.  The
right quadrilateral commutes by \cite[Lemma 5.7]{kqr:resind} (which gives
the analogous result for full crossed products), together with the
argument at the end of the proof of \cite[Theorem~5.6]{kqr:resind}, which
shows that the kernels of the regular representations match up.
Finally, the top and bottom triangles commute by
\thmref{eqvt-ekr-thm}.
\end{proof}

\section{Symmetric imprimitivity}
\label{symimp-sec}

In this section we show how to deduce the induced algebra results of
\cite[Section 4]{RaeIC} from the main symmetric imprimitivity theorem
of that paper, thus avoiding the need to repeat the arguments of
\cite[Section 1]{RaeIC}.  Suppose, therefore, that we are in the
setting of \cite[Section 4]{RaeIC}
(see also Section~1 of Appendix~\ref{imprim-chap}):
we have a locally compact space
${}_HP_K$ with commuting free and proper actions of locally compact
groups $H$ and $K$, a Morita equivalence ${}_DY_E$, commuting actions
$(\delta,\sigma)$ of $(H,K)$ on $D$ and $(\tau,\gamma)$ of $(H,K)$ on
$E$, and compatible actions $(\nu,\mu)$ of $(H,K)$ on ${}_DY_E$.  
We then have commuting actions
\[
L(\tau)\deq
\left(\begin{smallmatrix}
\delta&\nu
\\
\tilde\nu&\tau
\end{smallmatrix}\right):\H\to\aut L(Y)
\midtext{and}
L(\sigma)\deq
\left(\begin{smallmatrix}
\sigma&\mu
\\
\tilde\mu&\gamma
\end{smallmatrix}\right)\:K\to\aut L(Y)
\]
on the linking algebra \cite[Section 4]{com}.  The induced algebras
$\ind_H^P L(Y)$ and $\ind_K^P L(Y)$ carry actions $L(\alpha)$ of $K$
and $L(\beta)$ of $H$, respectively, and it follows from \cite[Theorem
1.1]{RaeIC} that $C_c(P,L(Y))$ can be given the structure of an
$\big((\ind_H^P L(Y))\times_{L(\alpha)}K\big)-\big((\ind_K^P
L(Y))\times_{L(\beta)}H\big)$ imprimitivity bimodule.

If we set
\[
\ind_H^P Y= \left\{
f\:P\to Y \left|
\begin{array}{l}
f(hp)=L(\tau)_h(f(p)) \text{ for } h\in H, \text{ and}
\\
Hp \mapsto \|\<f(p),f(p)\>_D\| \text{ vanishes at }\infty
\end{array}
\right. \right\},
\]
then $\ind_H^P L(Y)$ is naturally isomorphic to
\[
\left(\begin{matrix}
\ind_H^P D&\ind_H^P Y
\\
\ind_H^P\widetilde Y&\ind_H^P E
\end{matrix}\right)
= L(\ind_H^P Y).
\]
Under this identification, the action $L(\alpha)$ restricts in the
corners to the tensor product actions of $K$ on $\ind_H^P D$ and
$\ind_H^P E$, and hence the top left  corner in the decomposition
\[
\big(\ind_H^P L(Y)\big)\times_{L(\alpha)}K=
\left(\begin{matrix}
(\ind_H^P D)\times K&(\ind_H^P Y)\times K
\\
(\ind_H^P\widetilde Y)\times K&(\ind_H^P E)\times K
\end{matrix}\right)
\]
is the crossed product of $\ind_H^P D$ by the action $\alpha$ of
\cite[Section 4]{RaeIC}: with our conventions, $\alpha_k(f)(p)\deq
\sigma_k(f(pk))$.  Since the bottom right corner in the analogous
decomposition of $(\ind_K^P L(Y))\times_{L(\beta)} H$ is $(\ind_K^P
E)\times_\beta H$, we deduce immediately from the second part of
\lemref{link-lem} that the upper right corner $C_c(P,Y)$ in
$C_c(P,L(Y))$ completes to give an $\big((\ind_H^P D)\times_\alpha
K\big)-\big((\ind_K^P E)\times_\beta H\big)$ imprimitivity bimodule.
It is a straightforward matter to check that the module actions and
inner products are, modulo our change in conventions, the ones
described in \cite[page 384]{RaeIC}.  For example, for $f,g\in
C_c(P,Y)$, $\<f,g\>_{\ind E\times H}$ is the bottom right corner in
\[
\left\<\left(\begin{smallmatrix}
0&f
\\
0&0
\end{smallmatrix}\right),
\left(\begin{smallmatrix}
0&g
\\
0&0
\end{smallmatrix}\right)\right\>_{C_c(H,\ind L(Y))}=
\left(\begin{smallmatrix}
0&0
\\
0&r
\end{smallmatrix}\right),
\]
where $r$ is the function
\[
r(h,p)=\Delta(h)^{-1/2}\int_K\gamma_k\big(\<
f(h^{-1}p),\tau_h(g(h^{-1}pk))\>_E\big)\,dk.
\]
Thus we have shown how to deduce \cite[Theorem 4.1]{RaeIC} from
\cite[Theorem 1.1]{RaeIC}.

In retrospect, it was always relatively easy to deduce the existence
of the Morita equivalence in \cite[Theorem 4.1]{RaeIC} from
\cite[Theorem 1.1]{RaeIC}, by composing the Morita equivalence of
$\ind_H^P D\times K$ and $\ind_K^P D\times H$ given by \cite[Theorem
1.1]{RaeIC} with the equivalence of $\ind_K^P D\times H$ and $\ind_K^P
E\times H$ induced by ${}_DY_E$.  The third part of \lemref{link-lem}
says that the tensor product bimodule thus obtained is isomorphic to
the one we have just constructed.  Thus the main point of
\cite[Section 4]{RaeIC} is the specific nature of the bimodule.

\begin{rem}
\lemref{link-act-lem} and \lemref{link-coact-lem} would, in the case
$P=G$ studied in \cite{er:induced}, give an equivariant version of
\cite[Section 4]{RaeIC}.
\end{rem}

\appendix

%
%

\appendix
\chapter
{Crossed Products by Actions and Coactions}
\label{coactions-chap}

In this appendix we give an introduction to the theory of crossed products
by actions and coactions of groups on $C^*$-algebras. Since the theory
of coactions is much newer, and since there are at least three different
definitions of coactions of groups on $C^*$-algebras in the literature
(all somehow mixed together) we decided to present an almost 
self-contained exposition of coactions and their crossed products, including
all proofs for the basic constructions.  None of the results in this
appendix are new and the main sources are \cite{it, kat, lprs, qr:induced,
qui:full, qui:fullred, rae:full, rae:rep}.

Some general notation: If $X$ is a locally compact space and $E$ is a
normed vector space, then $C_c(X,E)$, $C_0(X,E)$, and
$C_b(X,E)$ denote the spaces of continuous $E$-valued 
compactly supported functions, functions which vanish at infinity,
and bounded functions, respectively.  
If $E=\bbC$ we simply write $C_c(X), C_0(X)$, and
$C_b(X)$, respectively. $M(B)$ denotes the multiplier algebra of a
$C^*$-algebra $B$ and $UM(B)$ denotes the group of unitary elements of
$M(B)$.  The \emph{strict topology} on $M(B)$ is the locally convex
topology generated by the seminorms $m\mapsto \|ma\|, \|am\|$, $a\in
A$.  Note that $M(A)$ is the strict completion of $A$, and we write
$M^{\beta}(A)$ if we consider $M(A)$ equipped with the strict
topology. For example, $C_c(X, M^{\beta}(A))$ 
will denote the strictly continuous functions
of $X$ into $M(A)$ with compact support.
A homomorphism $\phi\:A\to M(B)$ of a $C^*$-algebra $A$ into
$M(B)$ is called \emph{nondegenerate} if $\phi(A)B=B$.  We use the same
letter for a nondegenerate homomorphism $\phi\:A\to M(B)$ and its
unique (strictly continuous)
extension $M(A)\to M(B)$ \cite[Lemma 1.1]{lprs}.
Finally, if $\H$ is a Hilbert space, we shall always assume that
the inner product on $\H$ is conjugate linear in the first and
linear in the second variable.

\section{Tensor products}

Tensor products of $C^*$-algebras play a basic r\^ole in the theory of
crossed products by actions and coactions: in a certain sense, a crossed
product by a group action of $G$ on a $C^*$-algebra $A$
is just a skew tensor product of $A$ with $C^*(G)$, and a crossed
product by a coaction of $G$ on $A$ can be viewed as
a skew tensor product of $A$ with $C_0(G)$.

If $E$ and $F$ are complex vector spaces, then we denote by $E\odot F$
the \emph{algebraic tensor product} of $E$ and $F$.  If $A$ and $B$
are $C^*$-algebras, then $A\odot B$ becomes a $*$-algebra in the
canonical way and a $C^*$-cross norm on $A\odot B$ is a norm
$\|\cdot\|_{\nu}$ which satisfies $\|a\otimes b\|_{\nu}=\|a\|\|b\|$
for all elementary tensors $a\otimes b\in A\odot B$ and such that the
completion of $A\otimes_{\nu} B=\overline{A\odot B}^{\|\cdot\|_{\nu}}$
is a $C^*$-algebra.  $A\otimes_{\nu}B$ is then called the
\emph{$\nu$-tensor product} of $A$ and $B$.  We denote by
$k_A^{\nu}\:A\to M(A\otimes_{\nu}B)$ and $k_B^{\nu}\:B\to
M(A\otimes_{\nu}B)$ the canonical maps $k_A^{\nu}(a)=a\otimes 1$,
$k_B^{\nu}(b)=1\otimes b$.

Among the (possibly) many $C^*$-cross norms on $A\odot B$ there is a
maximal one and a minimal one.  The maximal norm $\|\cdot\|_{\max}$ is
characterized by the universal property that, whenever we have two
nondegenerate homomorphisms $\phi\:A\to M(D)$ and $\psi\:B\to M(D)$
with commuting ranges (\ie, $\phi(a)\psi(b)=\psi(b)\phi(a)$
for all $a\in A$, $b\in B$), then there exists a nondegenerate
homomorphism $\phi\otimes\psi\:A\otimes_{\max}B\to M(D)$
satisfying $(\phi\otimes\psi)\circ k_A^{\max}=\phi$ and
$(\phi\otimes\psi)\circ k_B^{\max}=\psi$.  Thus
\[
\biggl\|\sum_{i=1}^na_i\otimes b_i\biggr\|_{\max}=
\sup\biggl\{\biggl\|\sum_{i=1}^n\phi(a_i)\psi(b_i)\biggr\|\biggr\},
\]
where the supremum is taken over all commuting pairs of nondegenerate
homomorphisms of $A$ and $B$.  In particular, if $\|\cdot
\|_{\nu}$ is any other $C^*$-cross norm, then $k_A^{\nu}\otimes
k_B^{\nu}\:A\otimes_{\max}B\to A\otimes_{\nu}B$ is a surjection (since
it is the identity on $A\odot B$), and hence $\|\cdot\|_{\nu}$ is
dominated by $\|\cdot\|_{\max}$.

If $\pi\:A\to \B(\H)$ and $\rho\:B\to \B(\K)$ are
representations of $A$ and $B$ on the Hilbert spaces $\H$ and $\K$,
respectively, then there exists a representation
$\pi\otimes\rho\:A\odot B\to \B(\H\otimes \K)$ satisfying
$(\pi\otimes\rho)(a\otimes b)=\pi(a)\otimes\rho(b)$, and $\pi\otimes
\rho$ is faithful on $A\odot B$ if $\pi$ and $\rho$ are faithful.  The
\emph{minimal norm} on $A\odot B$ is defined by
\[
\biggl\|\sum_{i=1}^na_i\otimes b_i\biggr\|_{\min}=
\sup\biggl\{\biggl\|\sum_{i=1}^n\pi(a_i)
\otimes\rho(b_i)\biggr\|\biggr\},
\]
where the supremum is taken over all representations
$\pi$ and $\rho$ of $A$ and $B$, respectively.  It is a nontrivial
fact that $\|\cdot\|_{\min}$ is indeed smaller than any other
$C^*$-cross norm on $A\odot B$ (however, it is clear from the
definition that $\|\cdot\|_{\min}\leq \|\cdot\|_{\max}$).  We shall
simply write $\|\cdot\|$ for $\|\cdot \|_{\min}$ and $A\otimes B$ for
the \emph{minimal tensor product} $A\otimes_{\min} B$.  Note that it
follows from the minimality of $\|\cdot\|_{\min}$ that, whenever
$\pi\:A\to \B(\H)$ and $\rho\:B\to \B(\K)$ are faithful
representations of $A$ and $B$, respectively, then
$\pi\otimes\rho\:A\otimes B\to \B(\H\otimes\K)$ is a faithful
representation of $A\otimes B$.  Note also that
$\pi\otimes \rho$ is nondegenerate if and only if $\pi$ and $\rho$
are nondegenerate.
The following result will be used
frequently.

\begin{lem}
\label{lem-reptensor}
Let $\phi\:A\to M(C)$ and $\psi\:B\to M(D)$ be
homomorphisms.  Then there is a
homomorphism $\phi\otimes\psi\:A\otimes B\to M(C\otimes D)$
satisfying $(\phi\otimes\psi)(a\otimes b)=\phi(a)\otimes\psi(b)$.  If
$\phi$ and $\psi$ are nondegenerate \textup(resp. faithful\textup),
then so is $\phi\otimes\psi$.
\end{lem}

\begin{proof} Representing $C$ and $D$ faithfully on Hilbert spaces
turns $\phi$ and $\psi$ into
$*$-representations, and the result follows from the above-mentioned
properties of the minimal tensor product.
\end{proof}

\begin{rem}
Recall that a $C^*$-algebra $A$ is \emph{nuclear} if
$\|\cdot\|_{\max}=\|\cdot\|_{\min}$ on $A\odot B$ for any
$C^*$-algebra $B$, \ie, $A\otimes_{\max}B=A\otimes B$ for all $B$.
Basic examples of nuclear $C^*$-algebras are the commutative
$C^*$-algebras and the algebras $\K(\H)$ of compact operators on a
Hilbert space $\H$, but there are many others.
\end{rem}

In this paper, we often need to work with a certain
subalgebra $M_C(A\otimes C)$ of the multiplier algebra
$M(A\otimes C)$ of the minimal tensor product $A\otimes
C$.

\begin{defn}\label{defn-Cmult}
     Suppose that $A$ and $C$ are $C^*$-algebras.
     Then we define the {\em $C$-multiplier algebra} $M_C(A\otimes C)$
     of $A\otimes C$ as the set
     $$M_C(A\otimes C)=\{m\in M(A\otimes C)\mid m(1\otimes C)
     \cup (1\otimes C)m\subseteq
     A\otimes C\}.$$
     The {\em $C$-strict topology} on $M_C(A\otimes C)$ is 
     the locally convex topology generated
     by the seminorms $m\mapsto \|m(1\otimes c)\|,\|(1\otimes c)m\|$,
     $c\in C$.
\end{defn}

\begin{rem}\label{rem-Cmult}
     (1) It is straightforward to check that $M_C(A\otimes C)$
     is a closed $*$-subalgebra of $M(A\otimes C)$. We decided
     to call it the $C$-multiplier algebra, since it somehow
     consists of elements which are ``real'' multipliers
     only in the $C$-factor.
     Of course, there is a similar definition of  the
     $A$-multiplier algebra of $A\otimes C$.

     (2) If $C=C_0(V)$ for some locally compact space $V$,
     and if we identify $A\otimes C_0(V)$ with $C_0(V,A)$
     in the canonical way, then one can check, using
     the identification of $M(A\otimes C_0(V))$ with
     $C_b(V, M^{\beta}(A))$
     (see \cite[Corollary 3.4]{APT-MC}), that
     $M_{C_0(V)}(X\otimes C_0(V))$ can be identified
     with $C_b(V,A)$.

     (3) It is clear that a similar object can also be defined
     for the maximal tensor product $A\otimes_{\max}C$.
     However, since we only use minimal tensor products
     in the main body of this work, we stick to this case.
     Note that the $C$-multiplier algebra has appeared in several
     places in the literature in connection with coactions
     of groups and Hopf algebras (see, \eg, \cite{lprs,BS-CH},
     where it is denoted by $\widetilde{M}(A\otimes C)$).
\end{rem}

In what follows we gather some important
properties of the $C$-multiplier algebra.

\begin{prop}\label{prop-Cstrict}
     Let $A$ and $C$ be $C^*$-algebras.
     \begin{enumerate}
     \item The $C$-strict topology on $M_C(A\otimes C)$
     is stronger than the strict topology induced
     from $M(A\otimes C)$, and multiplication is separately
     $C$-strictly continuous on $M_C(A\otimes C)$.
     Also, the involution on $M_C(A\otimes C)$ is $C$-strictly
     continuous.
     \item $M_C(A\otimes C)$ is the $C$-strict completion of
     $A\otimes C$.
     \item We have $(1\otimes M(C))M_C(A\otimes C)\cup
     M_C(A\otimes C)(1\otimes M(C))\subseteq M_C(A\otimes C)$.
     \end{enumerate}
\end{prop}
\begin{proof} Let $(m_i)_{i\in I}$ be a net in $M_C(A\otimes C)$
     which converges $C$-strictly to $m$. If $z\in A\otimes C$
     we can factor $z=(1\otimes c)y$ for some $c\in C$, $y\in
     A\otimes C$, to conclude that
     $$m_iz=m_i(1\otimes c)y\to m(1\otimes c)y=mz,$$
     where convergence is in norm. A similar argument shows that
     $zm_i\to zm$.
     Separate $C$-strict continuity of multiplication and
     $C$-strict continuity of involution follows then from
     continuity with respect to the strict topology. Hence (i).

     For the proof of (ii),
     let $m\in M_C(A\otimes C)$ and let
     $(c_i)_{i\in I}$ be a bounded approximate unit of $C$.
     We claim that $(m(1\otimes c_i))_{i\in I}\subseteq A\otimes C$
     converges $C$-strictly to $m$.
     In fact, if $c\in C$, then $c_ic\to c$ in norm,
     which implies that $m(1\otimes c_i)(1\otimes c)\to m(1\otimes c)$
     in norm. On the other hand, one easily checks that
     $z(1\otimes c_i)\to z$ in norm for all $z\in A\otimes C$,
     from which it follows that $(1\otimes c)m(1\otimes c_i)
     \to (1\otimes c)m$ in norm.

     Suppose now that $(m_i)_{i\in I}$ is a
     $C$-strict Cauchy net in $M_C(A\otimes C)$. By (i) it follows
     that $(m_i)_{i\in I}$ is also a strict Cauchy net
     in $M(A\otimes C)$. Since $M(A\otimes C)$ is the strict
     completion of $A\otimes C$ we can find an $m\in M(A\otimes C)$
     such that $m_i\to m$ strictly.
     We claim that $m\in M_C(A\otimes C)$ and $m_i\to m$ $C$-strictly.
     For this we let $\eps>0$ and $c\in C$, and choose
     $i_0\in I$ such that $\|(m_i-m_j)(1\otimes c)\|\leq\eps$
     for all $i,j\geq i_0$. Since
     $m_j(1\otimes c)z\to m(1\otimes c)z$ in norm, it follows that
     $$\|(m_i-m)(1\otimes c)z\|=\lim_j\|(m_i-m_j)(1\otimes c)z\|\leq
     \eps$$
     for all $i\geq i_0$. Thus $m_i(1\otimes c)\to m(1\otimes c)$ in
     norm,
     and a similar argument shows that $(1\otimes c)m_i\to (1\otimes
     c)m$ in norm for all $c\in C$. Thus $m_i\to m$ $C$-strictly.
     Finally, since $m_i(1\otimes c), (1\otimes c)m_i\in A\otimes C$,
     it follows from the fact that $A\otimes C$ is norm closed
     in $M(A\otimes C)$ that $(1\otimes c)m, m(1\otimes c)\in A\otimes C$
     for all $c\in C$. Thus $m\in M_C(X\otimes C)$.
     This proves (ii).

     We omit the straightforward proof of (iii).
\end{proof}

\begin{prop}\label{prop-Cmultextend}
     Suppose that $A$, $B$, $C$, and $D$ are $C^*$-algebras.
     Let $\phi\:A\to M(B)$ be a possibly degenerate $*$-homomorphism
     and let $\psi\:C\to M(D)$ be a nondegenerate $*$-homomorphism.
     \begin{enumerate}
     \item There exists a unique
     $*$-homomorphism
     $$\overline{\phi\otimes\psi}\:M_C(A\otimes C)\to M(B\otimes D)$$
     which extends the homomorphism $\phi\otimes \psi\:A\otimes C
     \to M(B\otimes D)$ of \lemref{lem-reptensor}
     \item The homomorphism $\overline{\phi\otimes\psi}$
     of~\textup{(i)} is $C$-strict to strict continuous.
     \item If $\phi$ and $\psi$ are faithful, then so is
     $\overline{\phi\otimes\psi}$.
     \item If $\phi(A)\subseteq B$, then
     $\overline{\phi\otimes\psi}(M_C(A\otimes C))\subseteq M_D(B\otimes
     D)$ and $\overline{\phi\otimes\psi}$ is $C$-strict to $D$-strict
     continuous.
     \end{enumerate}
\end{prop}
\begin{proof}
     We show that $\phi\otimes \psi\:A\otimes C\to M(B\otimes D)$
     is $C$-strict to strict continuous.
     If $(a_i)_{i\in I}$ is a net in $A\otimes C$
     which converges $C$-strictly to $a\in A\otimes C$,
     and if $z\in B\otimes D$, we can factor $z=(1\otimes\psi(c))y$
     for some $c\in C$
     and $y\in B\otimes D$ (since $\psi$ is nondegenerate) to conclude that
     $$\phi\otimes\psi(a_i)z=\phi\otimes \psi(a_i(1\otimes c))y
     \to \phi\otimes\psi(a(1\otimes c))y=\phi\otimes \psi(a)z,$$
     where convergence is in norm. A similar argument shows that
     $z(\phi\otimes \psi(a_i))\to z(\phi\otimes \psi(a))$ in norm
     for all $z\in B\otimes D$.

     It follows that there exists a unique $C$-strict to strict
     continuous linear extension
     $\overline{\phi\otimes\psi}\:M_C(A\otimes C)\to M(B\otimes D)$
     which is automatically a $*$-homomorphism
     by the $C$-strict (resp. strict) continuity of involution
     and the separate $C$-strict (resp. strict) continuity
     of multiplication in both algebras.

    Assume now that
     ${\eta}\:M_C(A\otimes C)\to M(B\otimes D)$
     is another extension of $\phi\otimes\psi$ and let $m\in
     M_C(A\otimes C)$. If $z\in B\otimes D$, we can factor
     $z=(1\otimes\psi(c))y$ for some $c\in C$ and $y\in B\otimes D$
     to compute
  $$
{\eta}(m)z=\phi\otimes\psi(m(1\otimes c))y
=\overline{\phi\otimes\psi}(m)z,
$$
which implies ${\eta}=\overline{\phi\otimes\psi}$.
This finishes the proof of (i) and (ii).

If $\phi$ and $\psi$ are faithful, then so is
$\phi\otimes \psi\:A\otimes C\to M(B\otimes D)$ by
\lemref{lem-reptensor}. Thus, if $\overline{\phi\otimes\psi}(m)=0$
for $m\in M_C(A\otimes C)$, it follows that
$\phi\otimes \psi(m(A\otimes C))=\{0\}$, which
implies that $m(A\otimes C)=\{0\}$. Hence $m=0$. This proves (iii).

Finally, if $\phi(A)\subseteq B$, it is first clear that
$\phi\otimes\psi(A\otimes C)\subseteq M_D(B\otimes D)$
since $\phi\otimes\psi(a\otimes c)(1\otimes d)=
\phi(a)\otimes\psi(c)d\in B\otimes D$ for all
elementary tensors $a\otimes c\in A\odot  C$ and $d\in D$.
But this implies that $\phi\otimes\psi(A\otimes C)(1\otimes D)
\subseteq B\otimes D$, and applying the $*$-operation then
gives $(1\otimes D)\big(\phi\otimes\psi(A\otimes C)\big)\subseteq
B\otimes D$.
Factoring $D=\phi(C)D$, a similar argument as in the proof
of (i) then shows that
$\phi\otimes\psi\:A\otimes C\to M_D(B\otimes D)$ is $C$-strict to
$D$-strict continuous, which implies that
there exists a $C$-strict to $D$-strict continuous extension
$${\eta}\:M_C(A\otimes C)\to M_D(B\otimes D),$$
which is a $*$-homomorphism by $C$-strict (resp. $D$-strict)
continuity of the algebra operations.
But then the uniqueness clause of (i) implies that
${\eta}=\overline{\phi\otimes\psi}$.
\end{proof}

Some further information on $C$-multiplier algebras will be given
in \secref{sec-Cmult} of \chapref{hilbert-chap}.

\section{Actions and their crossed products}
\label{sec-actions}

If $G$ is a locally compact group, then $ds$ denotes
left Haar measure on $G$ and $\Delta$ its modular function,
\ie, we have $\int_G f(s)\, ds=\Delta(t)\int_G f(st)\, ds$ for all
$f\in C_c(G), t\in G$.
An \emph{action} of
$G$ on  a $C^*$-algebra $A$ is a strongly continuous
homomorphism  $\alpha$  of $G$ into the $*$-automorphism group
$\Aut(A)$ of $A$.
We also call the triple $(A,G,\alpha)$ an action.
If $(A,G,\alpha)$ is an action, then
$C_c(G,A)$ becomes a $*$-algebra with respect to the
convolution and involution defined by
\begin{equation}
f*g(s)=\int_G f(t)\alpha_t(g(t^{-1}s))\,dt
\midtext{and}
f^*(s)=\Delta(s^{-1})\alpha_s(f(s^{-1}))^*.
\end{equation}

A \emph{covariant homomorphism} of $(A,G,\alpha)$
into the multiplier algebra $M(D)$ of a $C^*$-algebra $D$
is a pair $(\pi, U)$, where
$\pi\:A\to M(D)$ is a nondegenerate homomorphism
and $U\:G\to UM(D)$ is a strictly continuous homomorphism satisfying
\[
\pi\circ\alpha_s=\ad U_s\circ\pi
\righttext{for all}s\in G,
\]
where for any unitary $v$ we let $\ad v$ denote the usual conjugation
automorphism $\ad v(b)=vbv^*$.
A \emph{covariant representation} of $(A,G,\alpha)$ on a
Hilbert space $\H$ is a covariant homomorphism into
$M(\K(\H))=\B(\H)$.  If $(\pi, U)$ is a covariant homomorphism into
$M(D)$, then the \emph{integrated form} $\pi\times U\:C_c(G,A)\to
M(D)$ is defined by
\begin{equation}
(\pi\times U)(f)=\int_G \pi(f(s))U_s\,ds.
\end{equation}
With the help of \cite[Lemma 7]{rae:rep} it is not hard to see that
$\pi\times U$ is well-defined and continuous with respect
to the inductive limit
topology on $C_c(G,A)$. Also, for each covariant representation
$(\pi,U)$ and each $f\in C_c(G,A)$ we have
\[
\|(\pi\times U)(f)\|\leq \int_G \|f(s)\|\, ds,
\]
from which it follows that for each $f\in C_c(G,A)$ we can
form
\[
\|f\|\deq \sup\{\|(\pi\times U)(f)\|\mid \text{$(\pi,U)$ is a
covariant representation
of $(A,G,\alpha)$}\}.
\]
One can check that $\|\cdot\|$ is a norm on $C_c(G,A)$. This leads to

\begin{defn}
\label{def-crossedaction}
Let $(A,G,\alpha)$ be an action. Then the completion
$A\times_{\alpha}G$ of $C_c(G,A)$ with respect to $\|\cdot\|$
is called the
(\emph{full}) \emph{crossed product} of
$(A,G,\alpha)$.
\end{defn}

\begin{rem}
\label{rem-universal}
(1) There is a canonical covariant homomorphism
$(i_A,i_G)$ of $(A,G,\alpha)$ into $M(A\times_{\alpha}G)$ given
by the formulas
\begin{align*}
(i_A(a) f)(s) &= a f(s)
&
(i_G(t) f)(s) &= \alpha_t(f(t^{-1}s))
\\
(f i_A(a))(s) &= f(s) \alpha_s(a)
&
(f i_G(t))(s)&= \Delta(t^{-1})f(st^{-1}),
\end{align*}
$f\in C_c(G,A)$. It follows that if
$f\in C_c(G,A)$, then
$(i_A\times i_G)(f)=f$ as elements of $A\times_{\alpha}G$.

(2)
The
triple $(A\times_{\alpha}G, i_A, i_G)$ enjoys the following
universal property: For any covariant homomorphism
  $(\pi, U)$ of $(A,G,\alpha)$ into $M(D)$ there exists a
unique nondegenerate homomorphism
$\pi\times U\:A\times_{\alpha}G\to M(D)$ such that
$(\pi\times U)\circ i_A=\pi$ and $(\pi\times U)\circ i_G =U$.

To see this observe that it follows from the definition of the
greatest $C^*$-norm on $C_c(G,A)$ that the integrated form $\pi\times
U\:C_c(G,A)\to M(D)$ is continuous with respect to this norm and
therefore extends uniquely to a homomorphism of
$A\times_{\alpha}G$.  Using an approximate identity of
$A\times_{\alpha}G$ which lies inside $C_c(G,A)$, it is not too hard
to check that the extension $\pi\times U\:A\times_{\alpha}G\to M(D)$
satisfies all properties mentioned above.%
\footnote{One can construct such an approximate identity as follows:
for each compact neighborhood $V$ of the identity $e\in G$ choose
$g_V\in C_c(G)^+$ such that $\supp g_V\subseteq V$ and $\int_G
g_V(s)\,ds=1$.  Then, if $\{a_i\}_i$ is an approximate identity of $A$,
$\{i_A(a_i) i_G(g_V)\}_{(V,i)}$ with ordering $(V,i)\leq
(W,j)\Leftrightarrow V\supseteq W$ and $i\leq j$ becomes an approximate
identity of $A\times_{\alpha}G$.}
Conversely, if $\rho\:A\times_{\alpha}G\to M(D)$ is a nondegenerate
homomorphism, then $(\rho\circ i_A, \rho\circ i_G)$ is a covariant
homomorphism of $(A,G,\alpha)$ such that
$\rho=(\rho\circ i_A)\times (\rho\circ i_G)$. Thus,
$(\pi,U)\mapsto \pi\times U$ gives a one-to-one correspondence between
the covariant homomorphisms of $(A,G,\alpha)$ and the
nondegenerate
homomorphisms
of $A\times_{\alpha}G$.

(3)
In \cite{rae:full} the universal properties
were used to \emph{define} the full crossed product  of $(A,G,\alpha)$
as a triple $(C, k_A, k_G)$ satisfying
\begin{enumerate}
\item
$(k_A,k_G)$ is a covariant homomorphism  of $(A,G,\alpha)$ into $M(C)$;

\item
$k_A(A)k_G(C^*(G))$ is dense in $C$;

\item
if $(\pi,U)$ is any covariant homomorphism of $(A,G,\alpha)$
into $M(D)$, then there exists a unique nondegenerate
homomorphism
$\phi_{\pi,U}\:C\to M(D)$ such that
$\phi_{\pi,U}\circ k_A=\pi$ and $\phi_{\pi,U}\circ k_G=U$.
\end{enumerate}
Clearly, $(A\times_{\alpha}G, i_A, i_G)$ is a crossed product in this
sense.  Moreover, if $(C, k_A, k_G)$ is any other triple satisfying
(i)--(iii), then $k_A\times k_G\:A\times_{\alpha}G\to C$
is an isomorphism with inverse $\phi_{i_A,i_G} \:C\to
A\times_{\alpha}G$, where $\phi_{i_A,i_G}$ is the homomorphism
associated to $(i_A,i_G)$ by (iii).

(4)
It is sometimes useful to be able to consider integrated forms
of pairs $(\pi, U)$ with $\pi\:A\to M(D)$ a $*$-homomorphism
and $U\:G\to UM(D)$ a strictly continuous homomorphism which satisfy
the covariance condition $\pi(\alpha_s(a))=U_s\pi(a)U_s^*$
for all $a\in A$, $s\in G$, but where $\pi\:A\to M(D)$ is
degenerate. We shall call such a pair a {\em degenerate covariant
homomorphism}.
As for nondegenerate homomorphisms we get a
$*$-homomorphism $\pi\times U\:C_c(G,A)\to M(D)$ by integration.
So the only problem is to see that $\pi\times U$ is
norm-decreasing, in order to obtain a unique
extension to $A\times_{\alpha}G$.
For this  represent $M(D)$ faithfully into $B(\H)$ for some
Hilbert space $\H$ and write $\H_1\deq \pi(A)\H$.
Then we get $U_s\H_1=U_s\pi(A)\H=\pi(A)U_s\H=\pi(A)\H=\H_1$,
from which we obtain a nondegenerate representation
$(\pi_1, U_1)$ of $(A,G,\alpha)$ into $\B(\H_1)$.
The desired result then follows from
$\|\pi\times U(f)\|=\|\pi_1\times U_1(f)\|$
for all $f\in C_c(G,A)$.
\end{rem}

The following (well-known) lemma serves as a first example
of the usefulness of the universal properties
of the full crossed product. It also indicates the conceptual similarity
of full crossed products with maximal tensor products of
$C^*$-algebras.

\begin{lem}
\label{lem-tensor}
  Let $(A,G,\alpha)$ be an action and let
$B$ be a $C^*$-algebra. Let $\id\otimes\alpha\:G\to \Aut(B\otimes_{\max}A)$
be the diagonal action of $G$ on $B\otimes_{\max}A$ and let
$k_B$ and $k_{A\times_{\alpha}G}$ denote the canonical maps
of $B$ and $A\times_{\alpha}G$ into $M(B\otimes_{\max}(A\times_{\alpha}G))$.
Further write $k_A=k_{A\times_{\alpha}G}\circ i_A$ and
$k_G=k_{A\times_{\alpha}G}\circ i_G$. Then
\[
(k_B\otimes k_A)\times k_G\:(B\otimes_{\max}A)\times_{\id\otimes\alpha}G
\to B\otimes_{\max}(A\times_{\alpha}G)
\]
is an isomorphism.
\end{lem}

\begin{proof} Just check (using the universal properties of the
crossed product and the maximal tensor product) that
$(B\otimes_{\max}(A\times_{\alpha}G), k_B\otimes k_A, k_G)$
satisfies conditions (i)--(iii) above.
\end{proof}

If $A=\bbC$ the
crossed product $\bbC\times_{\id} G$ is just the (full)
\emph{group $C^*$-algebra} $C^*(G)$ of $G$. We write
$u\:G\to UM(C^*(G))$ for the canonical map (whenever a name is
needed). Thus $C^*(G)$ enjoys the universal
property that for any
strictly continuous homomorphism $V\:G\to UM(D)$
there exists a unique nondegenerate
homomorphism $\bar V\:C^*(G)\to M(D)$ such that $\bar V\circ u=V$.
In what follows we shall make no notational distinction between
$V$ and its integrated form (denoted $\bar V$ above), \ie, we
simply write $V(z)$ for $z\in C^*(G)$.
The \emph{reduced group $C^*$-algebra} $C_r^*(G)$ is the image
of $C^*(G)$ under the \emph{left regular representation}
$\lambda\:G\to U(L^2(G))$ defined by $(\lambda_s\xi)(t)=\xi(s^{-1}t)$.
Notice that $\lambda\:C^*(G)\to C_r^*(G)$ is an isomorphism if and only if
$G$ is amenable \cite[Theorem 7.3.9]{PedCA}.

\begin{ex}
\label{ex-elementary}
One of the most interesting actions from the point of view of duality
theory is the action $(C_0(G), G,\tau)$, where $\tau\:G\to
\Aut(C_0(G))$ is given by left translation: $\tau_s(f)(t)=f(s^{-1}t)$.
If $M\:C_0(G)\to \B(L^2(G))$ denotes the representation of $C_0(G)$ as
multiplication operators on $L^2(G)$, then it is easily seen that
$(M,\lambda)$ is a covariant representation of $(C_0(G), G,\tau)$.  If
$f\in C_c(G\times G)\subseteq C_c(G, C_0(G))$, the integrated form of
$f$ can be written as
\[
\big((M\times \lambda)(f)\xi\big)(t)=\int_G f(s,t)\xi(s^{-1}t)\,ds
=\int_G \Delta(s^{-1})f(ts^{-1},t)\xi(s)\,ds.
\]
Thus, $(M\times \lambda)(f)$ is an integral operator with kernel
$k(s,t)=\Delta(s^{-1})f(ts^{-1},t)$ in $C_c(G\times G)$.  Since
the set of integral operators with kernels in $C_c(G\times G)$ is
dense in $\K(L^2(G))$ it follows that
$(M\times\lambda)(C_0(G)\times_{\tau}G)=\K(L^2(G))$.

It is a nontrivial result that any representation of
$(C_0(G), G,\tau)$ is unitarily equivalent to a representation of the
form $(M\otimes 1, \lambda\otimes 1)$ on $L^2(G)\otimes \H$ for some
Hilbert space $\H$. For $G=\RR^n$ this is a reformulation of the
uniqueness of the Heisenberg commutation relations, and for arbitrary
$G$ this is a special case of Mackey's imprimitivity theorem
(for a good treatment of this special case see
\cite{RieOU},
but we shall give an independent proof in \corref{cor-elementary} below).
It follows that $M\times \lambda\:C_0(G)\times_{\tau}G\to \K(L^2(G))$
is faithful and $C_0(G)\times_{\tau}G\cong \K(L^2(G))$.
\end{ex}

\begin{rem}
\label{rem-compact}
The first part of the above example shows in particular that
\[
\overline{M(C_0(G))\lambda(C^*(G))}
=\overline{\lambda(C^*(G))M(C_0(G))}
=\K(L^2(G)).
\]
A similar result is true for the \emph{right regular representation}
$\rho\:G\to U(L^2(G))$
{given by}
$(\rho_s\xi)(t)=\Delta(s)^{1/2}\xi(ts)$.
To see this, consider the self-adjoint unitary operator $U$ on $L^2(G)$
defined by $(U\xi)(s)=\Delta(s)^{-1/2}\xi(s^{-1})$ and observe that
$\ad U\circ M=M$ and $\ad U\circ\lambda=\rho$.
\end{rem}

If $(A,G,\alpha)$ and $(B,G,\beta)$ are two actions, then a
homomorphism
$\phi\:A\to M(B)$ is called \emph{$G$-equivariant} if
$\phi\circ\alpha_s=\beta_s\circ\phi$
(where we implicitly extend
each $\beta_s$ to an automorphism of $M(B)$).
If 
$(i_B, i_G)$ denotes the canonical maps of $(B,G)$ into
$M(B\times_{\beta}G)$, then we obtain a
\textup(possibly degenerate\textup)  homomorphism
\[
\phi \times G \deq
(i_B\circ \phi)\times i_G\:A\times_{\alpha}G\to M(B\times_{\beta}G).
\]
If $\phi\:A\to B$ is a $G$-equivariant isomorphism, then
$\phi \times G$ is an isomorphism between
$A\times_{\alpha}G$ and $B\times_{\beta}G$ with inverse given by
$\phi^{-1} \times G$.

\begin{ex}
\label{ex-extend}
Let $(A,G,\alpha)$ be an action and let $(C_0(G),G,\tau)$ be as in 
\exref{ex-elementary}. Let $\alpha\otimes\tau$ denote the diagonal action
of $G$ on $A\otimes C_0(G)=A\otimes_{\max}C_0(G)$.
We want to show that $(A\otimes C_0(G))\times_{\alpha\otimes\tau}G$
is canonically isomorphic to $A\otimes \K(L^2(G))$.
To see this define
$\phi\:A\otimes C_0(G)\to A\otimes C_0(G)$ by
$\phi(f)(s)=\alpha_{s^{-1}}(f(s))$ for $f\in C_0(G,A)\cong A\otimes C_0(G)$.
Then it is straightforward to check that $\phi$ is an
$(\alpha\otimes \tau)-(\id\otimes\tau)$ equivariant isomorphism.
Thus, $(A\otimes C_0(G))\times_{\alpha\otimes\tau}G$ is isomorphic to
$(A\otimes C_0(G))\times_{\id\otimes\tau}G$, which in turn is isomorphic
to $A\otimes_{\max} (C_0(G)\times_{\tau}G)\cong A\otimes\K(L^2(G))$
by \lemref{lem-tensor}, \exref{ex-elementary} and
the nuclearity of $\K(L^2(G))$.
\end{ex}

We are now going to define the reduced crossed product of an action
$(A,G,\alpha)$. Writing $\alpha(a)(s)\deq \alpha_{s^{-1}}(a)$,
we may view the action $\alpha$ as a
homomorphism of $A$ into $C_b(G,A)\subseteq M(A\otimes
C_0(G))$.
It is straightforward to check that
$(i_A^r,i_G^r)\deq \big((\id_A\otimes M)\circ \alpha, 1\otimes\lambda\big)$
is a covariant homomorphism of $(A,G,\alpha)$ into
$M(A\otimes \K(L^2(G)))$.

\begin{defn}
\label{def-red crossed}
Let $(A,G,\alpha)$ be an action and let $(i_A^r,i_G^r)=
\big((\id_A\otimes M)\circ \alpha, 1\otimes\lambda\big)$
be as above.
The \emph{reduced crossed product}
$A\times_{\alpha,r}G$ is the image
$(i_A^r\times i_G^r)(A\times_{\alpha}G)\subseteq M(A\otimes\K(L^2(G)))$.
We usually regard $(i_A^r,i_G^r)$ as maps from $(A,G)$ into
$M(A\times_{\alpha,r}G)$ and call them
the \emph{canonical maps}
of $(A,G)$ into $M(A\times_{\alpha,r}G)$.
\end{defn}

\begin{rem}\label{rem-Kmultcrossed}
     The reduced crossed product
     $A\times_{\alpha,r}G$ actually lies in the
     $\K(L^2(G))$-multiplier algebra
     $M_{\K(L^2(G))}(A\otimes \K(L^2(G)))$. To see this
    it is enough to show that
     $i_A^r(a)i_G^r(z)\in M_{\K(L^2(G))}(A\otimes \K(L^2(G)))$
     for all $a\in A$, $z\in C^*(G)$,
     and since $i_G^r(z)\in 1\otimes M(\K(L^2(G)))$
     we can use part (iii) of \propref{prop-Cstrict}
     to see that it actually suffices to show that
     $i_A^r(A)\subseteq M_{\K(L^2(G))}(A\otimes \K(L^2(G)))$.
     For this we first observe that $\alpha(A)\subseteq
     C_b(G,A)\cong M_{C_0(G)}(A\otimes C_0(G))$ (see
     part (2) of Remark \ref{rem-Cmult}).
     Then we use part (iv) of \propref{prop-Cmultextend}
     to see that
     $$i_A^r(A)=(\id_A\otimes M)\circ \alpha(A)\subseteq
     (\id_A\otimes M)\big((M_{C_0(G)}(A\otimes C_0(G))\big)
     \subseteq M_{\K(L^2(G))}(A\otimes\K(L^2(G))).$$
\end{rem}

Let $\pi\:A\to \B(\H)$ be a (possibly degenerate) representation
of $A$ on some Hilbert space $\H$. The
\emph{induced representation} $\Ind\pi\:A\times_{\alpha}G\to
\B(\H\otimes L^2(G))$
is defined as the integrated form of the covariant representation
$((\pi\otimes M)\circ \alpha, 1\otimes \lambda)$ on the
Hilbert space $\H\otimes L^2(G)\cong L^2(G,\H)$.
Thus we have
\[
\Ind\pi\deq (\pi\otimes \id)\circ (i_A^r\times i_G^r),
\]
which implies that $\Ind\pi$ factors through a representation
of $A\times_{\alpha,r}G$ on $\H\otimes L^2(G)$.
Note that by \propref{prop-Cmultextend}
the above composition makes perfect sense even if
$\pi$ is degenerate, since by the above remark
$A\times_{\alpha,r}G=i_A^r\times i_G^r(A\times_{\alpha}G)$
lies in the $\K(L^2(G))$-multiplier algebra of $A\otimes \K(L^2(G))$.
Since $\pi\otimes \id\: A\otimes \K(L^2(G))\to \B(\H\otimes L^2(G))$
is faithful if $\pi$ is, it follows that $\Ind\pi$
factors through a faithful representation of $A\times_{\alpha,r}G$
whenever $\pi$ is faithful. In this case we call
$\Ind\pi$ the \emph{regular representation  of $(A,G,\alpha)$ induced
from $\pi$}.

\begin{rem}
\label{rem-red crossed}
In the literature, the reduced crossed product
is often defined as the image $\Ind\pi(A\times_{\alpha}G)$
of any regular representation $\Ind\pi$ of $(A,G,\alpha)$.
Of course, this would determine $A\times_{\alpha,r}G$ only up
to isomorphism. However, it is often convenient to work with
regular representations and to identify $A\times_{\alpha,r}G$
with its image
$\Ind\pi(A\times_{\alpha,r}G)\subseteq \B(\H\otimes L^2(G))$.
Whenever we do this, we also identify $i_A^r$ with
$(\pi\otimes M)\circ \alpha$ and $i_G^r$ with
$1\otimes\lambda$ (although we will not do this outside of this
appendix).
\end{rem}

The following lemma is quite useful:

\begin{lem}\label{lem-equivdegenerate}
     Assume that $\alpha$ and $\beta$ are actions of $G$ on
     $C^*$-algebras $A$ and $B$, respectively. Assume further
     that $\phi\:A\to M(B)$ is a 
\textup(possibly degenerate\textup) $G$-equivariant
      $*$-homomorphism.
     Then there exists a $*$-homomorphism
     $$\phi\times_rG\deq (i_B^r\circ \phi)\times i_G^r\:
     A\times_{\alpha,r}G\to M(B\times_{\beta,r}G).$$
\end{lem}
\begin{proof} We have to check that the  homomorphism
     $(i_B^r\circ \phi)\times i_G^r\:A\times_{\alpha}G\to
     M(B\times_{\beta,r}G)$ factors through $A\times_{\alpha,r}G$.
     For this consider the composition
     $$
     \begin{CD}
	A\times_{\alpha} G@>i_A^r\times i_G^r>>
	M_{\K(L^2(G))}(A\otimes \K(L^2(G)) @>\phi\otimes \id>>
	M(B\otimes \K(L^2(G))).
      \end{CD}
      $$
      Since $\phi$ is $G$-equivariant, we have
      \begin{align*}
	 (\phi\otimes \id)\circ i_A^r&=
      (\phi\otimes \id)\circ(\id_A\otimes M)\circ \alpha\\
      &=(\id_B\otimes M)\circ (\phi\otimes \id)\circ \alpha
      =(\id_B\otimes M)\circ \beta\circ\phi.
      \end{align*}
      Since it is clear that $(\phi\otimes \id)\circ i_G^r=i_G^r$
      (in the appropriate sense), it follows that
      $$(\phi\otimes\id)\circ (i_A^r\times i_G^r)=(i_B^r\circ
      \phi)\times i_G^r$$
      and the result follows. If $\phi$ is faithful, then
      $\phi\otimes\id$ is faithful on
      $M_{\K(L^2(G))}(A\otimes\K(L^2(G))\supseteq
      A\times_{\alpha,r}G$ by \propref{prop-Cmultextend}.
      Hence $\phi\times_rG$ is faithful, too.
\end{proof}

\begin{rem}\label{rem-phiequiv}
     It is helpful to observe that the
     homomorphisms
     $$\phi\times G\:A\times_{\alpha}G\to M(B\times_{\beta}G)
     \quad\text{and}\quad
     \phi\times_rG\:A\times_{\alpha,r}G\to M(B\times_{\beta,r}G)
     $$ are given on elements of $C_c(G,A)$
     by $\phi\times G(f)=\phi\circ f\in C_c(G, M(B))$,
     acting on $C_c(G,B)$ via convolution.
\end{rem}

If $\pi\times V$ is a representation of $A\times_{\alpha}G$ on $\H$
and $U$ is a representation of $G$ on $\K$, then we write
$(\pi\times U)\otimes V\:A\times_{\alpha}G\to \B(\H\otimes \K)$ for the
integrated
form of the covariant representation $(\pi\otimes 1_{\K}, V\otimes U)$.
We shall frequently use:

\begin{lem}
\label{lem-regrep}
Let $\pi\:A\to \B(\H)$ be a representation of $A$,
and let $U$ be a representation of $G$ on some Hilbert space $\K$. Then
\begin{enumerate}
\item $(\Ind\pi)\otimes U$ is unitarily equivalent to 
$\Ind(\pi\otimes 1_{\K})$.
\item If $\K=\H$ and $(\pi,U)$ is covariant for $(A,G,\alpha)$,
then $\Ind\pi$ is unitarily equivalent to $(\pi\times U)\otimes \lambda$.
\end{enumerate}
\end{lem}
\begin{proof} For the proof of~(i), identify
$\H\otimes\K\otimes L^2(G)$ and $\H\otimes L^2(G)\otimes \H$ with
$L^2(G, \H\otimes \K)$ in the canonical way and check that the unitary
$W$ on $L^2(G, \H\otimes \K)$ defined
by $(W\xi)(s)=(1\otimes U_s)\xi(s)$ intertwines the representations.
Similarly, for the proof of (ii) check that
$W\in U(L^2(G,\H))$ defined by $(W\xi)(s)=U_s\xi(s)$ does the job.
\end{proof}

As a first application we get:

\begin{prop}
\label{prop-amenable}
If $G$ is amenable, then
$i_A^r\times i_G^r\:A\times_{\alpha}G\to A\times_{\alpha,r}G$
is an isomorphism for any action $\alpha\:G\to \Aut(A)$.
\end{prop}
\begin{proof}
Choose any faithful representation $\pi\times U$
of $A\times_{\alpha}G$ on some Hilbert space $\H$, and let
$1_G\:C^*_r(G)\cong C^*(G)\to \bbC$
denote the integrated form of the trivial representation of $G$.
We then have the identity
\[
(\id\otimes 1_G)\circ \bigl((\pi\times U)\otimes\lambda\bigr)
=\pi\times U
\]
under the canonical identification $\H\otimes \bbC\cong \H$.
In particular, it follows  that $(\pi\times U)\otimes\lambda$
is injective. But \lemref{lem-regrep} implies that
$(\pi\times U)\otimes\lambda$ is unitarily equivalent to a
regular representation of $(A,G,\alpha)$. Hence
$\ker (i_A^r\times i_G^r)=\ker (\pi\times U)\otimes\lambda=\{0\}$.
\end{proof}

Note that $\bbC\times_{\id}G=C^*(G)\cong C_r^*(G)=
\bbC\times_{\id,r}G$ via $i_{\bbC}^r\times i_G^r=\lambda$ if and
only if $G$ is amenable, so that the reduced crossed product does not
coincide with the full crossed product in general.
We should point out
that it follows from \exref{ex-elementary} that every
representation of $C_0(G)\times_{\tau}G$ is faithful, so
$C_0(G)\times_{\tau}G\cong C_0(G)\times_{\tau,r}G$ even if $G$ is
not amenable.\footnote{If the regular homomorphism $i_A^r \times
i_G^r \: A \times_\alpha G \to A \times_{\alpha,r} G$ is not faithful,
can the $C^*$-algebras $A \times_\alpha G$ and $A \times_{\alpha,r} G$
still be isomorphic? Presumably so, but we do not know an
example.}
Anyway, reduced crossed products can be viewed
in a certain sense
as  skew minimal tensor products of $A$ with the reduced group
algebra $C_r^*(G)$. For instance, parallel to \lemref{lem-tensor}
we have:

\begin{lem}
\label{lem-mintensor}
Let $(A,G,\alpha)$ be an action and let $B$ be a $C^*$-algebra.
Then there exists a canonical isomorphism
$(B\otimes A)\times_{\id\otimes\alpha,r}G\cong
B\otimes (A\times_{\alpha,r}G)$.
\end{lem}
\begin{proof} Let
$(k_B\otimes k_A)\times k_G\:(B\otimes_{\max}A)\times_{\id\otimes\alpha}G
\to B\otimes_{\max}(A\times_{\alpha}G)$ be the isomorphism of
\lemref{lem-tensor}. Then it is easy to check that
$(k_B\otimes k_A)\times k_G$
transforms the homomorphism
\[
(\id_B\otimes \id_A\otimes M)\circ
(\id_B\otimes\alpha)\times (1\otimes 1\otimes\lambda)\:
(B\otimes_{\max}A)\times_{\id\otimes\alpha}G\to
M(B\otimes A\otimes \K(L^2(G)))
\]
to the homomorphism
\[
\id_B\otimes
\big((\id_A\otimes M)\circ \alpha\times (1\otimes\lambda)\big)\:
B\otimes_{\max}(A\times_{\alpha}G)\to M(B\otimes A\otimes
\K(L^2(G))).
\]
But the image of the first
map is $(B\otimes A)\times_{\id\otimes \alpha,r}G$ and the image
of the second is $B\otimes(A\times_{\alpha,r}G)$.
\end{proof}

The combination of \lemref{lem-tensor}, \lemref{lem-mintensor}, and
\propref{prop-amenable}
implies the important fact that $A\times_{\alpha}G$ is nuclear
whenever $A$ is nuclear and $G$ is amenable: if $B$ is any
$C^*$-algebra, then
\begin{align*}
B\otimes_{\max}(A\times_{\alpha}G)
&\cong
(B\otimes_{\max}A)\times_{\id\otimes\alpha}G
\cong
(B\otimes A)\times_{\id\otimes\alpha}G
\\
&\cong (B\otimes A)\times_{\id\otimes\alpha,r}G
\cong B\otimes (A\times_{\alpha,r}G)
\cong B\otimes (A\times_{\alpha}G).
\end{align*}

\section{Coactions}
\label{sec-coaction}

If $(A,G,\alpha)$ is an action such that $G$ is abelian, then there
exists a natural action $\widehat{\alpha}$ of the dual group
$\widehat{G}$ of $G$ on
$A\times_{\alpha}G$ given for $f\in C_c(G,A)$ by 
\begin{equation}
\widehat{\alpha}_{\chi}(f)(s)=\chi(s)f(s).
\end{equation}
The famous Takesaki-Takai duality theorem
asserts that the double crossed product
$A\times_{\alpha}G\times_{\widehat\alpha}\widehat G$ is isomorphic
to $A\otimes\K(L^2(G))$, \ie, the double crossed product is stably
isomorphic to $A$.
If $G$ is nonabelian, then $\widehat G$ is not a group
and there is no dual action to talk about.

Coactions of groups on $C^*$-algebras were mainly introduced in order
to overcome this problem and to obtain a reasonable duality theory
for actions of nonabelian groups.
In order to define coactions of groups on $C^*$-algebras let us first
note that the group $C^*$-algebra $C^*(G)$ of $G$ carries a natural
\emph{comultiplication} given by the integrated form
$\delta_G\:C^*(G)\to M(C^*(G)\otimes C^*(G))$ of the strictly continuous
homomorphism $s\mapsto s\otimes s\in UM(C^*(G)\otimes C^*(G))$.
It follows directly from the definition that
$\delta_G$ satisfies the \emph{comultiplication identity}
\[
(\delta_G\otimes\id_G)\circ \delta_G
=(\id_G\otimes\delta_G)\circ \delta_G.
\]

\begin{defn}
\label{def-coaction}
A (\emph{full}) \emph{coaction} of a locally compact group $G$ on
a $C^*$-algebra $A$ is an injective and nondegenerate
homomorphism $\delta\:A\to M(A\otimes C^*(G))$ satisfying
\begin{enumerate}
\item
$\delta(A)(1\otimes C^*(G))\subseteq A\otimes C^*(G)$, and
\item
$(\delta\otimes \id_G)\circ \delta=(\id_A\otimes\delta_G)\circ
\delta$ as maps from $A$ into $M(A\otimes C^*(G)\otimes C^*(G))$ (the
\emph{coaction identity}).
\end{enumerate}
We also call the triple $(A,G,\delta)$ a coaction.
A coaction $\delta$ is called \emph{nondegenerate} if
\begin{enumerate}
\item[(iii)] $\overline{\delta(A)(1\otimes C^*(G))}=A\otimes C^*(G)$.
\end{enumerate}
\end{defn}
It might be helpful to realize that Condition (ii) above
simply means that the diagram
\begin{equation*}
\xymatrix{
{A}
\ar[r]^-{\delta}
\ar[d]_{\delta}
&{M(A\otimes C^*(G))}
\ar[d]^{\id_A\otimes\delta_G}
\\
{M(A\otimes C^*(G))}
\ar[r]_-{\delta\otimes\id_G}
&{M(A\otimes C^*(G)\otimes C^*(G))}
}
\end{equation*}
commutes.

\begin{rem}
\label{rem-coaction}
(1)
If we apply the $*$-operation to the inclusion in item (i) of
the definition, we obtain the inclusion $(1\otimes C^*(G))\delta(A)
\subseteq A\otimes C^*(G)$. This shows that condition (i) of the
definition is equivalent to the requirement that
$\delta(A)$ lies in the $C^*(G)$-multiplier algebra
$M_{C^*(G)}(A\otimes C^*(G))$ of $A\otimes C^*(G)$ (see
\defnref{defn-Cmult}). In order to simplify notation we
shall  write $M_G(A\otimes C^*(G))$ for $M_{C^*(G)}(A\otimes C^*(G))$
and will call it the {\em $G$-multiplier algebra} of $A\otimes C^*(G)$.

(2) The easiest example of a coaction is $\delta_G$ itself, which 
is a coaction of $G$
on $C^*(G)$. The notion of a coaction we use here is that of a
\emph{full} coaction of $G$ on $A$ in the sense of \cite{qui:fullred}.
The term \emph{full} refers to the fact that
we use $C^*(G)$ instead of the reduced group
algebra $C_r^*(G)$, which had for quite some time been
the standard setting in the literature
(\eg, see \cite{it, LanDT, kat, lprs}).
Full coactions were first introduced in
\cite{rae:rep}, using maximal tensor products instead
of the minimal tensor products we use here.
We shall discuss the (relatively small)
differences between the different approaches later
(see \secref{sec-other} below).

(3) Somehow irritatingly, the word \emph{nondegenerate}
 has two meanings in connection
with coactions. By definition, every coaction is a  nondegenerate
homomorphism of $A$ into $M(A\otimes C^*(G))$ (\ie,
$\delta(A)(A\otimes C^*(G))=A\otimes C^*(G)$), but being
nondegenerate \emph{as a coaction} is the apparently stronger condition that
$\overline{\delta(A)(1\otimes C^*(G))}=A\otimes C^*(G)$. The reason
for this terminology is that a nondegenerate coaction determines
a nondegenerate module action of the Fourier algebra $A(G)$ on $A$
(see \propref{prop-nondeg} below).
It is actually an
open question whether every coaction is automatically
nondegenerate, although it has been settled affirmatively for $G$ amenable
\cite[Lemma 3.8]{LanDT}, \cite[Proposition 6]{kat} and
$G$ discrete \cite{BS-CH}.
\end{rem}

\begin{ex}
\label{ex-abelian}
If $G$ is abelian, then there is a one-to-one correspondence between
coactions of $G$ and strongly continuous actions of the dual group
$\widehat G$.  To see this let us identify $C^*(G)$ with
$C_0(\widehat{G})$ via the Fourier transform $\mathcal F\:C^*(G)\to
C_0(\widehat G)$ given by $\mathcal F(x)(\chi)=\chi(x)$.  A brief
calculation shows that the comultiplication $\delta_G$ is then
translated to the formula
\[
\delta_G(f)(\chi,\mu)=f(\chi\mu)\in C_b(\widehat G\times \widehat G)
\]
for $f\in C_0(\widehat{G}), \chi,\mu\in \widehat G$.  If
$\alpha\:\widehat G\to \Aut(A)$ is an action, then $\alpha$ determines
an injective and nondegenerate homomorphism $\delta^{\alpha}\:A\to
C^b(\widehat G,A)\subseteq M(A\otimes C_0(\widehat{G}))$ by the
formula
\[
\delta^{\alpha}(a)(\chi)\deq \alpha_{\chi}(a)
\righttext{for} a\in A, \chi\in \widehat G.
\]
Condition (i) of the definition is,
in this setting, equivalent to
$\delta^{\alpha}$ taking values in the subalgebra $C^b(\widehat G,A)$
of $M(A\otimes C_0(\widehat{G}))$. The coaction identity (ii) follows from
a straightforward computation,
using multiplicativity of $\alpha$.
Thus $\delta^{\alpha}$ is a coaction
of $G$ on $A$.
Conversely, if $\delta\:A\to C_b(\widehat{G},A)$ is any
injective nondegenerate homomorphism which satisfies the coaction identity,
then we obtain an action  of $\widehat{G}$ on $A$ by putting
$\beta_{\chi}(b)\deq \delta(b)(\chi)$.
\end{ex}

Before we present some other important examples of coactions we need

\begin{lem}[{\cf~\cite[Remarks 2.2]{rae:rep}}]
\label{lem-injective}
Let $\delta\:A\to M(A\otimes C^*(G))$ be a nondegenerate homomorphism
which satisfies conditions \textup{(i)} and \textup{(ii)}
of \defnref{def-coaction}
above.
Let $1_G\:G\to \bbC$ denote the trivial representation of $G$ and let
us identify $A$ with $A\otimes\bbC$ in the canonical way.
Then $\delta\circ (\id_A\otimes 1_G)\circ \delta=\delta$.
In particular,
$\delta$ is injective \textup(hence a coaction\textup) if and only if
$(\id_A\otimes 1_G)\circ\delta\:A\to A$ is the identity on $A$.
\end{lem}
\begin{proof}
It follows from (i) that
$(\id_A\otimes 1_G)\big(\delta(a)(1\otimes z)\big)\in A$ for all
$a\in A$, $z\in C^*(G)$, and hence that
$(\id_A\otimes 1_G)\circ\delta(a)\in A$.
Since
$(\id_G\otimes 1_G)\circ \delta_G(s)=(\id_G\otimes 1_G)(s\otimes s)=s$
{for all} $s\in G$,
we have $(\id_G\otimes 1_G)\circ \delta_G=\id_G$.
This together with (ii) gives
\begin{align*}
\delta(a)&=(\id_A\otimes \id_G\otimes 1_G)\circ(\id_A\otimes\delta_G)\circ
\delta(a)
=(\id_A\otimes \id_G\otimes 1_G)\circ(\delta_A\otimes\id_G)\circ
\delta(a)\\
&=(\delta\otimes 1_G)\circ\delta(a)=\delta\big((\id_A\otimes 
1_G)(\delta(a))\big).
\end{align*}
This completes the proof.
\end{proof}

\begin{rem}
The map $\id_A\otimes 1_G\:M(A\otimes C^*(G))\to M(A)$ above is the first
appearance of a slice map. We shall see later that more general slice maps
play an important r\^ole in the theory.
\end{rem}

\begin{ex}
\label{ex-dualcoaction}
Let $(A,G,\alpha)$ be an action,
let $(i_A,i_G)$ denote the canonical maps from
$(A,G)$ into $M(A\times_{\alpha}G)$ and let $u$ denote the canonical
map from $G$ into $M(C^*(G))$. The \emph{dual coaction} of $G$ on
$A\times_{\alpha}G$
is defined as the integrated form
\[
\widehat{\alpha}\deq (i_A\otimes 1)\times (i_G\otimes u)\:
A\times_{\alpha}G\to M\big((A\times_{\alpha}G)\otimes C^*(G)\big).
\]
Let us check conditions~(i) and~(ii) of \defnref{def-coaction}.
For~(i) we consider the map
$\Psi\:C_c(G\times G,A)\to (A\times_{\alpha}G)\otimes C^*(G)$
given by
\[
\Psi(g)=\int_{G\times G}
(i_A\otimes 1)(f(s,t))\bigl(i_G(s)\otimes u(t)\bigr)\,d(s,t).
\]
If $z\in C_c(G,A)$ and $w\in C_c(G)$, let
$z\diamond w(s,t)=z(s)w(s^{-1}t)\in C_c(G\times G,A)$. Then
\begin{align*}
\widehat{\alpha}(z)(1\otimes u(w))
&=\int_{G\times G} (i_A\otimes 1)(z(s)w(t))(i_G(s)\otimes u(st))\,d(s,t)\\
&=\int_{G\times G} (i_A\otimes 1)(z(s)w(s^{-1}t))(i_G(s)\otimes u(t))\,d(s,t)\\
&=\Psi(z\diamond w)\in (A\times_{\alpha}G)\otimes C^*(G).
\end{align*}
Hence (i) follows from the fact that $C_c(G,A)$ and $C_c(G)$ are dense in
$A\times_{\alpha}G$ and $C^*(G)$, respectively. Now we check (ii).
Since $\widehat{\alpha}(i_A(a))=i_A(a)\otimes 1$ we first get
\[
(\widehat{\alpha}\otimes\id_G)\circ \widehat{\alpha}(i_A(a))
=i_A(a)\otimes 1\otimes 1=
(\id_{A\times_{\alpha}G}\otimes\delta_G)\circ \widehat{\alpha}(i_A(a))
\]
for all $a\in A$, and using $\widehat{\alpha}(i_G(s))=i_G(s)\otimes s$ and
$\delta_G(s)= s\otimes s\in M(C^*(G)\otimes C^*(G))$ we get
\[
(\widehat{\alpha}\otimes\id_G)\circ \widehat{\alpha}(i_G(s))
=i_G(s)\otimes s\otimes s
=(\id_{A\times_{\alpha}G}\otimes\delta_G)\circ \widehat{\alpha}(i_G(s))
\]
for $s\in G$. Hence (ii) follows from integration. Finally, injectivity follows
from \lemref{lem-injective} by observing that
$(\id_{A\times_{\alpha}G}\otimes 1_G)\circ \widehat\alpha= i_A\times i_G$.

Note that $\widehat\alpha$ is always nondegenerate.
This follows from inductive-limit density of
$C_c(G,A)\diamond C_c(G)$ in
$C_c(G\times G,A)$, continuity of $\Psi$ with respect to the inductive
limit topology on $C_c(G\times G,A)$, and density of
$\Psi\big(C_c(G\times G,A)\big)$ in
$(A\times_{\alpha}G)\otimes C^*(G)$ (since it contains
$C_c(G,A)\odot C_c(G)$).
\end{ex}

It is a good exercise for the reader to check that in case where $G$
is abelian, the dual coaction $\widehat{\alpha}$ corresponds to the
dual action of $\widehat{G}$ on $A\times_{\alpha}G$ under the
correspondence given in \exref{ex-abelian}.  The next example shows
that we also have a dual coaction on the reduced crossed product
$A\times_{\alpha,r}G$.

\begin{ex}
\label{ex-reduceddual}
Let $(i_A^r,i_G^r)$ denote the canonical maps of $(A,G)$ into
$M(A\times_{\alpha,r}G)$, and let $u\:G\to M(C^*(G))$ be as above.
Then \lemref{lem-regrep} implies that
$(i_A^r\otimes 1)\times(i_G^r\otimes u)$ factors
through a faithful representation $\widehat{\alpha}^n$
of $A\times_{\alpha,r}G$ into
$M\big((A\times_{\alpha,r}G)\otimes C^*(G)\big)$ (write $i_A^r\times
i_G^r=\Ind\pi$ for some regular representation $\Ind\pi$ and use part
(i) of \lemref{lem-regrep} to deduce that $(i_A^r\otimes
1)\times(i_G^r\otimes u)$ is injective).  Exactly the same arguments
as in the preceding example show that $\widehat{\alpha}^n$ is a
nondegenerate coaction of $G$ on $A\times_{\alpha,r}G$, which is
called the \emph{dual coaction} of $G$ on $A\times_{\alpha,r}G$
\end{ex}

It turns out later that $\widehat{\alpha}^n$ is the
\emph{normalization} of the coaction $\widehat\alpha$ of
\exref{ex-dualcoaction} (see \propref{prop-dualnormal} below).
This is the motivation for using the superscript ``$n$'' in
our notation.
However, we shall often skip the``$n$'' if no confusion
is possible---in particular, we shall always write $\widehat\alpha$
instead of $\widehat\alpha^n$ in the main body of this work.

The following example shows that coactions of $G$ restrict to
coactions of $G/N$ for any closed normal subgroup $N$ of $G$.  
In such a situation we shall always assume that Haar measures on
$G/N$ and $N$ are normalized so that $\int_Gf(s)\,ds=\int_{G/N}\int_N
f(sn)\, dn \, d{sN}$ for all $f\in C_c(G)$.

\begin{ex}
\label{ex-restrict}
Let $(A,G,\delta)$
be a coaction and
let $N$ be a closed normal subgroup of $G$.
Let $q\:G\to G/N$ denote the quotient map. Identifying $sN\in G/N$
with its image in $M(C^*(G/N))$ we may view $q$ as a unitary homomorphism
from $G$ into $UM(C^*(G/N))$ whose integrated form is the quotient map
$q\:C^*(G)\to C^*(G/N)$ (to see that it takes values in $C^*(G/N)$ and is
surjective one checks that it restricts on $C_c(G)$ to the surjective map
$C_c(G)\to C_c(G/N)$ given by $q(f)(sN)=\int_Nf(sn)\,dn$). Then
\[
\delta|\deq (\id_A\otimes q)\circ \delta\:A\to M(A\otimes C^*(G/N))
\]
is a coaction of $G/N$ on $A$, which is called the \emph{restriction}
of $\delta$ to $G/N$.
Note that Condition (i) of \defnref{def-coaction}
follows from
\begin{equation}
\label{eq-rest}
\delta|(a)(1\otimes q(z))=(\id_A\otimes q)(\delta(a))(1\otimes q(z)
=(\id_A\otimes q)\big(\delta(a)(1\otimes z)\big)
\end{equation}
for $z\in C^*(G)$, and the  coaction identity for $\delta|$
follows from applying $\id_A\otimes q\otimes q$ to both sides
of the coaction identity for $\delta$. Since $1_{G/N}\circ q=1_G$,
injectivity of $\delta|$ follows from \lemref{lem-injective}.
Moreover, it is a direct consequence of \eqeqref{eq-rest}
that $\delta|$ is nondegenerate if $\delta$ is.
\end{ex}

The next example shows that any coaction of a closed subgroup
$H$ of $G$ inflates to a coaction of $G$.

\begin{ex}
\label{ex-inflate}
Let $H$ be a closed subgroup of $G$. Let
$\phi \:C^*(H)\to M(C^*(G))$ denote the integrated form
of $u|_H\:H\to M(C^*(G))$, where $u\:G\to M(C^*(G))$ denotes the
canonical map. If $(A,H,\eps)$
is a coaction,
then $\Infl\eps\deq  (\id_A\otimes \phi )\circ\eps$ is a coaction
of $G$ on $A$: Condition (i) of the definition follows from
factoring $z\in C^*(G)$ as $z= \phi (v)w$ for some $v\in C^*(H)$
and $w\in C^*(G)$, and computing
\begin{equation}
\label{eq-inf}
\Infl\eps(a)(1\otimes z)
=(\id_A\otimes \phi )\big(\eps(a)(1\otimes v)\big)(1\otimes w)
\in A\otimes C^*(G).
\end{equation}
Condition (ii) follows from
$(\phi \otimes \phi )\circ\delta_H=\delta_G\circ \phi $,
and injectivity of $\Infl\eps$ is a
consequence of \lemref{lem-injective} and the fact that $1_G|_H=1_H$.
We call $\Infl\eps$ the coaction \emph{inflated} from $\eps$.
Note that (\ref{eq-inf}) also implies that $\Infl\eps$ is
nondegenerate if $\eps$ is nondegenerate.
\end{ex}

\section{Slice maps and nondegeneracy}
\label{sec-nondeg}

Slicing in tensor products is one of the basic tools in the theory
of coactions. In particular, slice maps are indispensable for
the study of covariant representations of coactions. Thus we
now provide the basic facts about slice maps which are needed
in this work.

If $B$ is a $C^*$-algebra and $\rho\:B\to \B(\H)$ is a nondegenerate
representation, then 
$b\mapsto \lk\rho(b)\xi,\eta\rk$ is a continuous
linear functional on $B$ for each $\xi,\eta\in\H$
(called a \emph{matrix coefficient} of
$\rho$), and every element of $B^*$ can
be written in this way.  We denote by $B^*_{\rho}$ the matrix
coefficients of $\rho$.  There are canonical left and right actions of
$M(B)$ on $B^*$ (respectively $B^*_{\rho}$) given by
\begin{equation}
\label{eq-act}
m\cdot f(b)= f(bm)
\midtext{and}
f\cdot m(b)=f(mb)
\end{equation}
for $m\in M(B)$, $f\in B^*$ (respectively $f\in B^*_{\rho}$).  For
$f(b)=\lk\rho(b)\xi, \eta\rk$ in $B^*_{\rho}$ we can factor $\eta=
\rho(c)\zeta$ and $\xi=\rho(c')\zeta'$ for some $c,c'\in B$ and
$\zeta, \zeta'\in \H$.  Thus, if $g, g'\in B^*_{\rho}$ are given by
$g(b)=\lk \rho(b)\xi, \zeta\rk$ and $g'(b)=\lk \rho(b)\zeta',\eta\rk$
we see that $f= g\cdot c^*= c'\cdot g'$.  The possibility of factoring
elements $f\in B^*_{\rho}$ as above turns out to be essential.

If $A$ and $B$ are $C^*$-algebras and $f\in B^*$,
then the \emph{slice map}
$S_f\:A\odot B\to A$ defined by
\[
S_f\biggl(\sum_{i=1}^na_i\otimes b_i\biggr)=\sum_{i=1}^na_if(b_i)
\]
extends to a bounded linear map of $A\otimes B$ into $A$
of norm $\|f\|$. We shall frequently use the following
lemma, which is essentially \cite[Lemma 15]{lprs}. We give
an alternative and (probably) more elementary proof here.

\begin{lem}[{\cf~\cite[Lemma 1.5]{lprs}}]
\label{lem-slice}
$S_f$ extends to a strictly continuous linear map $S_f\:M(A\otimes
B)\to M(A)$ satisfying
\begin{equation}
\label{eq-slice}
S_{b\cdot f}(m)a=S_f\big(m(a\otimes b)\big)
\midtext{and}
aS_{f\cdot b}(m)=S_f\big((a\otimes b)m\big)
\end{equation}
for all $a\in M(A)$, $f\in B^*$ and $b\in M(B)$.\footnote{Note that
our definition of the action of $B$ on $B^*$ differs from that used in
\cite{lprs}.} If $m\in M(A\otimes B)$ such that $S_f(m)=0$ for all
$f\in B^*$, then $m=0$.  Finally, if $\rho\:A\to M(D)$ is a
nondegenerate homomorphism then
\begin{equation}
\label{eq-homslice}
S_f\circ(\rho\otimes\id_B)=\rho\circ S_f.
\end{equation}
\end{lem}

\begin{proof}
It is straightforward to check that (\ref{eq-slice}) holds
for all $m\in A\otimes B$. To see that $S_f$ extends to $M(A\otimes B)$
write $f=c\cdot g\cdot c'$ for some $c,c'\in B$ and $g\in
B^*$. Then, for $m\in M(A\otimes B)$ and $a\in A$, define
\[
S_f(m)a= S_{g\cdot c'}(m(a\otimes c))
\midtext{and}
aS_f(m)=S_{c\cdot g}((a\otimes c')m).
\]
The maps $a\mapsto S_f(m)a$ and $a\mapsto aS_f(m)$ satisfy the
double-centralizer property:
\[
(aS_f(m))b
= S_{c\cdot g}((a\otimes c')m)b
=S_g\bigl((a\otimes c')m(b\otimes c)\bigr)
=a(S_f(m)b)
\]
for all $a,b\in A$, so $S_f(m)\in M(A)$. If $\{m_i\}_{i\in I}$ is a net
in $M(A\otimes B)$ which converges strictly to $m\in M(A\otimes B)$, then
\[
S_f(m_i)a=S_{g\cdot c'}(m_i(a\otimes c))\to S_{g\cdot c'}(m(a\otimes c))
=S_f(m)a
\]
for all $a\in A$ (since $S_f$ is norm continuous on $A\otimes B$).
Similarly, $aS_f(m_i)\to aS_f(m)$ for $a\in A$, which proves that
$S_f$ is strictly continuous.  Since $S_f$ extends the slice map on
$A\otimes B$ and $A\otimes B$ is strictly dense in $M(A\otimes B)$
this also shows that the definition of $S_f(m)$ does not depend on the
particular factorization $f= c\cdot g\cdot c'$.  But this implies that
(\ref{eq-slice}) holds for all $a\in A$, $b\in B$ and $m\in M(A\otimes
B)$.  Further, if $a\in M(A)$, $b\in M(B)$ and $f=c\cdot g$, then
\[
S_{b\cdot f}(m)ad
=S_{(bc)\cdot g}(m)ad
=S_{g\cdot c'}(m(ad\otimes bc))
=S_{c\cdot g}(m(a\otimes b))d
=S_f(m(a\otimes b))d
\]
for all $d\in A$, which implies that $S_{b\cdot f}(m)a=S_f(m(a\otimes
b))$.  Similarly, $aS_{f\cdot b}(m)=S_f((a\otimes b)m)$, so
(\ref{eq-slice}) holds for all $a\in M(A), b\in M(B)$, and $m\in
M(A\otimes B)$.

If $g\in A^*$, $f\in B^*$, then $g\otimes f\in (A\otimes B)^*$ extends
to a strictly continuous functional on $M(A\otimes B)$.  Moreover, the
set $\{g\otimes f\mid g\in A^*, f\in B^*\}$ separates the points of
$M(A\otimes B)$%
\footnote{To see this, choose faithful representations
$\pi$ and $\sigma$ of $A$ and $B$ and observe that $(g\otimes f)(m)=0$
for all $g\in A_{\pi}^*, f\in B_{\sigma}^*$ implies that
$(\pi\otimes\sigma)(m)=0$, and hence $m=0$ since $\pi\otimes \sigma$
is faithful on $A\otimes B$.}%
.  Since $m\mapsto g(S_f(m))$ and
$m\mapsto (g\otimes f)(m)$ are both strictly continuous and clearly
agree on $A\otimes B$, they also agree on $M(A\otimes B)$.  Thus, if
$S_f(m)=0$ for all $f\in B^*$, then $(g\otimes f)(m)=g(S_f(m))=0$ for
all $g\in A^*, f\in B^*$ which implies $m=0$.

The final assertion now follows because both sides of
\eqeqref{eq-homslice} are strictly continuous and coincide
on $A\otimes B$.
\end{proof}

In this work we are mainly interested in the case $B=C^*(G)$ for a
locally compact group $G$.  In this case we can identify $C^*(G)^*$
with the \emph{Fourier-Stieltjes algebra} $B(G)$ of $G$, which
consists of all bounded continuous functions on $G$ which can be
expressed as matrix coefficients of unitary representations of $G$,
\ie, a function $f\in B(G)$ is of the form $f(s)= \lk
V_s\xi,\eta\rk$, where $V$ is a unitary representation of $G$ on a
Hilbert space $\H$ and $\xi,\eta\in\H$.  Of course, if $x\in C^*(G)$
(or $M(C^*(G))$) and $f\in B(G)$ with $f(s)=\lk V_s\xi,\eta\rk$, then
$f(x)=\lk V(x)\xi,\eta\rk$.  In particular, we have $f(u(s))=f(s)$ if
$u\:G\to UM(C^*(G))$ denotes the canonical map.  Note that $B(G)$
becomes an involutive Banach algebra when equipped with the dual norm,
pointwise multiplication and the involution $f\mapsto \bar{f}$.

The \emph{Fourier algebra} $A(G)$ is the (closed) $*$-subalgebra of
$B(G)$ consisting  of all matrix coefficients of the
left regular representation  $\lambda$, \ie, $A(G)=C^*(G)^*_{\lambda}$.
The nontrivial fact that $A(G)$ is an algebra can be deduced from the fact
that $\lambda\otimes\lambda$ is unitarily equivalent to
$\lambda\otimes 1_{L^2(G)}$ (\lemref{lem-regrep}) and that
$\lambda\oplus\lambda$ is unitarily equivalent to
$\lambda\otimes 1_{\CC^2}$.
One can show that $A(G)$ is dense in $C_0(G)$ with respect to the
supremum-norm. Moreover, the subalgebra $A_c(G)$ of compactly supported
elements is dense in $A(G)$ (which follows from the fact that
the compactly supported elements are dense in $L^2(G)$) and for each
compact subset $C$ of $G$ there exists an element $f\in A_c(G)$
with $f|_C\equiv 1$.  For more details on $A(G)$ and $B(G)$ we
refer to \cite{eym}.
The following proposition indicates that any nondegenerate coaction
turns $A$ into a nondegenerate $A(G)$-module.

\begin{prop}
\label{prop-nondeg}
Let $(A,G,\delta)$
be a nondegenerate coaction.
Then for all $f \in B(G)$ the composition $\delta_f \deq  S_f \circ \delta$
maps $A$ back into itself%
, and
$\spn\{\delta_f(a)\mid a\in A, f\in A(G)\}$ is dense in $A$.
\end{prop}

\begin{proof}
Let $x\in C^*(G)$ and $f\in A(G)$ with $f(x)=1$.  Let $a\in A$ be
given.  Since $\delta(A)(1\otimes C^*(G))$ is dense in $A\otimes
C^*(G)$ we can approximate $a\otimes x$ by a finite sum
$\sum_{i=1}^n\delta(a_i)(1\otimes x_i)$ with $a_i\in A$ and $x_i\in
C^*(G)$.  Then, using \lemref{lem-slice}, we get
\begin{align*}
a
&=f(x)a
=S_f(a\otimes x)
\approx S_f\biggl(\sum_{i=1}^n\delta(a_i)(1\otimes x_i)\biggr)
=\sum_{i=1}^nS_{x_i\cdot f}(\delta(a_i))
= \sum_{i=1}^n\delta_{x_i\cdot f}(a_i).
\end{align*}
\end{proof}

Note that
by \cite[Corollary 1.5]{qui:fullred} the converse of the proposition above is
also true, \ie, $\delta$ is nondegenerate if and only if
$\spn\{\delta_f(a)\mid a\in A, f\in A(G)\}$ is dense in $A$.

\section
{Covariant representations and crossed products}

We are now going to introduce covariant representations of
coactions. Of course, if $G$ is abelian, these should coincide with
the covariant representations of the corresponding action
$(A,\widehat G,\alpha)$ as discussed in \secref{sec-coaction}.
So let us start our discussion with a covariant homomorphism
$(\pi, V)$ of an action $(A,\widehat G,\alpha)$ into $M(D)$,
where
$G$ an abelian group. Identifying $C^*(G)$ with $C_0(\widehat G)$ 
via Fourier transform
we may view the integrated form of $V$ as a homomorphism
$\mu\:C_0(\widehat G)\to M(D)$. Thus we may expect covariant
homomorphisms of coactions $(A,G,\delta)$ to be pairs
of nondegenerate maps $(\pi,\mu)$ with $\pi\:A\to M(D), 
\mu\:C_0(\widehat G)\to M(D)$.

Since $\delta^{\alpha}\:A\to C_b(\what G,A)\subseteq
M(A\otimes C_0(\widehat{G}))$ is given by the formula
$\delta^{\alpha}(a)(\chi)=\alpha_{\chi}(a)$, the covariance condition
on $(\pi, V)$ can be expressed as
\begin{equation}
\label{eq-cov}
(\pi\otimes \id)(\delta(a))= \widetilde{V}(\pi(a)\otimes 1)\widetilde{V}^*
\end{equation}
in $M(D\otimes C_0(\widehat{G}))$, where $\widetilde V\in
UM(D\otimes C_0(\widehat{G}))$
is given by the strictly continuous function
$\chi\mapsto V_\chi$.
Since the canonical embedding of $\widehat G$ into $M(C^*(\widehat
G))\cong C_b(G)$ maps a character $\chi\in \widehat G$ to the function
$s\mapsto \chi(s)$, and since the Fourier transform $\mathcal
F\: C^*(G)\to C_0(\widehat G)$ maps $s\in  M(C^*(G))$ to the function
$\chi\mapsto \chi(s)$ in $C_b(\widehat{G})$, it follows that
\[
\widetilde V=  (\mu\otimes \mathcal F)(w_G),
\]
where $w_G$ denotes the function $s\mapsto u(s)$ in $UM(C_0(G)\otimes
C^*(G))$.  Applying the
inverse Fourier transform (or more precisely $\id_D\otimes \mathcal F^{-1}$)
to both sides of \eqeqref{eq-cov}
we see that the covariance condition for $(\pi,V)$ translates into
the condition
\[
(\pi\otimes\id_G)\circ\delta^{\alpha}(a)=
(\mu\otimes \id_G)(w_G)(\pi(a)\otimes 1)(\mu\otimes\id_G)(w_G)^*
\]
in $M(D\otimes C^*(G))$.
Thus we are led to the following definition of covariant representations for
coactions.

\begin{defn}
\label{def-covariant}
Let $(A,G,\delta)$ be a coaction and let $w_G\in UM(C_0(G)\otimes C^*(G))$
be given by the map $s\mapsto u(s)\: G\to UM(C^*(G))$.
A \emph{covariant homomorphism} of $(A,G,\delta)$ into $M(D)$
is a pair $(\pi,\mu)$ of nondegenerate homomorphisms
of $(A,C_0(G))$ into $M(D)$ satisfying the covariance condition
\begin{equation}\label{eq-covariance-coact}
(\pi\otimes\id_G)\circ\delta(a)=
(\mu\otimes \id_G)(w_G)(\pi(a)\otimes 1)(\mu\otimes\id_G)(w_G)^*.
\end{equation}
If $D=\K(\H)$ for some Hilbert space $\H$, then $(\pi,\mu)$ is called
a \emph{covariant representation} of $(A,G,\delta)$ on $\H$.
\end{defn}

\begin{rem}\label{rem-coactdegenerate}
     It is sometimes useful to work with pairs $(\pi,\mu)$
     in which the homomorphism $\pi\:A\to M(D)$ is degenerate,
     but the pair $(\pi,\mu)$ satisfies all other
     conditions on a covariant representation (including
     nondegeneracy of $\mu\:C_0(G)\to M(D)$).
     We shall call such a pair a {\em degenerate} covariant
     representation of $(A,G,\delta)$.
     Note that the covariance condition
     \eqref{eq-covariance-coact} makes sense if $\pi$
     is degenerate, since $\delta(A)\in M_G(A\otimes C^*(G))$,
     and $\pi\otimes\id_G$ extends uniquely to
     $M_G(A\otimes C^*(G))$ by \propref{prop-Cmultextend}.

     However, if not explicitly
     stated otherwise, all covariant representations 
in this work are
     assumed to be nondegenerate as in the definition.
\end{rem}

When working with covariant representations, it is often necessary
to recover the nondegenerate homomorphism $\mu\:C_0(G)\to M(A)$
from the unitary $(\mu\otimes\id_G)(w_G)\in UM(A\otimes C^*(G))$
and to know some further properties of this unitary.
For notation: if $m\in M(A\otimes C)$ for some $C^*$-algebra $C$
then $m_{12}, m_{13}\in M(A\otimes C\otimes C)$ are defined by
$m_{12}=m\otimes 1$ and
$m_{13}=(\id_A\otimes\Sigma)(m\otimes 1)$,
where $\Sigma\:C\otimes C\to C\otimes C$ denotes the flip isomorphism
$c\otimes d\mapsto d\otimes c$.

\begin{prop}[{\cf~\cite[Lemma A1, Lemma A2]{qr:induced}}]
\label{slice-unitary}
Let $A$ be a $C^*$-algebra and let $\mu\:C_0(G)\to M(A)$ be a nondegenerate
homomorphism.
Then
\[
S_f\big( (\mu\otimes \id_G)(w_G) \big)=\mu(f)
\righttext{for all} f \in B(G) \subseteq C_b(G).
\]
Moreover, $w=(\mu\otimes \id_G)(w_G)$
satisfies the equation $w_{12}w_{13}=(\mu\otimes \delta_G)(w_G)=
(\id_A\otimes \delta_G)(w)$.
\end{prop}

\begin{proof} First, we have $S_f(w_G)=f$ for all $f\in B(G)$.
To see this let $s\in G$ and let $\chi_s\:C_b(G)\to \bbC$ denote evaluation
at $s$. Using \lemref{lem-slice} we get
\[
\chi_s(S_f(w_G))= S_f\big((\chi_s\otimes\id_G)(w_G)\big)=
S_f(1_{\CC}\otimes u(s))=f(s).
\]
Again using \lemref{lem-slice} this implies
\[
S_f \bigl( (\mu\otimes \id_G)(w_G) \bigr)=\mu(S_f(w_G))=\mu(f).
\]
 From the definition of $w_G$ it follows that
$(w_G)_{12}(w_G)_{13}\in M(C_0(G)\otimes C^*(G)\otimes C^*(G))$
is given by the map $s\mapsto u(s)\otimes u(s)$, which clearly
coincides with $(\id_{C_0(G)}\otimes\delta_G)(w_G)$.
Applying $\mu\otimes\id_G\otimes\id_G$ to both sides gives
$w_{12}w_{13}=(\mu\otimes\delta_G)(w_G)=
(\id_A\otimes \delta_G)(w)$ for $w=(\mu\otimes\id_G)(w_G)$.
\end{proof}

\begin{rem}
\label{rem-slice}
In fact, there is a certain converse to the above result:
If $w\in UM(A\otimes C^*(G))$ satisfies
$w_{12}w_{13}=(\id_A\otimes \delta_G)(w)$
in $M(A\otimes C^*(G)\otimes C^*(G))$, then
$w=(\mu\otimes\id_G)(w_G)$ for some nondegenerate homomorphism
$\mu\:C_0(G)\to M(A)$ (see \cite[Lemma A.1]{qr:induced}).
Of course, the above result then implies that
$\mu(f)=S_f(w)$ for all $f\in B(G)\subseteq C_b(G)$.
\end{rem}

\begin{prop}[{\cf~\cite[Lemma 2.10]{rae:rep}}]
\label{cov-image}
Let $(\pi,\mu)$ be a \textup(possibly degenerate\textup)
covariant homomorphism of $(A,G,\delta)$ into
$M(D)$ for some $C^*$-algebra $D$. Then
\[
C^*(\pi,\mu)\deq \overline{\pi(A)\mu(C_0(G))}
\]
is a $C^*$-algebra.
\end{prop}

\begin{proof} The result will follow as soon as it is clear
that an element $\mu(f)\pi(a)$ can be approximated
in norm by a finite sum $\sum_i \pi(a_i)\mu(f_i)$, where
$f,f_i\in A(G)$, $a, a_i\in A$.
So let $f\in A(G)$, and factor $f=g\cdot x$ for some
$x\in C^*(G)$ and $g\in A(G)$. The covariance identity together with
\lemref{lem-slice} and \propref{slice-unitary}
gives
\begin{align*}
\mu(f)\pi(a)
&=S_f\Bigl(\bigl( (\mu\otimes\id_G)(w_G) \bigr)(\pi(a)\otimes 1)\Bigr)
=S_f\Bigl(\bigl(\pi\otimes \id_G(\delta(a))\bigr)
\bigl( (\mu\otimes\id_G)(w_G) \bigr)\Bigr)\\
&=S_g\Bigl(\bigl(\pi\otimes\id_G((1\otimes x)\delta(a))\bigr)
\bigl( (\mu\otimes\id_G)(w_G) \bigr)\Bigr).
\end{align*}
Approximating $(1\otimes x)\delta(a)\in A\otimes C^*(G)$
by a finite sum $\sum_ia_i\otimes x_i$ then implies
\begin{align*}
\mu(f)\pi(a)&\approx
S_g\biggl(\Bigl(\sum_i\pi(a_i)\otimes x_i\Bigr)
\bigl( (\mu\otimes\id_G)(w_G) \bigr)\biggr)\\
&=\sum_i\pi(a_i)S_{g\cdot x_i}\bigl( (\mu\otimes\id_G)(w_G) \bigr)
=\sum_i\pi(a_i)\mu(g\cdot x_i),
\end{align*}
which completes the proof.
\end{proof}

The following proposition shows that covariant homomorphisms
do exist.

\begin{prop}[{\cf~\cite[Proposition 2.6]{rae:rep}}]
\label{prop-ind}
If $\pi$ is a nondegenerate homomorphism of $A$ into $M(D)$, then
$((\pi\otimes \lambda)\circ\delta, 1\otimes M)$ is a covariant homomorphism
of $(A,G,\delta)$ into $M(D\otimes \K(L^2(G)))$.
\end{prop}

\begin{proof}
We first establish the identity
\begin{equation}
\label{eq-ind}
(\lambda\otimes\id_G)\circ \delta_G=
\Ad(M\otimes \id_G)(w_G)\circ(\lambda \otimes 1).
\end{equation}
Applied to $s \in G$ (regarded as an element of $UM(C^*(G))$)
the left hand side becomes
$\lambda_s \otimes s$. To compute the right hand side at $s$,
assume that
$C^*(G)$ is represented faithfully on a Hilbert space $\H$.
Then for $\xi\in L^2(G,\H)$ and $t\in G$ we get
\begin{align*}
&\bigl( (M \otimes \id_G)(w_G) (\lambda_s \otimes 1)
(M \otimes \id_G)(w_G)^* \xi \bigr)(t)
= t \bigl( (\lambda_s \otimes 1)
(M \otimes \id_G)(w_G)^* \xi \bigr)(t)
\\&\quad= t \bigl( (M \otimes \id_G)(w_G)^* \xi \bigr)(s^{-1}t)
= t (t^{-1}s) \xi(s^{-1}t)
= \bigl( (\lambda_s \otimes s) \xi \bigr)(t).
\end{align*}
Thus (\ref{eq-ind}) follows from integration.
Together with the coaction identity this implies
\begin{align*}
&\bigl( (\pi \otimes \lambda) \circ \delta \otimes \id_G \bigr)
\circ \delta(a)
= (\pi \otimes \lambda \otimes \id_G) \circ
(\delta \otimes \id_G) \circ \delta(a)
\\&\quad= (\pi \otimes \lambda \otimes \id_G) \circ
(\id_A \otimes \delta_G) \circ \delta(a)
\\&\quad= (\pi \otimes \id_{\K} \otimes \id_G) \circ
\bigl( \id_A \otimes (\lambda \otimes \id_G) \circ \delta_G \bigr)
\circ \delta(a)
\\&\quad= (\pi \otimes \id_{\K }\otimes \id_G)
\Bigl( (1 \otimes M \otimes \id_G)(w_G)
\bigl( (\id_A \otimes \lambda)(\delta(a)) \otimes 1 \bigr)
(1 \otimes M \otimes \id_G)(w_G)^* \Bigr)
\\&\quad= (1 \otimes M \otimes \id_G)(w_G)
\bigl( (\pi \otimes \lambda)(\delta(a)) \otimes 1 \bigr)
(1 \otimes M \otimes \id_G)(w_G)^*,
\end{align*}
which completes the proof.
\end{proof}

\begin{defn}

\label{def-regular}
Let $\pi\:A\to M(D)$ be a nondegenerate homomorphism.
Then the covariant homomorphism
$\Ind\pi\deq \big((\pi\otimes \lambda)\circ\delta, 1\otimes M\big)$
of $(A,G,\delta)$ into $M(D\otimes\K(L^2(G)))$
is called
the homomorphism of $(A,G,\alpha)$ \emph{induced from $\pi$}.
If $\pi$ is
faithful, then $\Ind\pi$ is called a \emph{regular}
homomorphism of
$(A,G,\alpha)$.
\end{defn}

Somewhat surprisingly, it turns out that for coactions
there is no difference between full and reduced crossed products.
In fact we are now going to define the crossed product by a
coaction as the $C^*$-algebra generated by a certain
regular representation (similarly to the definition of the
reduced crossed product $A\times_{\alpha,r}G$ as given in
\defnref{def-red crossed}) and then show that this
crossed product enjoys universal properties similar to those
enjoyed by the full crossed product $A\times_{\alpha}G$ of
an action $\alpha\:G\to \Aut(A)$.
However, if we view the crossed product
$A\times_{\delta}G$ as a skew tensor product of $A$ with $C_0(G)$,
this might be less surprising in view of the fact that
any commutative $C^*$-algebra is nuclear.

\begin{defn}
\label{def-crossed}
Let $(A,G,\delta)$ be a coaction and let
$(j_A,j_G)=\big((\id_A\otimes\lambda)\circ\delta, 1\otimes M\big)$
be the regular homomorphism induced by $\id_A\:A\to A$.
Then $A\times_{\delta}G\deq C^*(j_A,j_G)$
  is called the
\emph{crossed product} of the coaction $(A,G,\delta)$.
The maps $(j_A, j_G)$, viewed as maps of $(A,C_0(G))$ into
$M(A\times_{\delta}G)$, are called the \emph{canonical maps}
of $(A, C_0(G))$ into $M(A\times_{\delta}G)$.
\end{defn}

\begin{rem}\label{rem-Kmultcoact}
    By definition, the crossed product $A\times_{\delta}G$
    is a subalgebra of $M(A\otimes \K(L^2(G)))$.
    It is important to observe that it actually lies in the
    $\K(L^2(G))$-multiplier algebra $M_{\K(L^2(G))}(A\otimes
    \K(L^2(G)))$. To see this just observe that
    $$j_A(A)=(\id_A\otimes \lambda)\circ \delta(A)
    \subseteq (\id_A\otimes \lambda)\big(M_G(A\otimes C^*(G))\big)
    \subseteq M_{\K(L^2(G))}(A\otimes
    \K(L^2(G))),$$
    where the last inclusion follows from \propref{prop-Cmultextend}.
    Since $j_G(C_0(G))=1\otimes M(C_0(G))\subseteq 1\otimes
    M(\K(L^2(G)))$, it follows from item (iii) of \propref{prop-Cstrict}
    that $A\times_\delta G=\overline{j_A(A)j_G(C_0(G))}\subseteq
    M_{\K(L^2(G))}(A\otimes \K(L^2(G)))$.
\end{rem}

\begin{thm}[{\cf~\cite[Theorem 4.1]{rae:rep}}]
\label{thm-crossed}
The triple $(A\times_{\delta}G, j_A, j_G)$ satisfies the following
universal property: if $(\pi,\mu)$ is any
covariant homomorphism
of $(A,G,\delta)$ into $M(D)$ for some $C^*$-algebra $D$, then
there exists a unique nondegenerate homomorphism
$\pi\times\mu\:A\times_{\delta}G\to M(D)$ such that
$(\pi\times\mu)\circ j_A=\pi$ and $(\pi\times \mu)\circ j_G=\mu$.

Moreover, if $(\pi,\mu)$ is a degenerate covariant
homomorphism, then there exists a unique
\textup(degenerate\textup) homomorphism
$\pi\times\mu\:A\times_{\delta}G\to M(D)$
such that $\pi\times\mu(j_A(a)j_G(f))=\pi(a)\mu(f)$.
\end{thm}

\begin{defn}
\label{def-integrated}
The homomorphism
  $\pi\times\mu$ of \thmref{thm-crossed}
is called the \emph{integrated form} of the covariant homomorphism
$(\pi,\mu)$ of $(A,G,\delta)$.
\end{defn}

\begin{rem}
\label{rem-co-universal0}
(1) It is not hard to check that if $\rho\:A\times_{\delta}G\to M(D)$
is a nondegenerate homomorphism of $A\times_{\delta}G$, then
$(\rho\circ j_A, \rho\circ j_G)$ is a covariant homomorphism. We then
have $\rho=(\rho\circ j_A)\times (\rho\circ j_G)$.
Thus, $(\pi,\mu)\mapsto \pi\times\mu$ gives a one-to-one correspondence
between the (nondegenerate) covariant homomorphisms of $(A,G,\delta)$ and the
nondegenerate homomorphisms of $A\times_{\delta}G$.

(2) As for actions, one could alternatively \emph{define} the crossed product
via universal properties, \ie, as a triple $(C, l_A, l_G)$ satisfying
\begin{enumerate}
\item $(l_A, l_G)$ is a covariant homomorphism
of $(A,G,\delta)$ and $C=C^*(l_A,l_G)$; and
\item  for every covariant homomorphism $(\pi, \mu)$ of $(A,G,\delta)$
into $M(D)$ there exists a unique nondegenerate homomorphism
$\pi\times \mu\:C\to M(D)$ such that
$\pi\times\mu\circ l_A=\pi$ and $\pi\times\mu\circ l_G=\mu$.
\end{enumerate}
Of course, for any such triple $(C, l_A, l_G)$, the $C^*$-algebra
$C$ is isomorphic to $A\times_{\delta}G$
via the integrated form
$l_A\times l_G\:A\times_{\delta}G\to C$ with inverse
map $j_A\times j_G\:C\to A\times_{\delta}G$, where $j_A\times j_G$
is the map associated
  to $(j_A, j_G)$ by
  (ii).

(3) If $\Ind\pi=((\pi\otimes\lambda)\circ \delta, 1\otimes M)$ is any
regular homomorphism of $(A,G,\delta)$
  induced by the
faithful homomorphism  $\pi\:A\to M(D)$, then its integrated form
(also denoted $\Ind\pi$) is a faithful homomorphism
of $A\times_{\delta}G$. For this one observes that
\[
\Ind\pi= (\pi\otimes \id_{\K})\circ \Ind\id_A= (\pi\otimes 
\id_{\K})\circ (j_A\times j_G),
\]
which is faithful since $\pi\otimes\id_\K$ is faithful.

(4) For abelian $G$, the reader can check without too much difficulty
(using the Plancherel isomorphism $L^2(G)\to L^2(\widehat G)$)
that regular representations of $(A,\widehat G,\alpha)$
correspond to regular representations of $(A, G,\delta^{\alpha})$,
which implies that
$A\times_{\alpha}\widehat G\cong A\times_{\delta^{\alpha}}G$.
\end{rem}

The proof of \thmref{thm-crossed} we give here is based on the
proof of \cite[Theorem 3.7]{lprs}.
We need the following lemma.

\begin{lem}
\label{lem-identities}
Let $(\pi,\mu)$ be a \textup(possibly degenerate\textup) covariant homomorphism of
$(A,G,\delta)$ into $M(E)$ for some $C^*$-algebra $E$,
let $W=(\mu\otimes\lambda)(w_G)\in M(E\otimes \K(L^2(G)))$
and let $\eps\:C_0(G)\to M(C_0(G)\otimes C_0(G))$ be the comultiplication
on $C_0(G)$, \ie, $\eps(f)(s,t)=f(st)$.
Then
\begin{equation}
\label{eq-redcov}
\Ad W^*\circ (\pi\otimes\lambda)\circ \delta=\pi \otimes 1
\midtext{and}
\Ad W^*\circ (1\otimes M)=(\mu\otimes M)\circ \eps.
\end{equation}
\end{lem}
\begin{proof}
The first equation follows from applying
$\id_E\otimes \lambda$ to both sides of the covariance condition
$\pi\otimes 1=\Ad \bigl( (\mu\otimes\id_G)(w_G) \bigr)^*
\circ (\pi\otimes \id_G)\circ \delta$.
For the proof of the second equation we first show that
\[
(M\otimes\lambda(w_G))^*(1\otimes M(f))(M\otimes\lambda(w_G))=(M\otimes
M)\big(\eps(f)\big)
\]
in $\B(L^2(G\times G))$.
For this write $W_G\deq M\otimes\lambda(w_G)$.
Then  $(W_G\xi)(s,t)=\xi(s, s^{-1}t)$, from which it follows that
\begin{align*}
\big(W_G^*&(1\otimes M(f))W_G\xi\big)(s,t)
=\big((1\otimes M(f))W_G\xi\big)(s,st)= f(st)\big(W_G\xi\big)(s,st)\\
&=f(st)\xi(s,t)=\big((M\otimes M)(\eps(f))\xi\big)(s,t).
\end{align*}
Since $M$ is faithful on $C_0(G)$ we get the equation
\[
(\id\otimes\lambda(w_G))^*(1\otimes M(f))(\id\otimes\lambda(w_G))=(\id\otimes
M)\big(\eps(f)\big)
\]
for all $f\in C_0(G)$, and applying $\mu\otimes \id_{\K}$
to
  both sides of
this equation gives the desired result.
\end{proof}

\begin{proof}[Proof of \thmref{thm-crossed}]
Let $(\pi,\mu)$ be a (possibly degenerate)
covariant homomorphism of $(A,G,\delta)$ into $M(D)$
and let $W=(\mu\otimes\lambda)(w_G)\in UM(D\otimes \K(L^2(G)))$. 
Further let $\phi\:A\times_{\delta}G\to M(D\otimes\K(L^2(G)))$
denote the restriction to $A\times_\delta G$ of
the homomorphism
\[
\Ad W^*\circ (\pi\otimes\id_{\K})\:M_{\K(L^2(G))}(A\otimes \K(L^2(G)))\to
M(D\otimes \K(L^2(G))).
\]
Using \lemref{lem-identities}, we get
\begin{align*}
\phi(j_A(a)j_G(f))&=
\phi\big((\id\otimes \lambda)(\delta(a))(1\otimes M_f)\big)\\
&=W^*((\pi\otimes \lambda)(\delta(a))WW^*(1\otimes M_f)W
=(\pi(a)\otimes 1)(\mu\otimes M)\big(\eps(f)\big)
\end{align*}
for $a\in A$, $f\in C_0(G)$. It follows
that $\phi(A\times_{\delta}G)$ lies in the image of
$M(D\otimes C_0(G))$ under the homomorphism
$\id_D\otimes M\: D\otimes C_0(G)\to M(D\otimes \K(L^2(G)))$.
Since $\id_D\otimes M$ is faithful, we therefore obtain a homomorphism
$\psi\:A\times_{\delta}G\to M(D\otimes C_0(G))$ satisfying
\[
\psi\big(j_A(a)j_G(f)\big)=
(\pi(a)\otimes 1)(\mu\otimes\id)(\eps(f)).
\]
Let $\chi_e\:C_0(G)\to \mathbb C$ denote evaluation at the identity $e\in G$
and define
\[
\pi\times \mu\deq  (\id_D\otimes \chi_e)\circ \psi\:A\times_{\delta}G
\to M(D).
\]
We then have
\begin{align*}
   \pi\times \mu(j_A(a)j_G(f))&=
   (\id_D\otimes \chi_e)\big((\pi(a)\otimes
   1)(\mu\otimes\id)(\eps(f))\big)\\
   &=(\id_D\otimes \chi_e)\big((\pi(a)\otimes 1)\big)
   (\id_D\otimes \chi_e)\big((\mu\otimes\id)(\eps(f))\big)\\
   &=\pi(a)(\mu\otimes\chi_e)(\eps(f))=\pi(a)\mu(f)
\end{align*}
for all $a\in A$, $f\in C_0(G)$, where the identity
$(\mu\otimes\chi_e)(\eps(f))\big)=\mu(f)$
follows from the equation
$(\id_{C_0(G)}\otimes\chi_e)(\eps(f))=f$ for $f\in C_0(G)$.
If $\pi$ is degenerate, then so is $\pi\times\mu$, since
$$\pi\times \mu(A\times_{\delta}G)D=\overline{\pi(A)\mu(C_0(G))}D
=\pi(A)(\mu(C_0(G))D)=\pi(A)D\neq D.$$
If $\pi$ is nondegenerate, then
as a composition of nondegenerate homomorphisms, $\pi\times\mu$
is nondegenerate, too. We then have
\[
(\pi\times \mu)\circ j_A(a)=(\id_D\otimes\chi_e)\big(\pi(a)\otimes 1)=\pi(a)
\]
for $a\in A$, and
\[
(\pi \times \mu) \circ j_G(f)
= (\mu \otimes \chi_e)(\eps(f))
= \mu(f),
\]
for all $f\in C_0(G)$.
\end{proof}

\begin{defn}\label{defn-equivariantcoact}
     Let $(A,\delta)$ and $(B,\epsilon)$ be coactions of $G$.
     A (possibly degenerate) homomorphism $\phi\:A\to M(B)$
     is called {\em $\delta-\epsilon$ equivariant} if
     $(\varphi\otimes \id_G)\circ \delta=\epsilon\circ \varphi$.
     (Note that the composition on the left hand side is well-defined
     by \propref{prop-Cmultextend} since $\delta(A)\subseteq
     M_G(A\otimes C^*(G))$.)
\end{defn}

\begin{lem}\label{lem-equivcrossed}
      Let $(A,\delta)$ and $(B,\epsilon)$ be coactions of $G$
      and let $\phi\:A\to M(B)$ be a $\delta-\epsilon$
     equivariant homomorphism. Then
     there is a well-defined homomorphism
     $$\phi\times G\deq (j_B\circ \phi)\times j_G\: A\times_{\delta}G
     \to M(B\times_{\epsilon}G),$$
     where $(j_B, j_G)\:(B, C_0(G))\to M(B\times_{\epsilon}G)$
     denotes the canonical pair of maps.
     Moreover, $\phi\times G$ is nondegenerate if \textup(and only if\textup) $\phi$
     is nondegenerate, and $\phi\times G$ is faithful
     if $\phi$ is faithful.
\end{lem}
\begin{proof} The covariance condition
     for $(j_B, j_G)$
     (except nondegeneracy) extends
     to $(j_B\circ \phi, j_G)$. Thus
     $\phi\times G$ is well-defined.
     If $\phi$ is nondegenerate, then so is $j_B\circ\phi$
     and hence $\phi\times G=(j_B\circ\phi)\times j_G$.
     So suppose finally that $\phi$ is faithful,
     and consider the composition
     $$
     \begin{CD}
	A\times_{\delta}G @>j_A\times j_G>>
	M_{\K(L^2(G))}(A\otimes\K(L^2(G)))
	@>\phi\otimes \id_{\K}>> M(B\otimes \K(L^2(G))).
     \end{CD}
     $$
     Using the definition of $j_A$, $j_B$ and $j_G$,
     and the $\delta-\epsilon$ equivariance of $\phi$
     we compute for all $a\in A$ and $f\in C_0(G)$:
     \begin{align*}
      \phi\otimes\id_{\K}(j_A(a)j_G(f))&=
      \phi\otimes\id_{\K}\left(\big((\id_A\otimes\lambda)\circ\delta(a)\big)
      (1\otimes M(f))\right)\\
      &=\big((\id_B\otimes\lambda)\circ (\phi\otimes\id_G)\circ
      \delta(a)\big)(1\otimes M(f))\\
      &=\big((\id_B\otimes\lambda)\circ
      \epsilon(\varphi(a))\big)(1\otimes M(f))\\
      &=j_B\circ \phi(a)j_G(f),
      \end{align*}
      which implies that $\phi\times G$ coincides with the
      above composition of maps. Thus, if $\phi$ is faithful,
      the composition $(\phi\otimes\id_\K)\circ (j_A\times j_G)$
      is faithful by \propref{prop-Cmultextend}.
\end{proof}


\section{Dual actions and decomposition coactions}
\label{sec-dec}


In \exref{ex-dualcoaction} we saw that for any action
$(A,G,\alpha)$ there is a dual coaction
$\widehat\alpha$ of $G$ on $A\times_{\alpha}G$.
Obversely, for any coaction $(A,G,\delta)$ there is a
dual action of $G$ on $A\times_{\delta}G$.

\begin{defn}
\label{def-dualact}
Let $A\times_{\delta}G=C^*(j_A, j_G)$ be the crossed product
of the coaction $(A,G,\delta)$. For $s\in G$ let $\sigma_s$ denote
right translation on $C_0(G)$, \ie, $\sigma_s(f)(t)=f(ts)$.
Then $\widehat{\delta}\:G\to \Aut(A\times_{\delta}G)$
defined by $\widehat{\delta}_s\deq j_A\times (j_G\circ \sigma_s)$
is called the \emph{dual action} of $G$ on $A\times_{\delta}G$.
\end{defn}

\begin{rem}
\label{rem-dualaction}
Note that $\widehat\delta_s$ is given on a typical element of the
form $j_A(a)j_G(f)$ by
\[
\widehat\delta_s(j_A(a)j_G(f))=j_A(a)j_G(\sigma_s(f)).
\]
Since right translation on $C_0(G)$ is continuous it follows
that $\widehat\delta\:G\to \Aut(A\times_{\delta}G)$ is strongly
continuous.
\end{rem}


In the main body of the paper we need to work with a certain
canonical coaction, the \emph{decomposition coaction}
$\delta\dec$ of $G$ on $A\times_{\delta|}G/N$, where
$\delta$ is a given coaction of $G$ on $A$ and $N$ is a closed normal
subgroup of $G$. 
(The reason for this terminology
is that,  by \cite[Theorem~3.1]{pr:twisted},
$A \times_\delta G$ can be decomposed into a twisted crossed
product of $A \times_{\delta|} G/N$ by $\delta\dec$.)
We are now going to describe this coaction.

\begin{lem}
\label{decom-coact}
Let $(A,G,\delta)$ and $N$ be as above. The formula
\begin{equation}\label{decom-eq}
\delta\dec(j_A(a)j_{G/N}(f))
= (j_A \otimes \id) \circ \delta(a)(j_{G/N}(f) \otimes 1)
\quad\text{ for } a\in A, f\in C_0(G/N)
\end{equation}
defines a  coaction of $G$ on $A
\times_{\delta|} G/N$, which is nondegenerate if $\delta$ is.
\end{lem}

\begin{proof}
First, it is easy to see that $(j_A \otimes \id) \circ
\delta$ and $j_{G/N} \otimes 1$ are nondegenerate homomorphisms of $A$
and $C_0(G)$, respectively, into $M((A \times_{\delta|} G/N) \otimes
C^*(G))$. We show that $((j_A \otimes \id) \circ
\delta, j_{G/N} \otimes 1)$ is a covariant pair: for $a \in A$ we have
\begin{align*}
&\ad (j_{G/N} \otimes 1 \otimes \id)(w_{G/N})
\bigl((j_A \otimes \id) \circ \delta(a) \otimes 1\bigr)
\\&\quad
= \ad (\id \otimes \sigma)
\bigl((j_{G/N} \otimes \id)(w_{G/N}) \otimes 1\bigr)
\bigl((j_A \otimes \id) \circ \delta(a) \otimes 1\bigr)
\\&\quad
= (\id \otimes \sigma) \circ
\ad \bigl((j_{G/N} \otimes \id)(w_{G/N}) \otimes 1\bigr)
\circ (\id \otimes \sigma)
\bigl((j_A \otimes \id) \circ \delta(a) \otimes 1\bigr)
\\&\quad
= (\id \otimes \sigma) \circ
\bigl( \ad (j_{G/N} \otimes \id)(w_{G/N}) \otimes \id \bigr)
\circ (j_A \otimes 1 \otimes \id) \circ \delta(a)
\\&\quad
= (\id \otimes \sigma) \circ
\bigl( \ad (j_{G/N} \otimes \id)(w_{G/N}) \circ (j_A \otimes 1)
\otimes \id \bigr) \circ \delta(a)
\\&\quad
=  (\id \otimes \sigma) \circ
\bigl( (j_A \otimes \id) \circ \delta| \otimes \id \bigr)
\circ \delta(a)
\\&\quad
= (j_A \otimes \id \otimes \id) \circ (\id \otimes \sigma)
\circ (\delta| \otimes \id) \circ \delta(a)
\\&\quad
=  (j_A \otimes \id \otimes \id) \circ
(\delta \otimes \id) \circ \delta|(a)
\\&\quad
= \bigl( (j_A \otimes \id) \circ \delta \otimes \id \bigr)
\circ \delta|(a).
\end{align*}

Thus there is a nondegenerate homomorphism $\delta\dec \: A
\times_{\delta|} G/N \to M((A \times_{\delta|} G/N) \otimes C^*(G))$.
This homomorphism is injective, because (letting $1_G$ denote the
trivial character on $G$)
$(\id \otimes 1_G) \circ (j_A \otimes \id) \circ \delta
= j_A \circ (\id \otimes 1_G) \circ \delta
= j_A$
and
$(\id \otimes 1_G) \circ (j_{G/N} \otimes 1)
= j_{G/N}$.
The coaction identity holds because
\begin{align*}
(\delta\dec \otimes \id) \circ \delta\dec \circ j_A
&= (\delta\dec \otimes \id) \circ (j_A \otimes \id) \circ \delta
= (\delta\dec \circ j_A \otimes \id) \circ \delta\\
&= \bigl((j_A \otimes \id) \circ \delta \otimes \id\bigr) \circ \delta
= (j_A \otimes \id \otimes \id) \circ (\delta \otimes \id) \circ \delta\\
&= (j_A \otimes \id \otimes \id) \circ (\id \otimes \delta_G) \circ \delta
= (\id \otimes \delta_G) \circ (j_A \otimes \id) \circ \delta
= (\id \otimes \delta_G) \circ \delta\dec \circ j_A
\end{align*}
and
\begin{align*}
&(\delta\dec \otimes \id) \circ \delta\dec \circ j_{G/N}
= (\delta\dec \otimes \id) \circ (j_{G/N} \otimes 1)
= (\delta\dec \circ j_{G/N})\otimes 1
\\&\quad
= j_{G/N} \otimes 1 \otimes 1
= (\id \otimes \delta_G) \circ (j_{G/N} \otimes 1)
= (\id \otimes \delta_G) \circ \delta\dec \circ j_{G/N}.
\end{align*}
Also, for $a \in A$, $f \in C_0(G/N)$, and $c \in C^*(G)$ we have
\begin{align*}
&\delta\dec\bigl( j_A(a) j_{G/N}(f) \bigr)
(1 \otimes c)
= (j_A \otimes \id) \circ \delta(a)
( j_{G/N}(f) \otimes c )
\\&\quad
= (j_A \otimes \id)
\bigl( \delta(a) (1 \otimes c) \bigr)
(j_{G/N}(f) \otimes 1)
\in \bigl( j_A(A) \otimes C^*(G) \bigr)
(j_{G/N}(f) \otimes 1)
\\&\quad
\subseteq (A \times_{\delta|} G/N) \otimes C^*(G).
\end{align*}
Hence $\delta\dec$ is a coaction.

Now assume $\delta$ is nondegenerate.
We must show $\delta\dec$ is nondegenerate, too. If
$a \in A$, $f \in C_0(G/N)$, and $c \in C^*(G)$ we have
\begin{align*}
\delta\dec\bigl( j_A(a) j_{G/N}(f) \bigr) (1 \otimes c)
&= (j_A \otimes \id) \circ \delta(a) (j_{G/N}(f) \otimes 1)
(1 \otimes c)
\\&= (j_A \otimes \id) \circ \delta(a) (1 \otimes c)
(j_{G/N}(f) \otimes 1).
\end{align*}
The latter elements densely span
$(j_A \otimes \id)\bigl( A \otimes C^*(G/N) \bigr)
\bigl( j_{G/N}(C_0(G)) \otimes 1 \bigr)$,
which of course densely spans $(A \times_{\delta|} G/N) \otimes
C^*(G)$.
\end{proof}

\section{Normal coactions and normalizations}
\label{sec-normal}

Normal coactions and normalizations of coactions will play
a fundamental r\^ole in this work. We start the discussion with:

\begin{defn}
\label{def-normal}
A coaction $(A,G,\delta)$
is \emph{normal} if
  $j_A\:A\to M(A\times_{\delta}G)$ is injective.
\end{defn}

\begin{rem}
\label{rem-normal}
Since for every covariant homomorphism $(\pi,\mu)$ of
$(A,G,\delta)$
the homomorphism $\pi$ of $A$
factors through $j_A$ (which follows from $(\pi\times\mu)\circ j_A=\pi$)
we see that a coaction $\delta$ is normal if and only
if it has a covariant homomorphism (or representation)
$(\pi,\mu)$ with $\pi$ faithful.
\end{rem}


\begin{ex}
\label{decom-normal}
Recall from the preceding section that if
  $(A,G,\delta)$
is a   coaction and
$N$ is a closed normal
subgroup of $G$, then
the decomposition coaction $\delta\dec$ of $G$ on $A
\times_{\delta|} G/N$  is defined on the
generators by
\[
\delta\dec(j_A(a)j_{G/N}(f))
= (j_A \otimes \id) \circ \delta(a)(j_{G/N}(f) \otimes 1).
\]
Suppose $\delta$ is normal. We will now show that $\delta\dec$ is
normal, too.
It suffices to show that there is a covariant homomorphism $(\pi,\mu)$
of $(A \times_{\delta|} G/N,G,\delta\dec)$ with $\pi$ faithful. We take
$\pi = j_A^G \times j_G|$ and $\mu = j_G$; as we mention in
\appxref{imprim-chap}, our hypotheses on $\delta$ guarantee that $\pi$
is a faithful nondegenerate homomorphism of $A \times_{\delta|} G/N$
into $M(A \times_\delta G)$. So, it remains to verify that the pair
$(\pi,\mu)$ is covariant: we have
\begin{align*}
&\ad (j_G \otimes \id)(w_G) \circ (j_A \times j_G|) \circ
(j_A \otimes 1)
= \ad (j_G \otimes \id)(w_G) \circ (j_A^{G} \otimes 1)
\\&\quad= (j_A^G \otimes \id) \circ \delta
= \bigl( (j_A^G \times j_G|) \circ j_A^{G/N} \otimes \id \bigr)
\circ \delta
\\&\quad= \bigl( (j_A^G \otimes j_G|) \otimes \id \bigr)
\circ (j_A^{G/N} \otimes \id) \circ \delta
= \bigl( (j_A^G \times j_G|) \otimes \id \bigr)
\circ \delta\dec \circ j_A,
\end{align*}
and
\begin{align*}
&\ad (j_G \otimes \id)(w_G) \circ (j_A \times j_G|) \circ
(j_{G/N} \otimes 1)
= \ad (j_G \otimes \id)(w_G) \circ (j_G| \otimes 1)
\\&\quad= j_G| \otimes 1
= (j_A \times j_G|) \circ j_{G/N} \otimes 1
\\&\quad= \bigl( (j_A \times j_G|) \otimes \id \bigr) \circ
(j_{G/N} \otimes 1)
= \bigl( (j_A \times j_G|) \otimes \id \bigr)
\circ \delta\dec \circ j_{G/N}.
\end{align*}
\end{ex}


In what follows next we show that for every coaction
$\delta$ there exists a canonically-defined normal coaction $\delta^n$
such that $\delta$ and $\delta^n$ have the same representation theory
and crossed products.

\begin{defn}
\label{defn-conj}
Let $(A,\delta)$ and $(B,\eps)$
be coactions of $G$.
If $\phi\:A\to B$ is a $\delta-\epsilon$ equivariant isomorphism,
then we say that $\phi$ is
a \emph{conjugacy} between $\delta$ and $\eps$.
We say that $\delta$ and $\eps$ are \emph{conjugate} if there exists a
conjugacy between $\delta$ and $\eps$.
\end{defn}

\begin{rem}
\label{rem-conj}
It is straightforward to check that conjugacy defines
an equivalence relation, and
that if $\phi\:A\to B$ is a conjugacy, then
the homomorphism $\phi\times G\:A\times_{\delta}G\to
B\times_{\epsilon}G$ of \lemref{lem-equivcrossed}
is an isomorphism with inverse $\phi^{-1}\times G$.
\end{rem}

\begin{lem}[{\cf~\cite[Proposition 2.4]{qui:full}}]
\label{lem-covcoact}
Let $(\pi,\mu)$ be a covariant homomorphism of the coaction
$(A,G,\delta)$
into $M(D)$. Then the map $\delta^{\mu} \: \pi(A) \to M \bigl( \pi(A)
\otimes C^*(G) \bigr)$ defined by
\[
\pi(a) \mapsto
\Ad \bigl( (\mu \otimes \id_G)(w_G) \bigr) \bigl( \pi(a) \otimes 1 \bigr)
\]
is a normal coaction of $G$ on $\pi(A)$ which is nondegenerate if
$\delta$ is  nondegenerate. Moreover,
$\pi\:A\to \pi(A)$ is $\delta-\delta^{\mu}$ equivariant.
\end{lem}
\begin{proof} The covariance condition for $(\pi,\mu)$ implies that
\begin{equation}
\label{eq-delta-mu}
\delta^{\mu}(\pi(a))= (\pi\otimes\id_G)\circ \delta(a).
\end{equation}
This shows that $\delta^{\mu}$ takes values in $M(\pi(A)\otimes C^*(G))$
and
\[
\delta^{\mu}(\pi(A))(1\otimes C^*(G))=\pi\otimes \id_G\big(\delta(A)(1\otimes
C^*(G))\big)\subseteq \pi(A)\otimes C^*(G).
\]
Of course we have equality if $\delta$ is nondegenerate.
The coaction identity for $\delta^{\mu}$ follows from applying
$\pi\otimes\id_G\otimes\id_G$ to
  both sides of the coaction identity of $\delta$.
Finally, observe that $(\id_{\pi(A)},\mu)$ is covariant for $\delta^{\mu}$,
so $\delta^{\mu}$ is normal by \remref{rem-normal}.
\end{proof}

Of course, by identifying $\pi(A)$ with $A/\ker\pi$
we may also view $\delta^{\mu}$ as a coaction on $A/\ker\pi$.

\begin{defn}
\label{defn-normal}
Let $(A,G,\delta)$ be a coaction and let
$A^n=j_A(A)\cong A/\ker j_A$.
Then $\delta^n\deq \delta^{j_G}\:A^n \to M(A^n\otimes C^*(G))$
is called the \emph{normalization} of $\delta$.
\end{defn}

\begin{rem}
\label{rem-normalrep}
It follows from \lemref{lem-covcoact} that $\delta^n$ is nondegenerate if
$\delta$ is nondegenerate. The converse is shown in
\cite[Proposition 2.5]{qui:fullred}.

Since $j_A\:A\to A^n$ is $\delta-\delta^n$ equivariant we see that
if $(\pi,\mu)$ is a covariant homomorphism of $(A^n, G,\delta^n)$,
then $(\pi\circ j_A, \mu)$ is a covariant homomorphism of $(A,G,\delta)$
and $C^*(\pi\circ j_A, \mu)=C^*(\pi,\mu)$.
Conversely, since $\ker j_A\subseteq \ker\rho$ for every
covariant homomorphism $(\rho,\nu)$ of $(A,G,\delta)$ (see \remref{rem-normal})
it follows that every covariant homomorphism of $(A,G,\delta)$ arises this way.
\end{rem}

The following proposition is now completely straightforward.

\begin{prop}[{\cf~\cite[Proposition 2.6]{qui:fullred}}]
\label{prop-normal}
For any coaction $(A,G,\delta)$, the map
\[
j_A \times G \: A\times_{\delta}G \to A^n\times_{\delta^n}G
\]
is a $\widehat{\delta}-\widehat{\delta^n}$ equivariant isomorphism.
\end{prop}

As an interesting corollary we get:

\begin{cor}[{\cf~\cite[Corollary 2.7]{qui:fullred}}]
\label{cor-reg}
Let $(A,G,\delta)$ be a coaction and let $\pi\:A\to M(D)$ be any
nondegenerate homomorphism of $A$ such that $\ker\pi\subseteq \ker j_A$.
Then
\[
\Ind\pi=(\pi\otimes\lambda)\circ \delta\times(1\otimes
M)\:A\times_{\delta}G
\to M(D\otimes \K(L^2(G)))
\]
is faithful.
\end{cor}
\begin{proof}
Let $\phi\:\pi(A)\to j_A(A)$ be such that $\phi\circ \pi=j_A$.
Since $\delta^n\circ j_A=(j_A\otimes\id_G)\circ\delta$ we have
\[
(\phi\circ\pi\otimes\lambda)\circ \delta=
(j_A\otimes \lambda)\circ \delta=(\id_{A^n}\otimes\lambda)\circ
\delta^n\circ j_A.
\]
Thus it follows from the proposition and \thmref{thm-crossed}
that $(\phi\otimes\id)\circ (\pi\otimes\lambda)\circ\delta\times(1\otimes M)
=(j_A\otimes \lambda)\circ \delta\times(1\otimes M)$ is faithful on
$A\times_{\delta}G= A^n\times_{\delta^n}G$.
\end{proof}

We can even strengthen the above result as follows:

\begin{prop}
\label{prop-faithful}
Let $(A,G,\delta)$ be a coaction and let $\pi\:A\to M(D)$ be a
nondegenerate homomorphism such that $\ker(\pi\otimes \lambda)\circ \delta
=\ker j_A$. Then $\Ind\pi\:A\times_{\delta}G\to M(D\otimes \K(L^2(G)))$
is faithful.
\end{prop}
\begin{proof} Let us assume that $D$ is represented faithfully on the Hilbert
space $\H$.
By \corref{cor-reg} we know that
\[
\Ind \big((\pi\otimes\lambda)\circ \delta\big)=
\big((\pi\otimes\lambda\otimes\lambda)\circ
(\delta\otimes\id_G)\circ \delta\big)\times(1\otimes 1\otimes M)
\]
is a faithful representation of $A\times_{\delta}G$ on $\H\otimes
L^2(G)\otimes L^2(G)$. Hence the result will follow if
we can show that $\Ind\big((\pi\otimes\lambda)\circ \delta\big)$ is unitarily
equivalent to $(\Ind\pi)\otimes 1$.
For this let $U\in U(L^2(G\times G))$ be defined by
$(U\xi)(s,t)=\Delta(t)^{-1/2}\xi(st^{-1},s)$. A quick calculation
shows that $U^*$ is given by $(U^*\xi)(s,t)=
\Delta(s^{-1}t)^{1/2}\xi(t,
s^{-1}t)$. Using this, another  straightforward computation
shows that
\begin{equation}
\label{eq-unitary1}
U\big((\lambda\otimes\lambda)\circ \delta_G(s)\big)U^*=
U(\lambda_s\otimes\lambda_s)U^*=\lambda_s\otimes 1
\end{equation}
for all $s\in G$ and
\begin{equation}
\label{eq-unitary2}
U(1\otimes M(f))U^*=M(f)\otimes 1
\end{equation}
for all $f\in C_0(G)$%
\footnote{Notice that $U$ is the composition
$U=V_G^*W_G^*$, where $(W_G\xi)(s,t)=\xi(s,s^{-1}t)$ and
$(V_G\xi)(s,t)=\Delta(t)^{1/2}\xi(st,t)$. $\Ad W_G^*$ moves $1\otimes M$
to $M\otimes M\circ \eps$ and $(\lambda\otimes \lambda)\circ 
\delta_G$ to $\lambda\otimes 1$,
where $\eps$ is the comultiplication on
$C_0(G)$, and $\Ad V_G^*$ moves $M\otimes M\circ \eps$ to $M\otimes 1$.}%
.  Using this and the
coaction identity for $\delta$, we get
\begin{align*}
(1\otimes U)&\big((\pi\otimes\lambda\otimes\lambda)\circ
(\delta\otimes\id_G)\circ \delta(a)\big)(1\otimes U^*)\\
&=(1\otimes U)\big((\pi\otimes\lambda\otimes\lambda)\circ
(\id_A\circ \delta_G)\circ \delta(a)\big)(1\otimes U^*)\\
&=(1\otimes U)\big((\pi\otimes((\lambda\otimes\lambda)\circ\delta_G))
\circ \delta(a)\big)(1\otimes U^*)\\
&\stackrel{(\ref{eq-unitary1})}{=}
(\pi\otimes(\lambda\otimes 1))\circ\delta(a)=\big((\pi\otimes\lambda)\circ
\delta(a)\big)\otimes 1
\end{align*}
for $a\in A$, and (\ref{eq-unitary2}) clearly implies that
\[
(1\otimes U)(1\otimes 1\otimes M(f))(1\otimes U^*)=
1\otimes M(f)\otimes 1
\]
for $f\in C_0(G)$. This completes the proof.
\end{proof}

The most enlightening example of normalizations is possibly given
in the case of dual coactions.

\begin{prop}
\label{prop-dualnormal}
Let $(A,G,\alpha)$ be an action and let $\widehat\alpha$ be the dual coaction
of $G$ on $A\times_{\alpha}G$. Then the normalization
$\widehat\alpha^n$ of $\widehat\alpha$ is the dual coaction of $G$ on
$A\times_{\alpha,r}G$. Hence $\widehat\alpha$ is normal if and only if
$i_A^r\times i_G^r\:A\times_{\alpha}G\to A\times_{\alpha,r}G$ is an
isomorphism. Moreover, the double crossed products
$A\times_{\alpha}G\times_{\widehat\alpha}G$ and
$A\times_{\alpha,r}G\times_{\widehat\alpha^n}G$ are canonically
isomorphic.
\end{prop}
\begin{proof}
Let $\pi\times U$ be any faithful representation
of $A\times_{\alpha}G$ on a Hilbert space $\H$.
Then \thmref{thm-crossed} implies that we can write
$j_{A\times_{\alpha}G}=(\pi\times U\otimes\lambda)\circ\widehat\alpha$.
By the definition of $\widehat\alpha$ it follows that
$j_{A\times_{\alpha}G}\circ i_A(a)=\pi(a)\otimes 1$ and
$j_{A\times_{\alpha}G}\circ i_G(s)=U_s \otimes \lambda_s$,
where $(i_A, i_G)$ denotes
the canonical map of $(A,G)$ into $M(A\times_{\alpha}G)$. Hence
$j_{A\times_{\alpha}G}=(\pi\otimes 1)\times (U\otimes \lambda)$, which
by
\lemref{lem-regrep} is unitarily equivalent to a regular representation of
$A\times_{\alpha}G$.  Thus, identifying $A\times_{\alpha,r}G$ with
$j_{A\times_{\alpha}G}(A\times_{\alpha}G)$, we get $i_A^r=\pi\otimes 1,
i_G^r=U\otimes\lambda$.

By definition (using $j_{A\times_{\alpha}G}=i_A^r\times i_G^r$), $\alpha^n$ is
determined by the equation
\[
\alpha^n\circ (i_A^r\times i_G^r)
=\bigl((i_A^r\times i_G^r)\otimes \id_G\bigr)\circ
\widehat\alpha.
\]
But the right hand side applied to the elements $i_A(a), i_G(s)\in
M(A\times_{\alpha}G)$ gives $i_A^r(a)\otimes 1$ and $i_G^r(s)\otimes s$,
respectively. Hence, $\widehat\alpha^n$ is precisely the dual coaction of $G$
on $A\times_{\alpha,r}G$ as defined in \exref{ex-reduceddual}.
The final assertion now follows from
\propref{prop-normal}.
\end{proof}

\begin{ex}
\label{ex-groupcrossed}
As an interesting application of the above results we
now show that $(\K(L^2(G)),\lambda, M)$ is a crossed product for
the coaction $(C^*(G),G,\delta^G)$.
By \propref{prop-dualnormal} we know that the normalization
$\delta_G^n$ is the coaction on $C_r^*(G)$ determined by the
map $\lambda_s\mapsto \lambda_s\otimes u(s)$ for $s\in G$
(note that $\delta_G$ is just the dual coaction on $C^*(G)=\bbC\times_{\id}G$).
Consider the trivial representation $1_G\:G\to\bbC$ of
$G$. Then $\Ind 1_G=(1_G\otimes\lambda)\circ\delta_G\times(1_{\bbC}\otimes M)
=\lambda \times M$.
  Thus
$\ker(1_G\otimes\lambda)\circ\delta_G=\ker \lambda$, which by
\propref{prop-dualnormal} coincides with $j_{C^*(G)}$,
and it follows from \propref{prop-faithful} that
$\lambda\times M$ is a faithful representation of $C^*(G)\times_{\delta_G}G$.
Since we already observed in \remref{rem-compact} that
$C^*(\lambda,M)=\lambda(C^*(G))M(C_0(G))=\K(L^2(G))$, the result follows.
\end{ex}

\exref{ex-groupcrossed} should be compared with
\exref{ex-elementary}, where we mentioned that $(\K(L^2(G)),M,\lambda)$ is a
crossed product for the action $(C_0(G),G,\tau)$. Note that for the latter
result we had to refer to a deep theorem of Mackey (the imprimitivity theorem),
while the result on $(C^*(G),G,\delta_G)$ is a natural
outcome of the theory of coactions. We shall see below
(\thmref{thm-isom}) that there is a natural isomorphism
between $C^*(G)\times_{\delta_G}G$ and $C_0(G)\times_{\tau}G$,
so the result for $(C_0(G),G,\tau)$ will actually be a consequence
of the above example.

\section
{The duality theorems of Imai-Takai and Katayama}

In this section we want to deduce the classical duality theorem
of Imai and Takai for actions and the duality theorem of
Katayama for coactions,
\ie, the analogues of the Takesaki-Takai duality theorem
for nonabelian groups. We start with a result which is actually more
general than what we need here, but which will  be important later.

Let $(A,G,\alpha)$ be an action and let $N$ be a closed normal
subgroup of $G$.  Recall from Examples~\ref{ex-dualcoaction} and
\ref{ex-restrict} that the restriction $\widehat\alpha|$ of the dual
coaction $\widehat\alpha$ to $G/N$ is the integrated form of the
covariant homomorphism $(i_A\otimes 1, i_G\otimes q)$ of
$(A,G,\alpha)$ on $M(A\times_{\alpha}G\otimes C^*(G/N))$, where
$q\:G\to G/N\subseteq M(C^*(G/N))$ denotes the quotient map.  In what
follows next we first want to show that
$(A\times_{\alpha}G)\times_{\widehat{\alpha}|}G/N$ is canonically
isomorphic to the crossed product $(A\otimes
C_0(G/N))\times_{\alpha\otimes\tau}G$, where $\tau\:G\to
\Aut(C_0(G/N))$ is given by $\tau_s(f)(tN)=f(s^{-1}tN)$, as in
\exref{ex-elementary}.

Recall that if $(\pi,\mu)$ is a pair of commuting nondegenerate
homomorphisms $\pi\:A\to M(D)$, $\mu\:C_0(G/N)\to M(D)$, then (since
$C_0(G/N)$ is nuclear) there is a nondegenerate
homomorphism
\[
\pi\otimes\mu\:A\otimes C_0(G/N)\to M(D)
\]
which is determined on elementary tensors by $\pi\otimes\mu(a\otimes
f)= \pi(a)\mu(f)$, and every nondegenerate homomorphism of $A\otimes
C_0(G)$ arises this way.  Since the canonical maps $A\to M(A\otimes
C_0(G/N))$, $C_0(G/N)\to M(A\otimes C_0(G/N))$ are $G$-equivariant, it
follows that a triple $(\pi,\mu,V)$ of nondegenerate homomorphisms
$\pi\:A\to M(D)$, $\mu\:C_0(G/N)\to M(D)$, $V\:C^*(G)\to M(D)$
determines a covariant homomorphism $(\pi\otimes\mu, V)$ of $(A\otimes
C_0(G/N), G,\alpha\otimes \tau)$ if and only if
\begin{enumerate}
\item $\pi$ and $\mu$ have commuting ranges in $M(D)$,
\item $(\pi, V)$ is covariant for $\alpha$, and
\item $(\mu, V)$ is covariant for $\tau$.
\end{enumerate}

\begin{prop}[{\cf~\cite[Examples 2.9]{rae:rep}}]
\label{prop-dual}
The assignment $(\pi\otimes \mu, V)\mapsto (\pi\times V, \mu)$ is a
one-to-one correspondence between the covariant homomorphisms of
$(A\otimes C_0(G/N), G,\alpha\otimes \tau)$ and the covariant
homomorphisms of $(A\times_{\alpha}G,G/N,\widehat\alpha|)$.  Moreover,
\[
  \bigl( (\pi\otimes\mu)\times V \bigr)
\big( (A\otimes C_0(G/N)\times_{\alpha\otimes\tau}G \big)
= \bigl( (\pi\times V)\times \mu \bigr)
\big( A\times_{\alpha}G\times_{\widehat{\alpha}|}G/N \big).
\]
\end{prop}

\begin{proof}
We only have to check that conditions (i)--(iii) above are
equivalent to the covariance condition for
$(\pi \times V,\mu)$
  with respect to $\widehat\alpha|$.
First note that,
  by \lemref{lem-slice},
the covariance condition for $(\pi\times V,\mu)$  is equivalent to
\setlength{\multlinegap}{1cm}
\begin{multline}
\label{eq-covariant}
S_f\bigl( (\pi \times V \otimes \id_{G/N}) \circ
\widehat{\alpha}|(c)((\mu\otimes\id_G)(w_{G/N})) \bigr)
\\
= S_f\bigl( (\mu \otimes \id_{G/N})(w_{G/N})
(\pi \times V(c) \otimes 1) \bigr)
\end{multline}
for all $c\in A\times_{\alpha}G$, $f\in B(G/N)$.
By integration it is enough to check this equation
on all elements $i_A(a), i_G(s)\in M(A\times_{\alpha}G)$.
For $c=i_A(a)$, using \lemref{lem-slice} and
\propref{slice-unitary}, the left-hand side of (\ref{eq-covariant}) becomes
\[
S_f \bigl( (\pi(a) \otimes 1) (\mu \otimes \id_G)(w_G) \bigr)
= \pi(a) S_f \bigl( (\mu \otimes \id_G)(w_G) \bigr)
= \pi(a) \mu(f),
\]
and a similar computation shows that the right hand side of
(\ref{eq-covariant}) becomes $\mu(f)\pi(a)$.
Thus, since $A(G/N)$ is norm-dense
in $C_0(G/N)$
and weak$^*$-dense in $B(G/N)$
we see that having
(\ref{eq-covariant}) hold for all $c=i_A(a)\in i_A(A)$ and $f\in B(G/N)$
is equivalent to (i).
For $c=i_G(s)$ the left hand side of
(\ref{eq-covariant}) becomes
\begin{align*}
S_f \bigl( (V_s\otimes q(s))(\mu\otimes\id_{G/N})(w_{G/N}) \bigr)
&
= V_s S_{f\cdot q(s)} \big( (\mu\otimes\id_{G/N})(w_{G/N}) \bigr)
\\
&
= V_s \mu(f\cdot q(s))
=V_s \mu(\tau_{s^{-1}}(f)),
\end{align*}
while the right hand side becomes $\mu(f)V_s$.
Thus, having (\ref{eq-covariant}) hold for all $c=i_G(s)$ and $f\in B(G/N)$
is equivalent to (iii).

The last assertion follows from the equality
\[
\pi(A)\mu(C_0(G/N))V(C^*(G))=\pi(A)V(C^*(G))\mu(C_0(G/N)).
\]
\end{proof}

Let us write $i_A\otimes i_{C_0(G/N)}$ for the
canonical map
$$i_{A\otimes C_0(G/N)}\:A\otimes C_0(G/N)\to M\big((A\otimes
C_0(G/N))\times_{\alpha\otimes\tau}G\big).$$
We define an action
$\beta\:G/N\to \Aut((A\otimes
C_0(G/N))\times_{\alpha\otimes\tau}G)$ by
\[
\beta_{sN}\big(i_A(a)i_{C_0(G/N)}(f)i_G(z)\big)=
i_A(a)i_{C_0(G/N)}(\sigma_{sN}(f))i_G(z),
\]
where $\sigma_{sN}$ denotes right translation by $sN$.
The following theorem is now a consequence of
\propref{prop-dual} and the universal
properties of the crossed products.

\begin{thm}
\label{thm-isom}
Let $(A,G,\alpha)$ be an action and let $N$ be a closed normal
subgroup of $G$.
Then
there is a $\beta-\widehat{\hat\alpha|}$
equivariant isomorphism
\[
(A \otimes C_0(G/N)) \times_{\alpha \otimes \tau} G
\cong
(A \times_\alpha G) \times_{\hat\alpha|} G/N
\]
taking $i_{A \otimes C_0(G/N)}(a \otimes f) i_G(c)$
{to}
$j_{G/N}(f) j_{A \times G}(i_A(a) i_G(c))$
for $a \in A$, $f \in C_0(G/N)$, and $c \in C^*(G)$.
\end{thm}

We also present a reduced version of \thmref{thm-isom}.

\begin{thm}
\label{thm-isom-red}
Let $(A,G,\alpha)$ be an action and let $N$ be a closed normal
subgroup of $G$.
Then
there is a $\beta-\widehat{\hat\alpha|}$
equivariant isomorphism
\[
(A \otimes C_0(G/N)) \times_{\alpha \otimes \tau,r} G
\cong
(A \times_{\alpha,r} G) \times_{\hat\alpha|} G/N
\]
taking $i^r_{A \otimes C_0(G/N)}(a \otimes f) i^r_G(c)$
{to}
$j_{G/N}(f) j_{A \times_r G}(i^r_A(a) i^r_G(c))$
for $a \in A$, $f \in C_0(G/N)$, and $c \in C^*(G)$.
\end{thm}

\begin{proof}
We show that, up to unitary equivalence, there are regular representations
of $(A\otimes C_0(G/N), G, \alpha\otimes\tau)$ and
$(A\times_{\alpha,r}G, G/N, \widehat\alpha|)$ which match up as
in  \propref{prop-dual}. This will give the result.

We start with a
faithful representation  $\pi\:A\to \B(\H)$. Let $M^{G/N}\:C_0(G/N)\to
\B(L^2(G/N))$ and  $M^G\:C_0(G) \to \B(L^2(G))$ denote the respective
representations by multiplication operators.
Consider the regular representation
\[
\Pi\deq \Ind(\pi\otimes M^{G/N})=
(\pi\otimes M^{G/N}\otimes M^G)\circ(\alpha\otimes\tau)\times
(1\otimes 1\otimes \lambda)
\]
of
$(A\otimes C_0(G/N))\times_{\alpha\otimes\tau}G$
on
$\H\otimes L^2(G/N)\otimes L^2(G)$
(and recall that in this formula we view $\alpha\otimes \tau$
as a homomorphism of $A\otimes C_0(G/N)$ to
$M(A\otimes C_0(G/N)\otimes C_0(G))$ as in
\defnref{def-red crossed}),
and the regular  representation
\[
\Lambda\deq \Ind(\Ind\pi)=
\big((\pi\otimes M^G)\circ\alpha\times
(1\otimes\lambda)\big)\circ\widehat\alpha|\times
\big( 1\otimes 1\otimes M^{G/N}\big)
\]
of $(A\times_{\alpha,r}G)\times_{\widehat\alpha|}G/N$ on
$\H\otimes L^2(G)\otimes L^2(G/N)$. Define a unitary operator
$U\:\H\otimes L^2(G/N)\otimes L^2(G)\to \H\otimes L^2(G)\otimes L^2(G/N)$
by 
\[
(U\xi)(t, sN)= \xi(t^{-1}sN, t).
\]
for $\xi\in L^2(G/N\times G,\H)=\H\otimes L^2(G/N)\otimes L^2(G)$.
Moreover, let's
write $i_A^r\otimes i_{C_0(G/N)}^r$ for
$i_{A\otimes C_0(G/N)}^r$ and $j_A\times j_G$ for $j_{A\times_rG}$.
We claim that
\begin{enumerate}
\item $\Ad U\circ \Pi\circ i_A^r= \Lambda\circ j_A$,
\item $\Ad U\circ \Pi\circ i_{C_0(G/N)}^r=\Lambda\circ j_{C_0(G/N)}$, and
\item $\Ad U\circ \Pi\circ i_G^r=\Lambda\circ j_G$.
\end{enumerate}
Since $\Ad U\circ \Pi$ and $\Lambda$ are faithful on the respective
(reduced) crossed products, the result will follow from
\propref{prop-dual}. For~(i),
for $a\in A$ and $\xi\in L^2(G\times G/N,\H)\cong
\H\otimes L^2(G)\otimes L^2(G/N)$  we compute
\begin{align*}
\big(\Ad U  &\circ \Pi(i_A^r(a))\xi\big)(t,sN)=
\big(U(\pi\otimes M^{G/N}\otimes M^G)\circ(\alpha\otimes\tau)(a\otimes
1)U^*\xi\big)(t,sN)\\
&=\big((\pi\otimes M^G\otimes M^{G/N})\circ(\alpha\otimes \tau)(a\otimes
1)U^*\xi\big)(t^{-1}sN,t)\\
&=\pi(\alpha_{t^{-1}}(a))\big(U^*\xi\big)(t^{-1}sN,t)
=\pi(\alpha_{t^{-1}}(a))\xi(t,sN)\\
&=\big((\pi\otimes M^G)\circ\alpha(a)\otimes 1)\xi\big)(t,sN)
=\big(\Lambda\circ j_A(a)\xi\big)(t,sN).
\end{align*}
For~(ii) 
(writing $(\pi\otimes M)\widetilde{\ }$ instead of
$(\pi\otimes M^{G/N}\otimes M^G)\circ (\alpha\otimes\tau)$),
for $f\in C_0(G/N)$ compute
\begin{align*}
\big(\Ad U&\circ\Pi(i_{C_0(G/N)}^r(f))\xi\big)(t,sN)=
\big(U(\pi\otimes M)\!\widetilde{\ }(1\otimes f)U^*\xi\big)(t,sN)\\
&=\big((\pi\otimes M)\!\widetilde{\ }(1\otimes f)U^*\xi\big)(t^{-1}sN,t)
=M(\tau_{t^{-1}}(f))\big(U^*\xi\big)(t^{-1}sN,t)\\
&=\tau_{t^{-1}}(f)(t^{-1}sN)\big(U\xi\big)(t^{-1}sN,t)
=f(sN)\xi(t,sN)
=\big((1\otimes 1\otimes M)(f)\xi\big)(t,sN)\\
&=\big(\Lambda\circ j_{C_0(G/N)}(f)\xi)(t,sN),
\end{align*}
and for~(iii), for $r\in G$ we compute
\begin{align*}
\big(\Ad U&\circ \Pi(i_G^r(r))\xi\big)(t,sN)
=\big(U(1\otimes 1\otimes\lambda(r))U^*\xi\big)(t,sN)\\
&=\big((1\otimes 1\otimes \lambda(r))U^*\xi\big)(t^{-1}sN,t)
=\big(U^*\xi\big)(t^{-1}sN, r^{-1}t)\\
&=\xi(r^{-1}t, r^{-1}sN)
=\big((1\otimes \lambda(r)\otimes \lambda^{G/N}(rN))\xi\big)(t,sN)
=\big(\Lambda\circ j_G(r)\xi\big)(t,sN).
\end{align*}
\end{proof}

In case $N=\{e\}$ and $A=\bbC$ we get, as a corollary of \thmref{thm-isom}
and  \exref{ex-groupcrossed},  the complete proof of the
statement of \exref{ex-elementary}.

\begin{cor}
\label{cor-elementary}
$(\K(L^2(G)), M,\lambda)$ is a crossed product for $(C_0(G), G,\tau)$.
\end{cor}

Another (more general) consequence of the special case $N=\{e\}$
together with \exref{ex-extend} is the Imai-Takai duality
theorem for actions.

\begin{thm}
\label{thm-ImaiTakai}
Let $(A,G,\alpha)$ be an action. Then
\[
(A\times_{\alpha}G)\times_{\widehat\alpha}G
\cong
A\otimes \K(L^2(G)),
\]
equivariantly for the double dual action
$\widehat{\widehat\alpha}$ and for the action $\alpha\otimes \Ad\rho$,
where $\rho$ denotes the right regular representation of $G$ on
$L^2(G)$.
\end{thm}

The only statement in the above theorem we haven't explicitly shown
so far is the statement about the action $\widehat{\widehat{\alpha}}$.
But this can be done by checking carefully what
the isomorphism of \exref{ex-extend} does to the action
$\beta$ of the theorem.
Since the dual coaction of $G$ on $A\times_{\alpha,r}G$
is the normalization of $\widehat\alpha$, we get
from \propref{prop-normal} a
reduced version of the above result

\begin{thm}
\label{thm-ImaiTakaired}
There exists an $\widehat{\widehat\alpha}-\alpha\otimes\Ad\rho$
equivariant isomorphism
\[
(A\times_{\alpha,r}G)\times_{\widehat\alpha}G
\cong A\otimes \K(L^2(G)).
\]
\end{thm}

We now turn our attention to the duality theorem of Katayama,
where we consider the double crossed product
$(A\times_{\delta}G)\times_{\widehat\delta}G$ of a coaction
$(A,G,\delta)$. Unfortunately the
situation is a bit more complicated in this case,
since the crossed product $A\times_{\delta}G$ doesn't
seem to see the full algebra $A$ but rather only the quotient
$j_A(A)$ of $A$. So in principle we cannot expect to recover
$A$ from the double crossed product. However, that this reasoning
is somewhat too naive follows from \cite[Theorem 3.7]{qui:full}
which shows that full duality does work for all dual coactions
on full crossed products, which are very often not normal
(see also the discussion in \cite[\S3]{EQ-IC}).

\begin{thm}[{\cf~\cite[Theorem 8]{kat}}]
\label{thm-kat}
Let $(A,G,\delta)$
be a nondegenerate coaction.
Then the reduced double crossed product
$A\times_{\delta}G\times_{\widehat\delta,r}G$ is isomorphic to
$(A/\ker j_A)\otimes \K(L^2(G)).$
\end{thm}

For the proof we need

\begin{lem}[{\cf~\cite[Theorem 5, Theorem 8]{kat}}]
\label{lem-AtensorK}
Let $\delta\:A\to M(A\otimes C^*(G))$ be a nondegenerate coaction,
and let $W_G=M\otimes\lambda(w_G)$.
Then
\begin{equation}
\label{eq-kat1}
(\id_A\otimes\lambda)\circ \delta(A)(1\otimes \K(L^2(G)))=
A \otimes \K(L^2(G))
\end{equation}
and
\begin{multline}
\label{eq-kat2}
(1\otimes W_G)\big((\id_A\otimes\lambda)\circ
\delta(A)\otimes 1\big)(1\otimes W_G^*)(1\otimes
1\otimes \K(L^2(G)))\\
=\big((\id_A\otimes\lambda)\circ \delta(A)\big)\otimes \K(L^2(G)).
\end{multline}
\end{lem}
\begin{proof} It follows from \remref{rem-compact} that
$\lambda(C^*(G))M(C_0(G))=\K(L^2(G))$. Since $\delta$ is nondegenerate
we have $(\id_A\otimes \lambda)\circ\delta(A)(1\otimes C_r^*(G))=
A\otimes C_r^*(G)$.
Multiplying both sides of the equation with $1\otimes M(C_0(G))$
from the right we get (\ref{eq-kat1}).

Since
$\Ad W_G\circ (\lambda\otimes 1)= (\lambda\otimes\lambda)\circ \delta_G$
(which follows from applying $\id\otimes \lambda$
to \eqeqref{eq-ind}),
we get
\begin{align*}
&(1\otimes W_G)\big((\id_A\otimes\lambda)\circ
\delta(A)\otimes 1\big)(1\otimes W_G^*)
=(\id_A \otimes\lambda\otimes \lambda)
\big((\id_A\otimes \delta_G)\circ \delta(A))
\\&\quad=(\id_A \otimes\lambda\otimes \lambda)
\big((\delta\otimes \id_G)\circ \delta(A))
=(\id_A \otimes\lambda\otimes \id_{\K})
\big((\delta\otimes \lambda)\circ \delta(A))
\end{align*}
Thus, the left hand side of (\ref{eq-kat2}) becomes
\begin{align*}
(\id_A \otimes\lambda&\otimes \id_{\K})\circ (\delta\otimes \id_{\K})
\big((\id_A\otimes \lambda)(\delta(A))(1\otimes\K(L^2(G)))\big)\\
&\stackrel{(\ref{eq-kat1})}{=}
(\id_A \otimes\lambda\otimes \id_{\K})\circ (\delta\otimes \id_{\K})
\big(A\otimes \K(L^2(G))\big)\\
&=\big((\delta\otimes\lambda)\circ \delta(A)\big)\otimes \K(L^2(G)).
\end{align*}
\end{proof}

\begin{proof}[Proof of \thmref{thm-kat}]
Let us assume that $A$ is faithfully represented on a Hilbert space
$\H$.  By \thmref{thm-crossed}, $\pi\deq (\id_A\otimes\lambda)\circ
\delta\times (1\otimes M)$ is a faithful representation of
$A\times_{\delta}G$ on $\H\otimes L^2(G)$, and then $\Ind\pi$ (the
regular representation of $(A\times_{\delta}G,G,\widehat\delta)$
induced from $\pi$) is a faithful representation of
$A\times_{\delta}G\times_{\widehat\delta,r}G$ on $\H\otimes
L^2(G)\otimes L^2(G)$.  Let us write $k_A, k_{C(G)}$ and $k_G$ for the
compositions $\Ind\pi\circ i_{A\times_{\delta}G}^r\circ j_A$,
$\Ind\pi\circ i_{A\times_{\delta}G}^r\circ j_G$ and $i_G^r$,
respectively, where $(i_{A\times_{\delta}G}^r, i_G^r)$ denote the
canonical maps of $(A\times_{\delta}G,G)$ into
$M(A\times_{\delta}G\times_{\widehat\delta,r}G)$.  By definition of
$\Ind\pi$ we have
\begin{gather*}
k_A(a)= \big((\id_A\otimes\lambda)\circ \delta(a)\big)\otimes 1,
\qquad
k_{C(G)}(f)= 1\otimes \big((M\otimes M)\circ \nu(f)\big)
\\
\text{and}\quad
k_G(s)=1\otimes 1\otimes \lambda_s,
\end{gather*}
where $\nu\:C_0(G)\to C_b(G\times G)$ is defined by $\nu(f)(s,t)=f(st^{-1})$.
We define $V\deq (1\otimes W_G)(1\otimes S)$, where $W_G=M\otimes\lambda(w_G)$ and
$S\:L^2(G\times G)\to L^2(G\times G)$ is the self-adjoint unitary operator
$(S\xi)(s,t)=\Delta(t)^{-1/2}\xi(s,t^{-1})$.
Since $S(\lambda\otimes 1)S=\lambda\otimes 1$ we have
$\Ad (1\otimes S)\circ k_A=k_A$. For $\xi\in L^2(G\times G)$, we compute
\begin{align*}
\big(W_GS&((M\otimes M)\circ \nu(f))SW_G^*\xi\big)(s,t)
=\Delta(t^{-1}s)^{1/2}\big((M\otimes M)\circ \nu(f))SW_G^*\xi\big)(s,t^{-1}s)\\
&=\Delta(t^{-1}s)^{1/2}f(t)\big(SW_G^*\xi\big)(s,t^{-1}s)=f(t)\xi(s,t)=
\big((1\otimes M(f))\xi\big)(s,t),
\end{align*}
which implies that $\Ad V\circ k_{C(G)}= 1\otimes 1\otimes M$,
and the identity $W_GS(1\otimes \lambda)SW_G^*=1\otimes\rho$ implies that
$\Ad V\circ k_G= 1\otimes 1\otimes \rho$,
where $\rho\:G\to U(L^2(G))$ denotes the right regular representation
of $G$.
Thus, since $C_0(G)\rho(C^*(G))=\K(L^2(G))$ by \remref{rem-compact},
it follows that
\begin{align*}
A\times_{\delta}G\times_{\widehat\delta,r}G
&\cong \Ad V(k_A(A))\Ad V(k_{C(G)}(C_0(G)))\Ad V(k_G(C^*(G)))
\\&=\big((1\otimes W_G)\big((\id_A\otimes \lambda)\circ \delta(A)\otimes
1 \big)(1\otimes W_G^*)\big)
\\&\hphantom{= \big(}\quad \cdot\big(1\otimes 1\otimes M(C_0(G))\big)
\big(1\otimes 1\otimes \rho(C^*(G))\big)
\\&=\Big((1\otimes W_G)
\big((\id_A\otimes \lambda)\circ \delta(A)\otimes 1\big)
(1\otimes W_G^*)\Big)
\big(1\otimes 1\otimes \K(L^2(G))\big)
\\&\stackrel{(\ref{eq-kat1})}{=}
\big((\id_A\otimes\lambda)\circ\delta(A)\big)\otimes \K(L^2(G)).
\end{align*}
Finally, $(\id_A\otimes\lambda)\circ \delta = j_A$ by
definition,
so the proof is complete.
\end{proof}

\section{Other definitions of coactions}
\label{sec-other}

As mentioned before, full coactions were first introduced by the fourth
author in \cite{rae:rep} using maximal tensor products instead of the
minimal ones we use here.  Although we shall see below that this
definition has some disadvantages, we should point out that
the approach of \cite{rae:rep} inspired much of the theory we use in
this work.

Let us denote by $\delta_G^{\max}\:C^*(G)\to M(C^*(G)\otimes_{\max}C^*(G))$
the integrated form of the unitary homomorphism
$s\mapsto u(s)\otimes_{\max}u(s)$. The fourth author defined a coaction
of $G$ on a $C^*$-algebra $A$ as a nondegenerate homomorphism
$\delta\:A\to M(A\otimes_{\max}C^*(G))$ satisfying
\begin{enumerate}
\item $\delta(a)(1\otimes z)\in A\otimes_{\max}C^*(G)$ for all $a\in A$, $z\in
C^*(G)$.
\item $(\delta\otimes \id_G)\circ \delta=(\id_A\otimes
\delta_G^{\max})\circ\delta$ as maps from $A$ to
$M(A\otimes_{\max}C^*(G)\otimes_{\max}C^*(G))$.
\end{enumerate}
Note that there is no assumption on injectivity here, but it
was pointed out in \cite[Remarks 2.2]{rae:rep} that if
$\delta$ is not injective, then it factors through an
injective coaction $\delta_1\:A/\ker\delta\to
M(A/\ker\delta\otimes_{\max}C^*(G))$.
Covariant representations of a full coaction $(A,G,\delta)$ (in this sense)
on a Hilbert space $\H$ are defined completely analogous to
\defnref{def-covariant}: They are pairs $(\pi,\mu)$ of
nondegenerate representations $\pi\:A\to\B(\H), \mu\:C_0(G)\to \B(\H)$
satisfying
\[
(\pi \otimes \id_G) \circ \delta
= \Ad (\mu \otimes \id_G)(w_G) \circ (\pi \otimes 1)
\]
in $M(\K(\H)\otimes C^*(G))$.

If $\delta\:A\to M(A\otimes_{\max}C^*(G))$ is an injective coaction as above,
then $\delta^{\min}=\Phi\circ \delta$ is a
coaction in our sense (\ie, in the sense of \cite{qui:fullred}), where
$\Phi\:A\otimes_{\max}C^*(G)\to A\otimes C^*(G)$ denotes the quotient map
(see \cite[Remarks 2.2 (4)]{rae:rep}). Moreover,
by \cite[Remarks 2.5]{rae:rep} a pair
$(\pi,\mu)$ is a covariant representation for $\delta$ if and only
if it is a covariant representation for $\delta^{\min}$, which
in particular shows (since both satisfy the same universal properties)
that the crossed products of $\delta$ and $\delta^{\min}$ are the same.

One of the main defects of the theory of full coactions with maximal
tensor products is the fact that, if $(\pi,\mu)$ is a covariant
representation of $(A,G,\delta)$, then $\delta^{\mu}(\pi(a)) = \Ad
\bigl( (\mu\otimes_{\max}\id_G)(w_G) \bigr)(\pi(a)\otimes 1)$ does not
define a coaction on $\pi(A)$ in general.  So one of the basic
features of the theory (see \secref{sec-normal}) is not available
anymore.

\begin{ex}[{\cf~\cite[Example 1.15]{qui:fullred}}]
\label{ex-full}
We consider the coaction $(C^*(G), G, \delta_G^{\max})$.
Since
$(\delta_G^{\max})^{\min}=\delta_G$ it follows from the above discussion that
$C^*(G)\times_{\delta_G^{\max}}G=C^*(G)\times_{\delta_G}G$.
Thus, \exref{ex-groupcrossed} implies that
$(\K(L^2(G)), \lambda,M)$ is a crossed product for 
$(C^*(G),G,\delta_G^{\max})$.
Now let $G$ be a discrete non-amenable group. We claim that
$\delta^M=\Ad(M\otimes_{\max}\id_G(w_G))\circ (\lambda(\cdot)\otimes 1)$
is not a full coaction on $\lambda(C^*(G))=C_r^*(G)$ in the sense
of \cite{rae:rep},
\ie, we show that it doesn't factor through $C_r^*(G)=C^*(G)/\ker\lambda$.
To see this, first observe that if it did, then it would be determined
by the map $\lambda_s\mapsto \lambda_s\otimes_{\max}u(s)$ for $s\in G$.
Let $\rho$ denote the right regular representation of $G$,
and define a representation $\pi$ of $C_r^*(G)\otimes_{\max}C^*(G)$
on $\ell^2(G)$ by $\pi(\lambda(x_i)\otimes y_i)=\lambda(x_i)\rho(y_i)$,
$x_i,y_i\in C^*(G)$ (since the left regular
representation commutes with  the right regular representation,
$\pi$ is well-defined). For $s\in G$ it follows that
$\pi\circ \delta^M(\lambda_s)=\lambda_s\rho_s$,
hence $\pi\circ \delta^M(\lambda_s)\Chi_{\{e\}}=\Chi_{\{e\}}$,
where $\Chi_{\{e\}}\in \ell^2(G)$ denotes the characteristic function
on the trivial element $e$. Thus,
$\ker \lambda\subseteq\ker \pi\circ \delta^M\circ \lambda\subseteq\ker 1_G$,
where $1_G$ is the trivial representation of $G$, which contradicts
the fact that $G$ is non-amenable.
\end{ex}

Note that the (non existent) coaction $\delta^M$ in the above
example corresponds to the dual coaction $\delta_G^n$
of $C^*(G)$ on the reduced group algebra $C_r^*(G)$ in our theory.
Thus there seems to be no natural definition of a dual coaction
of $C^*(G)$ on the reduced crossed product $A\times_{\alpha,r}G$
in the theory of \cite{rae:rep}. Of course this is another serious drawback.

We want to use this opportunity to point out a gap in
\cite[Example 1.15]{qui:fullred} (from where we actually have
extracted the above example). There it was stated that there exists
no full coaction $\delta$ of $G$ on $C_r^*(G)$ in the sense of
\cite{rae:rep}
such that $\delta^{\min}=\delta_G^n$, arguing that if such a coaction
did exist, it would have to be given by $\delta^M$ as in the example above.
But of course, there is no guarantee that, if $I$ is the kernel
of the quotient map $C_r^*(G)\otimes_{\max}C^*(G)\to C_r^*(G)\otimes C^*(G)$
in $M(C_r^*(G)\otimes_{\max}C^*(G))$, we could not modify
$\delta^M$ by adding suitable elements of $I$ so that
the new map would be a coaction of $G$ on $C_r^*(G)$ in the
sense of \cite{rae:rep}. However, we do not see how this can be done.

The other different approach (and the one which is actually most established
in the literature \cite{LanDT, kat, lprs}) is to define coactions of $G$ as
coactions of the reduced group $C^*$-algebra $C_r^*(G)$ with
comultiplication given by
$\delta_G^r\:\lambda_s\mapsto\lambda_s\otimes \lambda_s$.

\begin{defn}

\label{defn-red coact}
A \emph{reduced coaction} of $G$ on $A$ is an injective and nondegenerate
homomorphism
$\delta\:A\to M(A\otimes C_r^*(G))$ satisfying
\begin{enumerate}
\item $\delta(a)(1\otimes \lambda(x))\in A\otimes C_r^*(G)$ for all $x\in
C^*(G)$ and
\item $(\delta\otimes\id_G^r)\circ \delta=(\id_A\otimes \delta_G^r)\circ
\delta$ as maps from $A$ into $M(A\otimes C_r^*(G)\otimes C_r*(G))$.
\end{enumerate}
We also call the triple $(A,G,\delta)$ a reduced coaction.
\end{defn}
Again, the definition of covariant representations of reduced coactions
is completely analogous to the definition we use here:
They are pairs $(\pi,\mu)$ of nondegenerate representations of $(A,C_0(G))$
satisfying
\[
(\pi\otimes \id_G^r)\circ\delta
=\Ad \bigl( (\mu\otimes\lambda)(w_g) \bigr)
(\pi(\cdot)\otimes 1).
\]
If $\delta$ is a coaction in our sense, then
$\delta^r\deq (\id_A\otimes\lambda)\circ \delta$ factors through a
reduced coaction on $A^r\deq A/\ker j_A$ \cite[Corollary 3.4]{qui:fullred},
which is called the \emph{reduction} of $\delta$.
$\delta^r$ is nondegenerate (in the appropriate sense) if and only if $\delta$
is nondegenerate. In particular, we see that if $\delta$ is normal,
then $\delta^r$ coacts an $A$ itself and we have $(\delta^n)^r=\delta^r$.
Thus we see that taking the assignment $\delta\mapsto \delta^r$ is
not a one-to-one correspondence between (conjugacy classes of)
coactions in our
sense and (conjugacy classes of) reduced coactions. However, it is if we
restrict to nondegenerate normal coactions \cite[Theorem 4.7]{qui:fullred}.

The covariant homomorphisms of $\delta$ and $\delta^r$
coincide by \remref{rem-normalrep} and \cite[Proposition 3.7]{qui:fullred}.
To be more precise: If $(\pi,\mu)$ is a covariant homomorphism of
$(A^r,G,\delta^r)$, then $(\pi\circ j_A,\mu)$ is a covariant homomorphism
of $(A,G,\delta)$ and every covariant homomorphism of
$(A,G,\delta)$ arises this way.
In particular, it follows from the universal properties
that the crossed products of $(A^r,G,\delta^r)$ and $(A,G,\delta)$
are the same (more precisely,
$(A^r\times_{\delta^r}G, j_{A^r}\circ j_A, j_G)$ is a crossed product
for $(A,G,\delta)$).

The fact that the reductions of a coaction $\delta$ and its normalization
$\delta^n$ are the same shows that the theory of coactions we use here
is potentially richer than the theory of reduced coactions.
One major drawback in the theory of reduced coactions is given by
the fact that, in case where $N$ is a non-amenable normal subgroup
of $G$, there is no well-defined quotient map $C_r^*(G)\to C_r^*(G/N)$,
so it would be difficult to define the restriction
of a reduced coaction of $G$ to the quotient $G/N$ in the realm
of reduced coactions (see \exref{ex-reduceddual}).
Since restriction  of coactions is one of the most basic concepts
in this work, the use of reduced coactions would be inadequate for
our purposes. However, the above discussion clearly shows that
results on reduced coactions can
often be applied (with care)  to problems of coactions in our sense
(in particular to problems concerning the representation
theory and the crossed products of coactions $(A,G,\delta)$).
Of course, if $G$ is amenable, all the different notions of coactions
do coincide.

%
%

\chapter
{The Imprimitivity Theorems of Green and Mansfield}
\label{imprim-chap}

In this appendix we want to recall the main ingredients of
the imprimitivity theorems of Green and Mansfield.
The main idea behind those imprimitivity theorems is to show that
certain $C^*$-algebras which appear in the study of
crossed products by actions and coactions, and in particular in
the study of induced representations, are Morita equivalent in the
sense that there are canonical imprimitivity bimodules linking
those algebras together.

\section{Imprimitivity theorems for actions}

In what follows we recall the main results of \cite{RaeIC}.
We refer to \cite{RW} (and \chapref{hilbert-chap} of this
work)
for the necessary background on
(pre-)imprimitivity bimodules and Morita equivalence.

Consider a $C^*$-algebra $D$, two locally compact groups $K$ and $H$,
and a locally compact space $P$. Suppose that $K$ acts freely
and properly on the left of $P$, and that $H$ acts likewise on the right
such that these actions commute (\ie, $k(ph)=(kp)h$).
Suppose also that we have commuting actions $\beta$ of $K$ and
$\alpha$ of $H$ on $D$. For the left action of $K$ we define
the induced $C^*$-algebra $\Ind_K^P(D,\beta)$ (or just $\Ind_K^PD$,
if confusions seems unlikely) to be the
set of bounded continuous functions $f\:P\to D$ such that
$f(kp)=\beta_k(f(p))$
for all $k\in K$ and $p\in P$, and such that the function
$Kp\mapsto \|f(p)\|$ vanishes at infinity on $K\backslash P$.
For the right action of $H$ we define the \emph{induced $C^*$-algebra}
$\Ind_H^P(D,\alpha)$ (or just $\Ind_H^PD$) as the set of
bounded continuous functions $g\:P\to D$ such that
$g(ph)=\alpha_{h^{-1}}(g(p))$ for all $p\in P$ and $h\in H$, and such that
$pH\to \|g(p)\|$ vanishes at infinity on $P/H$. Note that it
follows from the properness of the actions of $K$ and $H$ on $P$
that the quotient spaces $K\backslash P$ and $P/H$ are locally compact
Hausdorff spaces.

The induced algebras $\Ind_K^PD$ and $\Ind_H^PD$
are $C^*$-algebras with pointwise operations, and carry actions
$\sigma\:H\to \Aut (\Ind_K^PD)$, $\rho\:K\to \Aut
(\Ind_H^PD)$
given by
\begin{equation}\label{eq-sym-action}
\sigma_h(f)(p)=\alpha_h(f(ph))
\quad\text{and}\quad
\rho_k(g)(p)=\beta_k(g(k^{-1}p)).
\end{equation}
Theorem 1.1 of \cite{RaeIC} states that $C_c(P,D)$ can be given a
$C_c(K,\Ind_H^PD)- C_c(H, \Ind_K^PD)$ pre-imprimitivity
bimodule structure
which completes to an imprimitivity bimodule, $Z$, for the
full crossed products
$\Ind_H^PD\times_{\rho}K$ and $\Ind_K^PD\times_{\sigma}H$.
The actions and inner products are given for
$b\in C_c(K,\Ind_H^PD)\subseteq \Ind_H^PD\times K$, $x,y\in C_c(P,D)$,
and $c\in C_c(H,\Ind_K^PD)\subseteq \Ind_K^PD\times H$ as follows:
\begin{equation}\label{eq-sym-pre}
\begin{split}
b\d x(p)&=\int_K b(k,p)\beta_k(x(k^{-1}p))\Delta_K(k)^{1/2}dk\\
x\d c(p)&=\int_H \alpha_h\bigl(x(ph)c(h^{-1}, ph)\bigr)
\Delta_H(h)^{-1/2}dh\\
_{C_c(K, \Ind_H^PD)}\<x,y\> (k,p)&= \Delta_K(k)^{-1/2}
\int_H\alpha_h\bigl(x(ph)\beta_k(y(k^{-1}ph)^*)\bigr)\,dh\\
\<x,y\>_{C_c(H,\Ind_K^PD)}(h,p)&= \Delta_H(h)^{-1/2}
\int_K\beta_k\bigl(x(k^{-1}p)^*\alpha_h(y(k^{-1}ph))\bigr)\,
dk.
\end{split}
\end{equation}

By some abuse of notation, we shall often regard an
element in $C_c(K, \ind_H^PD)$ 
as a continuous function of two variables via $b(k,p) \deq b(k)(p)$,
and similarly for $C_c(H,\ind_K^P D)$.

Note that the formulas above are not precisely those given in
\cite{RaeIC}. This results from the fact that we are working
with left \emph{and} right actions on $P$, while in
\cite{RaeIC} all actions were on the left. To convert to the
two-left-actions situation of \cite{RaeIC}, just put $hp\deq ph^{-1}$.
Note also that a similar symmetric imprimitivity theorem has
been deduced independently by Kasparov in \cite[Theorem 3.15]{kas},
which also gives a Morita equivalence for the reduced crossed
products. However, the construction of Kasparov's bimodule
is given by a composition of several other bimodules, and
does not really fit our needs in this work. Nevertheless, we do need to
know that the imprimitivity bimodule constructed above
factors through an imprimitivity bimodule for the reduced
crossed products, \ie, we want to have:

\begin{prop}\label{prop-sym-reduced}
If we view $C_c(K,\Ind_H^PD)$ and $C_c(H,\Ind_K^PD)$ as dense
subalgebras of the reduced crossed products
$\Ind_H^PD\times_{\rho,r}K$ and $\Ind_K^PD\times_{\sigma,r}H$,
respectively,
then the $C_c(K,\Ind_H^PD)-C_c(H,\Ind_K^PD)$ bimodule $C_c(P,D)$
with actions and inner products as given in Equation \ref{eq-sym-pre}
  completes to an
$(\Ind_H^PD\times_{\rho,r}K)-(\Ind_K^PD\times_{\sigma,r}H)$
imprimitivity bimodule.
\end{prop}
\begin{proof} The hard work for this proof has been done in
\cite{qs:regularity};
we show how the proof follows from \cite[Lemma 4.1]{qs:regularity}.
Using \cite[Proposition 3.24]{RW}, the content of
\cite[Lemma 4.1]{qs:regularity} can
be translated
into the   statement:

\emph{If $I\subseteq \Ind_K^PD\times H$ is the kernel of a regular
representation, say $\Ind\pi$, of $\Ind_K^PD\times H$, then the ideal
$Z$-$\Ind I$ of $\Ind_H^PD\times K$ induced from $I$ via $Z$ contains
the kernel, say $J$, of the regular representation\textup(s\textup) of
$\Ind_H^PD\times K$.  }

By symmetry, we also get $\widetilde{Z}$-$\Ind J\supseteq I$,
where $\widetilde Z$ denotes the conjugate bimodule.
Since by \cite[Corollary 3.31]{RW}
$Z$-$\Ind\circ \widetilde Z$-$\Ind$ is the
identity map on the set of closed ideals of
$\Ind_H^PD\times K$,
and since induction via $Z$
and $\widetilde Z$ preserves inclusion of ideals by
\cite[Theorem 3.22 and Proposition 3.24]{RW}, we get
\[
Z\dashind I\supseteq J
=Z\dashind\circ
\widetilde Z\dashind J\supseteq Z\dashind I,
\]
and hence $J= Z$-$\ind I$.
But it follows then from \cite[Proposition 3.25]{RW}
that $Z$ factors through an imprimitivity bimodule for the
quotients
$\Ind_H^PD\times_rK\cong (\Ind_H^PD\times K)/J$ and
$\Ind_K^PD\times_rH\cong (\Ind_K^PD\times H)/I$.
Applying the respective quotient maps on the dense
subspaces $C_c(K, \Ind_H^PD)$, $C_c(P,D)$ and $C_c(H,\Ind_K^PD)$
then completes the proof.
\end{proof}

We are now going to derive a one-sided version of the above
symmetric imprimitivity theorem. The full-crossed-product
version of this
result has first been deduced by Green (it follows from
\cite[Theorem 17]{gre:local}), so we shall call it {\em Green's
imprimitivity theorem for induced algebras}.

Let $G$ be a locally compact group
and let $\alpha\:H\to \Aut A$  be an action of a
closed subgroup $H$ of $G$.
We let $G$ act on itself by left translation, we let $H$ act
on $G$ by right translation, and we let $G$ act trivially on $A$.
Then, if we put $P=G$, $K=G$, and $D=A$ in the setting of the symmetric
imprimitivity theorem,
we obtain a Morita equivalence between the
crossed products $\Ind_H^GA\times_{\rho} G$ 
on the left and $\Ind_G^GA\times_{\sigma}H$ on the
right, and another between the respective reduced crossed products.
Moreover, it is easy to check that
$\Ind_G^GA\to A$
{defined by}
$f\mapsto f(e)$
is an isomorphism which transforms $\sigma$ to $\alpha$.
Hence,
the map $\Phi\:C_c(H, \Ind_G^GA)\to C_c(H,A)$
{defined by}
$\Phi(F)(h)=F(h,e)$ 
extends to an isomorphism
$\Ind_G^GA\times_{\sigma}H\cong A\times_{\alpha}H$ (and similarly
for the reduced crossed products).

In this special situation we shall always denote the
action $\rho\:G\to \Aut (\Ind_H^GA)$ of Equation \ref{eq-sym-action}
by $\Ind\alpha$. Since we start with the trivial action of $G$ on $A$,
$\Ind\alpha$
is given by the formula
\begin{equation}\label{ind-a-eq}
\Ind\alpha_s(f)(t)=f(s^{-1}t)
\end{equation}
for $f\in \Ind_H^G A$ and $s,t\in G$.
Using the map $\Phi\:C_c(H,\Ind_G^GA)\to C_c(H,A)$
described above, it follows  from
the formulas given in Equation \ref{eq-sym-pre} that
$C_c(G,A)$ becomes a $C_c(G,\Ind_H^GA)-C_c(H,A)$
pre-imprimitivity bimodule with actions and
inner products given by the formulas
\begin{equation}
\label{eq-ind-imp}
\begin{split}
b \d x(s)
&
= \int_G b(t,s) x(t^{-1}s) \Delta_G(t)^{1/2} dt
\\
x\d c(s)
&
= \int_H \alpha_h\bigl( x(sh) c(h^{-1}) \bigr) \Delta_H(h)^{-1/2} dh
\\
_{C_c(G,\Ind_H^GA)}\<x,y\>(s,t)
&
= \Delta_G(s)^{-1/2}
\int_H \alpha_h\bigl( x(th) y(s^{-1}th)^* \bigr) \,dh
\\
\<x,y\>_{C_c(H,A)}(h)
&
= \Delta_H(h)^{-1/2}
\int_G x(t^{-1})^* \alpha_h(y(t^{-1}h)) \,dt.
\end{split}
\end{equation}
More precisely, we obtain:

\begin{thm}
\label{thm-green-ind}
Regard $C_c(G,\Ind_H^GA)$ and $C_c(H,A)$
as dense subalgebras of the full crossed products
$\Ind_H^GA\times_{\Ind\alpha}G$ and $A\times_{\alpha}H$, respectively.
Then $C_c(G,A)$ with actions and inner products as given in
Equation \ref{eq-ind-imp} completes to
an  $(\Ind_H^GA\times_{\Ind\alpha}G)-(A\times_{\alpha}H)$
imprimitivity bimodule $V_H^G(A)^f$
\textup(the $f$ stands for ``full crossed
products''\textup).

Similarly, if we regard $C_c(G,\Ind_H^GA)$ and $C_c(H,A)$ as dense
subalgebras
of the respective reduced crossed products, $C_c(G,A)$ completes
to give a $(\Ind_H^GA\times_{\Ind\alpha,r}G)-(A\times_{\alpha,r}H)$
imprimitivity bimodule $V_H^G(A)^r$.
\end{thm}

We usually suppress the superscripts in the notation
of the bimodules, if confusion seems unlikely. In fact, in the main
body of our work we use \emph{reduced} crossed products almost
exclusively.

As an important special case of Green's imprimitivity theorem for
induced actions, we shall now derive Green's original imprimitivity
theorem, which we call \emph{Green's imprimitivity theorem for induced
representations}, since it was used by Green to develop a very strong
version of Mackey's imprimitivity theorem for induced representations
of locally compact groups.  Thus, the original motivation to develop
this theorem was to solve the following problem:

\emph{Suppose that $\alpha\:G\to \Aut A$ is an action, and $H$ is a
closed subgroup of $G$.  Assume further that $\pi\times V$ is a
representation of $A\times_{\alpha}G$ on a Hilbert space $\H$.  When
is this representation equivalent to a representation
$\Ind_H^G(\sigma\times W)$ induced from a representation $\sigma\times
W$ of $A\times_{\alpha|}H$?}

Of course, before one can solve this problem, one has to make clear
what the construction of the induced representation
$\Ind_H^G(\sigma\times W)$ should be.  Green's approach to this
problem was to use the modern bimodule techniques of Rieffel
\cite{RieIR} to define induction via a certain natural right-Hilbert
$(A\times_{\alpha}G)-(A\times_{\alpha|}H)$ bimodule.  We now derive
this bimodule as a special case of \thmref{thm-green-ind} above.

For this assume that $(A,G,\alpha)$ is an action and
$H$ is a
closed subgroup of $G$.
We observe that the map
$\phi\:\Ind_H^G(A,\alpha|)\to C_0(G/H,A)\cong A\otimes C_0(G/H)$
{defined by}
$\phi(f)(sH)=\alpha_s(f(s))$
is an $\Ind\alpha-(\alpha\otimes \tau)$
equivariant isomorphism,
where $\tau\:G\to \Aut C_0(G/N)$ is given by
$\tau_s(f)(tN)=f(s^{-1}tN)$, as in \exref{ex-elementary}.
Moreover, if we define
$\Psi\:C_c(G,A)\to C_c(G,A)$
{by}
$\Psi(x)(s)=\alpha_s(x(s))$,
then the formulas given in Equation \ref{eq-ind-imp} transform into
\begin{equation}
\label{eq-green-imp-thm}
\begin{split}
b \d x(s)
&
= \int_G b(t,sH) \alpha_t(x(t^{-1}s)) \Delta_G(t)^{1/2} dt
\\
x \d c(s)
&
= \int_H x(sh) \alpha_{sh}(c(h^{-1})) \Delta_H(h)^{-1/2} dh
\\
_{C_c(G,C_0(G/H,A))}\<x,y\>(s,tH)
&
= \Delta_G(s)^{-1/2}
\int_H x(th) \alpha_s(y(s^{-1}th)^*) \,dh
\\
\<x,y\>_{C_c(H,A)}(h)
&
= \Delta_H(h)^{-1/2}
\int_G \alpha_t\bigl( x(t^{-1})^* y(t^{-1}h) \bigr) \,dt,
\end{split}
\end{equation}
and as a direct consequence of \thmref{thm-green-ind} we obtain:

\begin{thm}
\label{thm-green-imp}
With the above formulas for the actions and inner products,
$C_c(G,A)$ completes to give
$((A \otimes C_0(G/H)) \times_{\alpha \otimes \tau} G)-(A
\times_{\alpha|} H)$
and
$((A \otimes C_0(G/H)) \times_{\alpha \otimes \tau,r} G)-(A
\times_{\alpha|,r} H)$
imprimitivity bimodules $X_H^G(A)^f$ and $X_H^G(A)^r$,
respectively.
\end{thm}

In practice, the dense subalgebra $C_c(G,C_0(G/H,A))$ of $(A \otimes
C_0(G/H)) \times_{\alpha \otimes \tau,r} G$ is actually replaced with
the smaller (but still dense) subalgebra $C_c(G \times G/H,A)$, and
$C_c(G,A)$ is regarded as a $C_c(G \times G/H,A) - C_c(H,A)$
pre-imprimitivity bimodule.  In
the (important) special case $H = \{e\}$, it takes a little effort to
remember that it is the second variable that comes from $C_0(G,A)$ (so
that for $b \in C_c(G \times G,A) \subseteq C_0(G,A) \times_r G$ we
have $b(s,t) = b(s)(t)$, where $b(s)$ is an element of the
$C^*$-algebra $C_0(G,A)$).

Let us briefly explain what this theorem has to do with the
problem mentioned above. For this let
$(k_A\otimes k_{C_0(G/H)}, k_G)$ denote the canonical maps of
$(A\otimes C_0(G/H), G)$ into $M\bigl((A\otimes
C_0(G/H))\times_{\alpha\otimes\tau}G\bigr)$.
Then we get a nondegenerate $*$-homomorphism
\[
k_A\times k_G\:A\times_{\alpha}G\to
M\bigl((A\times C_0(G/H))\times_{\alpha\otimes \tau}G\bigr).
\]
Identifying $M\bigl((A\otimes C_0(G/H))\times_{\alpha\otimes
\tau}G\bigr)$
with $\L(X_H^G(A)^f)$ (see, for example, \cite[Corollary
2.54]{RW}), this makes $X_H^G(A)^f$ into a right-Hilbert
$(A\times_{\alpha}G)-(A\times_{\alpha|}H)$ bimodule, which we can use
to induce representations from $A\times_{\alpha|}H$ to
$A\times_{\alpha}G$.  On the other side, since $X_H^G(A)^f$ is a
$((A\otimes
C_0(G/H))\times_{\alpha\otimes\tau}G)-(A\times_{\alpha|}H)$
imprimitivity bimodule, induction via $X_H^G(A)^f$ provides an
\emph{equivalence} between the representation spaces of
$A\times_{\alpha|}H$ and $(A\otimes
C_0(G/H))\times_{\alpha\otimes\tau}G$.  It follows from the definition
of the action of $A\times_{\alpha}G$ on $X_H^G(A)^f$, that if
$\sigma\times W$ is a representation of $A\times_{\alpha|}H$ and
$(\pi\otimes \mu)\times V$ is the representation of $(A\otimes
C_0(G/H))\times_{\alpha\otimes\tau}G$ induced from $\sigma\times W$
via $X_H^G(A)^f$, then $\pi\times V$ is the representation of
$A\times_{\alpha}G$ induced from $A\times_{\alpha|}H$ via
$X_H^G(A)^f$.  Thus we arrive at the following answer to the
above-stated problem:

\begin{thm}[Mackey-Green]
\label{Mackey-thm}
Let $\pi\times V$ be a representation of $A\times_{\alpha}G$
on a Hilbert space $\H$.
Then $\pi\times V$ is induced from a representation
$\sigma\times W$ of $A\times_{\alpha|}H$ 
\textup(via $X_H^G(A)^f$\textup)
if and only if there exists a nondegenerate representation
$\mu\:C_0(G/H)\to \B(\H)$ such that $\mu$ and $\pi$ have
commuting images in $\B(\H)$ and such that $(\pi\otimes \mu, V)$
is a covariant representation of $(A\otimes
C_0(G/H))\times_{\alpha\otimes\tau}G$.
\end{thm}

We finally want to give an interpretation of the above imprimitivity
theorems in terms of duality theory. For this we specialize
even further
to the case where $N=H$ is \emph{normal} in $G$. If $\alpha$ is
an action of $G$, we can form the dual coaction
$\widehat\alpha$ of $G$ on the crossed product $A\times_{\alpha}G$,
and the dual coaction $\widehat\alpha^n$ of $G$ on $A\times_{\alpha,r}G$.
As described in \exref{ex-restrict}, we may restrict the
coactions $\widehat\alpha$ and $\widehat\alpha^n$ to coactions
$\widehat\alpha|$ and $\widehat\alpha^n|$ of $G/N$, respectively.
\thmref{thm-isom} provides
a canonical isomorphism
\[
(A\otimes C_0(G/N))\times_{\alpha\otimes
\tau}G\cong (A\times_{\alpha}G)\times_{\widehat\alpha|}G/N
\]
and \thmref{thm-isom-red} provides an isomorphism
\[
(A\otimes C_0(G/N))\times_{\alpha\otimes\tau,r}G
\cong
(A\times_{\alpha,r}G)\times_{\alpha^n|}G/N.
\]
Thus, replacing stable isomorphism by Morita equivalence, the
above-derived Morita equivalences should be regarded as a
generalization of the Imai-Takai duality theorems (see
\thmref{thm-ImaiTakai} and \thmref{thm-ImaiTakaired}):

\begin{thm}\label{thm-green-dual}
Suppose that $(A,G,\alpha)$ is an action and $N$ is a closed normal
subgroup of $G$.  Then the above constructions provide
\textup(natural\textup)
$$((A\times_{\alpha}G)\times_{\widehat\alpha|}G/N)
-(A\times_{\alpha|}N)\quad\text{and}\quad
((A\times_{\alpha,r}G)\times_{\widehat\alpha^n|}G/N)
-(A\times_{\alpha|,r}N)$$
imprimitivity bimodules, still denoted $X_N^G(A)^f$ and $X_N^G(A)^r$,
respectively.
\end{thm}
The main purpose of this paper is to study how ``natural'' these
bimodules are, and to study various actions and coactions on these
bimodules, which fit with certain canonical actions and coactions on
the algebras involved.

\section{Mansfield's imprimitivity bimodule}\label{sec-appmansfield}

Starting with a coaction $(A,G,\delta)$, in \cite{ManIR2} Mansfield
provided a dual mirror to \thmref{thm-green-dual}, at least if $N$ is
an \emph{amenable} closed normal subgroup of $G$.  Later, in
\cite{kq:imprimitivity}, Mansfield's result was generalized to
arbitrary closed normal subgroups $N$ of $G$, provided the coaction
$\delta$ satisfies certain extra conditions.  Since those conditions are
always satisfied if $\delta$ is \emph{normal}, we may formulate:

\begin{thm}
[{\cf~\cite[Theorem 27]{ManIR2},\cite[Theorem 3.3]{kq:imprimitivity}}]
\label{thm-mans-imp}
Assume that $(A,G,\delta)$ is a nondegenerate normal coaction,
and let $N$ be a closed normal subgroup of $G$.  There exists an
$(A\times_{\delta}G)\times_{\hat\delta,r}N - A\times_{\delta|}G/N$
imprimitivity bimodule $Y_{G/N}^G(A)$.
\end{thm}

Of course, in the same way as \thmref{thm-green-dual} should be viewed
as a generalization of the Imai-Takai duality theorem, this theorem
should be viewed as a generalization of Katayama's duality theorem
(see \thmref{thm-kat}).

Let us be a bit more precise about how the above-stated result follows
from \cite[Theorem 3.3]{kq:imprimitivity}.  For this let $(j_A, j_G)$
denote the canonical maps of $(A, C_0(G))$ into
$M(A\times_{\delta}G)$.  Restricting $j_G$ to $C_0(G/N)\subseteq
M(C_0(G))$ gives a homomorphism $j_G|\:C_0(G/N)\to
M(A\times_{\delta}G)$, and one can check without too much pain that
the pair $(j_A, j_G|)$ is then a covariant homomorphism of $(A, G/N,
\delta|)$ into $M(A\times_{\delta}G)$.  Now, \cite[Theorem
3.3]{kq:imprimitivity} states that the constructions of Mansfield
provide a $(A\times_{\delta}G\times_{\hat\delta|,r}N)-(j_A\times
j_G|)(A\times_{\delta|}G/N)$ imprimitivity bimodule $Y_{G/N}^G$.  On
the other hand, if $\delta$ is normal, then it follows from
\cite[Lemma 3.2]{kq:imprimitivity} that $\delta|$ is normal, too, and
that $j_A\times j_G|\:A\times_{\delta|}G/N\to M(A\times_{\delta}G)$ is
faithful.  Thus in this situation we may identify $(j_A\times
j_G|)(A\times_{\delta|}G/N)$ with $A\times_{\delta|}G/N$.  The
following example shows that \thmref{thm-mans-imp} does not hold for
arbitrary coactions.

\begin{ex}
\label{ex-mans}
Consider the coaction $\delta_G$ of $G$ on $C^*(G)$, and let $N=G$.
The restriction of $\delta_G$ to $G/G$ is of course the trivial
coaction of the trivial group, so we get
$C^*(G)\times_{\delta_G|}G/G=C^*(G)$.  On the other hand, we have
$C^*(G)\times_{\delta_G}G\cong C_r^*(G)\times_{\delta_G^n}G\cong
\K(L^2(G))$ by \exref{ex-groupcrossed} and \propref{prop-dualnormal}.
Applying the Imai-Takai duality theorem (see
\thmref{thm-ImaiTakaired}) to the trivial action of $G$ on $\bbC$, we
see that the dual action $\what{\delta_G}$ corresponds to the unitary
action $\Ad\rho$ of $G$ on $\K(L^2(G))$.  But this implies that
$C^*(G)\times_{\delta_G}G\times_{\what{\delta_G},r}G \cong
\K(L^2(G))\times_{\Ad\rho,r}G\cong \K(L^2(G))\otimes C_r^*(G)$.  Thus,
if it were true for $\delta_G$, \thmref{thm-mans-imp} would provide a
Morita equivalence between $C^*(G)$ and $C_r^*(G)$.  Such an
equivalence should rarely exist for any non-amenable group, and it
certainly does not exist if $G=\FF_2$, the free group on two
generators, since $C_r^*(\FF_2)$ is simple but $C^*(\FF_2)$ is not.
\end{ex}

Note that it is this kind of problem which makes us stick to reduced
crossed products and normal coactions in the body of the paper.

In the rest of this section we want to briefly recall the construction
of the bimodule $Y_{G/N}^G(A)$ of \thmref{thm-mans-imp}.  As in the
other imprimitivity theorems we have to work with certain dense
subalgebras, but, unfortunately, the constructions which have to be
done here are much more complicated than the corresponding
constructions in the action case (however, see \cite{ekr} for an
easier construction in the case of a dual coaction).

For a compact subset $E\subseteq G$ let $C_E(G)=\{f\in C_c(G)\mid \supp
f\subseteq E\}$, and let $A_c(G)=A(G)\cap C_c(G)$.  Let
$\phi\:C_c(G)\to C_c(G/N)$ denote the surjection
\[
\phi(f)(sN)=\int_N f(sn)\,dn,
\]
with the usual convention on the Haar measures such that
$\int_{G/N}\circ \int_N=\int_G$.  Recall from \propref{prop-nondeg}
that if $u\in A(G)$, then $\delta_u\:A\to A$ denotes the composition
$S_u \circ \delta$, where $S_u\:M(A\otimes C^*(G))\to M(A)$ is the
slice map corresponding to $u$.

\begin{defn}\label{defn-DN}
For fixed $u\in A_c(G)$
and compact $E\subseteq G$ we put
\[
\D_{(u,E,N)}=\overline{j_A(\delta_u(A))j_G|(\phi(C_E(G)))}
\subseteq M(A\times_{\delta}G),
\]
and we define
\[
\D_N=\bigcup\{\D_{(u,E,N)}\mid u\in A_c(G),E\subseteq G\text{ compact}\}.
\]
Moreover, we write $\D_{(u,E)}$ for $\D_{(u,E,\{e\})}$
and $\D$ for $\D_{\{e\}}$.
\end{defn}

It follows from Mansfield's computations in \cite{ManIR2} that $\D$ is
a dense $*$-sub\-al\-ge\-bra of $A\times_{\delta}G$ and that $\D_N$ is a
dense $*$-subalgebra of $(j_A\times j_G|)(A\times_{\delta|}G/N)$.
Thus, if $\delta$ is normal, and if we identify $A\times_{\delta|}G/N$
with its image in $M(A\times_{\delta}G)$, we can regard $\D_N$ as a
dense $*$-subalgebra of $A\times_{\delta|}G/N$.  It is also a
consequence of Mansfield's computations that there is a linear map
$\Psi\:\D\to \D_N$ such that
\[
\Psi\bigl(j_A(a)j_G(f)\bigr)=j_A(a)j_G|(\phi(f))
\quad\text{ for }
a\in \delta_{A_c(G)}(A), f\in C_c(G).
\]
Indeed it follows from \cite[Lemma 18]{ManIR2} that
for all $x,y\in \D$ the maps
$s\mapsto \hat\delta_s(x)y$ and $s\mapsto y\hat\delta_s(x)$ are
norm continuous with compact supports and that,
as an element of $M(A\times_{\delta}G)$, $\Psi(x)$
is characterized by
\begin{equation}
\label{Psi-Int}
\Psi(x)y=\int_N\hat\delta_n(x)y\,dn\quad\text{and}\quad
y\Psi(x)=\int_Ny\hat\delta_n(x)\,dn
\righttext{for}y\in \D.
\end{equation}
So, although $n\mapsto \hat\delta_n(x)$
rarely has compact support, we shall often write
\[
\Psi(x)=\int_N\hat\delta_n(x)\,dn.
\]
Actually, the following proposition shows that the latter expression
exists as a \emph{weak} integral:

\begin{prop}\label{prop-integral}
For each $x\in \D$ and $\omega\in (A\times_{\delta}G)^*$ the integral
$\int_N\omega(\hat\delta_n(x))\,dn$ exists, and the element
$\Psi(x)\in M(A\times_{\delta}G)$ is uniquely determined by the
equations
\[
\omega(\Psi(x))=\int_N\omega(\hat\delta_n(x))\,dn.
\]
\end{prop}

\begin{proof}
It follows from
\cite[Lemma~9]{ManIR2} that
$\c D$ is closed under multiplication on either side by $j_G(C_b(G))$,
and in particular by $j_G(C_c(G))$.  Hence, \cite[Lemmas~3.5, 3.8,
and 3.10 and Corollary~3.6]{QuiLD} (which are based upon
\cite[Section~2]{op:connes1} and
\cite[note added in proof]{op:connes2}) tell us that for each $x \in \c D$
there
is a unique element $\int_N\hat\delta_n(x)\,dn$ of $M(A\times G)$ such
that
\begin{equation*}
\omega\left(\int_N\hat\delta_n(x)\,dn\right)
=\int_N\omega(\hat\delta_n(x))\,dn
\righttext{for} \omega\in(A\times G)^*.
\end{equation*}
To see that this integral coincides with $\Psi(x)$, we use
Equation~\ref{Psi-Int} to see that for all $y\in \D$ and
all $\omega\in (A\times_{\delta}G)^*$ we have
\begin{align*}
\omega(\Psi(x)y)
&= \omega\left(\int_N \hat\delta_n(x) y \,dn\right)
= \int_N \omega\bigl(\hat\delta_n(x) y\bigr) \,dn
= \int_N y \d \omega\bigl(\hat\delta_n(x)\bigr) \,dn
\\&= y \d \omega\left(\int_N \hat\delta_n(x) \,dn\right)
= \omega\left(\int_N \hat\delta_n(x) \,dn \,y\right)
\end{align*}
and similarly
$\omega(y\Psi(x))=\omega\left(y\int_N\hat\delta_n(x)\,dn\right)$.
Since $\D$ is dense in $A\times_{\delta}G$, this finishes
the proof.
\end{proof}

If we view $C_c(N,\D)\subseteq C_c(N, A\times_{\delta}G)$ as a dense
subalgebra of  $A\times_{\delta}G\times_{\hat\delta|,r}N$,
then $\D$ becomes a  $C_c(N,\D)-\D_N$ pre-imprimitivity bimodule
with module actions and inner products given for
$g\in C_c(N,\D)$, $x,y\in \D$, and $c\in \D_N$ by
\begin{equation}\label{eq-actions-mans}
\begin{split}
g\d x&=\int_N g(n) \hat\delta_n(x)\Delta(n)^{1/2}dn\\
x\d c&= xc\quad\quad\quad\text{(product in $M(A\times_{\delta}G)$)}\\
_{C_c(N,\D)}\<x,y\>(n)&= x\hat\delta_n(y^*)\Delta(n)^{-1/2}\\
\<x,y\>_{\D_N}&=\Psi(x^*y)
\quad \left( =\int_N\hat\delta_n(x^*y)\,dn \right).
\end{split}
\end{equation}
The
$(A\times_{\delta}G\times_{\hat\delta|,r}N) -(A\times_{\delta|}G/N)$
imprimitivity bimodule $Y_{G/N}^G(A)$ of \thmref{thm-mans-imp} is the
completion of the $C_c(N,\D)-\D_N$ pre-imprimitivity bimodule $\D$.

%
%

\chapter{Function Spaces}
\label{indlim-chap}

In this appendix we develop machinery allowing us to work with
certain function spaces in the various multiplier algebras of
crossed products and multiplier bimodules.

\section
{The spaces $C_c(T,\c X)$ for locally convex spaces $\c X$}

We start with some general remarks on inductive limit topologies on
functions with compact support taking values in locally convex spaces.
Assume that $T$ is a locally compact space and $\c X$ is a Hausdorff
locally convex vector space over $\b C$.  Let $\c P$ be a set of
continuous seminorms on $\c X$ which generate the topology on $\c X$.
Let $C_c(T,\c X)$ denote the space of all continuous functions from
$T$ to $\c X$ with compact support.  If $K\subseteq T$ is compact,
then $C_K(T,\c X)$ denotes the space of all $f\in C_c(T,\c X)$ with
$\supp f\subseteq K$.  For each $p\in \c P$ we define a seminorm
$p_K\:C_K(T,\c X)\to \b R$ by
\[
p_K(f)=\sup_{t\in T}p(f(t)),
\]
and we topologize $C_K(T,\c X)$ via the family of seminorms
$\c P_K=\{p_K\mid p\in \c P\}$. It is clear that $f_i\to 0$
in $C_K(T,\c X)$ if and only if
$p\circ f_i\to 0$ uniformly for all $p\in \c P$.
It is also clear that if $L\subseteq T$ is a compact subset
containing $K$, then the inclusion $C_K(T,\c X)\to C_L(T,\c X)$
is a homeomorphism onto its image. Let $\c C$ denote the set of
all compact subsets of $T$ ordered by inclusion.
Then $C_c(T,\c X)=\bigcup_{K\in \c C}C_K(T,\c X)$ and we may
equip
$C_c(T,\c X)$ with the inductive limit topology
(see \cite[II section 6]{sch:top}
for the precise definitions).
The following useful properties of the
inductive limit topology on $C_c(T,\c X)$ shall be used without reference:
\begin{enumerate}
\item For all $K\in \c C$ the inclusion
$C_K(T,\c X)\to C_c(T,\c X)$ is a homeomorphism
onto its image.
\item If $\c Y$ is any locally convex space and
$\Phi\:C_c(T,\c X)\to \c Y$
is a linear map, then $\Phi$ is continuous if and only if its
restriction to $C_K(T,\c X)$ is continuous for all $K\in\c C$.
\item
If $\c Y$ is any locally convex space and $\Psi \: \c X \to \c Y$ is
a continuous linear map, then the linear map $f \mapsto \Psi \circ f
\: C_c(T,\c X) \to C_c(T,\c Y)$ is continuous.
\end{enumerate}

\begin{prop}\label{propproduct}
Suppose that $S$ and $T$ are locally compact spaces and $\c X$ is a
locally convex space. Then the natural embedding
$\Phi\:C_c(S\times T,\c X)\to C_c(S,C_c(T,\c X))$
given by  $\Phi(f)(s)=f(s,\cdot)$
is continuous.
\end{prop}
\begin{proof}
We only have to check that the restriction of $\Phi$ to
$C_{L\times K}(S\times T,\c X)$ for any
compact subset of the form $L\times K$ of $S\times T$ is continuous.
The image of $C_{L\times K}(S\times T,\c X)$ clearly lies in
$C_L(S, C_K(T,\c X))$, and if $q\:C_c(T,\c X)\to \b R$ is any
continuous
seminorm, its restriction to $C_K(T,\c X)$ is dominated by a finite sum
of seminorms of the form $p_K$ as above. So the result follows from
the fact that $p_K(\Phi(f))=p_{L\times K}(f)$ for all
$f\in C_{L\times K}(S\times T,\c X)$.
\end{proof}

In the sequel we need to know that $C_c(T,\c X)$ is complete
whenever $\c X$ is. Unfortunately, it seems to be not
clear whether this is true in general, but the following result shows
that it is true whenever $T$ is a locally compact group.

\begin{prop}\label{propcomplete}
Suppose that $T$ is a disjoint union of a collection
$(T_\lambda)_{\lambda\in \Lambda}$ of $\sigma$-compact open subsets of
$T$.  Then $C_c(T,\c X)$ is complete whenever $\c X$ is.  In
particular, if $T$ is a locally compact group, or a quotient of a
locally compact group by a closed subgroup, then $C_c(T,\c X)$ is
complete whenever $\c X$ is.
\end{prop}
\begin{proof}
If $T$ is a disjoint union of a family of open subsets $T_{\lambda}$,
$\lambda\in \Lambda$, then the natural map $C_c(T,\c X)\to
\bigoplus_{\lambda\in \Lambda}C_c(T_{\lambda},\c X)$ is an
isomorphism.  Hence by \cite[page~55]{sch:top} $C_c(T,\c X)$ is complete
if and only if $C_c(T_{\lambda},\c X)$ is complete for all $\lambda\in
\Lambda$.  Thus we have reduced to the case where $T$ is
$\sigma$-compact.  But in this case there exists an increasing
sequence $(K_n)_{n\in \b N}$ such that $T=\cup_{n\in \b N}K_n$ and
$C_c(T,\c X)$ is the inductive limit of the sequence $C_{K_n}(T,\c
X)$.  Using \cite[page~59]{sch:top} it follows that $C_c(T,\c X)$ is
complete if $C_{K_n}(T,\c X)$ is complete for all $n\in \b N$.  But
standard arguments show that each $C_{K_n}(T,\c X)$ is complete
whenever $\c X$ is.

If $T$ is a locally compact group and $V$ is a compact neighborhood of
the identity in $T$, then the subgroup generated by $V$ is open and
$\sigma$-compact.  Thus $T$ is a disjoint union of open
$\sigma$-compact subsets and the last assertion follows from the
above.  Finally, if $T$ is a quotient of a locally compact group by a
closed subgroup, the result follows from the fact that the quotient
map is open and continuous.
\end{proof}

We need to integrate functions in $C_c(T,\c X)$ with respect to a
Radon measure $\mu$ on $T$.  For $f\in C_c(T,\c X)$ the integral
$\int_T f(t) d{\mu}(t)$ is defined as the unique element $y\in \c X$
(if it exists) such that
\[
x'(y)=\int_Tx'(f(t))\,d{\mu}(t)
\]
for every continuous linear functional $x'$ on $\c X$.  By
\cite[Theorem 3.27]{RudFA} such a $y$ always exists if the convex hull
of $f(T)$ has compact closure in $\c X$.  But if $\c X$ is complete,
this latter property of $f(T)$ follows from \cite[Theorem 3.4]{RudFA}
and the fact that the closure of a set $E\subseteq \c X$ is compact if
and only if $E$ is totally bounded \cite[page~198]{KelGT}.  Moreover,
if $p$ is any continuous seminorm on $\c X$, then similar arguments as
used in \cite[Theorem 3.29]{RudFA} show that
\[
p\left(\int_Tf(t)\,d{\mu}(t)\right)\leq \int_Tp(f(t))\,d{\mu}(t)
\]
for all $f\in C_c(T,\c X)$.

We are now collecting some important properties of integration:

\begin{prop}\label{propint}
Assume that $\c X$ is a complete locally convex space, $\mu$ is a
Radon measure on $T$ and $S$ is a locally compact space which
satisfies the conditions of \propref{propcomplete}.  Then
\begin{enumerate}
\item The linear map
$\int\:C_c(T,\c X)\to \c X$ given by
$f\mapsto \int_T f(t) d{\mu}(t)$
is continuous\textup;
\item  $\int_T\:C_c(S\times T, \c X)\to C_c(S,\c X)$ given by
$f\mapsto \int_T f(\cdot,t)\,d{\mu}(t) $
is well-defined and continuous. Moreover,
\[
\left(\int_Tf(\cdot, t)\,d{\mu}(t)\right)(s)=\int_T f(s,t)d{\mu}(t)
\]
for all $s\in S$.
\end{enumerate}
\end{prop}
\begin{proof} Since $\int$ is linear, it is enough to check continuity
on $C_K(T,\c X)$ for all compact $K\subseteq T$.  So let $f_i\to 0$ in
$C_K(T,\c X)$ and let $p$ be a continuous seminorm on $\c X$.  By
definition of the topology on $C_K(T,\c X)$ it follows that $p\circ
f_i\to 0$ uniformly, which implies that
\[
p\left(\int_T f(t)\,d{\mu}(t)\right)\leq \int_Tp(f(t))\,d{\mu}(t)\to 0.
\]
This proves (i).
In order to see (ii) note that by \propref{propproduct}
the map $C_c(S\times T,\c X)\to C_c(T,C_c(S,\c X))$ which maps
$f\in C_c(S\times T,\c X)$ to the function $t\mapsto f(\cdot, t)$ is
continuous. By the assumption on $S$ it follows from
\propref{propcomplete} that $C_c(S,\c X)$ is complete.
Therefore the integral $\int_Tf(\cdot,t)\, d{\mu}(t)$ makes sense
and takes its value in $C_c(S,\c X)$. The continuity of $\int_T$ follows
then from (i) and  \propref{propproduct}.
Finally, since evaluation at $s\in S$ is continuous
on $C_c(S,\c X)$ the last equality follows from the definition of
the integral.
\end{proof}

\section
{Functions in multiplier algebras and multiplier bimodules}

In what follows let $A$ be a $C^*$-algebra and let $M(A)$ denote the
multiplier algebra of $A$. We shall always write $M^{\beta}(A)$ if we
want to consider $M(A)$  with the strict topology.
While some of the results we present here could be proven in much
greater generality (involving things like separately continuous
bilinear maps among complete locally convex spaces), we choose to
emphasize the $C^*$-techniques
rather than follow a more Bourbaki-esque approach.

\begin{lem}\label{lemmulti}
Let $A$ be a $C^*$-algebra and $T$ a locally compact space.  For maps
$f,g\:T\to M(A)$ let $L_f(g)=fg$ and $R_f(g)=gf$ denote
left and right pointwise multiplication of $g$ by $f$.
\begin{enumerate}
\item If $f\:T\to A$ is norm continuous, then $L_f(g),R_f(g)\in
C_c(T,A)$ for every $g\in C_c(T,M^{\beta}(A))$, and the maps
$L_f,R_f\:C_c(T,M^{\beta}(A))\to C_c(T, A)$ are continuous.
\item If $f\:T\to M^{\beta}(A)$ is \textup(strictly\textup)
continuous, then $L_f(g),R_f(g)\in C_c(T,A)$ for every $g\in C_c(T,A)$,
and $L_f,R_f\:C_c(T,A)\to C_c(T, A)$ are continuous.
\item If $f\:T\to M^{\beta}(A)$ is \textup(strictly\textup)
continuous, then $L_f(g),R_f(g)\in C_c(T,M^{\beta}(A))$ for every
$g\in C_c(T,M^{\beta}(A))$, and $L_f,R_f\:C_c(T,M^{\beta}(A))\to C_c(T,
M^{\beta}(A))$ are continuous.
\end{enumerate}
\end{lem}
\begin{proof}
Since the natural pairings $M^{\beta}(A)\times A\to A$ and $A\times
M^{\beta}(A)\to A$ are continuous on bounded sets and multiplication
on $M^{\beta}(A)$ is strictly continuous on bounded sets, it follows
that all linear maps in (i)--(iii) are well-defined.  In order to
prove continuity in (i) we have to prove continuity on
$C_K(T,M^{\beta}(A))$ for all compact sets $K\subseteq T$.  Let
$h=f|_K$.  By Cohen's factorization theorem we may write
$h=a_1h_1=h_2a_2$ with $h_1,h_2\in C(K,A)$ and $a_1,a_2\in A$, where,
for example, $(a_1h_1)(t)=a_1h_1(t)$.  If $g_i\to 0$ in $C_K(T,
M^{\beta}(A))$, then by definition of the topology on
$C_K(T,M^{\beta}(A))$ it follows that $g_ia_1\to 0$ and $a_2g_i\to 0$
uniformly on $K$.  But by continuity of multiplication in $C(K,A)$
this implies that $ g_i=h_2a_2g_i\to 0$ and $g_if=g_ia_1h_1\to 0$
uniformly on $K$.  This proves (i).

Suppose now that $f$ is as in (ii).  If $K\subseteq T$ is compact,
$f|_K$ is a bounded strictly continuous function on $K$ with values in
$M(A)$ and hence an element in the multiplier algebra of $C(K,A)$.
Thus, if $g_i\to 0$ in $C_K(T,A)$, then $g_if\to 0$
and $fg_i\to 0$ in $C(K,A)$ and hence in $C_c(T,A)$.

Having proved (i) and (ii), (iii) follows from the fact that if
$g_i\to 0$ in $C_K(T,M^{\beta}(A))$, then $g_ifa\to 0$ uniformly on
$K$ by (i) for all $a\in A$ since $fa\in C(T,A)$, and $ag_if\to 0$ by
(ii) since $ag_i\to 0$ in $C_K(T,A)$.  Similar arguments show that
$afg_i\to 0$ and $fg_ia\to 0$ uniformly on $K$.
\end{proof}

We shall need the following easy consequence of Lemma \ref{lemmulti}.

\begin{cor}\label{cordense}
If $A$ is a $C^*$-algebra, then $C_c(T,A)$ is dense in
$C_c(T, M^{\beta}(A))$.
\end{cor}
\begin{proof}
Let $\{e_i\}_{i\in I}$ be an approximate unit in $A$ with $\|e_i\|=1$
for all $i\in I$.  We claim that $\{fe_i\}_{i\in I}\subseteq C_c(T,A)$
converges to $f$ in $C_c(T, M^{\beta}(A))$ for all $f\in
C_c(T,M^{\beta}(A))$.  Since all $fe_i$ have support in a single
compact subset $K\subseteq T$ it is enough to show that if $a \in A$
then $(fe_ia)_{i\in I}$ and $(afe_i)_{i\in I}$ converge to $fa$ and
$af$ uniformly on $K$, respectively.  For this let $\psi\in C_c(T)$
with $\psi|_K\equiv 1$.  Then convergence of $(fe_ia)_{i\in I}$ to
$fa$ follows from Lemma \ref{lemmulti} and the fact that the functions
$\psi_i(t)\deq\psi(t)e_ia$ converge to $\psi a$ in $C_c(T,A)$, and
convergence of $(afe_i)_{i\in I}$ follows from Lemma \ref{lemmulti}
and the fact that $af\in C_c(T,A)$ and $\tilde{\psi}_i(t)\deq\psi(t)e_i$
converges to $\psi 1$ in $C_c(T, M^{\beta}(A))$.
\end{proof}

\begin{prop}\label{propconv}
Let $(A,\,G,\,\alpha)$ be an action.
\begin{enumerate}
\item
The formulas
\[
f*g(s)=\int_Gf(t)\alpha_t(g(t^{-1}s))\,dt\quad\text{and}\quad
f^*(s)=\Delta_G(s^{-1})\alpha_s(f(s^{-1}))^*
\]
for convolution and involution on $C_c(G,A)$ extend to
$C_c(G,M^{\beta}(A))$ in such a way that $C_c(G,M^{\beta}(A))$ becomes
a locally convex $*$-algebra with separately continuous
multiplication.
\item
The pairing $C_c(G,M^{\beta}(A))\times C_c(G,M^{\beta}(A))\to C_c(G,
M^{\beta}(A))$ given by convolution restricts to separately continuous
pairings
\[
C_c(G,M^{\beta}(A))\times C_c(G, A)\to C_c(G, A)
\]
and
\[
C_c(G, A)\times C_c(G,M^{\beta}(A))\to C_c(G,  A).
\]
\item
If $(\pi, u)$ is a covariant homomorphism of $(A,G,\alpha)$ into
$M(D)$ for some $C^*$-algebra $D$, then the formula
\[
(\pi\times u)(f)=\int_G\pi(f(s))v_s\,ds
\righttext{for}f\in C_c(G,M^{\beta}(A))
\]
determines the unique continuous $*$-homomorphism $\pi\times u\:C_c(G,
M^{\beta}(A))\to M^{\beta}(D)$ which extends the usual integrated form
$\pi\times u$ on $C_c(G,A)$.  Moreover, $\pi\times u$ is faithful on
$C_c(G, M^{\beta}(A))$ if and only if it is faithful on $C_c(G,A)$.
\end{enumerate}
\end{prop}

\begin{proof}
Let $f,g\in C_c(G,M^{\beta}(A))$. Then
$(s,t)\mapsto f(t)\alpha_t(g(t^{-1}s))\in M(A)$
is clearly strictly continuous with compact support, so that
by \propref{propint} the formula
\[
f*g(s)=\int_G f(t)\alpha_t(g(t^{-1}s))\,dt
\]
makes sense and delivers a function $f*g\in C_c(G, M^{\beta}(A))$.
Since the $*$-operation is continuous in $M^{\beta}(A)$ it is also clear
that involution carries over to a continuous involution on
$C_c(G,M^{\beta}(A))$.

We claim that convolution on $C_c(G,M^{\beta}(A))$ is separately
continuous.  For this let $g\in C_c(G,M^{\beta}(A))$ be fixed.  Then
by linearity of $f\mapsto f*g$ it is enough to show continuity on
$C_K(G, M^{\beta}(A))$ for any compact subset $K\subseteq G$.  Choose
a compact subset $C\subseteq G\times G$ such that $(s,t)\mapsto
f(t)\alpha_t(g(t^{-1}s))$ has support in $C$ for all $f\in
C_K(G,M^{\beta}(A))$, and choose $\psi\in C_c(G\times G)$ such that
$\psi|_C\equiv 1$.  Then $f\mapsto \psi f$ with $\psi f(s,t)=\psi(s,t)
f(t)$ is a continuous mapping from $C_K(G, M^{\beta}(A))$ into
$C_c(G\times G, M^{\beta}(A))$; hence it follows from Lemma
\ref{lemmulti} that $f\mapsto \psi_f$, where
\[
\psi_f(s,t)=\psi(s,t)f(t)\alpha_t(g(t^{-1}s))=f(t)\alpha_t(g(t^{-1}s))
\]
is also a continuous map from $C_K(G, M^{\beta}(A))$ into $C_c(G\times
G,M^{\beta}(A))$.  Now the continuity of $f\mapsto f*g$ follows from
\propref{propint}.  In order to see that $g\mapsto f*g$ is also
continuous we use the same trick: chose a compact subset $C\subseteq
G\times G$ such that $(s,t)\mapsto f(t)\alpha_t(g(t^{-1}s))$ has
support in $C$ for all $g\in C_K(G,M^{\beta}(A))$, and then look at
the map $g\mapsto g\psi$ with $\psi$ as above, and
$g\psi(s,t)=\alpha_t(g(t^{-1}s))\psi(s,t)$.  This proves (i).

Arguments similar to the above together with part (i) of Lemma
\ref{lemmulti} show that $f*g$ and $g*f$ are in $C_c(G,A)$ for all $f\in
C_c(G,M^{\beta}(A))$ and $g\in C_c(G,A)$, and that the resulting
pairings $C_c(G,M^{\beta}(A))\times C_c(G,A)\to C_c(G,A)$ and
$C_c(G,A)\times C_c(G,M^{\beta}(A))\to C_c(G,A)$ are separately
continuous.  This proves (ii).

In order to verify (iii) we first observe that the integrand
$\psi_f(s) \deq \pi(f(s))v_s$ is the pointwise product of $\pi \circ f$
in $C_c(G,M^\beta(D))$ and the strictly continuous function $u$, so by
\lemref{lemmulti} the map $f \mapsto \psi_f \: C_c(G,M^\beta(A)) \to
C_c(G,M^\beta(D))$ is continuous.  Thus by \propref{propint} the
linear map $\pi \times u \: C_c(G,M^\beta(A)) \to M^\beta(D)$ is
continuous.
Since $C_c(G,A)$ is dense in $C_c(G, M^{\beta}(A))$, it follows that
this extension is unique, and the separate continuity of
multiplication and the continuity of involution imply that this
extension is a $*$-homomorphism.

Assume now that $\pi\times u$ is faithful on $C_c(G,A)$.  If $f\in
C_c(G, M^{\beta}(A))$ with $(\pi\times u)(f)=0$, it follows that
$0=\pi(a)(\pi\times u)(f)=(\pi\times u)(af)$ for all $a\in A$.  Thus
$af(s)=0$ for all $a\in A, s\in G$, which implies that $f(s)=0$ for
all $s\in G$.
\end{proof}

Recall that the canonical embedding of $C_c(G,A)$ into
$A\times_{\alpha}G$ (respectively $A\times_{\alpha,r}G$) is given by
the integrated form of the canonical maps $(i_A, i_G)$ (respectively
$(i_A^r, i_G^r)$) of $(A,G)$ into $M(A\times_{\alpha}G)$ (respectively
$M(A\times_{\alpha,r}G)$).  Thus, applying part (iii) of the above
proposition to $i_A\times i_G$ (respectively $i_A^r\times i_G^r$),
gives:

\begin{cor}
\label{corconv}
There are unique continuous injective $*$-homomorphisms
\[
i_A\times i_G\:C_c(G, M^{\beta}(A))\to M^{\beta}(A\times_{\alpha}G)
\quad\text{and}
\]
\[
i_A^r\times i_G^r\:C_c(G, M^{\beta}(A))\to M^{\beta}(A\times_{\alpha,r}G)
\]
extending the usual embeddings of $C_c(G,A)$ into $A\times_{\alpha}G$
and $A\times_{\alpha,r}G$, respectively.  Moreover, if we view $C_c(G,
M^{\beta}(A))$ as a subalgebra of $M(A\times_{\alpha}G)$ via the above
embedding, and if $(\pi,u)$ is a covariant homomorphism of
$(A,G,\alpha)$ into $M(D)$, then the integrated form $\pi\times
u\:C_c(G, M^{\beta}(A))\to M^{\beta}(D)$ of \propref{propconv}
\textup{(iii)} coincides with the restriction to $C_c(G,
M^{\beta}(A))$ of the usual extension of $\pi\times u$ to
$M(A\times_{\alpha}G)$.  If $\pi\times u$ factors through
$A\times_{\alpha,r}G$, a similar result is true for the reduced
crossed product.
\end{cor}

\begin{proof} The first part is a direct consequence of
\propref{propconv} (iii) and the second part follows
from the uniqueness of the continuous extension of
$\pi\times u\:C_c(G,A)\to M^{\beta}(D)$ to all of $C_c(G, M^{\beta}(A))$.
\end{proof}

In what follows we want to show how the above approach may be used to
get similar embeddings of the appropriate function spaces into various
multiplier bimodules.  The first and easiest example of this kind of
result is given by the situation where we have a right-Hilbert
bimodule action ${}_{(A,\alpha)}(X,\gamma)_{(B,\beta)}$ of $G$.  Then
there is a right-Hilbert $(A\times_{\alpha,r}G)-(B\times_{\beta,r}G)$
bimodule $X\times_{\gamma,r}G$ (and similarly for the full crossed
products), where $X\times_{\gamma,r}G$ is the completion of the
pre-right-Hilbert $C_c(G,A)-C_c(G,B)$ bimodule $C_c(G,X)$, with
operations given by \eqeqref{X x G} (see \chapref{functors-chap} for
more details).

In the following we denote by $M^{\beta}(X)$ the multiplier bimodule
of a right-Hilbert bimodule $_AX_B$ equipped with the strict
topology (see \defnref{def-strict}).

\begin{lem}\label{lemmor}
Suppose that $(X,\gamma)$ is as above.
\begin{enumerate}
\item All pairings among $C_c(G,A)$, $C_c(G,X)$, and $C_c(G,B)$ as given
by \eqeqref{X x G} extend uniquely to separately continuous pairings among
$C_c(G,M^{\beta}(A))$, $C_c(G,M^{\beta}(X))$, and
$C_c(G,M^{\beta}(B))$, given by the same formulas.
\item The natural
bimodule homomorphism
\[
_{C_c(G,A)}C_c(G,X)_{C_c(G,B)}\to X\times_{\gamma,r}G
\] extends
uniquely to a continuous bimodule homomorphism
\[
_{C_c(G,M^{\beta}(A))}C_c(G,M^{\beta}(X))_{C_c(G,M^{\beta}(B))}
\to M^{\beta}(X\times_{\gamma,r}G)
\] which also preserves
the inner products. Similar results hold for the full crossed products.
\end{enumerate}
\end{lem}

\begin{proof}
First assume $X$ is a right-partial $A - B$ imprimitivity
bimodule
(see Definition \ref{defn-partial}).
Let $L = L(X)$ be the linking algebra, with associated action $\nu$
of $G$. We have
\[
L \times_{\nu,r} G =
\begin{pmatrix}
A \times_{\eta,r}G & X \times_{\gamma,r} G
\\
* & B \times_{\beta,r} G
\end{pmatrix}
\]
by \lemref{lem-linkingaction},
and it follows
from \propref{multsoflink} that
\[
M(L \times_{\nu,r} G) =
\begin{pmatrix}
M(A \times_{\eta,r}G) & M(X \times_{\gamma,r} G)
\\
* & M(B \times_{\beta,r} G)
\end{pmatrix},
\]
where the strict topology on $M(L\times G)$ agrees with the product of
the strict topologies on the corners by
\propref{multsoflink}.
Since convolution in $C_c(G,
L)$ and $C_c(G,M^{\beta}(L))$ are given by the above pairings among
the corners, the desired conclusions follow in this case from the
continuous embedding of $C_c(G, M^{\beta}(L))=\left(
\begin{smallmatrix}
C_c(G,M^{\beta}(A)) & C_c(G,M^{\beta}(X)) \\
* & C_c(G,M^{\beta}(B))
\end{smallmatrix}
\right)$ into $M^{\beta}(L\times_{\nu,r}G)$
as in \corref{corconv}.

For the general case, let $K=\c K_B(X)$ be the algebra of compact
operators on $X$,
with associated action $\mu$.  Then from the above we get a pairing
$C_c(G,M^\beta(K))\times C_c(G,M^\beta(X))\to C_c(G,M^\beta(X))$
satisfying the required properties.  Letting $\phi\:A\to M(K)$ be the
associated nondegenerate homomorphism, we get a continuous
homomorphism
\[
\phi \times_r G \: C_c(G,M^\beta(A)) \to C_c(G,M^\beta(K)).
\]
It follows that the pairing
\[
C_c(G,M^\beta(A)) \times C_c(G,M^\beta(X))
\to C_c(G,M^\beta(X))
\]
has the required properties also.  This proves the result for the case
of reduced crossed products; exactly the same arguments (replacing all
reduced crossed products by full crossed products) give the
corresponding result for full crossed products.
\end{proof}

Another result of this kind is the embedding of $C_c(G,M^{\beta}(A))$
into the multiplier bimodule $M(V_H^G(A))$ of the
$(\Ind_H^GA\times_{\Ind\alpha,r}G)-(A\times_{\alpha,r}H)$
imprimitivity bimodule $V_H^G(A)$ of \thmref{thm-green-ind}, where
$(A,\alpha)$ is an action of a closed subgroup $H$ of $G$.  Recall
that $V_H^G(A)$ is the completion of the $C_c(G,\Ind_H^GA)-C_c(H,A)$
pre-imprimitivity bimodule $C_c(G,A)$ with operations given by
\eqeqref{eq-ind-imp}.
Note that the embedding $C_c(G,A) \hookrightarrow V_H^G(A)$ is
continuous from the inductive limit topology to the norm topology.

It follows from \cite[Corollary 3.2]{qr:induced} that we can view
elements in $C_c(G, M^{\beta}(\Ind_H^GA))$ as continuous functions of
$G\times G$ into $M^{\beta}(A)$.  We only have to check that for $f\in
C_c(G,M^{\beta}(\Ind_H^GA))$ the functions $(s,t)\mapsto af(s,t)$
and $(s,t)\mapsto f(s,t)a$ 
of $G\times G$ into $A$ are continuous for every $a\in A$,
which follows from the fact that for each $t\in G$ and $a\in A$ there
exists a $g\in \Ind_H^GA$ such that $g(t)=a$.

\begin{prop}\label{propgreen}
Let $G$ and $(A,\,H,\,\alpha)$ be as above.
\begin{enumerate}
\item All pairings among $C_c(G,\Ind_H^GA)$, $C_c(G,A)$, and $C_c(H,A)$
given in \eqeqref{eq-ind-imp}
extend uniquely to separately continuous pairings among
$C_c(G, M^{\beta}(\Ind_H^GA))$,
$C_c(G,M^{\beta}(A))$, and
$C_c(H, M^{\beta}(A))$,
given by the same formulas.
\item The canonical bimodule homomorphism
\[
_{C_c(G,\Ind_H^GA)}C_c(G,A)_{C_c(H,A)}\to
   V_H^G(A)
\]
extends uniquely to a continuous bimodule homomorphism
\[
_{C_c(G,M^{\beta}(\Ind_H^GA))}C_c(G,M^{\beta}(A))_{C_c(H,M^{\beta}(A))}
\to
M^{\beta}(V_H^G(A))
\]
which also preserves the inner products.
\end{enumerate}
\end{prop}

\begin{proof}
Arguments similar to those used in \propref{propconv} show that
the formulas for the pairings extend uniquely to
separately continuous
pairings among
$C_c(G, M^{\beta}(\Ind_H^GA))$, $C_c(G,M^{\beta}(A))$, and
$C_c(H, M^{\beta}(A))$.

Note that by \corref{corconv} we have continuous embeddings $C_c(G,
M^{\beta}(\Ind_H^GA))\to M^{\beta}(\Ind_H^GA\times_{\Ind\alpha,r}G)$
and $C_c(H,M^{\beta}(A))\to M^{\beta}(A\times_{\alpha,r}H)$ which
uniquely extend the canonical embeddings of $C_c(G,\Ind_H^GA)$ and
$C_c(H,A)$ into $\Ind_H^GA\times_{\Ind\alpha,r}G$ and
$A\times_{\alpha,r}H$, respectively.

We claim that the embedding $C_c(G,A)\to V_H^G(A)$ is continuous with
respect to the topology on $C_c(G,A)$ inherited from
$C_c(G,M^{\beta}(A))$ and the strict topology on $V_H^G(A)$.  Since
$C_c(G,A)$ is dense in $C_c(G,M^{\beta}(A))$ we then obtain a unique
continuous extension $C_c(G,M^{\beta}(A))\to M^{\beta}(V_H^G(A))$,
which by the separate continuity of the above pairings preserves all
actions and inner products.

The claim can be shown by factoring elements $b\in
A\times_{\alpha,r}H$ as $i_A^r(a)c$ for some $c\in
\mbox{$A\times_{\alpha,r}H$}$ and $a\in A$ and by factoring elements $d\in
\Ind_H^GA\times_{\Ind\alpha,r}G$ as $e (i_{\Ind A}^r(g))$, for some
$e\in \Ind_H^GA\times_{\Ind\alpha,r}G$ and $g\in \Ind_H^GA$. Then, for $x
\in C_c(G,M^\beta(A))$ the function $x \cdot i_A^r(a)$ given by $x
\cdot i_A^r(a)(s) = x(s) \alpha_s(a)$ is the pointwise product of $x$
with the norm-continuous function $s \mapsto \alpha_s(a)$, hence by
\lemref{lemmulti} the linear map
$x \mapsto x \cdot i_A^r(a)$ of 
$C_c(G,M^\beta(A))$ into $C_c(G,A)$
is continuous.  Since $C_c(G,A) \hookrightarrow V_H^G(A)$ is
continuous, it follows that the map $C_c(G,A) \to V_H^G(A)$ given by
$x\mapsto (x\cdot i_A^r(a))\cdot c=x\cdot b$ is continuous with
respect to the topology on $C_c(G,A)$ inherited from
$C_c(G,M^{\beta}(A))$.  Similarly, $i_{\Ind A}^r(g)\cdot x(s)= g(s)
x(s)$, which is an element of $C_c(G,A)$, and this implies that
$x\mapsto e\cdot(i_{\Ind A}^r(g)\cdot x)=d\cdot x$ is continuous for
all $d$.
\end{proof}

Throughout this work we freely use strictly continuous
multiplier-valued functions inside multiplier bimodules.  In every
case the computations can be justified by closely following the lines
of Propositions~\ref{propconv}--\ref{propgreen}, so we omit any
further details.


\backmatter


%
%

\providecommand{\bysame}{\leavevmode\hbox to3em{\hrulefill}\thinspace}

\bibliographystyle{amsplain}


\printindex

\end{document}